\documentclass[reqno,11pt]{amsart}
\usepackage{amsfonts}
\usepackage{a4wide}
\usepackage{color}
\usepackage{mathrsfs}
\usepackage{mathtools}
\usepackage{amsmath}
\usepackage{amssymb}
\usepackage{bbm}
\usepackage{esint}
\usepackage{nicefrac}
\numberwithin{equation}{section}
\usepackage[colorlinks,citecolor=green,linkcolor=blue]{hyperref}

\usepackage[latin1]{inputenc}
\input xy
\xyoption{all}

\newtheorem{theorem}{Theorem}[section]
\newtheorem{Ca}[theorem]{Corollary}
\newtheorem{Th}[theorem]{Theorem}
\newtheorem{Lm}[theorem]{Lemma}
\newtheorem{Prop}[theorem]{Proposition}
\newtheorem{Def}[theorem]{Definition}

\newtheorem{Remark}[theorem]{Remark}
\newtheorem{Problem}[theorem]{Problem}

\title{Traces of Sobolev spaces to irregular subsets of metric measure spaces}
\author{Alexander I. Tyulenev}
\address{Steklov Mathematical Institute of Russian Academy of Sciences}
\email{tyulenev-math@yandex.ru, tyulenev@mi-ras.ru}
\begin{document}
\allowdisplaybreaks
\begin{abstract}
Given $p \in (1,\infty)$, let $(\operatorname{X},\operatorname{d},\mu)$ be a metric measure space with uniformly locally doubling measure $\mu$ supporting a weak local $(1,p)$-Poincar\'e inequality. For each $\theta \in [0,p)$, we characterize the trace space of the
Sobolev $W^{1}_{p}(\operatorname{X})$-space to lower codimension-$\theta$ content regular closed subsets $S \subset \operatorname{X}$.
In particular, if
the space $(\operatorname{X},\operatorname{d},\mu)$ is Ahlfors $Q$-regular for some $Q \geq 1$ and $p \in (Q,\infty)$, then we get an intrinsic description of the trace-space of the Sobolev $W^{1}_{p}(\operatorname{X})$-space to arbitrary closed nonempty set
$S \subset \operatorname{X}$.
\end{abstract}
\maketitle
\tableofcontents
\markright{Traces of Sobolev Spaces}
\section*{Introduction}
\noindent

The theory of Sobolev spaces on metric measure spaces $\operatorname{X}=(\operatorname{X},\operatorname{d},\mu)$ is an important
rapidly growing area of modern Geometric Analysis. Since no additional regularity structure on $\operatorname{X}$ is a priori assumed, it is not surprising that
the most studies available so far are related to the \textit{first-order} Sobolev spaces $W^{1}_{p}(\operatorname{X})$, $p \in (1,\infty)$.
We refer to the recent beautiful monograph \cite{HKST} and the lecture notes \cite{GP20}
containing an exhaustive treatment of the subject.
However, some natural questions concerning the spaces $W^{1}_{p}(\operatorname{X})$ remain open.
One of the most difficult and exciting among them is the so-called \textit{trace problem}, i.e., the problem of \textit{sharp intrinsic description of the trace-space} of the
space $W^{1}_{p}(\operatorname{X})$, $p \in (1,\infty)$, to different closed sets $S \subset \operatorname{X}$. In all previously known studies, this problem was considered
under some extra regularity assumptions on $S$. In the present paper we introduce a \textit{sufficiently broad} class of closed sets and
solve the corresponding trace problem for sets from that class.

In order to pose the problem precisely, we recall several concepts from analysis on metric measure spaces.
First of all, by \textit{a metric measure space} (an m.m.s., for short) we always mean a triple $\operatorname{X}=(\operatorname{X},\operatorname{d},\mu)$, where $(\operatorname{X},\operatorname{d})$
is a \textit{complete separable metric space} and $\mu$ is a \textit{Borel regular} measure giving \textit{finite positive mass}
to all balls $B_{r}(x)$ centered at $x \in \operatorname{X}$ with radii $r \in (0,\infty)$.
Furthermore, we will deal with m.m.s. $\operatorname{X}=(\operatorname{X},\operatorname{d},\mu)$ that are \textit{$q$-admissible} for some $q \in (1,\infty)$ (see Section 2.1 for the details). This means that:
\begin{itemize}
\item[\(\rm A)\)] the measure $\mu$ has the uniformly locally doubling property;

\item[\(\rm B)\)] the space $\operatorname{X}$ supports a weak local $(1,q)$-Poincar\'e inequality.
\end{itemize}
Given an m.m.s. $\operatorname{X}=(\operatorname{X},\operatorname{d},\mu)$ and a parameter $p \in (1,\infty)$,
there are at least five approaches to the definition of the first-order Sobolev spaces $W_{p}^{1}(\operatorname{X})$
commonly used in the modern literature \cite{AGS1, AGS2, Ch, Haj, KS, Shan1}.
It is remarkable that
in the case when $\operatorname{X}$ is $p$-admissible, all of them \textit{are equivalent} in an appropriate sense (see Section 2.2 for the details).
In this paper, we focus on the approach proposed by J.~Cheeger \cite{Ch} but in the equivalent modern form used in \cite{GT}.
This approach is more suitable for our framework.

\textbf{Statement of the problem.}
Recall the notion of $p$-capacity $C_{p}$ (see Subsection 1.4 in \cite{BB} for the details).
It is well known that, given $p \in (1,\infty)$ and a $p$-admissible m.m.s. $\operatorname{X}$, for every element $F \in W^{1}_{p}(\operatorname{X})$ there is a Borel representative $\overline{F}$ which
has Lebesgue points everywhere except a set of $p$-capacity zero. Any such a representative will be called
a \textit{$p$-sharp representative} of $F$. Given a set $S \subset \operatorname{X}$ of positive $p$-capacity, we define the \textit{$p$-sharp trace}
of any element $F \in W_{p}^{1}(\operatorname{X})$ to $S$ as the equivalence class (modulo $p$-capacity zero) consisting
of all pointwise restrictions of the $p$-sharp representatives of $F$ to the set $S$ and denote it by $F|_{S}$. In what follows, we will \textit{not distinguish}
between $F|_{S}$ and the pointwise restriction of any $p$-sharp representative of $F$ to $S$.
We define the \textit{$p$-sharp trace space} $W_{p}^{1}(\operatorname{X})|_{S}$ as the linear space of $p$-sharp
traces $F|_{S}$ of all elements $F \in W_{p}^{1}(\operatorname{X})$. We equip this space with the usual quotient space norm.
We also introduce the $p$-sharp trace operator $\operatorname{Tr}|_{S}:W_{p}^{1}(\operatorname{X}) \to W_{p}^{1}(\operatorname{X})|_{S}$,
which acts on $F \in W_{p}^{1}(\operatorname{X})$ by $\operatorname{Tr}|_{S}(F):=F|_{S}$. Finally, we say that $F \in W_{p}^{1}(\operatorname{X})$ is a \textit{p-sharp extension} of a given function $f: S \to \mathbb{R}$ provided that $f=F|_{S}$. The first problem we consider in the present paper can be formulated as follows.

\begin{Problem}[{\textbf{The sharp trace problem}}]
\label{TraceProblem}
Let $p \in (1,\infty)$ and let $\operatorname{X}=(\operatorname{X},\operatorname{d},\mu)$ be a $p$-admissible metric measure space. Let $S \subset \operatorname{X}$ be a closed nonempty set with $C_{p}(S) > 0$.
\begin{itemize}
\item[\((\textbf{Q}1)\)] Given a Borel function $f: S \to \mathbb{R}$, find necessary and sufficient conditions for the existence of
a $p$-sharp extension $F \in W_{p}^{1}(\operatorname{X})$ of $f$.

\item[\((\textbf{Q}2)\)] Using only geometry of the set $S$ and values of a function $f \in W_{p}^{1}(\operatorname{X})|_{S}$, compute the trace norm
$\|f|W_{p}^{1}(\operatorname{X})|_{S}\|$ up to some universal constants.

\item[\((\textbf{Q}3)\)] Does there exist a bounded linear operator $\operatorname{Ext}_{S,p}: W_{p}^{1}(\operatorname{X})|_{S} \to W_{p}^{1}(\operatorname{X})$, called a
$p$-sharp extension operator,
such that $\operatorname{Tr}|_{S} \circ \operatorname{Ext}_{S,p} = \operatorname{Id}$ on $W_{p}^{1}(\operatorname{X})|_{S}$?
\end{itemize}
\end{Problem}

In many particular cases, the concepts of the $p$-sharp trace space and the $p$-sharp extension \textit{should be relaxed} in an appropriate sense.
For example, if a given set $S \subset \operatorname{X}$ has ``a constant Hausdorff dimension'', then it is natural to use the corresponding
Hausdorff measure rather than the $C_{p}$-capacity to describe ``negligible sets''. For example, in \cite{GKS,GS,JonWal,Maly,MNS,Ry,SakSot,Shv1}
the corresponding notions of traces of Sobolev-type functions were introduced with the help of the corresponding Hausdorff-type measures rather than capacities.
However, the situation becomes more involved if we deal with a set $S$ composed of infinitely many ``pieces of different dimensions''. Clearly, in this case,
the use of a single Hausdorff-type measure in the definition of traces of Sobolev spaces can lead to ill-posed trace problems.
At the same time, the use of $C_{p}$-capacities can be unnatural.

Inspired by the above observations, we introduce a \textit{more flexible concept} of traces of Sobolev functions.
Let $\operatorname{X}=(\operatorname{X},\operatorname{d},\mu)$ be an m.m.s. and let $S \subset \operatorname{X}$ be a closed nonempty set.
Given  a Borel regular locally finite measure $\mathfrak{m}$ on $\operatorname{X}$, by $L_{0}(\mathfrak{m})$ we denote
the space of all $\mathfrak{m}$-equivalence classes of Borel functions $f:\operatorname{supp}\mathfrak{m} \to \mathbb{R}$.
Assume that $\operatorname{supp}\mathfrak{m}=S$ and the measure $\mathfrak{m}$ is absolutely continuous with respect to $C_{p}$, i.e.,
for each Borel set $E \subset S$, the equality $C_{p}(E)=0$ implies the equality $\mathfrak{m}(E)=0$.
We define the $\mathfrak{m}$-trace $F|_{S}^{\mathfrak{m}}$ of any element $F \in W_{p}^{1}(\operatorname{X})$ to $S$ as the
$\mathfrak{m}$-equivalence class of the $p$-sharp trace $F|_{S}$.
By $W_{p}^{1}(\operatorname{X})|^{\mathfrak{m}}_{S}$ we denote the linear space of $\mathfrak{m}$-traces of all $F \in W_{p}^{1}(\operatorname{X})$ equipped with the corresponding quotient space norm.
We also introduce the $\mathfrak{m}$-trace operator $\operatorname{Tr}|_{S}^{\mathfrak{m}}:W_{p}^{1}(\operatorname{X}) \to W_{p}^{1}(\operatorname{X})|_{S}^{\mathfrak{m}}$
which acts on $W_{p}^{1}(\operatorname{X})$ by $\operatorname{Tr}|_{S}^{\mathfrak{m}}(F):=F|_{S}^{\mathfrak{m}}$. We say that $F \in W_{p}^{1}(\operatorname{X})$ is an \textit{$\mathfrak{m}$-extension} of a given function $f: S \to \mathbb{R}$ provided that $f=F|^{\mathfrak{m}}_{S}$.
The second problem considered in the present paper can be formulated as follows.

\begin{Problem}[{\textbf{The $\mathfrak{m}$}-\textbf{trace problem}}]
\label{MeasureTraceProblem}
Let $p \in (1,\infty)$ and let $\operatorname{X}=(\operatorname{X},\operatorname{d},\mu)$ be a $p$-admissible metric measure space.
Let $\mathfrak{m}$ be a positive locally finite Borel regular measure on $\operatorname{X}$ that is absolutely continuous with respect to $C_{p}$ and let $S=\operatorname{supp}\mathfrak{m}$.

\begin{itemize}
\item[\((\textbf{MQ}1)\)] Given  $f \in L_{0}(\mathfrak{m})$, find necessary and sufficient conditions for the existence of
an $\mathfrak{m}$-extension $F \in W_{p}^{1}(\operatorname{X})$ of $f$.

\item[\((\textbf{MQ}2)\)] Using only geometry of the set $S$, properties of $\mathfrak{m}$, and values of $f \in W_{p}^{1}(\operatorname{X})|^{\mathfrak{m}}_{S}$, compute the trace norm
$\|f|W_{p}^{1}(\operatorname{X})|^{\mathfrak{m}}_{S}\|$ up to some universal constants.

\item[\((\textbf{MQ}3)\)] Does there exist a bounded linear operator $\operatorname{Ext}_{\mathfrak{m},p}: W_{p}^{1}(\operatorname{X})|^{\mathfrak{m}}_{S} \to W_{p}^{1}(\operatorname{X})$, called an $\mathfrak{m}$-extension
operator, such that $\operatorname{Tr}|^{\mathfrak{m}}_{S} \circ \operatorname{Ext}_{\mathfrak{m},p} = \operatorname{Id}$ on $W_{p}^{1}(\operatorname{X})|^{\mathfrak{m}}_{S}$?
\end{itemize}
 \end{Problem}

\textbf{Previously known results.}
As far as we know, Problem \ref{TraceProblem} was considered
only in the case $\operatorname{X}=(\mathbb{R}^{n},\|\cdot\|_{2},\mathcal{L}^{n})$. Furthermore, this problem is still unsolved in full generality.
Below we briefly recall the most powerful particular results available in the literature.

\begin{itemize}
\item[\((\textbf{R}.1.1)\)] In \cite{Shv1, Shv3} the case $p > n$ was fully covered,
i.e., Problem \ref{TraceProblem} was solved without any additional regularity assumptions on $S$.

\item[\((\textbf{R}.1.2)\)] In the case $p \in (1,n]$, for each $d \in (n-p,n]$, Problem \ref{TraceProblem} was solved for any closed \textit{$d$-lower content regular}
(or equivalently \textit{$d$-thick}) set $S \subset \mathbb{R}^{n}$ \cite{TV1}.

\item[\((\textbf{R}.1.3)\)]
Very recently, a weakened version of Problem \ref{TraceProblem} was solved by the author without any additional regularity assumption on $S$ \cite{T1, T3}.

\end{itemize}

Now we briefly describe the available results concerning Problem \ref{MeasureTraceProblem}. Let $\operatorname{X}=(\operatorname{X},\operatorname{d},\mu)$ be a metric  measure space.
Throughout the paper, we use the symbol
$B_{r}(x)$ to denote the closed ball centered at $x \in \operatorname{X}$ with radius $r > 0$, i.e., $B_{r}(x):=\{y \in \operatorname{X}:\operatorname{d}(x,y) \le r\}$.
Since in general metric spaces the behavior of $\mu(B_{r}(x))$ is not so transparent, it was observed in \cite{Maly, MNS, GS, GKS} that codimensional
analogs of the Hausdorff contents $\mathcal{H}_{\theta,\delta}$ (see Section 2.1 for the corresponding definition) are more suitable in this case.
Following \cite{Maly, MNS}, given an m.m.s. $\operatorname{X}$
and a parameter $\theta \geq 0$, we say that a closed $S \subset \operatorname{X}$ is \textit{Ahlfors--David codimension-$\theta$ regular} provided that there exist
constants $c_{\theta,1}(S), c_{\theta,2}(S) > 0$ such that
\begin{equation}
\label{eqq.ADR}
c_{\theta,1}(S) \frac{\mu(B_{r}(x))}{r^{\theta}} \le \mathcal{H}_{\theta}(B_{r}(x) \cap S) \le c_{\theta,2}(S) \frac{\mu(B_{r}(x))}{r^{\theta}} \quad \hbox{for all} \quad (x,r) \in S \times (0,1].
\end{equation}
The class of all  Ahlfors--David codimension-$\theta$ regular sets will be denoted by $\mathcal{ADR}_{\theta}(\operatorname{X})$.

\begin{itemize}
\item[\((\textbf{R}.2.1)\)] In \cite{Shv1}
traces of Calderon--Sobolev spaces and Hajlasz--Sobolev spaces to sets $S \in \mathcal{ADR}_{0}(\operatorname{X})$ were considered.
It was assumed in \cite{Shv1} that the measure $\mu$
is globally doubling and, in addition, satisfies the \textit{reverse doubling property.}

\item[\((\textbf{R}.2.2)\)] In \cite{SakSot} traces of Besov, Lizorkin--Triebel, and Hajlasz--Sobolev functions to
porous Ahlfors--David regular closed subsets of $\operatorname{X}$ were discussed.
The paper \cite{SakSot} also shows possible relaxation of Ahlfors--David regularity by
replacing it with the Ahlfors--David codimension-$\theta$ regularity.

\item[\((\textbf{R}.2.3)\)] In \cite{Maly}, given $\theta > 0$ and a uniform domain $\Omega \subset \operatorname{X}$ whose boundary $\partial\Omega$
satisfies the corresponding Ahlfors--David codimension $\theta$ regularity condition, an exact description of traces
of Sobolev $W_{p}^{1}(\Omega)$-spaces to $\partial \Omega$ was obtained. Furthermore, very recently,
the same problem was considered for homogeneous Sobolev-type spaces $D^{1}_{p}(\Omega)$ \cite{GS}.

\item[\((\textbf{R}.2.4)\)] Very recently, an analog of Problem \ref{MeasureTraceProblem} for Banach-valued Sobolev mappings were studied in \cite{GBIZ} for the case $S \in \mathcal{ADR}_{0}(\operatorname{X})$.

\end{itemize}

\textbf{Aims of the paper.}
An analysis of the results mentioned in $(\textbf{R}.2.1)$--$(\textbf{R}.2.4)$ shows that Problem \ref{MeasureTraceProblem} was considered for subsets $S \subset \operatorname{X}$
satisfying Ahlfors--David-type regularity conditions. In particular, methods and tools available so far were found to be inapplicable even
in the case when $S=S_{1} \cup S_{2}$, where $S_{i} \in \mathcal{ADR}_{\theta_{i}}(\operatorname{X})$, $i=1,2$ with $\theta_{1} \neq \theta_{2}$
such that $S_{1} \cap S_{2} \neq \emptyset$ and $\mathcal{H}_{\max\{\theta_{1},\theta_{2}\}}(S_{1} \cap S_{2}) = 0$. This elementary obstacle shows that the classes $\mathcal{ADR}_{\theta}(\operatorname{X})$, $\theta \geq 0$ are too narrow
to built a fruitful framework for the trace problems.
Hence, it is natural to introduce a relaxation of the Ahlfors--David regularity condition \eqref{eqq.ADR} replacing the Hausdorff measure by
the corresponding Hausdorff content.
We say that a set $S \subset \operatorname{X}$
is \textit{lower codimension-$\theta$ content regular} if there exists a constant $\lambda_{S} \in (0,1]$ such that
\begin{equation}
\label{eqq.LCR}
\lambda_{\theta}(S) \frac{\mu(B_{r}(x))}{r^{\theta}} \le  \mathcal{H}_{\theta,r}(B_{r}(x) \cap S) \le \frac{\mu(B_{r}(x))}{r^{\theta}} \quad \hbox{for all} \quad (x,r) \in S \times (0,1].
\end{equation}
While the upper bound in \eqref{eqq.LCR} holds trivially, we formulate our definition in the present form for the reader convenience.
By $\mathcal{LCR}_{\theta}(\operatorname{X})$ we denote the class of all lower codimension-$\theta$ content regular subsets of $\operatorname{X}$.
This notion is a natural extension of the concept of $d$-thick sets introduced by V.~Rychkov \cite{Ry} to the case of general metric measure spaces.
Indeed, in the case $\operatorname{X}=(\mathbb{R}^{n},\|\cdot\|_{2},\mathcal{L}^{n})$ and $d \in [0,n]$, a given set $S \subset \mathbb{R}^{n}$ is $d$-thick in the sense of V.~Rychkov if and only if $S \in \mathcal{LCR}_{n-d}(\mathbb{R}^{n})$. Very recently, some interesting geometric properties of $d$-thick sets in $\mathbb{R}^{n}$ were actively studied \cite{AzSh, AzVil, Ty5}.
One can show that
$\mathcal{ADR}_{\theta}(\operatorname{X}) \subset \mathcal{LCR}_{\theta}(\operatorname{X})$ for each $\theta \geq 0$, but the inclusion is strict
in general (see Section 4 for the details). The class $\mathcal{LCR}_{\theta}(\operatorname{X})$ is very broad. For example, any
path-connected set $\Gamma \subset \mathbb{R}^{n}$ containing at least two distinct points belongs to the class $\mathcal{LCR}_{n-1}(\mathbb{R}^{n})$. Furthermore, if
an m.m.s.\ space $\operatorname{X}$ is Ahlfors $Q$-regular for some $Q > 0$, then, for each $\theta \geq Q$, every nonempty set $S \subset \operatorname{X}$ belongs
to $\mathcal{LCR}_{\theta}(\operatorname{X})$.

The \textit{aim of this paper} is to solve Problems \ref{TraceProblem} and \ref{MeasureTraceProblem}
for closed sets $S \in \mathcal{LCR}(\operatorname{X})$. Our results will
cover all the previously known results \cite{Shv1, Shv2, SakSot, TV1}. Furthermore,
we provide an illustrative Example 11.4 in which Problem \ref{MeasureTraceProblem} is solved for a set composed of two Ahlfors--David regular sets
of different codimensions. Note that even this elementary example \textit{was beyond the scope} of the previously known studies.
Finally, as a particular case of our main results, given $Q \geq 1$ and an Ahlfors $Q$-regular m.m.s.\ $\operatorname{X}$,
for each $p > Q$, we present in Example 11.5 a solution to Problem \ref{TraceProblem} for an \textit{arbitrary closed nonempty} set $S \subset \operatorname{X}$.
This example gives a natural generalization of Theorem 1.2 from \cite{Shv2}.


\textbf{Statements of the main results.}
In order to formulate the main results of the present paper, we introduce the keystone tools.

Given an m.m.s.\ $\operatorname{X}=(\operatorname{X},\operatorname{d},\mu)$ and a parameter $\theta \geq 0$, we say that a sequence of
measures $\{\mathfrak{m}_{k}\}:=\{\mathfrak{m}_{k}\}_{k=0}^{\infty}$ is $\theta$-regular if there exists $\epsilon=\epsilon(\{\mathfrak{m}_{k}\}) \in (0,1)$ such that
the following conditions are satisfied:
\begin{itemize}
\item[\((\textbf{M}1)\)] there exists a closed nonempty set $S \subset \operatorname{X}$ such that
\begin{equation}
\label{M.1}
\operatorname{supp}\mathfrak{m}_{k}=S \quad \text{for all} \quad k \in \mathbb{N}_{0};
\end{equation}

\item[\((\textbf{M}2)\)]
there exists a constant $C_{1} > 0$ such that for each $k \in \mathbb{N}_{0}$
\begin{equation}
\label{M.2}
\mathfrak{m}_{k}(B_{r}(x)) \le C_{1} \frac{\mu(B_{r}(x))}{r^{\theta}} \quad \hbox{for every} \quad x \in \operatorname{X} \quad \hbox{and every} \quad r \in (0,\epsilon^{k}];
\end{equation}

\item[\((\textbf{M}3)\)]  there exists a constant $C_{2} > 0$ such that for each $k \in \mathbb{N}_{0}$
\begin{equation}
\label{M.3}
\mathfrak{m}_{k}(B_{r}(x)) \geq
C_{2}\frac{\mu(B_{r}(x))}{r^{\theta}} \quad \hbox{for every} \quad x \in S \quad \hbox{and every} \quad r \in [\epsilon^{k},1];
\end{equation}

\item[\((\textbf{M}4)\)]  $\mathfrak{m}_{k}=w_{k}\mathfrak{m}_{0}$ with $w_{k} \in L_{\infty}(\mathfrak{m}_{0})$ for every $k \in \mathbb{N}_{0}$ and, furthermore, there exists a constant
$C_{3} > 0$ such that for all $k,j \in \mathbb{N}_{0}$
\begin{equation}
\label{M.4}
\frac{\epsilon^{\theta j}}{C_{3}} \le \frac{w_{k}(x)}{w_{k+j}(x)} \le C_{3} \quad \text{for} \quad \mathfrak{m}_{0}-\text{a.e. } x \in S.
\end{equation}
\end{itemize}
The class of all $\theta$-regular sequences of measures $\{\mathfrak{m}_{k}\}$ satisfying \eqref{M.1} will be denoted by $\mathfrak{M}_{\theta}(S)$.
Furthermore, we say that a sequence $\{\mathfrak{m}_{k}\} \in \mathfrak{M}_{\theta}(S)$ is strongly $\theta$-regular if
\begin{itemize}
\item[\((\textbf{M}5)\)] for each Borel set $E \subset S$,
\begin{equation}
\label{M5}
\varlimsup\limits_{k \to \infty}\frac{\mathfrak{m}_{k}(B_{\epsilon^{k}}(\underline{x}) \cap E)}{\mathfrak{m}_{k}(B_{\epsilon^{k}}(\underline{x}))}  > 0 \quad \text{for} \quad \mathfrak{m}_{0}-\text{a.e. } \underline{x} \in E.
\end{equation}
\end{itemize}
The class of all strongly $\theta$-regular sequences of measures $\{\mathfrak{m}_{k}\}$ satisfying \eqref{M.1} will be denoted by
$\mathfrak{M}^{str}_{\theta}(S)$.
Condition $(\textbf{M}5)$ can be considered as a multiweight generalization of the famous $A_{\infty}$-condition (compare with Section 5.7 in \cite{Stein}).
It is clear that, given $\theta \geq 0$, we always have the inclusion $\mathfrak{M}_{\theta}^{str}(S) \subset \mathfrak{M}_{\theta}(S)$.
The question of coincidence of $\mathfrak{M}_{\theta}^{str}(S)$ and $\mathfrak{M}_{\theta}(S)$ is rather subtle and will be discussed in Section 5.2 of the present paper.

\textit{The first main result} of this paper looks like an auxiliary statement. Nevertheless, this result is new and we believe that it can be of independent interest. It can be considered
as a natural and far-reaching generalization of a simple characterization of Ahlfors--David regular sets in $\mathbb{R}^{n}$ (see Definition 1.1 and Theorem 1 in Section 1 of the monograph \cite{JonWal}).
\begin{Th}
\label{Th.regularsequence}
Let $\theta \geq 0$ and let $S \subset \operatorname{X}$ be a closed nonempty set. If $S \in \mathcal{LCR}_{\theta}(\operatorname{X})$, then
$\mathfrak{M}^{str}_{\theta}(S) \neq \emptyset$. If $\mathfrak{M}_{\theta}(S) \neq \emptyset$, then $S \in \mathcal{LCR}_{\theta}(\operatorname{X})$.
\end{Th}

For an exposition of the subsequent results it will be convenient to fix the following data:

\begin{itemize}

\item[\((\textbf{D}1)\)]
a parameter $p \in (1,\infty)$ and a $p$-admissible metric measure space $\operatorname{X}=(\operatorname{X},\operatorname{d},\mu)$;

\item[\((\textbf{D}2)\)]  a parameter $\theta \in [0,p)$ and a closed set $S \in \mathcal{LCR}_{\theta}(\operatorname{X})$;

\item[\((\textbf{D}3)\)]  a sequence of measures $\{\mathfrak{m}_{k}\} \in \mathfrak{M}_{\theta}(S)$ with parameter $\epsilon=\epsilon(\{\mathfrak{m}_{k}\}) \in (0,\frac{1}{10}]$.
\end{itemize}

Given $r > 0$, we introduce an \textit{important notation} by letting $k(r):= \max\{k \in \mathbb{Z}:r \le \epsilon^{k}\}$.
Now we introduce several \textit{keystone functionals}, which will be the main tools for different characterizations of the trace spaces of $W_{p}^{1}(\operatorname{X})$.
Given $p \in [0,+\infty)$, we will use the notation
$L_{p}(\{\mathfrak{m}_{k}\}):=\cap_{k=0}^{\infty} L_{p}(\mathfrak{m}_{k})$, $L^{loc}_{p}(\{\mathfrak{m}_{k}\}):=\cap_{k=0}^{\infty} L^{loc}_{p}(\mathfrak{m}_{k})$.
Given a nonzero Borel locally finite measure $\mathfrak{m}$ on $\operatorname{X}$ and an element $f \in L_{1}^{loc}(\mathfrak{m})$,
for every Borel set $G$ with $\mathfrak{m}(G) \in (0,+\infty)$ we put
$$
\mathcal{E}_{\mathfrak{m}}(f,G):=\inf_{c \in \mathbb{R}}\frac{1}{\mathfrak{m}(G)}\int\limits_{G}|f(x)-c|\,d\mathfrak{m}(x).
$$
For each $r > 0$, we put
\begin{equation}
\label{eqq.tricky_average}
\widetilde{\mathcal{E}}_{\mathfrak{m}}(f,B_{r}(x)):=
\begin{cases}
\mathcal{E}_{\mathfrak{m}}(f,B_{2r}(x)), \quad \text{if} \quad B_{r}(x) \cap \operatorname{supp}\mathfrak{m} \neq \emptyset;\\
0, \quad \text{if} \quad B_{r}(x) \cap \operatorname{supp}\mathfrak{m} = \emptyset.
\end{cases}
\end{equation}
Now, given $f \in L_{1}^{loc}(\{\mathfrak{m}_{k}\})$, we define the \textit{$\{\mathfrak{m}_{k}\}$-Calder\'on maximal function} by the formula
\begin{equation}
\notag
f^{\sharp}_{\{\mathfrak{m}_{k}\}}(x):=\sup\limits_{r \in (0,1]}\frac{1}{r}\widetilde{\mathcal{E}}_{\mathfrak{m}_{k(r)}}(f,B_{r}(x)), \quad x \in \operatorname{X}.
\end{equation}
Furthermore, given $p \in (1,\infty)$, we consider the \textit{Calder\'on functional} on the space $L_{1}^{loc}(\{\mathfrak{m}_{k}\})$ (with values in $[0,+\infty]$) by letting, for each $f \in L_{1}^{loc}(\{\mathfrak{m}_{k}\})$,
\begin{equation}
\label{eq.main1}
\mathcal{CN}_{p,\{\mathfrak{m}_{k}\}}(f):=\|f^{\sharp}_{\{\mathfrak{m}_{k}\}}|L_{p}(\mu)\|.
\end{equation}
Note that if $\operatorname{X}=S=\mathbb{R}^{n}$ and $\mathfrak{m}_{k}:=\mathcal{L}^{n}$, $k \in \mathbb{N}_{0}$, then the
mapping $f^{\sharp}_{\{\mathfrak{m}_{k}\}}$ coincides with the classical maximal function $f^{\sharp}$ introduced by Calder\'on \cite{Calderon}. Furthermore,
in \cite{Calderon} Calder\'on proved that, for $p \in (1,\infty]$, a function $f$ lies in $W_{p}^{1}(\mathbb{R}^{n})$ if and only
if $f$ and $f^{\sharp}$ are both in $L_{p}(\mathbb{R}^{n})$. This justifies the name of our functional $\mathcal{CN}_{p,\{\mathfrak{m}_{k}\}}$.

Given $p \in (1,\infty)$ and $c > 1$, we also introduce the \textit{Brudnyi--Shvartsman functional} on $L_{1}^{loc}(\{\mathfrak{m}_{k}\})$ (with values in $[0,+\infty]$) by letting, for each $f \in L_{1}^{loc}(\{\mathfrak{m}_{k}\})$,
\begin{equation}
\begin{split}
\label{eq.main2}
\mathcal{BSN}_{p,\{\mathfrak{m}_{k}\},c}(f):=
\sup&\Bigl(\sum\limits_{i=1}^{N} \frac{\mu(B_{r_{i}}(x_{i}))}{r^{p}_{i}}\Bigl(\widetilde{\mathcal{E}}_{\mathfrak{m}_{k(r_{i})}}(f,B_{cr_{i}}(x_{i}))\Bigr)^{p}\Bigr)^{\frac{1}{p}},
\end{split}
\end{equation}
where the supremum is taken over all finite families of closed balls $\{B_{r_{i}}(x_{i})\}_{i=1}^{N}$ such that:

\begin{itemize}
\item[\((\textbf{F}1)\)] $B_{r_{i}}(x_{i}) \cap B_{r_{j}}(x_{j}) = \emptyset$ if $i \neq j$;

\item[\((\textbf{F}2)\)] $\max\{r_{i}:i=1,...,N\} \le 1$;

\item[\((\textbf{F}3)\)] $B_{cr_{i}}(x_{i}) \cap S \neq \emptyset$ for all $i \in \{1,...,N\}$.
\end{itemize}
Note that if $\operatorname{X}=S=\mathbb{R}^{n}$ and $\mathfrak{m}_{k}:=\mathcal{L}^{n}$, $k \in \mathbb{N}_{0}$, then the functional $\mathcal{BSN}_{p,\{\mathfrak{m}_{k}\},c}$ is very close in spirit to that  
used by Brudnyi \cite{Brud} to characterize
the Sobolev-type spaces on $\mathbb{R}^{n}$. In the case when $\operatorname{X}=\mathbb{R}^{n}$, $p > n$, and $S \subset \mathbb{R}^{n}$ is an arbitrary closed nonempty set,
our functional is also very close in spirit to the corresponding functionals used by Shvartsman in \cite{Shv2,Shv3}. These observations justify the name of our functional.

Given a Borel set $S \subset \operatorname{X}$ and a parameter $\sigma \in (0,1]$, we say that a ball $B$ is \textit{$(S,\sigma)$-porous} if
there is a ball $B' \subset B \setminus S$ such that $r(B') \geq \sigma r(B)$.
Furthermore, given $r \in (0,1]$, we put
\begin{equation}
\label{eqq.notation_porous}
S_{r}(\sigma):=\{x \in S:B_{r}(x) \text{ is } (S,\sigma)\text{-porous}\}.
\end{equation}
We say that $S$ is $\sigma$-porous if $S=S_{r}(\sigma)$ for all $r \in (0,1]$. Porous sets arise naturally in many areas of modern geometric analysis (see, for example, the survey \cite{Shmerkin}).
In the classical Euclidean settings, the porosity properties of lower content regular sets were studied in \cite{Ty5}.

We define a natural analog of the Besov seminorm. More precisely,
given $p \in (1,\infty)$ and $\sigma \in (0,1]$, we introduce the \textit{Besov functional} on $L_{1}^{loc}(\{\mathfrak{m}_{k}\})$ by letting,
for each $f \in L_{1}^{loc}(\{\mathfrak{m}_{k}\})$,
\begin{equation}
\label{eq.main3}
\mathcal{BN}_{p,\{\mathfrak{m}_{k}\},\sigma}(f):=\|f^{\sharp}_{\{\mathfrak{m}_{k}\}}|L_{p}(S,\mu)\|+\Bigl(\sum\limits_{k=1}^{\infty}
\epsilon^{k(\theta-p)}\int\limits_{S_{\epsilon^{k}}(\sigma)}\Bigl(\mathcal{E}_{\mathfrak{m}_{k}}(f,B_{\epsilon^{k}}(x))\Bigr)^{p}\,d\mathfrak{m}_{k}(x)\Bigr)^{\frac{1}{p}}.
\end{equation}
If the space $\operatorname{X}$ is Ahlfors $Q$-regular for some $Q > 0$, $S \subset \operatorname{X}$
is a closed Ahlfors $\theta$-regular set for some $\theta \in (0,Q)$, and $\mathfrak{m}_{k}=\mathcal{H}_{Q-\theta}\lfloor_{S}$, $k \in \mathbb{N}_{0}$, then
the functional $\mathcal{BN}_{p,\{\mathfrak{m}_{k}\},\sigma}$ coincides with the corresponding Besov seminorm \cite{SakSot}. This justifies the name of
our functional.

\textit{The second main result} of the present paper gives answers to questions $(\textbf{MQ}1)$ and $(\textbf{MQ}2)$
in Problem \ref{MeasureTraceProblem}. Namely, we present several equivalent characterizations of the trace space.
It is important that condition \eqref{M.4} implies that $W_{p}^{1}(\operatorname{X})|^{\mathfrak{m}_{0}}_{S}=W_{p}^{1}(\operatorname{X})|^{\mathfrak{m}_{k}}_{S}$
for all $k \in \mathbb{N}_{0}$.

\begin{Th}
\label{Th.SecondMain}
If $\{\mathfrak{m}_{k}\} \in \mathfrak{M}^{str}_{\theta}(\operatorname{X})$, $c \geq \frac{3}{\epsilon}$ and $\sigma \in (0,\frac{\epsilon^{2}}{4c})$, then, given $f \in L_{1}^{loc}(\{\mathfrak{m}_{k}\})$, the following conditions are equivalent:

\begin{itemize}
\item[\((i)\)] $f \in W_{p}^{1}(\operatorname{X})|_{S}^{\mathfrak{m}_{0}}$;

\item[\((ii)\)]  $\operatorname{CN}_{p,\{\mathfrak{m}_{k}\}}(f):=\|f|L_{p}(\mathfrak{m}_{0})\|+\mathcal{CN}_{p,\{\mathfrak{m}_{k}\}}(f) < +\infty$;

\item[\((iii)\)] $\operatorname{BSN}_{p,\{\mathfrak{m}_{k}\},c}(f):=\|f|L_{p}(\mathfrak{m}_{0})\|+\mathcal{BSN}_{p,\{\mathfrak{m}_{k}\},c}(f) < +\infty$;

\item[\((iv)\)] $\operatorname{BN}_{p,\{\mathfrak{m}_{k}\},\sigma}(f):=\|f|L_{p}(\mathfrak{m}_{0})\|+\mathcal{BN}_{p,\{\mathfrak{m}_{k}\},\sigma}(f) < +\infty$.


\end{itemize}
Furthermore, for each $c \geq \frac{3}{\epsilon}$ and $\sigma \in (0,\frac{\epsilon^{2}}{4c})$, for every $f \in L_{1}^{loc}(\{\mathfrak{m}_{k}\})$,
\begin{equation}
\label{eq.717}
\begin{split}
&\|f|W_{p}^{1}(\operatorname{X})|_{S}^{\mathfrak{m}_{0}}\| \approx \operatorname{CN}_{p,\{\mathfrak{m}_{k}\}}(f) \approx \operatorname{BSN}_{p,\{\mathfrak{m}_{k}\},c}(f) \approx
\operatorname{BN}_{p,\{\mathfrak{m}_{k}\},\sigma}(f)
\end{split}
\end{equation}
with the equivalence constants independent of $f$.
\end{Th}

In Section 11 we will show the equivalence between $(i)$ and $(iv)$ in Theorem \ref{Th.SecondMain} implies Theorem 1.5 from \cite{SakSot} as a particular case.
Furthermore, in the Euclidean settings the equivalence between $(i)$ and $(iv)$ strengthens the author's joint result \cite{TV1}.

\textit{The third main result} gives answers to questions $(\textbf{Q}1)$ and $(\textbf{Q}2)$ in Problem \ref{TraceProblem}.
\begin{Th}
\label{Th.FourthMain}
A Borel function $f: S \to \mathbb{R}$ belongs to $W_{p}^{1}(\operatorname{X})|_{S}$
if and only if the following conditions hold:
\begin{itemize}
\item[\((A)\)] the $\mathfrak{m}_{0}$-equivalence class $[f]_{\mathfrak{m}_{0}}$ of $f$ belongs to $W_{p}^{1}(\operatorname{X})|^{\mathfrak{m}_{0}}_{S}$;

\item[\((B)\)] there exists a set $\underline{S}_{f} \subset S$ with $C_{p}(S \setminus \underline{S}_{f})=0$ such that
\begin{equation}
\label{eq.delicate}
\lim\limits_{k \to \infty}\fint\limits_{B_{\epsilon^{k}}(x)}|f(\underline{x})-f(y)|\,d\mathfrak{m}_{k}(y)=0 \quad \text{for all} \quad \underline{x} \in \underline{S}_{f}.
\end{equation}
\end{itemize}
Furthermore, for each $c \geq \frac{3}{\epsilon}$ and $\sigma \in (0,\frac{\epsilon^{2}}{4c})$,
for every $f \in L_{1}^{loc}(\{\mathfrak{m}_{k}\})$,
\begin{equation}
\begin{split}
\label{eq.main2'}
&\|f|W_{p}^{1}(\operatorname{X})|_{S}\| \approx \operatorname{CN}_{p,\{\mathfrak{m}_{k}\}}(f) \approx \operatorname{BSN}_{p,\{\mathfrak{m}_{k}\},c}(f) \approx
\operatorname{BN}_{p,\{\mathfrak{m}_{k}\},\sigma}(f)
\end{split}
\end{equation}
with the equivalence constants independent of $f$.
\end{Th}

Note that in contrast to Theorem \ref{Th.SecondMain} condition (B)
in Theorem \ref{Th.FourthMain} is delicate and important. Roughly speaking, given $f \in L_{p}(\mathfrak{m}_{0})$, the finiteness of functionals  \eqref{eq.main1}, \eqref{eq.main2} and \eqref{eq.main3} is not
sufficient for constructing a $p$-sharp extension of $f$. On the other hand, the additional condition (B)
allows one to relax restrictions on the sequence of measures. Indeed, we do not require
that $\{\mathfrak{m}_{k}\} \in \mathfrak{M}^{str}_{\theta}(\operatorname{X})$ in Theorem \ref{Th.FourthMain}. We show in Example 11.5 that if $\operatorname{X}$ is geodesic and Ahlfros $Q$-regular for some $Q \geq 1$, then
Theorem \ref{Th.FourthMain} gives an exact intrinsic description of the $p$-sharp trace space $W_{p}^{1}(\operatorname{X})|_{S}$ to \textit{arbitrary closed} nonempty set $S \subset \operatorname{X}$.

\textit{The third main result} of the present paper gives answers to questions $(\textbf{Q}3)$ and $(\textbf{MQ}3)$ of Problem \ref{TraceProblem} and Problem \ref{MeasureTraceProblem}, respectively.
Furthermore, it clarifies a deep connection between Problems \ref{TraceProblem} and \ref{MeasureTraceProblem}. That connection is given by the existence
of a canonical isomorphism between a priori different trace spaces $W_{p}^{1}(\operatorname{X})|_{S}$ and $W_{p}^{1}(\operatorname{X})|^{\mathfrak{m}_{0}}_{S}$. This fact sheds light to the main reason why the
concept of the $p$-sharp trace space was not popular in the literature.
As usual, given normed linear spaces $E_{1}=(E_{1},\|\cdot\|_{1})$ and $E_{2}=(E_{2},\|\cdot\|_{2})$, by $\mathcal{L}(E_{1},E_{2})$ we denote the linear space of all bounded linear mappings from $E_{1}$ to $E_{2}$.

\begin{Th}
\label{Th.FifthMain}
Let $\{\mathfrak{m}_{k}\} \in \mathfrak{M}^{str}_{\theta}(S)$. Then the following assertions holds:
\begin{itemize}
\item[\((1)\)] There exists an $\mathfrak{m}_{0}$-extension operator $\operatorname{Ext}:=\operatorname{Ext}_{S,\{\mathfrak{m}_{k}\}} \in \mathcal{L}(W_{p}^{1}(\operatorname{X})|^{\mathfrak{m}_{0}}_{S},W_{p}^{1}(\operatorname{X}))$;

\item[\((2)\)] There exists a $p$-sharp extension operator $\overline{\operatorname{Ext}}:=\overline{\operatorname{Ext}}_{S,\{\mathfrak{m}_{k}\},p} \in \mathcal{L}(W_{p}^{1}(\operatorname{X})|_{S},W_{p}^{1}(\operatorname{X}))$;

\item[\((3)\)] The canonical imbedding $\operatorname{I}_{\mathfrak{m}_{0}}:W_{p}^{1}(\operatorname{X})|_{S} \to W_{p}^{1}(\operatorname{X})|^{\mathfrak{m}_{0}}_{S}$ that takes $f \in W_{p}^{1}(\operatorname{X})|_{S}$ and returns the $\mathfrak{m}_{0}$-equivalence class $[f]_{\mathfrak{m}_{0}}$ of $f$
is an isometrical isomorphism;

\item[\((4)\)] We have the following commutative diagram

\xymatrix{
&W_{p}^{1}(\operatorname{X}) \ar[d]^{\operatorname{Tr}|_{S}} \ar@<1ex>[r]^{\operatorname{Id}}
&W_{p}^{1}(\operatorname{X}) \ar[d]^{\operatorname{Tr}|^{\mathfrak{m}_{0}}_{S}} \ar@<1ex>[l]^{\operatorname{Id}} \\
&W_{p}^{1}(\operatorname{X})|_{S} \ar@<1ex>[r]^{\operatorname{I}_{\mathfrak{m}_{0}}} \ar@<1ex>[u]^{\overline{\operatorname{Ext}}} & W_{p}^{1}(\operatorname{X})|^{\mathfrak{m}_{0}}_{S}
\ar@<1ex>[l]^{\operatorname{I}^{-1}_{\mathfrak{m}_{0}}} \ar@<1ex>[u]^{\operatorname{Ext}}.}
\end{itemize}

Furthermore, there exists a constant $\overline{C} > 0$
such that the operator norms of $\operatorname{Ext}$ and $\overline{\operatorname{Ext}}$ are bounded from above by $\overline{C}$.
\end{Th}

\textbf{The keystone innovations.} Note that even in the case $\operatorname{X}=(\mathbb{R}^{n},\|\cdot\|_{2},\mathcal{L}^{n})$ our results are new. Indeed,
characterizations via Brudnyi--Shvartsman-type functionals were never considered
in the literature for $p \in (1,n]$ (the case $p > n$ was considered in \cite{Shv2}).
The keystone innovations of the present paper can be summarized as follows:

\begin{itemize}
\item[\(\bullet\)] In contrast with the classical Whitney's method used in the previous investigations, we built a new extension operator by constructing, for a given $f \in L_{1}^{loc}(\{\mathfrak{m}_{k}\})$, a special
approximating sequence of functions and obtain the resulting extension as the weak limit of that sequence.

\item[\(\bullet\)]
In contrast with the previously known studies related to Problem \ref{MeasureTraceProblem}, we use the so-called ``vertical approach'' to Sobolev spaces
introduced by J.~Cheeger \cite{Ch}. This gives a natural symbiosis with our new extension operator and  leads to the Brudnyi--Shavrtsman-type characterization of the trace space.

\item[\(\bullet\)] We introduce the new concept of $\mathfrak{m}$-traces of Sobolev $W^{1}_{p}(\operatorname{X})$-spaces and investigate its relations to the $p$-sharp traces of $W^{1}_{p}(\operatorname{X})$-spaces.

\item[\(\bullet\)] We introduce the class $\mathcal{LCR}_{\theta}(\operatorname{X})$ which is a natural generalization
of the Rychkov class of $d$-thick sets to metric measure settings.

\item[\(\bullet\)] We introduce the new class of measures $\mathfrak{M}_{\theta}^{str}(\operatorname{X})$.
This allows one to obtain characterizations of the $\mathfrak{m}_{0}$-trace space using only finiteness of the corresponding functionals.

\item[\(\bullet\)] We introduce the Brudnyi--Shvartsman-type functional in metric measure settings.

\end{itemize}

\textbf{Organization of the paper.} The paper is organized as follows:
\begin{itemize}
\item[\(\bullet\)] In Section 2 we collect some background results about metric measure spaces and
Sobolev functions defined on such spaces that
will be needed for our analysis.

\item[\(\bullet\)] In Section 3 we introduce weakly noncollapsed measures and show that they possess some sort of asymptotically
doubling properties, which will be very important in proving the existence of strongly $\theta$-regular sequences of measures in Section 5.

\item[\(\bullet\)] Section 4 is devoted to some elementary properties of lower codimension-$\theta$ content regular subsets of metric measure spaces.
We also present some examples.

\item[\(\bullet\)] Section 5 is a ``technical basis'' of the paper. We prove Theorem \ref{Th.regularsequence} and study in details various properties of $\theta$-regular sequences of measures.
Furthermore, we present elementary examples of sets $S$ for which one can easily construct explicit examples of the corresponding sequences of measures.

\item[\(\bullet\)] Section 6 is devoted to investigations of some delicate pointwise properties of functions.
This section will be crucial in proving that our new extension operator is a right-inverse of the corresponding trace operator.

\item[\(\bullet\)] In Section 7 we construct our new extension operator.

\item[\(\bullet\)]  Sections 8 and 9 contain a technical foundation for the proofs of the so-called reverse and direct trace theorems respectively.

\item[\(\bullet\)] In Section 10 we prove the main results of the paper, i.e., Theorems \ref{Th.SecondMain}, \ref{Th.FourthMain} and \ref{Th.FifthMain}.

\item[\(\bullet\)] We conclude our paper with Section 11 where we show that the most part of the available results
are merely particular cases of our main results. On the other hand, we present simple examples which do not fall into the scope
of the previous investigations.
\end{itemize}

{\bf Acknowledgements.}
I am grateful to Nageswari Shanmugalingam for warm conversations that
inspired me to write this paper. I thank Igor Verbitsky for valuable remarks concerning trace theorems for Riesz potentials
and Wolf-type inequalities. I thank Pavel Shvartsman for fruitful discussions.
Finally, I am grateful to my student Roman Oleinik who found some typos in the preliminary version of the paper.

\section{Preliminaries}

This preliminary section is meant to introduce some basic material and to set the
terminology that we shall adopt in the paper.

\subsection{Geometric analysis background}

Given a metric space $\operatorname{X}=(\operatorname{X},\operatorname{d})$ and a set $E \subset \operatorname{X}$, by $\operatorname{int}E$,
$\operatorname{cl}E$ and $\partial E$ we denote the \textit{interior} of $E$, the \textit{closure} of $E$ and the \textit{boundary} of $E$ in the metric topology of $\operatorname{X}$, respectively.
\textit{If no otherwise stated, all the balls in $\operatorname{X}$ are assumed to be closed.} More precisely, we put $B_{r}(x):=\{y \in \operatorname{X}:\operatorname{d}(x,y) \le r\}$
for $(x,r) \in \operatorname{X} \times [0,+\infty).$
Given a ball $B=B_{r}(x)$ in $\operatorname{X}$ and a parameter $\lambda > 0$, we will use \textit{notation} $\lambda B:=B_{\lambda r}(x)$.

Given a metric space $\operatorname{X}=(\operatorname{X},\operatorname{d})$ and a Borel set $E \subset \operatorname{X}$, by $\mathfrak{B}(E)$, we denote the set of all Borel functions
$f: E \to [-\infty,+\infty]$. By $C(\operatorname{X})$ and $C_{c}(\operatorname{X})$
we denote the space of all continuous and compactly supported continuous functions, respectively. We equip these spaces
with the usual sup-norm. By $\operatorname{LIP}(\operatorname{X})$ we denote the class of all Lipschitz functions.
For each
$f:\operatorname{X} \to \mathbb{R}$, we define the local Lipschitz constant $\operatorname{lip}f: \operatorname{X} \to [0,+\infty]$ of $f$ as
\begin{equation}
\label{eqq.local_Lipschitz_constant}
\operatorname{lip}f(x):=
\begin{cases}
\varlimsup\limits_{y \to x}\frac{|f(y)-f(x)|}{\operatorname{d}(x,y)}, \qquad \text{$x$ is an accumulation point};\\
0, \qquad \text{$x$ is an isolated point}.
\end{cases}
\end{equation}
It is clear that
for each $f \in \operatorname{LIP}(\operatorname{X})$ the local Lipschitz constant $\operatorname{lip}f$ is finite everywhere
on $\operatorname{X}$ and belongs to $\mathfrak{B}(\operatorname{X})$. Below we summarize elementary properties of the local Lipschitz constants of Lipschitz functions. 

\begin{Prop}
\label{Prop.lipprop}
Given a metric space $\operatorname{X}=(\operatorname{X},\operatorname{d})$, the following properties hold:

\begin{itemize}
\item[\((1)\)] if $f \equiv c$ on $\operatorname{X}$ for some $c \in \mathbb{R}$, then $\operatorname{lip}f \equiv 0$ on $\operatorname{X}$;

\item[\((2)\)] If $f_{1},f_{2} \in \operatorname{LIP}(\operatorname{X})$ then
\begin{equation}
\notag
\operatorname{lip}(f_{1}+f_{2})(x) \le \operatorname{lip}f_{1}(x)+\operatorname{lip}f_{2}(x) \quad \text{for all} \quad x \in \operatorname{X};
\end{equation}

\item[\((3)\)] $\operatorname{lip}(f+c) \equiv \operatorname{lip}f$ for each $f \in \operatorname{LIP}(\operatorname{X})$ and any $c \in \mathbb{R}$.

\end{itemize}

\end{Prop}

The following fact is well known (see, for example, Corollary 1.6 in \cite{GP20}).
\begin{Prop}
\label{Prop.GP}
If $\operatorname{X}=(\operatorname{X},\operatorname{d})$ is a compact metric space, then the space $C(\operatorname{X})$ is separable.
\end{Prop}

Given a separable metric space $\operatorname{X}=(\operatorname{X},\operatorname{d})$ and a number $\epsilon \in (0,1)$,
for each $k \in \mathbb{Z}$, we denote by $Z_{k}(\operatorname{X},\epsilon)$ an arbitrary maximal $\epsilon^{k}$-separated set in $\operatorname{X}$,
i.e.,
\begin{equation}
\label{eqq.centers_of_dyadic_cubes}
Z_{k}(\operatorname{X},\epsilon)=\{z_{k,\alpha}:\alpha \in \mathcal{A}_{k}(\operatorname{X},\epsilon)\},
\end{equation}
where $\mathcal{A}_{k}(\operatorname{X},\epsilon)$ is an at most countable index set.
Furthermore,  we set
\begin{equation}
\label{eq.lattice}
\mathcal{B}_{k}(\operatorname{X},\epsilon):=\{B_{\epsilon^{k}}(z_{k,\alpha}):\alpha \in \mathcal{A}_{k}(\operatorname{X},\epsilon)\}.
\end{equation}
Finally,  we put
\begin{equation}
\label{eqq.centers_of_dyadic_cubes_2}
Z(\operatorname{X},\epsilon):=\bigcup\limits_{k \in \mathbb{Z}}Z_{k}(\operatorname{X},\epsilon) \quad \text{and} \quad
Z^{\underline{k}}(\operatorname{X},\epsilon):=\bigcup\limits_{k \geq \underline{k}}Z_{k}(\operatorname{X},\epsilon), \quad \underline{k} \in \mathbb{Z}.
\end{equation}

Given a complete separable metric space $\operatorname{X}=(\operatorname{X},\operatorname{d})$, we say that $\mathfrak{m}$ is a \textit{measure on $\operatorname{X}$} if
$\mathfrak{m}$ is a \textit{Borel regular nonzero locally finite outer measure} on $\operatorname{X}$, i.e., $\operatorname{supp}\mathfrak{m} \neq \emptyset$ and $\mathfrak{m}(B_{r}(x)) < +\infty$ for all pairs $(x,r) \in \operatorname{X} \times [0,+\infty)$.
Given a Borel set $E \subset \operatorname{X}$ and a measure $\mathfrak{m}$ on $\operatorname{X}$,
we define the restriction $\mathfrak{m}\lfloor_{E}$ of the measure $\mathfrak{m}$ to the set $E$, as usual, by the formula
\begin{equation}
\label{eq.restriction_measure}
\mathfrak{m}\lfloor_{E}(F):=\mathfrak{m}(F \cap E) \quad \text{for all Borel sets} \quad F \subset \operatorname{X}.
\end{equation}
Sometimes it will be convenient to work with the so-called \textit{weighted measures}. More precisely, if $\mathfrak{m}$ is a measure on $\operatorname{X}$, we say
that $\gamma \in \mathfrak{B}(\operatorname{X})$ is an $\mathfrak{m}$-weight if $\gamma(x) > 0$ for $\mathfrak{m}$-a.e.\ $x \in \operatorname{X}$. In this case, the symbol
$\gamma \mathfrak{m}$ will be used for the measure on $\operatorname{X}$ defined by the formula
\begin{equation}
\label{eqq.weighted_measure}
\gamma \mathfrak{m}(E):=\int\limits_{E}\gamma(x)\,d\mathfrak{m}(x) \quad \text{for every Borel set} \quad E \subset \operatorname{X}.
\end{equation}

Given a locally compact separable metric space $\operatorname{X}=(\operatorname{X},\operatorname{d})$,
following \cite{AFP}, we say that a sequence of measures $\{\mathfrak{m}_{k}\}:=\{\mathfrak{m}_{k}\}_{k=0}^{\infty}$ on $\operatorname{X}$ \textit{locally weakly converges} to a measure $\mathfrak{m}$
on $\operatorname{X}$ and write $\mathfrak{m}_{k} \rightharpoonup \mathfrak{m}$, $k \to \infty$ if
\begin{equation}
\notag
\lim\limits_{k \to \infty}\int\limits_{\operatorname{X}}\varphi(x)\,d\mathfrak{m}_{k}(x)=\int\limits_{\operatorname{X}}\varphi(x)\,d\mathfrak{m}(x) \quad \text{for every} \quad \varphi \in C_{c}(\operatorname{X}).
\end{equation}
The following fact is well known. For a detailed proof, see, for example, Corollary 1.60 in \cite{AFP}.
\begin{Lm}
\label{Lm.weakconv}
Let $\operatorname{X}=(\operatorname{X},\operatorname{d})$ be a locally compact separable metric space and let
$\{\mathfrak{m}_{k}\}:=\{\mathfrak{m}_{k}\}_{k=0}^{\infty}$ be a sequence of measures on $\operatorname{X}$ such that
\begin{equation}
\label{eq.compact}
\sup_{k \in \mathbb{N}_{0}}\mathfrak{m}_{k}(B) < +\infty \quad \text{for every ball} \quad B \subset \operatorname{X}.
\end{equation}
Then there is a locally weakly convergent subsequence $\{\mathfrak{m}_{k_{l}}\}:=\{\mathfrak{m}_{k_{l}}\}_{l=1}^{\infty}$ of the sequence $\{\mathfrak{m}_{k}\}$.
\end{Lm}

We also recall some standard properties of locally weakly convergent sequences of measures.
\begin{Prop}
\label{Prop.weakconv}
Let $(\operatorname{X},\operatorname{d})$ be a locally compact separable metric space. If
a sequence of measures $\{\mathfrak{m}_{k}\}_{k=0}^{\infty}$ on $\operatorname{X}$ locally weakly converges to a measure $\mathfrak{m}$ on $\operatorname{X}$, then, for every open set $G \subset \operatorname{X}$
and for every compact set $F \subset \operatorname{X}$,
\begin{equation}
\begin{split}
&\varliminf\limits_{k \to \infty}\mathfrak{m}_{k}(G) \geq \mathfrak{m}(G), \qquad \varlimsup\limits_{k \to \infty}\mathfrak{m}_{k}(F) \le \mathfrak{m}(F).
\end{split}
\end{equation}
\end{Prop}

Throughout the paper, by a \textit{metric measure space} (an m.m.s., for short) we always mean a triple $(\operatorname{X},\operatorname{d},\mu)$, where
$(\operatorname{X},\operatorname{d})$ is a \textit{complete separable metric space} and $\mu$ is a measure on $\operatorname{X}$ such that
$\operatorname{supp}\mu = \operatorname{X}$.

\begin{Remark}
In what follows, given an m.m.s.\ $\operatorname{X}=(\operatorname{X},\operatorname{d},\mu)$, by a measure on $\operatorname{X}$ we always mean a measure on
the metric space $(\operatorname{X},\operatorname{d})$
\end{Remark}

Given a complete separable metric space $\operatorname{X}=(\operatorname{X},\operatorname{d})$ and a measure $\mathfrak{m}$ on $\operatorname{X}$, we will assume that the collection of
$\mathfrak{m}$-measurable sets is the completion of the Borel $\sigma$-algebra with respect
to $\mathfrak{m}$. Furthermore, given a function $f \in \mathfrak{B}(\operatorname{X})$,
we put
\begin{equation}
\label{eq.m_class_equivalence}
[f]_{\mathfrak{m}}:=\{\widetilde{f}:\operatorname{X} \to [-\infty,+\infty]: \widetilde{f}(x)=f(x) \text{ for $\mathfrak{m}$-a.e. } x \in \operatorname{supp}\mathfrak{m}\}.
\end{equation}
Finally, we put $L_{0}(\mathfrak{m}):=\{[f]_{\mathfrak{m}}: f \in \mathfrak{B}(\operatorname{X})\}$ and introduce the following mapping.

\begin{Def}
\label{Def.measure_identity}
Let $\operatorname{X}=(\operatorname{X},\operatorname{d})$ be a complete separable metric space. Given a measure $\mathfrak{m}$ on $\operatorname{X}$,
we define the mapping $\operatorname{I}_{\mathfrak{m}}: \mathfrak{B}(\operatorname{X}) \to L_{0}(\mathfrak{m})$ by letting $\operatorname{I}_{\mathfrak{m}}(f):=[f]_{\mathfrak{m}}$
for all $f \in \mathfrak{B}(\operatorname{X})$.
\end{Def}
Given $p \in (0,\infty)$, by $L_{p}(\mathfrak{m})$ ($L^{loc}_{p}(\mathfrak{m})$) we denote the linear space of
$\mathfrak{m}$-equivalence classes $[f]_{\mathfrak{m}}$ of all (locally) $p$-integrable with respect to $\mathfrak{m}$ functions $f \in \mathfrak{B}(\operatorname{X})$.
By $L_{\infty}(\mathfrak{m})$ ($L^{loc}_{\infty}(\mathfrak{m})$) we denote the linear space of $\mathfrak{m}$-equivalence classes of all (locally) bounded on $\operatorname{supp}\mathfrak{m}$ Borel functions.
Without loss of generality, we will assume that, for each $[f]_{\mathfrak{m}} \in L^{loc}_{p}(\mathfrak{m})$, the functions $f \in [f]_{\mathfrak{m}}$ \textit{are defined only on} $\operatorname{supp}\mathfrak{m}$, i.e., $L^{loc}_{p}(\mathfrak{m}):=L^{loc}_{p}(\operatorname{supp}\mathfrak{m},\mathfrak{m})$.
Given $p \in [0,\infty]$ and a sequence of measures $\{\mathfrak{m}_{k}\}_{k=0}^{\infty}$ on $\operatorname{X}$, we put $L_{p}(\{\mathfrak{m}_{k}\}):=\cap_{k=0}^{\infty}L_{p}(\mathfrak{m}_{k})$
and $L^{loc}_{p}(\{\mathfrak{m}_{k}\}):=\cap_{k=0}^{\infty}L^{loc}_{p}(\mathfrak{m}_{k})$.
Given an m.m.s.\ $\operatorname{X}=(\operatorname{X},\operatorname{d},\mu)$, a parameter $p \in [0,\infty]$, and a Borel set $S \subset \operatorname{X}$,
we use \textit{notation} $L_{p}(S):=L_{p}(\mu\lfloor_{S})$.

\begin{Remark}
Typically, given a complete separable metric space $\operatorname{X}=(\operatorname{X},\operatorname{d})$, a measure $\mathfrak{m}$ on $\operatorname{X}$, and a parameter $p \in [0,\infty]$,
it will be convenient to identify, for each $[f]_{\mathfrak{m}} \in L_{p}(\mathfrak{m})$, the functions $f \in [f]_{\mathfrak{m}}$.
We shall follow this path whenever our statements will depend only on the equivalence class without further mention, provided that it is clear from the context. In this
case we use the symbol $f$ instead of $[f]_{\mathfrak{m}}$. But we shall not consider functions agreeing $\mathfrak{m}$-a.e. to be identical if we are concerned with fine
properties of the single function.
\end{Remark}

Given a metric space $\operatorname{X}=(\operatorname{X},\operatorname{d})$ and a family of sets $\mathcal{G} \subset 2^{\operatorname{X}}$, by $\mathcal{M}(\mathcal{G})$ we denote its covering multiplicity, i.e.,
the minimal $M' \in \mathbb{N}_{0}$ such that every point $x \in \operatorname{X}$ belongs to at most $M'$ sets from $\mathcal{G}$.
The following proposition is elementary, we omit the proof.

\begin{Prop}
\label{Prop.covering_multiplicity}
Let $\mathfrak{m}$ be a measure on a complete separable metric space $\operatorname{X}=(\operatorname{X},\operatorname{d})$.
Let $\mathcal{G} \subset 2^{\operatorname{X}}$ be an at most countable family of sets with $\mathcal{M}(\mathcal{G}) < +\infty$. Then
\begin{equation}
\sum\limits_{G \in \mathcal{G}}\int\limits_{G}|f(x)|\,d\mathfrak{m}(x) \le \mathcal{M}(\mathcal{G})\int\limits_{\operatorname{G}}|f(x)|\,d\mathfrak{m}(x) \quad \text{for every} \quad f \in L_{1}(\mathfrak{m}\lfloor_{\operatorname{G}}),
\end{equation}
where $\operatorname{G}=\cup\{G:G \in \mathcal{G}\}$.
\end{Prop}

Given a complete separable metric space $\operatorname{X}=(\operatorname{X},\operatorname{d})$ and a measure $\mathfrak{m}$ on $\operatorname{X}$, for each $f \in L_{1}^{loc}(\mathfrak{m})$, and every Borel set $G \subset \operatorname{X}$ with $\mathfrak{m}(G) < +\infty$, we put
\begin{equation}
\label{eqq.notation_average}
f_{G,\mathfrak{m}}:=\fint\limits_{G}f(x)\,d\mathfrak{m}(x):=
\begin{cases}
\frac{1}{\mathfrak{m}(G)}\int\limits_{G}f(x)\,d\mathfrak{m}(x), \quad \mathfrak{m}(G) > 0;\\
0, \quad \mathfrak{m}(G) = 0.
\end{cases}
\end{equation}
Furthermore, we put
\begin{equation}
\label{eqq.best_approximation_constant}
\mathcal{E}_{\mathfrak{m}}(f,G):=\inf_{c \in \mathbb{R}}\fint_{G}|f(x)-c|\,d\mathfrak{m}(x).
\end{equation}

Sometimes we will use the following rough upper estimate $\mathcal{E}_{\mathfrak{m}}(f,G)$, which is an easy consequence of
Remark \ref{Rem.best_approx}, H\"older's inequality for sums and H\"older's inequality for integrals.
\begin{Prop}
\label{Prop.rough_estimate}
If $p \in [1,\infty)$, then $\Bigl(\mathcal{E}_{\mathfrak{m}}(f,G)\Bigr)^{p} \le 2^{p}\fint_{G}|f(x)|^{p}\,d\mathfrak{m}(x).$
\end{Prop}

In order to built a fruitful theory we shall work with measures satisfying some restrictions.

\begin{Def}
\label{Def.unformly_locally_doubling}
Given a complete separable metric space $\operatorname{X}=(\operatorname{X},\operatorname{d})$, we say that
a measure $\mathfrak{m}$ on $\operatorname{X}$ has a uniformly locally doubling property if for each $R > 0$,
\begin{equation}
\label{eq.doubling}
C_{\mathfrak{m}}(R):=\sup\limits_{r (0,R]}\sup\limits_{x \in \operatorname{X}}\frac{\mathfrak{m}(B_{2r}(x))}{\mathfrak{m}(B_{r}(x))} < +\infty.
\end{equation}
\end{Def}

\begin{Remark}
\label{Rem.best_approx}
Given a measure $\mathfrak{m}$ on $\operatorname{X}$ and a constant $c \in \mathbb{R}$, we have
\begin{equation}
\begin{split}
\notag
&\mathcal{E}_{\mathfrak{m}}(f,G) \le \fint\limits_{G}\bigl|f(x)-f_{G,\mathfrak{m}}|\,d\mathfrak{m}(x) \le \fint\limits_{G}\fint\limits_{G}|f(x)-f(y)|\,d\mathfrak{m}(x)\,d\mathfrak{m}(y) \le 2\mathcal{E}_{\mathfrak{m}}(f,G).
\end{split}
\end{equation}
Furthermore, if $\mathfrak{m}$ has a uniformly locally doubling property, then it follows easily from the above inequality that for each $R > 0$,  $c \geq 1$ there is $C > 0$ such that for every $(x,r) \in \operatorname{X} \times (0,R]$
\begin{equation}
|f_{B_{r}(x'),\mathfrak{m}}-f_{B_{cr}(x),\mathfrak{m}}| \le C \mathcal{E}_{\mathfrak{m}}(f,B_{cr}(x)) \quad \text{for each}  \quad B_{r}(x') \subset B_{cr}(x).
\end{equation}
\end{Remark}

Given an m.m.s.\ $\operatorname{X}=(\operatorname{X},\operatorname{d},\mu)$, it is well known that the global doubling property of the measure $\mu$ implies the globally metric doubling property of the space $(\operatorname{X},\operatorname{d})$ (see, for example, page 102 in \cite{HKST}). Similarly, we have the following result (we put $[c]:=\max\{k \in \mathbb{Z}:k \le c\}$).
\begin{Prop}
\label{Prop.metric_doubling}
Let $\operatorname{X}=(\operatorname{X},\operatorname{d},\mu)$ be a metric measure space. If $\mu$ is a uniformly locally doubling measure, then, for each $R > 0$ and $c \geq 1$,
any closed ball $B=B_{cR}(x)$ contains at most $N_{\mu}(R,c):=[(C_{\mu}((c+1)R))^{\log_{2}(2c)+1}]+1$ disjoint closed balls of radii $R$.
\end{Prop}

\begin{proof}
If $B'=B_{R}(x') \subset B$ then $B \subset B_{2cR}(x')$. Applying \eqref{eq.doubling} $[\log_{2}(2c)]+1$ times we have
\begin{equation}
\notag
\mu(B) \le (C_{\mu}((c+1)R))^{[\log_{2}(2c)]+1}\mu(B') \le (C_{\mu}((c+1)R))^{\log_{2}(2c)+1}\mu(B').
\end{equation}
If $\mathcal{B}$ is a disjoint family of closed balls with
radii $R$ lying in $B$, then $\sum\{\mu(B'):B' \in \mathcal{B}\} \le \mu(B)$. Hence,
\begin{equation}
\notag
\#\mathcal{B}\frac{\mu(B)}{(C_{\mu}((c+1)R))^{\log_{2}(2c)+1}} \le \sum\{\mu(B'):B' \in \mathcal{B}\} \le \mu(B).
\end{equation}
Consequently, we have $\#\mathcal{B} \le N_{\mu}(R,c)$.
\end{proof}

Given a number $\epsilon \in (0,1)$ and a family of closed balls $\mathcal{B} \subset 2^{\operatorname{X}}$, we put, for each $k \in \mathbb{Z}$,
\begin{equation}
\label{eq.strata}
\mathcal{B}(k,\epsilon):=\{B \in \mathcal{B}: r(B) \in (\epsilon^{k+1},\epsilon^{k}]\}.
\end{equation}

\begin{Prop}
\label{Prop.finite_intersection}
Let $\operatorname{X}=(\operatorname{X},\operatorname{d},\mu)$ be a metric measure space. If $\mu$ is a uniformly locally doubling measure, then for each $c \geq 1$, $\epsilon \in (0,1)$, and any disjoint family of closed balls $\mathcal{B}$,
\begin{equation}
\notag
\mathcal{M}(\{cB: B \in \mathcal{B}(k,\epsilon)\}) \le N_{\mu}(\epsilon^{k+1},\frac{2c}{\epsilon}) \quad \text{for every} \quad k \in \mathbb{Z},
\end{equation}
where $N_{\mu}(\epsilon^{k+1},\frac{2c}{\epsilon})$ is the same as in Proposition \ref{Prop.metric_doubling}.
\end{Prop}

\begin{proof}
We fix $c \geq 1$, $\epsilon \in (0,1)$, a disjoint family of closed balls $\mathcal{B}$, and a number $k \in \mathbb{Z}$. Consider the family
$\widetilde{\mathcal{B}}(k,\epsilon)$ consisting of all closed balls whose centers are exactly the same as in the family $\mathcal{B}(k,\epsilon)$
but with radii $\epsilon^{k+1}$.
Given a point $x \in \operatorname{X}$,
if $x \in cB$ for some $B \in \mathcal{B}(k,\epsilon)$, then $B \subset B_{2c\epsilon^{k}}(x)$.
Since the family $\widetilde{\mathcal{B}}(k,\epsilon)$ is disjoint, by Proposition \ref{Prop.metric_doubling},
\begin{equation}
\begin{split}
&\mathcal{M}(\{cB: B \in \mathcal{B}(k,\epsilon)\}) \le \sup\limits_{x \in \operatorname{X}}\sum\limits_{B \in \mathcal{B}(k,\epsilon)}\chi_{cB}(x)\\
&\le \sup\limits_{x \in \operatorname{X}}\#\{B \in \widetilde{\mathcal{B}}(k,\epsilon): B \subset B_{2c\epsilon^{k}}(x)\} \le N_{\mu}(\epsilon^{k+1},\frac{2c}{\epsilon}).
\end{split}
\end{equation}
The proof is complete.
\end{proof}

The following proposition which is an easy consequence of Proposition \ref{Prop.metric_doubling} is also well known. 
\begin{Prop}
\label{Prop.loccomp}
Let $(\operatorname{X},\operatorname{d},\mu)$ be a metric measure space with uniformly locally doubling measure $\mu$. Then every closed ball
$B=B_{r}(x)$ is a compact subset of $\operatorname{X}$.
\end{Prop}

We recall notation \eqref{eqq.centers_of_dyadic_cubes} and \eqref{eqq.centers_of_dyadic_cubes_2}.

\begin{Def}
\label{Def.partial_order}
Let $\operatorname{X}=(\operatorname{X},\operatorname{d})$ be a complete separable metric space and $\epsilon \in (0,1)$. We say that a partial order $\preceq$
on $Z(\operatorname{X},\epsilon)$ is admissible if the following properties hold:

\begin{itemize}
\item[\((\rm PO1)\)] if $z_{k,\alpha} \preceq z_{l,\beta}$ for some $k,l \in \mathbb{Z}$, then $k \geq l$;

\item[\((\rm PO2)\)] for any $l \le k$ and $z_{k,\alpha} \in Z_{k}(\operatorname{X},\epsilon)$ there is a unique $z_{l,\beta} \in Z_{l}(\operatorname{X},\epsilon)$ such that
$z_{k,\alpha} \preceq z_{l,\beta}$;

\item[\((\rm PO3)\)] if $k \in \mathbb{Z}$ and $z_{k,\alpha} \preceq z_{k-1,\beta}$, then $\operatorname{d}(z_{k,\alpha},z_{k-1,\beta}) < \epsilon^{k-1}$;

\item[\((\rm PO4)\)] if $k \in \mathbb{Z}$ and  $\operatorname{d}(z_{k,\alpha},z_{k-1,\beta}) < \frac{\epsilon^{k-1}}{2}$, then $z_{k,\alpha} \preceq z_{k-1,\beta}$.

\end{itemize}

\end{Def}

The following proposition was proved in \cite{Christ}.
\begin{Prop}
\label{Prop.partial_order}
Given a complete separable metric space $\operatorname{X}=(\operatorname{X},\operatorname{d})$ and a parameter $\epsilon \in (0,1)$,  there exists at least one
admissible partial order on the set $Z(\operatorname{X},\epsilon)$.
\end{Prop}

According to one beautiful result of Christ \cite{Christ}, given a metric measure space $\operatorname{X}=(\operatorname{X},\operatorname{d},\mu)$,
if the measure $\mu$ is globally doubling, then there exists a natural
analog of the Euclidean dyadic cubes available in $\mathbb{R}^{n}$. However, an analysis of the arguments
used in \cite{Christ} shows that in fact the uniformly locally doubling property of $\mu$
is sufficient to establish the following result.

\begin{Prop}
\label{Prop.mms_dyadic}
Let $\operatorname{X}=(\operatorname{X},\operatorname{d})$ be a complete separable metric space. Let
$\epsilon \in (0,\frac{1}{10}]$ and let $\preceq$ be an admissible partial order
on the set $Z(\operatorname{X},\epsilon)$. Given $a \in (0,\frac{1}{8}]$, for each $k \in \mathbb{Z}$, and every $\alpha \in \mathcal{A}_{k}(\operatorname{X},\epsilon)$ we define a generalized dyadic cube $Q_{k,\alpha}$ in the space $\operatorname{X}$ by the equality
\begin{equation}
\label{eqq.generalised_dyadic_cube}
Q_{k,\alpha}:=\bigcup\limits_{z_{j,\beta} \preceq z_{k,\alpha}}\operatorname{int}B_{a\epsilon^{k}}(z_{j,\beta}).
\end{equation}
Then the family $\{Q_{k,\alpha}\}:=\{Q_{k,\alpha}:k \in \mathbb{Z}, \alpha \in \mathcal{A}_{k}(\operatorname{X},\epsilon)\}$
satisfies the following properties:
\begin{itemize}
\item[\((\rm DQ1)\)] for each $k \in \mathbb{Z}$, the set $\mathcal{A}_{k}(\operatorname{X},\epsilon)$ is at most countable and $\operatorname{X} \setminus \cup_{\alpha \in \mathcal{A}_{k}(\operatorname{X},\epsilon)}\operatorname{cl}Q_{k,\alpha}=\emptyset$;

\item[\((\rm DQ2)\)]  if $j \geq k$ then either $Q_{j,\beta} \subset Q_{k,\alpha}$ or $Q_{j,\beta} \cap Q_{k,\alpha} = \emptyset$;

\item[\((\rm DQ3)\)] if $l < k$ and $\alpha \in \mathcal{A}_{k}(\operatorname{X},\epsilon)$ there is a unique $\beta \in \mathcal{A}_{l}(\operatorname{X},\epsilon)$ such that $Q_{k,\alpha} \subset Q_{l,\beta}$;

\item[\((\rm DQ4)\)] $B_{\frac{\epsilon^{k}}{4}}(z_{k,\alpha}) \subset Q_{k,\alpha} \subset B_{2\epsilon^{k}}(z_{k,\alpha})$ for each $k \in \mathbb{Z}$ and any $\alpha \in \mathcal{A}_{k}(\operatorname{X},\epsilon)$.

\end{itemize}

If, in addition, $\mu$ is a uniformly locally doubling measure on $\operatorname{X}$ with $\operatorname{supp}\mu = \operatorname{X}$, then

\begin{itemize}
\item[\((\rm DQ5)\)] $\mu(\partial Q_{k,\alpha}) = 0$ for each $k \in \mathbb{Z}$ and every $\alpha \in \mathcal{A}_{k}(\operatorname{X},\epsilon)$.
\end{itemize}
\end{Prop}

Given $\epsilon \in (0,1)$ and $r > 0$, we will use the \textit{following important notation}
\begin{equation}
\label{eqq.important_index_notation}
k(r):=k_{\epsilon}(r):=\max\{k \in \mathbb{Z}:r \le \epsilon^{k}\}.
\end{equation}

\begin{Prop}
\label{Prop.114}
Let $\operatorname{X}=(\operatorname{X},\operatorname{d},\mu)$ be a metric measure space with uniformly locally doubling measure $\mu$. Let $\epsilon \in (0,1)$, $\underline{k} \in \mathbb{Z}$, and
$\{Q_{k,\alpha}\}$ be a collection of generalized dyadic cubes.  Given $c \geq 1$, there exists a constant $C_{D}(c,\underline{k}) > 0$ depending only on $C_{\mu}((c+1+\frac{4}{\epsilon})\epsilon^{\underline{k}})$, $\epsilon$, $c$ and $\underline{k}$
such that, for each $x \in \operatorname{X}$ and any $r \in (0,\epsilon^{\underline{k}}]$,
\begin{equation}
\#\{\alpha \in \mathcal{A}_{k(r)}(\operatorname{X},\epsilon): \operatorname{cl}Q_{k(r),\alpha} \cap B_{cr}(x) \neq \emptyset\} \le C_{D}(c,\underline{k}).
\end{equation}
\end{Prop}

\begin{proof}
Note that if $\operatorname{cl}Q_{k(r),\alpha} \cap B_{cr}(x)$, then  by $(\operatorname{DQ}4)$ of Proposition \ref{Prop.mms_dyadic} we have the inclusion
$\operatorname{cl}Q_{k(r),\alpha} \subset (c+\frac{4\epsilon^{k(r)}}{r})B_{r}(x) \subset (c+\frac{4}{\epsilon})B_{r}(x)$. On the other hand,
the closed balls $\frac{1}{2}B_{\epsilon^{k(r)}}(z_{k(r),\alpha})$, $\alpha \in \mathcal{A}_{k(r)}(\operatorname{X},\epsilon)$ are disjoint. As a result, an application of Proposition \ref{Prop.metric_doubling}
proves the claim.
\end{proof}

Let $\operatorname{X}=(\operatorname{X},\operatorname{d},\mu)$ be a metric measure space. In cases when $\mu$ is a uniformly locally doubling measure, given $q \in (1,\infty)$, $\alpha \geq 0$,
we shall occasionally consider the \textit{local fractional maximal function} $M^{R}_{q,\alpha}(f)$ of a function
$f \in L_{1}^{loc}(\operatorname{X})$ defined, for any given $R > 0$, as
\begin{equation}
\label{eq.Maximal_function}
M^{R}_{q,\alpha}(f)(x):=\sup\limits_{r \in (0,R]}r^{\alpha}\Bigl(\fint\limits_{B_{r}(x)}|f(y)|^{q}\,d\mu(y)\Bigr)^{\frac{1}{q}}.
\end{equation}
It is well known that the doubling condition coupled with Vitali's 5B-covering lemma (see Section 3.3 in \cite{HKST}) gives the
following property.
\begin{Prop}
\label{Prop.maximal_function}
Let $(\operatorname{X},\operatorname{d},\mu)$ be a metric measure space with a uniformly locally doubling measure $\mu$. Let $p \in (1,\infty)$ and $q \in (1,p)$. Then
for every $R > 0$ there is a constant $C > 0$ depending  only on $p$,$q$ and $C_{\mu}(R)$ such
that
\begin{equation}
\|M^{R}_{q,0}(f)|L_{p}(\operatorname{X})\| \le C\|f|L_{p}(\operatorname{X})\| \quad \text{for all} \quad f \in L_{p}(\operatorname{X}).
\end{equation}
\end{Prop}

Recall that given $q \in [1,\infty)$, a metric measure space $\operatorname{X}=(\operatorname{X},\operatorname{d},\mu)$ is
said to \textit{support a weak local $(1,q)$-Poincar\'e
inequality} if for any $R > 0$, there are constants $C=C(R) > 0$, $\lambda=\lambda(R) \geq 1$ such that, for any Lipschitz
function $f : \operatorname{X} \to \mathbb{R}$ (we use notation \eqref{eqq.best_approximation_constant}),
\begin{equation}
\label{eq.Poincare}
\mathcal{E}_{\mu}(f,B_{r}(x)) \le C r \Bigl(\fint\limits_{B_{\lambda r}(x)}(\operatorname{lip}f(y))^{q}\,d\mu(y)\Bigr)^{\frac{1}{q}} \qquad \text{for all} \quad (x,r) \in \operatorname{X} \times (0,R].
\end{equation}

\begin{Remark}
Notice that in the literature this inequality is typically required to hold for
continuous functions and upper gradients: our formulation is equivalent to that one, see \cite{AGS2}.
\end{Remark}
%

In this paper we will always work with a special class of metric measure spaces which are commonly used in the modern Geometric Analysis.
\begin{Def}
\label{Def.admissiblespace}
Given $q \in [1,\infty)$, we say that a metric measure space $\operatorname{X}=(\operatorname{X},\operatorname{d},\mu)$ is $q$-admissible and write $\operatorname{X} \in \mathfrak{A}_{q}$
if $\mu$ is uniformly locally doubling and $\operatorname{X}$ supports a weak local $(1,q)$-Poincar\'e
inequality.
\end{Def}
We will occasionally use the following powerful result by Keith and Zhong (see Ch.12 in \cite{HKST} for a detailed proof and historical remarks).
\begin{Prop}
\label{KeithZhong}
Let $p \in (1,\infty)$ and $\operatorname{X} \in \mathfrak{A}_{p}$. Then there is $q \in [1,p)$
such that $\operatorname{X} \in \mathfrak{A}_{q}$.
\end{Prop}
Given a complete separable metric space $\operatorname{X}=(\operatorname{X},\operatorname{d})$, we say that a measure $\mathfrak{m}$ on $\operatorname{X}$ satisfies the \textit{uniformly locally reverse doubling property} if for each $R \in (0,\frac{\operatorname{diam}\operatorname{X}}{2})$
there is a constant $\operatorname{c}_{\mathfrak{m}}(R) \in (0,1)$
such that
\begin{equation}
\label{eqq.reverse_doubling}
\mathfrak{m}(B_{r}(x)) \le \operatorname{c}_{\mathfrak{m}}(R)\mathfrak{m}(B_{2r}(x)) \quad \text{for all} \quad (x,r) \in \operatorname{X}\times (0,R].
\end{equation}

Given a complete separable metric space $\operatorname{X}=(\operatorname{X},\operatorname{d})$, we say that a measure $\mathfrak{m}$ on $\operatorname{X}$ has a \textit{relative lower volume decay of order $Q > 0$} up to a scale $R > 0$ if there is a constant $C(Q,R) > 0$
such that, for any closed balls $B', B$ in $\operatorname{X}$ with $B' \subset B$ and $r(B) < R$,
\begin{equation}
\label{eqq.lower_decay_exponent}
\Bigl(\frac{r(B')}{r(B)}\Bigr)^{Q} \le C(Q,R) \frac{\mathfrak{m}(B')}{\mathfrak{m}(B)}.
\end{equation}

Now we summarize the basic properties of $q$-admissible m.m.s.\ that will be used below.

\begin{Prop}
\label{Prop.reverse_doubling}
Let $\operatorname{X} \in \mathfrak{A}_{p}$ for some $p \in [1,\infty)$. Then the following properties hold:

\begin{itemize}
\item[\((1)\)] for each $x \in \operatorname{X}$ and any $r \in (0,\infty)$ the open ball $\operatorname{int}B_{r}(x)$ is connected in the induced topology;

\item[\((2)\)] for each $R > 0$ there is $Q(R) > 0$ such that $\mu$ has the relative volume decay property of order $Q(R)$ up to the scale $R$;

\item[\((3)\)] the measure $\mu$ satisfies the uniformly locally reverse doubling property.

\end{itemize}

\end{Prop}

\begin{proof}
In fact it is sufficient to repeat almost verbatim the arguments from the proof of Proposition 8.1.6, Lemma 8.1.13 and Remark 8.1.15 in \cite{HKST} with obvious modifications
from global to local properties.
\end{proof}

Based on the above proposition we can formulate the following definition.

\begin{Def}
\label{Def.lower_decay_exponent}
Let $\operatorname{X}=(\operatorname{X},\operatorname{d},\mu)$ be a metric measure space with uniformly locally doubling measure $\mu$.
For each $R > 0$, we let $\operatorname{Q}_{\mu}(R)$ denote the set of all $Q > 0$ for each of which $\mu$ has
a relative lower volume decay of order $Q$ up to the scale $R$. \textit{The lower decay exponent} $\underline{Q}_{\mu}(R)$ of the measure $\mu$ is defined as
$\underline{Q}_{\mu}(R):=\inf\{Q:Q \in \operatorname{Q}_{\mu}(R)\}$.
\end{Def}

It is well known that in the Euclidean space $\mathbb{R}^{n}$ the $d$-Hausdorff measures provide a useful tool in measuring some $\mathcal{L}^{n}$-negligible sets.
Unfortunately, given an m.m.s.\ $(\operatorname{X},\operatorname{d},\mu)$ with locally uniformly doubling measure $\mu$ and a parameter $R > 0$,
the constant $C_{\mu}(R)$ in \eqref{eq.doubling}
can be essentially larger than $(\operatorname{c}_{\mu}(R))^{-1}$ in \eqref{eqq.reverse_doubling}.
Hence, the dependence of $\mu(B_{r}(x))$ on $r$ is not a power of $r$ in general. As a result, it is natural to consider
\textit{codimensional substitutions} for the usual Hausdorff contents and measures.
More precisely, following \cite{GS,GKS,Maly,MNS,SakSot}, given an m.m.s.\ $\operatorname{X}=(\operatorname{X},\operatorname{d},\mu)$ with locally uniformly doubling measure $\mu$ and a parameter $\theta \geq 0$, for each set $E \subset \operatorname{X}$ and any $\delta \in (0,\infty]$, we put
\begin{equation}
\label{eq.cod.content}
\mathcal{H}_{\theta,\delta}(E):=\inf\{\sum \frac{\mu(B_{r_{i}}(x_{i}))}{(r_{i})^{\theta}}:E \subset \bigcup B_{r_{i}}(x_{i}) \text{ and } r_{i} < \delta\},
\end{equation}
where the infimum is taken over all at most countable coverings of $E$ by balls $\{B_{r_{i}}(x_{i})\}$ with radii $r_{i} < \delta$.
Given $\delta > 0$, the mapping $\mathcal{H}_{\theta,\delta}:2^{\operatorname{X}} \to [0,+\infty]$ is called the \textit{codimension-$\theta$ Hausdorff content} at the scale $\delta$. We define the \textit{codimension-$\theta$ Hausdorff measure} by the equality
\begin{equation}
\label{eqq.definition_Hausdorff_measure}
\mathcal{H}_{\theta}(E):=\lim_{\delta \to 0}\mathcal{H}_{\theta,\delta}(E).
\end{equation}

\begin{Remark}
\label{Rem.Hausdorff_is_a_Borel_measure}
If there exists $R > 0$ such that $\theta \in [0,\underline{Q}_{\mu}(R))$, then it easily follows from Definition \ref{Def.lower_decay_exponent}
that $\mathcal{H}_{\theta}(\emptyset)=0$.
Hence, by Theorem 4.2 in \cite{Mat},
$\mathcal{H}_{\theta}$ is a Borel regular measure with $\operatorname{supp}\mathcal{H}_{\theta} = \operatorname{X}$. Typically, it is not locally finite, and hence, $\mathcal{H}_{\theta}$
is not a measure on $\operatorname{X}$ in accordance with the terminology adapted in this paper.
\end{Remark}

There is a special class of m.m.s.\, for which the behavior of $\mu(B_{r}(x))$ is, roughly speaking, expressed by $r^{Q}$ for some $Q > 0$. The detailed discussion of such spaces is beyond the scope
of this paper. We mention only  \cite{Ch,JJKR,SakSot} for some interesting results related to such spaces.

\begin{Def}
\label{Def.Ahlfors_regular_space}
Given $Q > 0$, we say that a metric measure space $\operatorname{X}=(\operatorname{X},\operatorname{d},\mu)$ is Ahlfors $Q$-regular if there exist constants $c_{\mu,1}$, $c_{\mu,2} > 0$ such that
\begin{equation}
\notag
c_{\mu,1}r^{Q} \le \mu(B_{r}(x)) \le c_{\mu,2}r^{Q} \quad \text{for all} \quad (x,r) \in \operatorname{X} \times [0,\operatorname{diam}\operatorname{X}).
\end{equation}
\end{Def}

\subsection{Sobolev calculus on metric measure spaces}
As mentioned in the introduction, given an m.m.s.\ $\operatorname{X}=(\operatorname{X},\operatorname{d},\mu)$ and a parameter $p \in (1,\infty)$,
there are at least five different approaches to the definition of Sobolev-type spaces on $\operatorname{X}$.
In the literature, the corresponding spaces are as follows: Korevaar--Schoen--Sobolev spaces $KS^{1}_{p}(\operatorname{X})$ \cite{KS,GT},
Hajlasz--Sobolev spaces $M^{1}_{p}(\operatorname{X})$ \cite{Haj}, Cheeger--Sobolev spaces
$Ch_{p}^{1}(\operatorname{X})$ \cite{Ch}, Newtonian--Sobolev spaces $N^{1}_{p}(\operatorname{X})$ \cite{Shan1}, and Sobolev
spaces $W^{1}_{p}(\operatorname{X})$ \cite{AGS1, AGS2, ACDM15}. The reader can also find some useful information related to these spaces
in Chapter 10 of \cite{HKST} and in the lecture notes \cite{GP20}.

\begin{Remark}
\label{Rem.Sobolev_approaches}
Given an arbitrary m.m.s.\ $\operatorname{X}=(\operatorname{X},\operatorname{d},\mu)$ and a parameter $p \in (1,\infty)$, there are canonical isometric isomorphisms
between $Ch_{p}^{1}(\operatorname{X})$, $N_{p}^{1}(\operatorname{X})$ and $W_{p}^{1}(\operatorname{X})$ \cite{AGS2}. Furthermore, if $\operatorname{X} \in \mathfrak{A}_{p}$ for some $p \in (1,\infty)$,
a combination of results from \cite{AGS2} and \cite{GT} gives the equalities
$KS_{p}^{1}(\operatorname{X})=M_{p}^{1}(\operatorname{X})=Ch_{p}^{1}(\operatorname{X})=N_{p}^{1}(\operatorname{X})=W_{p}^{1}(\operatorname{X})$ as linear spaces, the corresponding norms
being equivalent. Since in all main results of this paper we will always
assume that $\operatorname{X} \in \mathfrak{A}_{p}$, we identify different Sobolev-type spaces and we can use the symbol $W_{p}^{1}(\operatorname{X})$ for each of them.
\end{Remark}

Keeping in mind Remark \ref{Rem.Sobolev_approaches} we recall the approach of J.~Cheeger to Sobolev spaces.
\begin{Def}
\label{Def.Cheeger_Sobolev}
Given $p \in (1,\infty)$, the Sobolev space $W^{1}_{p}(\operatorname{X})$
is a linear space consisting of all $F \in L_{p}(\operatorname{X})$ with $\operatorname{Ch}_{p}(F) < +\infty$,
where $\operatorname{Ch}_{p}(F)$ is \textit{a Cheeger energy of $F$} defined by
\begin{equation}
\notag
\operatorname{Ch}_{p}(F):=\inf\{\varliminf\limits_{n \to \infty}\int\limits_{\operatorname{X}}(\operatorname{lip}F_{n})^{p}\,d\mu: \{F_{n}\} \subset \operatorname{LIP}(\operatorname{X}),\ F_{n} \to F \text{ in } L_{p}(\operatorname{X})\}.
\end{equation}
The space $W^{1}_{p}(\operatorname{X})$ is normed by $\|F|W_{p}^{1}(\operatorname{X})\|:=\|F|L_{p}(\operatorname{X})\|+(\operatorname{Ch}_{p}(F))^{\frac{1}{p}}.$
\end{Def}
\begin{Remark}
\label{Rem.minimal_upper_gradient}
It is well known that for each $p \in (1,\infty)$ and $F \in W^{1}_{p}(\operatorname{X})$ there is a well-defined nonnegative function $|DF|_{p} \in L_{p}(\operatorname{X})$, called minimal $p$-weak
upper gradient, which, if $\operatorname{X}$ is a smooth space, coincides $\mu$-a.e.\ with the modulus
of the distributional differential of $F$. Furthermore, $\operatorname{Ch}_{p}(F)=\||DF|_{p}|L_{p}(\operatorname{X})\|$ \cite{AGS2}.
However, in contrast with the classical settings the minimal $p$-weak upper gradient may depend on $p$ \cite{DiMarinoSpeight15}.
\end{Remark}

The following assertion will be important in proving some key estimates of Section 10.

\begin{Prop}
\label{Prop.Sobolev_Poincare}
Let $R > 0$, $q \in (1,\infty)$, $p \geq q$ and  $\operatorname{X} \in \mathfrak{A}_{q}$. Then for each $F \in W_{p}^{1}(\operatorname{X})$,
\begin{equation}
\label{eq.Sobolev_Poincare}
\mathcal{E}_{\mu}(F,B_{r}(x))
\le Cr\Bigl(\fint\limits_{B_{\lambda r}(x)}(|DF|_{p})^{q}\,d\mu(y)\Bigr)^{\frac{1}{q}} \quad \text{for all} \quad  (x,r) \in \operatorname{X} \times (0,R],
\end{equation}
where $C=C(R)$ and $\lambda=\lambda(R)$ are the same constants as in \eqref{eq.Poincare}.
\end{Prop}

\begin{proof}
According to the main results of \cite{AGS2}, given $F \in W_{p}^{1}(\operatorname{X})$, there is a sequence $\{F_{n}\} \subset \operatorname{LIP}(\operatorname{X})$
such that $F_{n} \to F$, $n \to \infty$ in $L_{p}(\operatorname{X})$-sense and $\operatorname{lip}F_{n} \to |DF|_{p}$, $n \to \infty$ in $L_{p}(\operatorname{X})$-sense.
Hence, taking into account that $W_{p}^{1}(\operatorname{X}) \subset W_{q}^{1,loc}(\operatorname{X})$, we use \eqref{eq.Poincare} and pass to the limit as $n \to \infty$.
This gives \eqref{eq.Sobolev_Poincare} and completes the proof.
\end{proof}

It was shown in \cite{ACDM15} that under mild assumptions on an m.m.s.\ $\operatorname{X}=(\operatorname{X},\operatorname{d},\mu)$, the
Sobolev space $W_{p}^{1}(\operatorname{X})$ is reflexive for every $p \in (1,\infty)$. In particular, we have the following result.
\begin{Prop}
\label{Prop.reflexivity}
Let $p \in (1,\infty)$ and $\operatorname{X} \in \mathfrak{A}_{p}$. Then the Sobolev space $W_{p}^{1}(\operatorname{X})$ is reflexive.
\end{Prop}
\begin{Remark}
In fact, it was assumed in \cite{ACDM15} that the metric space $(\operatorname{X},\operatorname{d})$ is globally metrically doubling. However, a careful analysis of the proof shows that
the uniformly locally doubling property of the measure $\mu$ is sufficient.
\end{Remark}
%

\subsection{Traces of Sobolev spaces}

We assume that the reader is familiar with the notion and basic properties of the so-called Sobolev $p$-capacities $C_{p}$, $p \in (1,\infty)$ (see Sections 7.2 and 9.2 in \cite{HKST} and Section 1.4 in \cite{BB} for details).
In fact, the main properties of $p$-capacities sufficient for our purposes are contained in the following proposition.

\begin{Prop}
\label{Prop.Capacity}
Let $p \in (1,\infty)$ and $\operatorname{X} \in \mathfrak{A}_{p}$.
Then the following properties hold:

\begin{itemize}
\item[\(\rm (1)\)] for each $F \in W_{p}^{1}(\operatorname{X})$, there is a set $E_{F}$ with $C_{p}(E_{F})=0$ such that
\begin{equation}
\overline{F}(x):=\varlimsup\limits_{r \to 0}\fint\limits_{B_{r}(x)}F(y)\,d\mu(y) \in \mathbb{R} \quad \text{for all} \quad x \in \operatorname{X} \setminus E_{F},
\end{equation}
and, furthermore, each $x \in \operatorname{X} \setminus E_{F}$ is a $\mu$-Lebesgue point of $\overline{F}$;

\item[\(\rm (2)\)] if $\theta \in [0,p)$, then $C_{p}(E)=0$ implies $\mathcal{H}_{\theta}(E)=0$ for any Borel set $E \subset \operatorname{X}$.
\end{itemize}
\end{Prop}

\begin{proof}
To prove $(1)$ one should repeat almost verbatim the arguments from the proof of Theorem 9.2.8 in \cite{HKST} and note that the additional requirement $Q \geq 1$
was used only in the end of the proof to establish higher order integrability.

Property (2) was proved in the recent paper \cite{GKS} (see Proposition 3.11 therein).
\end{proof}

Given a metric measure space $(\operatorname{X},\operatorname{d},\mu)$ and a parameter $p \in (1,\infty)$, a measure $\mathfrak{m}$ on $\operatorname{X}$ is
said to be \textit{absolutely continuous with respect to the $p$-capacity} if, for any Borel set $E \subset \operatorname{X}$, the equality $C_{p}(E)=0$ implies
the equality $\mathfrak{m}(E)=0$.

\begin{Def}
\label{Def.sharptrace}
Let $p \in (1,\infty)$ and $\operatorname{X} \in \mathfrak{A}_{p}$. Let $S \subset \operatorname{X}$ be a Borel set with $C_{p}(S) > 0$.
Given an element $F \in W_{p}^{1}(\operatorname{X})$, we define the $p$-sharp trace $F|_{S}$ of $F$ to the set $S$ as the class of all
Borel functions $f: S \to \mathbb{R}$ such that $f(x)=\overline{F}(x)$ everywhere on $S$ except a set of $p$-capacity zero. Furthermore,
we define the $p$-sharp trace space by
\begin{equation}
W_{p}^{1}(\operatorname{X})|_{S}:=\{f: S \to \mathbb{R}: f=F|_{S} \text{ for some } F \in W_{p}^{1}(\operatorname{X})\}
\end{equation}
and equip this space with the usual quotient space norm, i.e.,
\begin{equation}
\notag
\|f|W_{p}^{1}(\operatorname{X})|_{S}\|:=\inf\{\|F|W_{p}^{1}(\operatorname{X})\|:f=F|_{S}\}, \quad f \in W_{p}^{1}(\operatorname{X})|_{S}.
\end{equation}
\end{Def}
As we have already mentioned in the introduction, sometimes one should work with a relaxed version of the $p$-sharp trace space
of the Sobolev $W_{p}^{1}(\operatorname{X})$-space. This motivates us to introduce the following concept.
\begin{Def}
\label{Def.m_trace_space}
Let $p \in (1,\infty)$ and $\operatorname{X} \in \mathfrak{A}_{p}$. Let $S \subset \operatorname{X}$ be a Borel set with $C_{p}(S) > 0$.
Let $\mathfrak{m}$ be a nonzero measure on $\operatorname{X}$
which is absolutely continuous with respect to $C_{p}$ such that $S \subset \operatorname{supp}\mathfrak{m}$.
Given an element $F \in W_{p}^{1}(\operatorname{X})$, we define the $\mathfrak{m}$-trace $F|^{\mathfrak{m}}_{S}$ of $F$ to the set $S$ as
the $\mathfrak{m}$-equivalence class of the $p$-sharp trace, i.e., $F|^{\mathfrak{m}}_{S}:=[F|_{S}]_{\mathfrak{m}}$.
Furthermore, we define the $\mathfrak{m}$-trace space as the image of the $p$-sharp trace space in $L_{0}(\mathfrak{m})$ under
the mapping $\operatorname{I}_{\mathfrak{m}}$ equipped with the corresponding quotient space norm, i.e., we put $W_{p}^{1}(\operatorname{X})|^{\mathfrak{m}}_{S}:=\operatorname{I}_{\mathfrak{m}}(W_{p}^{1}(\operatorname{X})|_{S})$ and
\begin{equation}
\|f|W_{p}^{1}(\operatorname{X})|^{\mathfrak{m}}_{S}\|:=\inf\{\|F|W_{p}^{1}(\operatorname{X})\|:f=F|^{\mathfrak{m}}_{S}\}, \quad f \in W_{p}^{1}(\operatorname{X})|^{\mathfrak{m}}_{S}.
\end{equation}
\end{Def}
Having at our disposal different notions of trace spaces, it is natural to define the corresponding trace and extension operators.

\begin{Def}
\label{Def.m_trace_operator}
Let $p \in (1,\infty)$ and $\operatorname{X} \in \mathfrak{A}_{p}$. Let $S \subset \operatorname{X}$ be a Borel set with $C_{p}(S) > 0$.
Let $\mathfrak{m}$ be a nonzero measure on $\operatorname{X}$
which is absolutely continuous with respect to $C_{p}$ such that $S \subset \operatorname{supp}\mathfrak{m}$.
We define the $p$-sharp trace operator by the formula
\begin{equation}
\operatorname{Tr}|_{S}(F)=F|_{S}, \quad F \in W_{p}^{1}(\operatorname{X}).
\end{equation}
Furthermore, we define the $\mathfrak{m}$-trace operator by the equality $\operatorname{Tr}|^{\mathfrak{m}}_{S}:= \operatorname{I}_{\mathfrak{m}} \circ \operatorname{Tr}|_{S}.$

Finally, we say that a mapping $\operatorname{Ext}_{S,p}:W_{p}^{1}(\operatorname{X})|_{S} \to W_{p}^{1}(\operatorname{X})$ is a $p$-sharp extension operator if it
is the right inverse of $\operatorname{Tr}|_{S}$, and we say that a mapping $\operatorname{Ext}_{S,\mathfrak{m}}:W_{p}^{1}(\operatorname{X})|^{\mathfrak{m}}_{S} \to W_{p}^{1}(\operatorname{X})$ is an $\mathfrak{m}$-extension operator if it
is the right inverse of $\operatorname{Tr}|^{\mathfrak{m}}_{S}$.
\end{Def}

\begin{Remark}
In view of Definitions \ref{Def.sharptrace}, \ref{Def.m_trace_space} and \ref{Def.m_trace_operator}, it is clear that the $p$-sharp trace operator and the
$\mathfrak{m}$-trace operator are linear and bounded.
\end{Remark}
%
%
%
%
%
%
%

%
\section{Relaxation of the doubling property}

Given a metric measure space $\operatorname{X}=(\operatorname{X},\operatorname{d},\mu)$, in forthcoming sections we will frequently work with measures $\mathfrak{m}$ on $\operatorname{X}$
that fail to satisfy the uniformly locally doubling property \eqref{eq.doubling}.
However, in some cases it will be sufficient to have a some sort of uniformly doubling property along only a sequence
of balls. This motivates us to introduce the following concept.
\begin{Def}
\label{Def.doubling.sequence}
Given a measure $\mathfrak{m}$ on a metric measure space $\operatorname{X}=(\operatorname{X},\operatorname{d},\mu)$, we say that $\mathfrak{m}$
satisfies the uniformly weak asymptotically doubling property if, for each constant $c > 0$,
\begin{equation}
\label{eq.doubling.sequence}
\underline{\operatorname{C}}_{\mathfrak{m}}(c):=\lim\limits_{R \to +0}\sup\limits_{x \in \operatorname{supp}\mathfrak{m}}\inf\limits_{r \in (0,R]}\frac{\mathfrak{m}(B_{cr}(x))}{\mathfrak{m}(B_{r}(x))} < +\infty.
\end{equation}
\end{Def}

\begin{Remark}
The word ``weak'' in the above definition was used because in \cite{HKST} one can find a notion of uniformly asymptotically doubling property, which means that
\begin{equation}
\notag
\overline{\operatorname{C}}_{\mathfrak{m}}(c):=\lim\limits_{R \to+ 0}\sup\limits_{x \in \operatorname{supp}\mathfrak{m}}\sup\limits_{r \in (0,R]}\frac{\mathfrak{m}(B_{cr}(x))}{\mathfrak{m}(B_{r}(x))} < +\infty.
\end{equation}
\end{Remark}

Typically, in the present paper we will deal with measures that can not degenerate too fast.
\begin{Def}
\label{Def.weaknoncol}
We say that a measure $\mathfrak{m}$ on a metric measure space $\operatorname{X}=(\operatorname{X},\operatorname{d},\mu)$ is weakly noncollapsed if
\begin{equation}
\label{eq.weaknoncol}
C^{\mathfrak{m}}_{\mu}:=\inf\limits_{x \in \operatorname{supp}\mathfrak{m}}\varliminf\limits_{r \to 0}\frac{\mathfrak{m}(B_{r}(x))}{\mu(B_{r}(x))}>0.
\end{equation}
\end{Def}

It is well-known that, given a measure $\mathfrak{m}$ on the Euclidean space $(\mathbb{R}^{n},\|\cdot\|_{2})$, there are
a lot of ``doubling balls''. This fact was mentioned in \cite{Tolsa} without a proof.
We are grateful to D.~M.~Stolyarov who kindly shared with us the key idea of that proof.
Using a similar idea, we establish the following simple result, which will be quite important in what follows.
We recall property (2) of Proposition \ref{Prop.reverse_doubling} and Definition \ref{Def.lower_decay_exponent}.

\begin{Lm}
\label{Lm.doublingsequence}
Given a metric measure space $\operatorname{X}=(\operatorname{X},\operatorname{d},\mu)$, assume that the measure $\mu$
is uniformly locally doubling. If a measure $\mathfrak{m}$ on $\operatorname{X}$ is weakly noncollapsed, then
it satisfies the uniformly weak asymptotically doubling property. Furthermore, for each $c > 1$ and $Q \in \operatorname{Q}_{\mu}(1)$,
\begin{equation}
\underline{\operatorname{C}}_{\mathfrak{m}}(c) \le 2^{([c]+1)Q}.
\end{equation}
\end{Lm}
\begin{proof}
We fix $c > 1$ and assume for the contrary that $\underline{\operatorname{C}}_{\mathfrak{m}}(c) > 2^{([c]+1)Q}.$
We fix arbitrary numbers $\underline{k} \geq [c]+1$ and $M \in (2^{\underline{k}Q},\underline{\operatorname{C}}_{\mathfrak{m}}(c)).$
It is clear that there exist a point $\underline{x} \in \operatorname{supp}\mathfrak{m}$ and
a number $\overline{r}=\overline{r}(M,c) \in (0,1)$ such that
\begin{equation}
\label{eqq.1.6}
\frac{\mathfrak{m}(B_{cr}(\underline{x}))}{\mathfrak{m}(B_{r}(\underline{x}))} > M \quad \text{for all} \quad r \in (0,\overline{r}].
\end{equation}
By Definition \ref{Def.weaknoncol} we clearly have
\begin{equation}
\label{eqq.1.7}
\varlimsup\limits_{r \to 0}\frac{\mu(B_{r}(\underline{x}))}{\mathfrak{m}(B_{r}(\underline{x}))} \le \frac{1}{C^{\mathfrak{m}}_{\mu}}.
\end{equation}
Hence, combining \eqref{eqq.lower_decay_exponent} with \eqref{eqq.1.6} and \eqref{eqq.1.7} we get, for all large enough $i \in \mathbb{N}$,
\begin{equation}
\begin{split}
\notag
&2^{-Q\underline{k}i} \le C(Q,1) \frac{\mu(B_{\frac{\overline{r}}{2^{i\underline{k}}}}(\underline{x}))}{\mu(B_{\overline{r}}(\underline{x}))}\\
&= C(Q,1) \frac{\mu(B_{\frac{\overline{r}}{2^{i\underline{k}}}}(\underline{x}))}{\mathfrak{m}(B_{\frac{\overline{r}}{2^{i\underline{k}}}}(\underline{x}))}\frac{\mathfrak{m}(B_{\overline{r}}(\underline{x}))}{\mu(B_{\overline{r}}(\underline{x}))}
\frac{\mathfrak{m}(B_{\frac{\overline{r}}{2^{i\underline{k}}}}(\underline{x}))}{\mathfrak{m}(B_{\overline{r}}(\underline{x}))} \le \frac{2 C(Q,1)}{C^{\mathfrak{m}}_{\mu}} \frac{\mathfrak{m}(B_{\overline{r}}(\underline{x}))}{\mu(B_{\overline{r}}(\underline{x}))} \Bigl(\frac{1}{M}\Bigr)^{i}.
\end{split}
\end{equation}
Letting $i \to \infty$ we get a contradiction with the choice of $M$.
\end{proof}

Let $\operatorname{X}=(\operatorname{X},\operatorname{d},\mu)$ be a metric measure space. We recall notation \eqref{eqq.important_index_notation}.
Given a sequence of measures $\{\mathfrak{m}_{k}\}:=\{\mathfrak{m}_{k}\}_{k=0}^{\infty}$ on $\operatorname{X}$, a parameter $\epsilon \in (0,1)$, and a Borel set $E \subset \bigcap_{k=0}^{\infty} \operatorname{supp}\mathfrak{m}_{k}$,
for each $x \in E$ we introduce the \textit{lower} and the \textit{upper} $(\{\mathfrak{m}_{k}\},\epsilon)$-densities of $E$ at $x$ by letting
\begin{equation}
\label{eq.lower_and_upper_densities}
\underline{D}^{\{\mathfrak{m}_{k}\}}_{E}(x,\epsilon):=\varliminf\limits_{r \to 0}\frac{\mathfrak{m}_{k_{\epsilon}(r)}(B_{r}(x)\cap E)}{\mathfrak{m}_{k_{\epsilon}(r)}(B_{r}(x))}, \quad
\overline{D}^{\{\mathfrak{m}_{k}\}}_{E}(x,\epsilon):=\varlimsup\limits_{r \to 0}\frac{\mathfrak{m}_{k_{\epsilon}(r)}(B_{r}(x)\cap E)}{\mathfrak{m}_{k_{\epsilon}(r)}(B_{r}(x))}.
\end{equation}
We say that $x \in E$ is an \textit{$(\{\mathfrak{m}_{k}\},\epsilon)$-density point of $E$} if $\underline{D}^{\{\mathfrak{m}_{k}\}}_{E}(x,\epsilon)=\overline{D}^{\{\mathfrak{m}_{k}\}}_{E}(x,\epsilon)=1$.
It is clear that if there is a measure $\mathfrak{m}$ on $\operatorname{X}$ such that $\mathfrak{m}_{k}=\mathfrak{m}$ for all $k \in \mathbb{N}_{0}$,
then we get the standard lower and upper $\mathfrak{m}$-densities of $E$ at $x$, which will be denoted by $\underline{D}^{\mathfrak{m}}_{E}(x)$
and $\overline{D}^{\mathfrak{m}}_{E}(x)$, respectively (in this case the parameter $\epsilon$ is irrelevant and we omit it from our notation).

It is well known that if $\mathfrak{m}$ is a locally uniformly doubling measure on $\operatorname{X}$, then $\mathfrak{m}$-almost every point $x \in E$
is an $\mathfrak{m}$-density point of $E$. Unfortunately, this is
not the case if $\mathfrak{m}$ fails to satisfy the locally uniformly doubling property. However, the following result holds.
\begin{Lm}
\label{prop.density}
Let $\operatorname{X}=(\operatorname{X},\operatorname{d},\mu)$ be a metric measure space with uniformly locally doubling measure $\mu$.
Let $\mathfrak{m}$ be a weakly noncollapsed measure on $\operatorname{X}$.
Then, given a Borel set $E \subset \operatorname{X}$ and a parameter $c \geq 1$, for $\mathfrak{m}$-almost every point $x \in E$, there is a decreasing to zero sequence $\{r_{l}(x)\}$
such that
\begin{equation}
\label{doubling.sequence}
\varlimsup\limits_{l \to \infty} \frac{\mathfrak{m}(B_{\max\{c,5\}r_{l}(x)}(x))}{\mathfrak{m}(B_{r_{l}(x)}(x))} \le N \quad \text{and} \quad \varliminf\limits_{l \to \infty} \frac{\mathfrak{m}(B_{r_{l}(x)}(x) \cap E)}{\mathfrak{m}(B_{r_{l}(x)}(x))} \geq \frac{1}{2N},
\end{equation}
where $N=\underline{\operatorname{C}}_{\mathfrak{m}}(5\max\{c,5\})$. In particular, $\overline{D}^{\mathfrak{m}}_{E}(x) > 0$ for $\mathfrak{m}$-a.e. $x \in E$.
\end{Lm}

\begin{proof}
By Lemma \ref{Lm.doublingsequence}, for every point $x \in E$, there exists
a sequence $r_{l}(x) \downarrow 0$ satisfying
\begin{equation}
\frac{\mathfrak{m}(B_{\max\{c,5\} r_{l}(x)}(x))}{\mathfrak{m}(B_{\frac{r_{l}(x)}{5}}(x))} \le \underline{\operatorname{C}}_{\mathfrak{m}}(5\max\{c,5\})=N \quad \text{for all} \quad l \in \mathbb{N}.
\end{equation}
Given $n \in \mathbb{N}$, we consider the set
\begin{equation}
G_{n}:=\Bigl\{x \in E: \varliminf\limits_{l \to \infty} \frac{\mathfrak{m}(B_{r_{l}(x)}(x) \cap E)}{\mathfrak{m}(B_{r_{l}(x)}(x))} < \frac{1}{n}\Bigr\}.
\end{equation}
We show that $\mathfrak{m}(G_{n})=0$ for all $n \in \mathbb{N} \cap (2N,+\infty)$.
Without loss of generality we may assume that all sets $G_{n}$, $n \in \mathbb{N}$ are bounded.
Hence, in the rest of the proof we may assume that $\mathfrak{m}(G_{n}) < +\infty$ for all $n \in \mathbb{N}$.
Applying the $5B$-covering lemma (see p. 60 in \cite{HKST} for details) and taking into account the Borel regularity of the measure $\mathfrak{m}$, we get, for each $n \in \mathbb{N}$,
a family of closed balls $\mathcal{B}_{n}=\{B_{r_{l_{i}}(x_{i})}(x_{i})\}$ such that:

\begin{itemize}
\item[\((1)\)] the family $\widetilde{\mathcal{B}}_{n}:=\{\frac{1}{5}B:B \in \mathcal{B}_{n}\}$ is disjoint;

\item[\((2)\)] $G_{n} \subset \bigcup \{B:B \in \mathcal{B}_{n}\} \subset U_{\varepsilon_{n}}(G_{n})$ for some $\varepsilon_{n} > 0$;

\item[\((3)\)] $|\mathfrak{m}(U_{\varepsilon_{n}}(G_{n}))-\mathfrak{m}(G_{n})| < \frac{1}{2n}$;

\item[\((4)\)] $\mathfrak{m}(B) \le \frac{3}{2}N \mathfrak{m}(\frac{1}{5}B)$ for all  $B \in \mathcal{B}_{n}$;

\item[\((5)\)]  $\mathfrak{m}(B \cap E) < \frac{1}{n}\mathfrak{m}(B)$ for all  $B \in \mathcal{B}_{n}$.


\end{itemize}
We fix an arbitrary $n > 2N$ and assume that $\mathfrak{m}(G_{n}) > 0$ (note that if $G_{n}$ is not $\mathfrak{m}$-measurable, we consider $\mathfrak{m}$ as an outer measure).
Hence, taking $\varepsilon > 0$ small enough, we deduce from the above properties (1)--(5)
\begin{equation}
\notag
\begin{split}
&\mathfrak{m}(G_{n}) \le \sum\{\mathfrak{m}(B \cap G_{n}):B \in \mathcal{B}_{n}\} \le \sum\{\mathfrak{m}(B \cap E):B \in \mathcal{B}_{n}\}\\
&\le \frac{3N}{2n}\sum\{\mathfrak{m}(\frac{1}{5}B):B \in \mathcal{B}_{n}\}
\le  \frac{3N}{2n}\mathfrak{m}(U_{\varepsilon_{n}}(G_{n})) \le  \frac{2N}{n}\mathfrak{m}(G_{n}).
\end{split}
\end{equation}
This contradicts the assumption $\mathfrak{m}(G_{n}) > 0$.

As a result, we get $\mathfrak{m}(G_{n})=0$ for every $n > 2N$ and complete the proof.

\end{proof}

Now we introduce a new concept, which can be looked upon as a natural generalization of the notion of a Lebesgue point of locally integrable functions.
This concept will be extremely useful in the analysis of traces of Sobolev functions. We recall notation \eqref{eqq.important_index_notation}.
\begin{Def}
\label{Def.multiweight_Lebesgue}
Let $\operatorname{X}=(\operatorname{X},\operatorname{d},\mu)$ be a metric measure space. Let $\{\mathfrak{m}_{k}\}=\{\mathfrak{m}_{k}\}_{k=0}^{\infty}$ be a sequence of measures on $\operatorname{X}$ and $\epsilon \in (0,1)$. Given
$f \in L_{1}^{loc}(\{\mathfrak{m}_{k}\})$, we say that $\underline{x} \in \cap_{k=0}^{\infty}\operatorname{supp}\mathfrak{m}_{k}$ is an $(\{\mathfrak{m}_{k}\},\epsilon)$-Lebesgue point of $f$ if
\begin{equation}
\lim\limits_{r \to 0}\fint\limits_{B_{r}(\underline{x})}|f(\underline{x})-f(y)|\,d\mathfrak{m}_{k_{\epsilon}(r)}(y) = 0.
\end{equation}

By $\mathfrak{R}_{\{\mathfrak{m}_{k}\},\epsilon}(f)$ we denote the set of all $(\{\mathfrak{m}_{k}\},\epsilon)$-Lebesgue points of $f$.
\end{Def}

If there is a measure $\mathfrak{m}$ on $\operatorname{X}$  such that $\mathfrak{m}_{k}=\mathfrak{m}$ for all $k \in \mathbb{N}_{0}$, an $(\{\mathfrak{m}_{k}\},\epsilon)$-Lebesgue point of $f$
will be called an \textit{$\mathfrak{m}$-Lebesgue point of $f$} (in this case, the parameter $\epsilon$ is irrelevant, and we omit it from the notation).

\section{Lower codimension-$\theta$ content regular sets}

\textit{Throughout this section, we fix a metric measure space $\operatorname{X}=(\operatorname{X},\operatorname{d},\mu)$ with
uniformly locally doubling measure $\mu$.} We also recall that all the balls are assumed to be closed.

The following concept was actively used in \cite{GS, Maly, MNS}, where problems similar to Problem \ref{MeasureTraceProblem} were studied. We recall \eqref{eqq.definition_Hausdorff_measure}.

\begin{Def}
\label{Def.Ahlfors_David_regular}
Given $\theta \geq 0$, a set $S \subset \operatorname{X}$ is said to be codimension-$\theta$ Ahlfors--David regular if there exist constants $c_{\theta,1}(S), c_{\theta,2}(S) > 0$ such that
\begin{equation}
\label{eqq.32''}
c_{\theta,1}(S) \frac{\mu(B_{r}(x))}{r^{\theta}} \le \mathcal{H}_{\theta}(B_{r}(x) \cap S) \le  c_{\theta,2}(S) \frac{\mu(B_{r}(x))}{r^{\theta}} \quad \text{for all} \quad (x,r) \in S \times (0,1].
\end{equation}
The class of all codimension-$\theta$ Ahlfors--David regular sets is denoted by $\mathcal{ADR}_{\theta}(\operatorname{X})$.
\end{Def}

The following proposition shows that the scale $1$ in \eqref{eqq.32''} is not crucial. The proof is quite simple and follows easily from Proposition \ref{Prop.metric_doubling}.
The details are left to the reader.

\begin{Prop}
\label{Prop.Ahlfors_David_different_scales}
Let $\theta \geq 0$ and $S \in \mathcal{ADR}_{\theta}(\operatorname{X})$. Then, for every $R \geq 1$, there exist constants $c_{\theta,1}(S,R) > 0$ and $c_{\theta,2}(S,R) > 0$
such that
\begin{equation}
\label{eqq.32'''}
c_{\theta,1}(S,R) \frac{\mu(B_{r}(x))}{r^{\theta}} \le \mathcal{H}_{\theta}(B_{r}(x) \cap S) \le  c_{\theta,2}(S,R) \frac{\mu(B_{r}(x))}{r^{\theta}} \quad \text{for all} \quad (x,r) \in S \times (0,R].
\end{equation}
\end{Prop}

\begin{Remark}
In the case $\theta=0$, sets $S \in \mathcal{ADR}_{0}(\operatorname{X})$ were called regular sets in \cite{Shv1}.
\end{Remark}

Now we introduce a natural generalization of the class $\mathcal{ADR}_{\theta}(\operatorname{X})$. We recall \eqref{eq.cod.content}.

\begin{Def}
\label{Def.content}
Given $\theta \geq 0$, we say that a set $S \subset \operatorname{X}$ is lower codimension-$\theta$ content regular if
there exists a constant $\lambda_{\theta}(S) > 0$ such that
\begin{equation}
\label{eq.thick}
\lambda_{\theta}(S) \frac{\mu(B_{r}(x))}{r^{\theta}} \le \mathcal{H}_{\theta,r}(B_{r}(x)\cap S) \le \frac{\mu(B_{r}(x))}{r^{\theta}} \quad \text{for all} \quad (x,r) \in S \times (0,1].
\end{equation}
The class of all codimension-$\theta$ lower content regular sets will be denoted by $\mathcal{LCR}_{\theta}(\operatorname{X})$.
\end{Def}

\begin{Remark}
It is clear that the right-hand side of \eqref{eq.thick} holds automatically. The essence of the matter is contained in the left-hand side of \eqref{eq.thick}.
This justifies the term ``lower content regular''. On the other hand, we use the two-sided estimate in Definition \ref{Def.content} to make a clear comparison
of the classes $\mathcal{ADR}_{\theta}(\operatorname{X})$ and $\mathcal{LCR}_{\theta}(\operatorname{X})$, respectively.
\end{Remark}

\begin{Remark}
Let $n \in \mathbb{N}$ and $\operatorname{X}=(\mathbb{R}^{n},\|\cdot\|,\mathcal{L}^{n})$. It is easy to see that, given $\theta \in [0,n]$, a set $S$ lies in $\mathcal{LCR}_{\theta}(\operatorname{X})$
if and only if
\begin{equation}
\label{eq.thick'}
\mathcal{H}_{\theta,\infty}(B_{r}(x)\cap S) \geq  \lambda_{\theta}(S) \frac{\mathcal{L}^{n}(B_{r}(x))}{r^{\theta}} \quad \text{for all} \quad r \in (0,1].
\end{equation}
In other words, in the classical Euclidean settings one can replace $\mathcal{H}_{\theta,r}$ by $\mathcal{H}_{\theta,\infty}$.
Hence, a set $S$ lies in $\mathcal{LCR}_{\theta}(\mathbb{R}^{n})$ if and only if it is an $(n-\theta)$-thick set in the sense of Rychkov \cite{Ry}.
Furthermore, $d$-thick sets, $d \in [0,n]$, were actively studied in \cite{AzSh, AzVil}, where they were called $d$-lower content regular.
The choice of the set function $\mathcal{H}_{\theta,r}$ in Definition \ref{Def.content} is motivated by a possible gap between the constants $C_{\mu}(1)$ and $(\operatorname{c}_{\mu}(1))^{-1}$ in \eqref{eq.doubling}
and \eqref{eqq.reverse_doubling}, respectively.
\end{Remark}

The following lemma was proved in \cite{TV1} in the particular case $\operatorname{X}=(\mathbb{R}^{n},\|\cdot\|_{2},\mathcal{L}^{n})$.
The proof in the general case is similar. We present the details for completeness of our exposition.

\begin{Lm}
\label{Lm.Ahlfors_implies_thick}
Given $\theta \geq 0$, $\mathcal{ADR}_{\theta}(\operatorname{X}) \subset \mathcal{LCR}_{\theta}(\operatorname{X})$.
\end{Lm}

\begin{proof}
Fix $\theta \geq 0$ and $S \in \mathcal{ADR}_{\theta}(\operatorname{X})$. Assume that $S \neq \emptyset$. Consider an arbitrary closed ball $B_{r}(x)$ with $x \in S$ and $r \in (0,1]$.
Let $\mathcal{B}$ be a countable family of closed balls such that $B_{r}(x) \cap S \subset \cup\{B:B \in \mathcal{B}\}$, $r(B) < r$ for all $B \in \mathcal{B}$, and
\begin{equation}
\label{eqq.34'}
\sum\{\frac{\mu(B)}{(r(B))^{\theta}}:B \in \mathcal{B}\} \le 2\mathcal{H}_{\theta,r}(B_{r}(x) \cap S).
\end{equation}
Without loss of generality we may assume that, for each ball $B \in \mathcal{B}$, $B \cap S \neq \emptyset$.
For each $B \in \mathcal{B}$ we chose an arbitrary point $x_{B} \in B \cap S$ and consider a ball $\widetilde{B}$ centered
at $x_{B}$ with radius $2r(B)$. Clearly, $B \cap S \subset \widetilde{B} \cap S$ and $\widetilde{B} \subset 4B$ for all $B \in \mathcal{B}$.
Hence, using \eqref{eqq.34'}, the subadditivity property of $\mathcal{H}_{\theta}$ and Proposition \ref{Prop.Ahlfors_David_different_scales},
we obtain the required estimate
\begin{equation}
\notag
\begin{split}
&2\mathcal{H}_{\theta,r}(B_{r}(x) \cap S) \geq \frac{1}{(C_{\mu}(2))^{2}}\sum\{\frac{\mu(\widetilde{B})}{(r(\widetilde{B}))^{\theta}}:B \in \mathcal{B}\}\\
&\geq \frac{1}{c_{\theta,2}(S,2)(C_{\mu}(2))^{2}}\sum\{\mathcal{H}_{\theta}(\widetilde{B} \cap S):B \in \mathcal{B}\} \\
&\geq \frac{1}{c_{\theta,2}(S,2)(C_{\mu}(2))^{2}}\mathcal{H}_{\theta}(B_{r}(x) \cap S) \geq \frac{c_{\theta,1}(S,1)}{c_{\theta,2}(S,2)(C_{\mu}(2))^{2}}\frac{\mu(B_{r}(x))}{r^{\theta}}.
\end{split}
\end{equation}
The proof is complete.
\end{proof}

The following example demonstrates that if $\operatorname{X}$ is regular enough, then the classes $\mathcal{LCR}_{\theta}(\operatorname{X})$, $\theta \geq 0$
are broad enough. We recall Definition \ref{Def.Ahlfors_regular_space}.

\textbf{Example 4.7.}
Assume that the space $\operatorname{X}$ is Ahlfors $Q$-regular for some $Q > 0$.
We fix $\theta \in [\max\{0,Q-1\},Q)$ and show that any path-connected set $S \subset \operatorname{X}$ belongs to the class $\mathcal{LCR}_{\theta}(\operatorname{X})$.
Indeed, we fix $x \in S$, $r \in (0,1]$ and consider two cases.
In the fist case, $S \subset B_{r}(x)$. Hence,
\begin{equation}
\label{eqq.38''}
\mathcal{H}_{\theta,r}(B_{r}(x) \cap S) = \mathcal{H}_{\theta,r}(S) \geq \mathcal{H}_{\theta,1}(S)r^{Q-\theta} \geq \frac{\mathcal{H}_{\theta,1}(S)}{c_{\mu,2}}\frac{\mu(B_{r}(x))}{r^{\theta}}.
\end{equation}
In the second case, there is a point $y \in S \setminus B_{r}(x)$. Hence, there exists a curve $\gamma_{x,y}$ joining $x$ and $y$.
Let $\mathcal{B}$ be an at most countable family of closed balls such that $B_{r}(x) \cap S \subset \cup\{B:B \in \mathcal{B}\}$ and $r(B) < r$ for all $B \in \mathcal{B}$.
Consider the family $\underline{\mathcal{B}}:=\{\operatorname{int}(2B):B \in \mathcal{B}\}$.
By Proposition \ref{Prop.loccomp}, the set $\gamma_{x,y} \cap B_{r}(x)$ is compact.  Hence, since the set $\gamma_{x,y}$ is connected, there is a finite family
$\{B_{i}\}_{i=1}^{N} \subset \mathcal{B}$, $N \in \mathbb{N}$ such that $\gamma_{x,y} \subset \cup_{i=1}^{N}\operatorname{int}2B_{i}$ and
$\gamma_{x,y}\cap \operatorname{int}2B_{i}\cap \operatorname{int}2B_{i+1} \neq \emptyset$ for all $i \in \{1,...,N-1\}$.
Hence, by the triangle inequality, we deduce the crucial estimate
\begin{equation}
\notag
r \le \operatorname{diam}(\gamma_{x,y}) \le \sum_{i=1}^{N}\operatorname{diam}(\gamma_{x,y}\cap\operatorname{int}2B_{i}) \le \sum_{i=1}^{N}\operatorname{diam}(\operatorname{int}2B_{i})\le \sum_{i=1}^{N}4r(B_{i}).
\end{equation}
As a result, since $Q-\theta \in (0,1]$, we get
\begin{equation}
\label{eqq.39''}
\sum\limits_{B \in \mathcal{B}}\frac{\mu(B)}{(r(B))^{\theta}} \geq c_{\mu,1}\sum\limits_{B \in \underline{\mathcal{B}}}(r(B))^{Q-\theta} \geq
c_{\mu,1}\Bigl(\sum\limits_{i=1}^{N}r(B_{i})\Bigr)^{Q-\theta} \geq c_{\mu,1}\Bigl(\frac{r}{4}\Bigr)^{Q-\theta}.
\end{equation}
Taking the infimum in \eqref{eqq.39''} over all families $\mathcal{B}$ we obtain
\begin{equation}
\label{eqq.310''}
\mathcal{H}_{\theta,r}(B_{r}(x) \cap S) \geq c_{\mu,1}\Bigl(\frac{r}{4}\Bigr)^{Q-\theta} \geq \frac{c_{\mu,1}}{4^{Q-\theta}c_{\mu,2}}\frac{\mu(B_{r}(x))}{r^{\theta}}.
\end{equation}
Finally, a combination of \eqref{eqq.38''} and \eqref{eqq.310''} proves the required claim.
\hfill$\Box$

\begin{Remark}
Even in the case $\operatorname{X}=(\mathbb{R}^{2},\|\cdot\|,\mathcal{L}^{2})$,
it is clear that generic path-connected sets $S \subset \operatorname{X}$ may fail to satisfy the Ahlfors--David codimension-$1$
regularity condition.
Furthermore, it was shown in \cite{T2} that the corresponding examples can be obtained as graphs of locally Lipschitz functions.
Coupled with Lemma \ref{Lm.Ahlfors_implies_thick}, this shows that, given $\theta > 0$, the family $\mathcal{ADR}_{\theta}(\operatorname{X})$ can be a very poor subfamily of $\mathcal{LCR}_{\theta}(\operatorname{X})$
in general.
\end{Remark}

In the sequel, we will use the following result from \cite{GKS} (see Lemma 3.10 and the discussion after the lemma).

\begin{Prop}
\label{Prop.Evans}
Let $f \in L_{1}^{loc}(\operatorname{X})$, suppose $t > 0$ and define
\begin{equation}
\notag
\Lambda_{t}:=\{x \in \operatorname{X}: \varlimsup_{r \to 0}r^{t}\fint\limits_{B_{r}(x)}|f(y)|\,d\mu(y) > 0\}.
\end{equation}
Then $\mathcal{H}_{t}(\Lambda_{t})=0$.
\end{Prop}
\textbf{Example 4.9.} Assume that $\operatorname{X}$ is Ahlfors $Q$-regular for some $Q > 0$. Then each nonempty set $S \subset \operatorname{X}$
belongs to $\mathcal{LCR}_{\theta}(\operatorname{X})$ for every $\theta \geq Q$. Indeed, we fix $x \in S$, $r \in (0,1]$,
and an at most countable family of balls $\{B_{i}\}$ with radii $r_{i} < r$ that cover $B_{r}(x) \cap S$. Hence,
\begin{equation}
\notag
\sum \frac{\mu(B_{i})}{(r_{i})^{\theta}} \geq c_{\mu,1} \sum (r_{i})^{Q-\theta} \geq c_{\mu,1} r^{Q-\theta} \geq \frac{c_{\mu,1}}{c_{\mu,2}}\frac{\mu(B_{r}(x))}{r^{\theta}}.
\end{equation}
Taking the infimum  over all such coverings, we get the required claim.

\section{$\theta$-regular sequences of measures}
\textit{Throughout this section we fix $p \in (1,\infty)$ and an m.m.s.\ $\operatorname{X}=(\operatorname{X},\operatorname{d},\mu) \in \mathfrak{A}_{p}$.}
We recall the definitions of the classes $\mathfrak{M}_{\theta}(\operatorname{X})$ and $\mathfrak{M}^{str}_{\theta}(\operatorname{X})$ given in the introduction.
We also recall notation $B_{r}(x)$ given in Section 2.
It is clear that there are smallest constants for which conditions $(\textbf{M}2)$ and $(\textbf{M}4)$
hold. We denote them by $C_{\{\mathfrak{m}_{k}\},1}$ and $C_{\{\mathfrak{m}_{k}\},3}$, respectively. Similarly, there is a largest
constant for which $(\textbf{M}3)$ holds. We denote it by $C_{\{\mathfrak{m}_{k}\},2}$. We will use the symbol
$\mathcal{C}_{\{\mathfrak{m}_{k}\}}$ to denote the family of these constants, i.e., $\mathcal{C}_{\{\mathfrak{m}_{k}\}}:=\{C_{\{\mathfrak{m}_{k}\},1},C_{\{\mathfrak{m}_{k}\},2},C_{\{\mathfrak{m}_{k}\},3}\}$.
\subsection{Elementary properties}
We recall notation \eqref{eq.lower_and_upper_densities}.

The following fact is an immediate consequence of \eqref{M.2} and the definition of the set functions $\mathcal{H}_{\theta}$, $\theta \geq 0$. We omit an elementary proof.

\begin{Prop}
\label{Prop.absolute_continuity_measure}
Let $\theta \geq 0$ and $\{\mathfrak{m}_{k}\} \in \mathfrak{M}_{\theta}(S)$ for some closed set $S \subset \operatorname{X}$. Then, for each $k \in \mathbb{N}_{0}$,
the measure $\mathfrak{m}_{k}$ is absolutely continuous with respect to $\mathcal{H}_{\theta}$. Furthermore, for any Borel set $E \subset S$, for each $k \in \mathbb{N}_{0}$,
we have $\mathfrak{m}_{k}(E) \le C_{\{\mathfrak{m}_{k}\},1}\mathcal{H}_{\theta}(E)$.
\end{Prop}

Given a sequence $\{\mathfrak{m}_{k}\} \in \mathfrak{M}_{\theta}(S)$, it is naturally to ask whether the measures $\mathfrak{m}_{k}$, $k \in \mathbb{N}_{0}$
have the doubling properties? Unfortunately, this is not the case in general. Nevertheless, we have the following important result (we put $B_{k}(x):=B_{\epsilon^{k}}(x)$
for all $x \in \operatorname{X}$ and $k \in \mathbb{Z}$).

\begin{Th}
\label{Th.doubling_type}
Let $\theta \geq 0$, $S \in \mathcal{LCR}_{\theta}(\operatorname{X})$ and $\{\mathfrak{m}_{k}\} \in \mathfrak{M}_{\theta}(S)$. Then
for each $c \geq 1$, there is a constant $C > 0$ depending on $c$, $\mathcal{C}_{\{\mathfrak{m}_{k}\}}$ and $C_{\mu}(c)$ such that, for
each $k \in \mathbb{N}_{0}$,
\begin{equation}
\label{eq.doubling_type}
\frac{1}{C}\mathfrak{m}_{k}(B_{k}(y))\le \mathfrak{m}_{k}(\frac{1}{c}B_{k}(y)) \le \mathfrak{m}_{k}(cB_{k}(y)) \le C  \mathfrak{m}_{k}(B_{k}(y)) \quad \text{for all} \quad y \in S.
\end{equation}
\end{Th}
\begin{proof}
We set $\underline{k}:=\min\{k \in \mathbb{Z}:\epsilon^{k} < \frac{1}{c}\}$ and consider the upper and the lower bounds in \eqref{eq.doubling_type} separately.

To prove the upper bound we consider two cases. In the case $k \in \{0,...,\underline{k}\}$, a combination of \eqref{M.2}, \eqref{M.4} and the uniformly locally doubling property of $\mu$
gives the existence of a constant $C > 0$ such that (we recall \eqref{eq.lattice} and take into account Proposition \ref{Prop.finite_intersection})
\begin{equation}
\begin{split}
\label{eq.232}
&\mathfrak{m}_{k}(cB_{k}(y)) \le C \mathfrak{m}_{0}(cB_{k}(y)) \le C \sum\limits_{B \in \mathcal{B}_{k}(\operatorname{X},\epsilon)}\mathfrak{m}_{0}(B \cap cB_{k}(y))\\
&\le C\sum\limits_{\substack{B \in \mathcal{B}_{k}(\operatorname{X},\epsilon)\\ B \cap cB_{k}(y) \neq \emptyset}}\mu(B) \le C \mu((c+2)B_{k}(y)) \le C \mu(B_{k}(y)).
\end{split}
\end{equation}
In the case $k > \underline{k}$, using \eqref{M.2}, \eqref{M.3}, \eqref{M.4}, and the uniformly locally doubling property of $\mu$, we obtain
\begin{equation}
\label{eq.233}
\begin{split}
&\mathfrak{m}_{k}(cB_{k}(y)) \le  C\mathfrak{m}_{k-\underline{k}}(cB_{k}(y)) \le
C\mathfrak{m}_{k-\underline{k}}(B_{k-\underline{k}}(y))\\
&\le C\frac{\mu(B_{k-\underline{k}}(y))}{\epsilon^{k-\underline{k}}} \le C \frac{\mu(B_{k}(y))}{\epsilon^{k}} \le C \mathfrak{m}_{k}(B_{k}(y)).
\end{split}
\end{equation}
Combining \eqref{eq.232} and \eqref{eq.233} we deduce the required upper bound in \eqref{eq.doubling_type}.

Now we fix $k \in \mathbb{N}_{0}$.
To prove the lower bound in \eqref{eq.doubling_type}, an appeal to  \eqref{M.2}, \eqref{M.3}, \eqref{M.4}  and the uniformly locally doubling property of $\mu$,
gives us the required estimate
\begin{equation}
\begin{split}
&\mathfrak{m}_{k}(\frac{1}{c}B_{k}(y)) \geq C \mathfrak{m}_{k+\underline{k}}(\frac{1}{c}B_{k}(y)) \geq C \mathfrak{m}_{k+\underline{k}}(B_{k+\underline{k}}(y))\\
&\geq C \frac{\mu(B_{k+\underline{k}}(y))}{\epsilon^{k+\underline{k}}} \geq C \frac{\mu(B_{k}(y))}{\epsilon^{k}} \geq C \mathfrak{m}_{k}(B_{k}(y)).
\end{split}
\end{equation}

The proof is complete.
\end{proof}

The above theorem leads to the following useful corollary (we put $B_{k}(x):=B_{\epsilon^{k}}(x)$).

\begin{Prop}
\label{Prop.4.1}
Let $\theta \geq 0$, $S \in \mathcal{LCR}_{\theta}(\operatorname{X})$ and $\{\mathfrak{m}_{k}\} \in \mathfrak{M}_{\theta}(\operatorname{X})$. Then, for each $c \geq 1$,
there exists a constant $C > 0$ such that, for each $k \in \mathbb{N}_{0}$,
\begin{equation}
\int\limits_{cB_{k}(z)}\frac{1}{\mathfrak{m}_{k}(cB_{k}(y))}\,d\mathfrak{m}_{k}(y) \le C \quad \text{for all} \quad z \in S.
\end{equation}
\end{Prop}

\begin{proof}
Fix $k \in \mathbb{N}_{0}$ and $z \in S$.
Note that $2cB_{k}(y) \supset cB_{k}(z)$ for all $y \in cB_{k}(z)$. Hence, by Theorem \ref{Th.doubling_type}, 
\begin{equation}
\notag
\sup\limits_{y \in cB_{k}(z)}\frac{1}{\mathfrak{m}_{k}(cB_{k}(y))} \le \sup\limits_{y \in cB_{k}(z)}\frac{C}{\mathfrak{m}_{k}(2cB_{k}(y))} \le \frac{C}{\mathfrak{m}_{k}(cB_{k}(z))}
\end{equation}
Consequently,
\begin{equation}
\int\limits_{cB_{k}(z)}\frac{1}{\mathfrak{m}_{k}(cB_{k}(y))}\,d\mathfrak{m}_{k}(y) \le C \int\limits_{cB_{k}(z)}\frac{1}{\mathfrak{m}_{k}(cB_{k}(z))}\,d\mathfrak{m}_{k}(y) = C.
\end{equation}
The proof is complete.
\end{proof}

\subsection{Comparison of different classes of measures}

Now we formulate, for a given closed set $S \subset \operatorname{X}$, a simple condition that is sufficient
for the equality of classes $\mathfrak{M}_{\theta}^{str}(S)$ and $\mathfrak{M}_{\theta}(S)$. We recall Definition \ref{Def.lower_decay_exponent}, Remark \ref{Rem.Hausdorff_is_a_Borel_measure}. We also recall notation \eqref{eqq.important_index_notation} and put $k(r):=k_{\epsilon}(r)$.

\begin{Th}
\label{Th.coincideness_of_classes_of_measures}
Let $\theta \in [0,\underline{Q}_{\mu}(1))$ and let $S \subset \operatorname{X}$ be such that $\mathcal{H}_{\theta}(S) \in (0,+\infty)$. Then
\begin{equation}
\notag
\mathfrak{M}_{\theta}^{str}(S) = \mathfrak{M}_{\theta}(S).
\end{equation}
\end{Th}

\begin{proof}
Clearly, it is sufficient to show that $\mathfrak{M}_{\theta}(S) \subset \mathfrak{M}^{str}_{\theta}(S)$.
Assume that $\mathfrak{M}_{\theta}(S) \neq \emptyset$ and
fix an arbitrary $\{\mathfrak{m}_{k}\} \in \mathfrak{M}_{\theta}(S)$. We also fix a Borel set $E \subset S$ and verify \eqref{M5}.
We put $N:=\{x \in E: \overline{D}^{\{\mathfrak{m}_{k}\}}_{E}(x,\epsilon) = 0\}.$
Since $\theta \in [0,\underline{Q}_{\mu}(1))$, by Remark \ref{Rem.Hausdorff_is_a_Borel_measure}
the set function $\mathcal{H}_{\theta}\lfloor_{S}$ is a finite measure on $\operatorname{X}$. Furthermore, by Proposition \ref{Prop.absolute_continuity_measure} $\mathfrak{m}_{0}(S) < +\infty$.
If $\mathfrak{m}_{0}(N) > 0$, then using Egorov's theorem, given $\varepsilon > 0$, we
find a compact set $K_{\varepsilon} \subset N$ and a number $\delta(\varepsilon) > 0$ such that $\mathfrak{m}_{0}(N \setminus K_{\varepsilon}) < \varepsilon$ and
\begin{equation}
\label{eqq.46''}
\sup\limits_{x \in K_{\varepsilon}}\sup\limits_{r < \delta(\varepsilon)}\frac{\mathfrak{m}_{k(r)}(E \cap B_{r}(x))}{\mathfrak{m}_{k(r)}(B_{r}(x))} < \varepsilon.
\end{equation}
By the assumption of the lemma we have $\mathcal{H}_{\theta}(S) < +\infty$.
Hence, we find an arbitrary at most countable covering of $K_{\varepsilon}$ by balls $\{B_{j}\}_{j=1}^{N}$, $N \in \mathbb{N} \cup \{\infty\}$ with $r_{j}:=r(B_{j}) \le \frac{\delta(\varepsilon)}{4}$ such that
\begin{equation}
\label{eqq.47''}
\sum_{j} \frac{\mu(B_{j})}{(r_{j})^{\theta}} \le 2\mathcal{H}_{\theta,\delta(\varepsilon)}(K_{\varepsilon}).
\end{equation}
Without loss of generality we may assume that, for each $B_{j}$, $K_{\varepsilon} \cap B_{j} \neq \emptyset$.
For each $j$ we fix a point $x_{j} \in B_{j} \cap K_{\varepsilon}$. We obviously have $B_{j} \subset B_{2r_{j}}(x_{j}) \subset 3B_{j}$.
Hence, combining \eqref{M.2}, \eqref{M.4} with \eqref{eqq.46''}, \eqref{eqq.47''} and taking into account
the uniformly locally doubling property of $\mu$, we have
\begin{equation}
\begin{split}
&\mathfrak{m}_{0}(K_{\varepsilon}) \le \sum\limits_{j} \mathfrak{m}_{0}(E \cap B_{j}) \le C\sum_{j} \mathfrak{m}_{k(r_{j})}(E \cap B_{2r_{j}}(x_{j}))
< \varepsilon C \sum_{j} \mathfrak{m}_{k(r_{j})}(B_{2r_{j}}(x_{j}))\\
&< \varepsilon C \sum_{j} \frac{\mu(B_{2r_{j}}(x_{j}))}{(r_{j})^{\theta}} < \varepsilon C \sum_{j} \frac{\mu(3B_{j})}{(r_{j})^{\theta}} \le \varepsilon C \mathcal{H}_{\theta,\delta(\varepsilon)}(K_{\varepsilon}) \le
\varepsilon C \mathcal{H}_{\theta}(S).
\end{split}
\end{equation}
Hence, for all small enough $\varepsilon > 0$ we have $\mathfrak{m}_{0}(K_{\varepsilon})=0$. Since $\mathfrak{m}_{0}(N \setminus K_{\varepsilon}) < \varepsilon$
and $\varepsilon > 0$ can be chosen arbitrarily, we get $\mathfrak{m}_{0}(N)=0$, completing the proof.
\end{proof}

In the proof of the next lemma we built a simple example that exhibits a delicate difference between the classes
$\mathfrak{M}_{\theta}^{str}(S)$ and $\mathfrak{M}_{\theta}(S)$. Despite of its simplicity, the corresponding constructions are typical
and reflect the essence of the matter. One can built similar examples in higher dimensions and
even in some nice classes of metric measure spaces. However, the corresponding machinery will be much
less transparent.

\begin{Th}
\label{Th.different_classes_of_measures}
Let $\operatorname{X}=(\mathbb{R}^{2},\|\cdot\|_{2},\mathcal{L}^{2})$ and $S:=\{(x_{1},x_{2})\in\mathbb{R}^{2}:x_{1} \in [0,1], x_{2}=0\}$. Then, for each $\theta \in (1,2)$, there exists
a sequence of measures $\{\mathfrak{m}_{k}\} \in \mathfrak{M}_{\theta}(S) \setminus \mathfrak{M}^{str}_{\theta}(S)$.
\end{Th}

\begin{proof}
We fix an arbitrary $\theta \in (1,2)$ and put $c_{1}(\theta):=2\sum_{k=1}^{\infty}\frac{1}{k^{\theta}}$, $c_{2}(\theta):=\min_{j \in \mathbb{N}_{0}}\frac{2^{j}}{(1+j)^{\theta-1}}$.
It will be convenient to split the proof into several steps.

\textit{Step 1.} Let $E$ denote
the closed Cantor-type set built inductively as follows.
At the first step we put $E_{1}:=[0,1] \setminus (\frac{1}{2}-(2c_{1}(\theta))^{-1},\frac{1}{2}+(2c_{1}(\theta))^{-1})$ and $U_{1}:=(\frac{1}{2}-(2c_{1}(\theta))^{-1},\frac{1}{2}+(2c_{1}(\theta))^{-1})$.
Suppose that for some $k \in \mathbb{N}$
we have already built closed sets $E_{1} \supset ... \supset E_{k}$ and open sets $U_{1},...,U_{k}$ such that
$E_{i} \cup (\cup_{j=1}^{i}U_{j}) = [0,1]$ and $\mathcal{L}^{1}(U_{i}) = \frac{1}{c_{1}(\theta)i^{\theta}}$ for all $i \in \{1,...,k\}.$
Furthermore, given $i \in \{1,...,k\}$, each $E_{i}$ is a disjoint union of $2^{i}$ closed intervals $I_{i,l}$ and each $U_{i}$ is a disjoint union of $2^{i-1}$ open intervals $J_{i,l}$.
We delete from the middle of each closed interval $I_{k,l}$ an open interval of length $\frac{1}{c_{1}(\theta)2^{k}k^{\theta}}$ and consider the union of the remaining closed sets. We obtain the set $E_{k+1}$
and put $U_{k+1}:=E_{k} \setminus E_{k+1}$. As a result, we obtain the sequence $\{E_{k}\}_{k=1}^{\infty}$ of closed sets and the sequence $\{U_{k}\}_{k=1}^{\infty}$ of open sets.
Furthermore, for each $k \in \mathbb{N}$, we let $\mathcal{I}_{k}$ and $\mathcal{J}_{k}$ denote the corresponding families of closed and open intervals respectively. More precisely,
$E_{k}=\cup\{I: I \in \mathcal{I}_{k}\}$ and $U_{k}=\cup\{J: J \in \mathcal{J}_{k}\}$ for all $k \in \mathbb{N}$.
Now we put $E:=\cap_{n=0}^{\infty}E_{n}$ and define the weight functions $\omega_{k} \in L_{1}([0,1])$, $k \in \mathbb{N}_{0}$, by (we put $U_{0}:=\emptyset$)
\begin{equation}
\begin{split}
\label{eqq.49''}
&\omega_{k}(x):=\chi_{E}(x)+\sum\limits_{i=0}^{k}2^{(\theta-1)i}i^{\theta-1}\chi_{U_{i}}(x)
+2^{(\theta-1)k}(k+1)^{\theta-1}\sum\limits_{i=k+1}^{\infty}\chi_{U_{i}}(x), \quad x \in [0,1].
\end{split}
\end{equation}
Finally, we recall \eqref{eqq.weighted_measure} and put $\mathfrak{m}_{k}:=\omega_{k}\mathcal{H}^{1}\lfloor_{S}$, $k \in \mathbb{N}_{0}$ (here $\mathcal{H}^{1}$ is the usual $1$-dimensional Hausdorff measure).
We put $\epsilon = \frac{1}{2}$ and claim that $\{\mathfrak{m}_{k}\}:=\{\mathfrak{m}_{k}\}_{k \in \mathbb{N}_{0}} \in \mathfrak{M}_{\theta}(S) \setminus \mathfrak{M}^{str}_{\theta}(S)$.
This will be shown in the next steps.

\textit{Step 2.} Note that $\operatorname{supp}\mathfrak{m}_{k}=S$ for all $k \in \mathbb{N}_{0}$. This verifies (\textbf{M.1}).

\textit{Step 3.}
By \eqref{eqq.49''} it is easy to see that, for each $k \in \mathbb{N}_{0}$ and every $j \in \mathbb{N}_{0}$,
\begin{equation}
\label{eqq.511''}
\frac{c_{2}(\theta)}{2^{\theta j}} \le \frac{1}{2^{(\theta-1)j}(1+j)^{\theta-1}} \le \frac{(k+1)^{\theta-1}}{2^{(\theta-1)j}(k+1+j)^{\theta-1}} \le \frac{w_{k}(x)}{w_{k+j}(x)} \le 1, \quad  x \in [0,1].
\end{equation}
This proves that condition (\textbf{M.4}) is satisfied with $C_{3}=\max\{1,(c_{2}(\theta))^{-1}\}$.

\textit{Step 4.}
To verify (\textbf{M.2}) we proceed as follows. We fix arbitrary $k \in \mathbb{N}_{0}$, $j \geq k$ and $Q \in \mathcal{D}_{j}$ (by $\mathcal{D}_{j}$ here and below we denote the family of closed dyadic intervals of length $2^{-j}$). Given $i \in \mathbb{N}$, there are two cases to consider.

\textit{In the first case}, $(c_{1}(\theta))^{-1}2^{-i}i^{-\theta}< 2^{-j}$. Since $\theta > 1$, we obviously have
\begin{equation}
\begin{split}
\label{eqq.additional_estimate}
&\frac{2^{(\theta-1)i}i^{\theta-1}}{2^{i}i^{\theta}}=\frac{2^{(\theta-1)i}i^{\theta-1}}{2^{(\theta-1+2-\theta)i}i^{\theta(\theta-1+2-\theta)}} \le \frac{1}{2^{(2-\theta)i}}\frac{1}{i^{\theta(2-\theta)}} \le \frac{(c_{1}(\theta))^{2-\theta}}{2^{(2-\theta)j}}.
\end{split}
\end{equation}
Consequently, given $J \in \mathcal{J}_{i}$,
by \eqref{eqq.49''}, \eqref{eqq.additional_estimate} we obtain (since $(c_{1}(\theta))^{1-\theta} \le 1$)
\begin{equation}
\label{eqq.413''}
\frac{1}{2}\mathfrak{m}_{k}(Q \cap J) \le \frac{1}{2}\mathfrak{m}_{k}(J) \le
\begin{cases}
\frac{2^{(\theta-1)k}(k+1)^{\theta-1}}{c_{1}(\theta)2^{i}i^{\theta}} \le \frac{2^{(\theta-1)i}i^{\theta-1}}{c_{1}(\theta)2^{i}i^{\theta}} \le \frac{1}{2^{(2-\theta)j}}, \quad  \quad i > k;\\
\frac{2^{(\theta-1)i}i^{\theta-1}}{c_{1}(\theta)2^{i}i^{\theta}} \le \frac{1}{2^{(2-\theta)j}}, \quad  \quad i \le k.\\
\end{cases}
\end{equation}

\textit{In the second case}, $(c_{1}(\theta))^{-1}2^{-i}i^{-\theta} \geq 2^{-j}$. Since $\theta >1$, we clearly have
\begin{equation}
\label{eqq.additional_estimate_2}
\frac{2^{(\theta-1)i}i^{\theta-1}}{2^{j}} \le \frac{2^{(\theta-1)i}i^{\theta-1}}{2^{(2-\theta)j}2^{(\theta-1) j}} \le \frac{2^{(\theta-1)i}i^{\theta-1}}{(c_{1}(\theta))^{\theta-1}2^{(2-\theta)j}2^{(\theta-1)i}i^{\theta(\theta-1)}}
\le \frac{(c_{1}(\theta))^{1-\theta}}{2^{(2-\theta)j}}.
\end{equation}
Consequently, given $J \in \mathcal{J}_{i}$, by \eqref{eqq.49''}, \eqref{eqq.additional_estimate_2} we obtain (since $(c_{1}(\theta))^{1-\theta} \le 1$)
\begin{equation}
\label{eqq.414''}
\frac{1}{2}\mathfrak{m}_{k}(Q \cap J) \le
\begin{cases}
\frac{2^{(\theta-1)k}(k+1)^{\theta-1}}{2^{j}} \le \frac{2^{(\theta-1)i}i^{\theta-1}}{2^{j}} \le \frac{1}{2^{(2-\theta)j}}, \quad  \quad i > k;\\
\frac{2^{(\theta-1)i}i^{\theta-1}}{2^{j}} \le \frac{1}{2^{(2-\theta)j}}, \quad  \quad i \le k.
\end{cases}
\end{equation}

As a result, combining \eqref{eqq.413''} and \eqref{eqq.414''}, we have
\begin{equation}
\begin{split}
\label{eqq.415''}
\mathfrak{m}_{k}(Q \cap J) \le \frac{1}{2^{(2-\theta)j}} \quad \text{for each} \quad i \in \mathbb{N}, \quad \text{for all} \quad Q \in \mathcal{D}_{j}\text{ and }J \in \mathcal{J}_{i}.
\end{split}
\end{equation}

We fix a closed interval $I \in \mathcal{I}_{j}$ and note that $I \cap U_{i} = \emptyset$ for all $i \in \{1,...,j\}$. Hence, taking
into account that for each $i \geq j+1$ the set $U_{i} \cap I$ is composed of at most $2^{i-j}$ open intervals of length $(c_{1}(\theta))^{-1}2^{-i}i^{-\theta}$,
we have $c_{1}(\theta)\sum_{i=j+1}^{\infty}\mathcal{L}^{1}(U_{i} \cap I) \le 2^{-j}\sum_{i=j+1}^{\infty}i^{-\theta} \le \frac{2}{\theta-1}2^{-j}j^{1-\theta}$. Using this observation, and keeping in mind that
$\theta \in (1,2)$ and $c_{1}(\theta) \geq 2$, we have by \eqref{eqq.49''}
\begin{equation}
\begin{split}
\label{eqq.416''}
&\mathfrak{m}_{k}(Q \cap I) \le \mathfrak{m}_{k}(I) \le 2^{-j}+2^{(\theta-1)k}(k+1)^{\theta-1}\sum\limits_{i=j+1}^{\infty}\mathcal{L}^{1}(U_{i} \cap I)\\
&\le \frac{1}{2^{(2-\theta)j}}+ 2\frac{(c_{1}(\theta))^{-1}}{\theta-1}\frac{2^{(\theta-1)k}(k+1)^{\theta-1}}{2^{j}j^{\theta-1}} \le \frac{1}{2^{(2-\theta)j}}+\frac{1}{\theta-1}\frac{(j+1)^{\theta-1}}{j^{\theta-1}}\frac{1}{2^{(2-\theta)j}}\\
&\le \frac{2^{\theta}}{\theta-1}\frac{1}{2^{(2-\theta)j}} \quad \text{for each} \quad Q \in \mathcal{D}_{j}.
\end{split}
\end{equation}

By the construction of $U_{i}$, we have $\sum_{i=1}^{l}\mathcal{L}^{1}(U_{i}) < 1/2$ for all $l \in \mathbb{N}$. Hence, it is easy to see that
each closed interval $I \in \mathcal{I}_{j}$ has length greater than $2^{-j-1}$.
Consequently, it is easy to see that
$Q$ can intersect at most 3 different closed intervals
from $\mathcal{I}_{j}$ and at most 2 different open intervals from $\cup_{i=1}^{j}\mathcal{J}_{i}$. Hence, given $x \in \mathbb{R}^{2}$ and $r \in (2^{-j-1},2^{-j}]$,
combining the above observations with \eqref{eqq.415''}, \eqref{eqq.416''} we get
\begin{equation}
\label{eqq.417'''}
\mathfrak{m}_{k}(B_{r}(x)) \le \sum\limits_{\substack{Q \in \mathcal{D}_{j}}}\mathfrak{m}_{k}(Q \cap B_{r}(x)) \le \frac{2^{\theta}}{\theta-1}\frac{15}{2^{(2-\theta)j}}.
\end{equation}
Hence, we see that $\{\mathfrak{m}_{k}\}$ satisfies condition (\textbf{M.2}).

\textit{Step 5.}  We fix $\underline{x} \in E$ and $k \in \mathbb{N}_{0}$.
By the construction of $E$, there exists $I_{k}(\underline{x}) \in \mathcal{I}_{k}$ such that $\underline{x} \in I_{k}(\underline{x})$. Hence, taking into
account that, for each $i \geq k+1$, the set $U_{i} \cap I_{k}(\underline{x})$ consists of $2^{i-k}$ disjoint intervals of length $(c_{1}(\theta))^{-1}2^{-i}i^{-\theta}$ we get
\begin{equation}
\label{eqq.418'''}
\begin{split}
&\mathfrak{m}_{k}(B_{2^{-k}}(\underline{x})) \geq \mathfrak{m}_{k}(I_{k}(\underline{x})) \geq 2^{(\theta-1)k}(k+1)^{\theta-1}\sum\limits_{i=k+1}^{\infty}\mathcal{L}^{1}(U_{i} \cap I_{k}(\underline{x}))\\
&\geq \frac{2^{(\theta-1)k}(k+1)^{\theta-1}}{c_{1}(\theta)2^{k}}\sum\limits_{i=k+1}^{\infty}\frac{1}{i^{\theta}} \geq \frac{1}{c_{1}(\theta)(\theta-1)}\frac{1}{2^{k(2-\theta)}}.
\end{split}
\end{equation}
This observation in combination with \eqref{eqq.511''} easily implies condition (\textbf{M.3}).

\textit{Step 6.} By \eqref{eqq.49''} it is easy to see that
$\mathfrak{m}_{k}(E \cap B_{2^{-k}}(x)) \le 2^{-k+1}$ for all $x \in E$ and all $k \in \mathbb{N}_{0}$. Since $\theta > 1$, the observation above in combination with \eqref{eqq.418'''} gives
\begin{equation}
\frac{\mathfrak{m}_{k}(B_{2^{-k}}(x) \cap E)}{\mathfrak{m}_{k}(B_{2^{-k}}(x))} \le \frac{2c_{1}(\theta)(\theta-1)}{2^{k(\theta-1)}} \to 0, \quad k \to \infty \quad \text{for all} \quad x \in E.
\end{equation}
This proves that $\{\mathfrak{m}_{k}\} \notin \mathfrak{M}_{\theta}^{str}(S)$.

The proof is complete.
\end{proof}

\subsection{Proof of Theorem \ref{Th.regularsequence}}
We start by proving that, given $\theta \geq 0$, one can construct $\theta$-regular sequences of measures only for
lower codimension-$\theta$ content regular sets $S$.

\begin{Th}
\label{Th.necessary_regular_sequence}
Let $S \subset \operatorname{X}$ be a closed nonempty set. If $\theta \geq 0$ is such that $\mathfrak{M}_{\theta}(S) \neq \emptyset$, then
$S \in \mathcal{LCR}_{\theta}(\operatorname{X})$.
\end{Th}

\begin{proof}
Let $\{\mathfrak{m}_{k}\} \in \mathfrak{M}_{\theta}(S)$ and $\epsilon=\epsilon(\{\mathfrak{m}_{k}\}) \in (0,1)$. Given $r \in (0,1]$ and $x \in S$, let
$\mathcal{B}=\{B_{j}\}_{j \in \mathbb{N}}=\{B_{r_{j}}(x_{j})\}_{j \in \mathbb{N}}$
be a sequence of closed balls such that $r_{j} < r$ for all $j \in \mathbb{N}$, $B_{j} \cap S \neq \emptyset$ for all $j \in \mathbb{N}$, $B_{r}(x) \cap S \subset \cup_{j \in \mathbb{N}}B_{j}$ and
\begin{equation}
\label{eqq.426'''}
\sum\limits_{j}\frac{\mu(B_{j})}{(r_{j})^{\theta}} \le 2\mathcal{H}_{\theta,r}(B_{r}(x) \cap S).
\end{equation}
We recall notation \eqref{eqq.important_index_notation}, put $k_{j}:=k(r_{j})$, $j \in \mathbb{N}$, and fix, for every $j \in \mathbb{N}$, a ball $\widetilde{B}_{j}$ centered at some point $x_{j} \in S \cap B_{j}$ such that  $r(\widetilde{B}_{j})=2r_{j}$.
It is clear that $B_{j} \subset \widetilde{B}_{j} \subset 4B_{j}$
for all $j \in \mathbb{N}$. Hence, using the uniformly locally doubling property of the measure $\mu$, applying \eqref{M.2}, and then taking into account Theorem \ref{Th.doubling_type}, we obtain
\begin{equation}
\label{eqq.427'''}
\begin{split}
&\sum\limits_{j}\frac{\mu(B_{j})}{(r_{j})^{\theta}} \geq C\sum\limits_{j}\frac{\mu(4B_{j})}{(4r_{j})^{\theta}} \geq C \sum\limits_{j}\frac{\mu(\frac{1}{2}\widetilde{B}_{j})}{(\frac{r_{j}}{2})^{\theta}} \geq C\sum\limits_{j}\mathfrak{m}_{k_{j}}(\frac{1}{2}\widetilde{B}_{j}) \geq C\sum\limits_{j}\mathfrak{m}_{k_{j}}(\widetilde{B}_{j}).
\end{split}
\end{equation}
Now we combine \eqref{eqq.426'''} with \eqref{eqq.427'''}, take into account that $k(r) \le k_{j}$ for all $j \in \mathbb{N}$, and use \eqref{M.4}. This gives
\begin{equation}
\notag
\begin{split}
\mathcal{H}_{\theta,r}(B_{r}(x) \cap S) \geq  C\sum\limits_{j}\mathfrak{m}_{k(r)}(\widetilde{B}_{j}) \geq C \mathfrak{m}_{k(r)}(B_{r}(x) \cap S).
\end{split}
\end{equation}
Since $r \in (\epsilon^{k(r)+1},\epsilon^{k(r)}]$, by Theorem \ref{Th.doubling_type} and \eqref{M.3} we continue the previous estimate and get
\begin{equation}
\label{eqq.428'''}
\begin{split}
\mathcal{H}_{\theta,r}(B_{r}(x) \cap S) \geq C \mathfrak{m}_{k(r)}(B_{\epsilon^{k(r)}}(x) \cap S) \geq C\frac{\mu(B_{\epsilon^{k(r)}}(x))}{\epsilon^{k(r)\theta}} \geq C \frac{\mu(B_{r}(x))}{r^{\theta}}.
\end{split}
\end{equation}
Since $x \in S$ and $r \in (0,1]$ were chosen arbitrarily, the theorem follows from Definition \ref{Def.content}.
\end{proof}

The following result gives conditions on a given set $S \subset \operatorname{X}$ sufficient for the existence of a strongly $\theta$-regular sequence of measures
whose supports coincide with $S$.

\begin{Th}
\label{Th.sufficiency_regular_sequence}
Let $\theta \geq 0$. If $S \in \mathcal{LCR}_{\theta}(\operatorname{X})$, then $\mathfrak{M}^{str}_{\theta}(S) \neq \emptyset$.
\end{Th}

\begin{proof}

We fix an arbitrary $\epsilon \in (0,\frac{1}{10}]$ and recall notation \eqref{eqq.centers_of_dyadic_cubes}--\eqref{eqq.centers_of_dyadic_cubes_2}. Since $S$ is a closed subset of a complete separable
metric space $(\operatorname{X},\operatorname{d})$, the space $\operatorname{S}:=(S,\operatorname{d}|_{S})$ is a complete separable metric space (here $\operatorname{d}|_{S}$ is a
restriction of the metric $\operatorname{d}$ to the set $S$).

\textit{Step 1.}
We recall Definition \ref{Def.partial_order}, Proposition \ref{Prop.partial_order}  and fix an admissible partial order on $Z(\operatorname{S},\epsilon)$. Given $k \in \mathbb{N}_{0}$ and $z_{k,\alpha} \in Z_{k}(\operatorname{S},\epsilon)$, we put
\begin{equation}
\label{eqq.pseudo_dyadic_1}
\widetilde{Q}_{k,\alpha}:=\cup \{\operatorname{int}B_{\frac{\epsilon^{j}}{8}}(z_{j,\beta}):z_{j,\beta} \preceq z_{k,\alpha}\}.
\end{equation}
Note that $\widetilde{Q}_{k,\alpha}$, $k \in \mathbb{N}_{0}$, $\alpha \in \mathcal{A}_{k}(\operatorname{S},\epsilon)$,
are open subsets in $\operatorname{X}$. However, they are neither generalized dyadic cubes in $\operatorname{X}$ nor generalized dyadic cubes in $\operatorname{S}$.
At the same time, by \eqref{eqq.generalised_dyadic_cube} $\widetilde{Q}_{k,\alpha} \cap S$ is a generalized dyadic cube in the space $\operatorname{S}$ for each $k \in \mathbb{N}_{0}$ and any $\alpha \in \mathcal{A}_{k}(\operatorname{S},\epsilon)$. The only reason for introduction of such specific sets $\widetilde{Q}_{k,\alpha}$ is that the ``centers'' of such ``quasicubes'' belong
to the set $S$. This fact will be crucial at step 8 below.

Since $\epsilon \in (0,\frac{1}{10}]$ it is easy to see from (PO3) of Definition \ref{Def.partial_order} and \eqref{eqq.pseudo_dyadic_1} that
\begin{equation}
\label{eqq.pseudo_dyadic_4}
\widetilde{Q}_{k,\alpha} \subset B_{2\epsilon^{k}}(z_{k,\alpha})\quad \text{for each} \quad k \in \mathbb{N}_{0} \quad \text{and any} \quad \alpha \in \mathcal{A}_{k}(\operatorname{S},\epsilon).
\end{equation}
Repeating almost verbatim the arguments from the proof of Lemma 15 in \cite{Christ}, we have
\begin{equation}
\label{eqq.pseudo_dyadic_2}
\text{if} \quad l \geq k \quad \text{then either} \quad \widetilde{Q}_{l,\beta} \subset \widetilde{Q}_{k,\alpha} \quad \text{or} \quad \widetilde{Q}_{l,\beta} \cap \widetilde{Q}_{k,\alpha} = \emptyset.
\end{equation}
Furthermore, by \eqref{eqq.pseudo_dyadic_1} and \eqref{eqq.pseudo_dyadic_2} we clearly have
\begin{equation}
\label{eqq.pseudo_dyadic_3}
\sum\limits_{z_{k+1,\beta} \preceq z_{k,\alpha}}\mu(\widetilde{Q}_{k+1,\beta}) \le \mu(\widetilde{Q}_{k,\alpha}) \quad \text{for each} \quad k \in \mathbb{N}_{0} \quad \text{and any} \quad \alpha \in \mathcal{A}_{k}(\operatorname{S},\epsilon).
\end{equation}
Finally, let $B=B_{r}(x)$ be an arbitrary closed ball with $x \in \operatorname{X}$, $r \in (0,1]$. Let $c \geq 1$ be such that $B_{cr}(x) \cap S \neq \emptyset$.
The same arguments as in the proof of Proposition \ref{Prop.114} give
\begin{equation}
\label{eqq.intersection'}
\begin{split}
&\#\{\alpha \in \mathcal{A}_{k(r)}(\operatorname{S},\epsilon): \operatorname{cl}(\widetilde{Q}_{k(r),\alpha} \cap S) \cap cB \neq \emptyset\} \le C_{D}(c,0).
\end{split}
\end{equation}

\textit{Step 2.}
For each $j \in \mathbb{N}_{0}$ and any $\beta \in \mathcal{A}_{j}(\operatorname{S},\epsilon)$ we put $h_{j,\beta}:=\frac{\mu(\widetilde{Q}_{j,\beta})}{\epsilon^{j\theta}}$.
Now, for each $j \in \mathbb{N}_{0}$, we define the measure $\mathfrak{m}^{j,j}$ on $S$ by the formula (by $\delta_{x}$ we denote the Dirac measure concentrated at $x \in S$)
\begin{equation}
\label{weight2}
\mathfrak{m}^{j,j}:=\sum\limits_{z_{j,\beta} \in Z_{j}(\operatorname{S},\epsilon)}h_{j,\beta}\delta_{z_{j,\beta}}.
\end{equation}
Given $j \in \mathbb{N}$, we modify the measure $\mathfrak{m}^{j,j}$ in the following way. If $\alpha \in \mathcal{A}_{j-1}(\operatorname{S},\epsilon)$ is such that
$\mathfrak{m}^{j,j}(\widetilde{Q}_{j-1,\alpha} \cap S)=\mathfrak{m}^{j,j}(\{z_{j,\beta}:z_{j,\beta} \preceq z_{j-1,\alpha}\}) > h_{j-1,\alpha}$, then
we reduce the mass of $\mathfrak{m}^{j,j}$ uniformly on $\{z_{j,\beta}:z_{j,\beta} \preceq z_{j-1,\alpha}\}$ until
it becomes $h_{j-1,\alpha}$. On the other hand, if $\alpha \in \mathcal{A}_{j-1}(\operatorname{S},\epsilon)$ is such that
$\mathfrak{m}^{j,j}(\widetilde{Q}_{j-1,\alpha} \cap S) \le h_{j-1,\alpha}$, then
we leave $\mathfrak{m}^{j,j}$ on $\{z_{j,\beta}:z_{j,\beta} \preceq z_{j-1,\alpha}\}$ unchanged.
In this way we clearly obtain a new measure $\mathfrak{m}^{j,j-1}$. We repeat this procedure for $\mathfrak{m}^{j,j-1}$ obtaining $\mathfrak{m}^{j,j-2}$,
and after $j$ steps we obtain $\mathfrak{m}^{j,0}$.
Given $j \in \mathbb{N}_{0}$ and $k \le j$, it follows from this construction that
\begin{equation}
\label{eq.uppergrowth}
\mathfrak{m}^{j,k}(\widetilde{Q}_{i,\beta} \cap S) \le h_{i,\beta}
\quad \text{for each} \quad i \in \{k,...,j\} \quad \text{for all} \quad \beta \in \mathcal{A}_{i}(\operatorname{S},\epsilon).
\end{equation}

By \eqref{eqq.pseudo_dyadic_3} it is clear that
\begin{equation}
\label{eqq.important_lower_bound}
\underline{M}:=\inf\limits_{k \in \mathbb{N}_{0}}\inf\limits_{z_{k,\alpha} \in Z_{k}(\operatorname{S},\epsilon)} \Bigl(\sum_{z_{k+1,\beta} \preceq z_{k,\alpha}}h_{k+1,\beta}\Bigr)^{-1} h_{k,\alpha} \geq \epsilon^{\theta}.
\end{equation}
Note that by the construction, for each $k \in \mathbb{N}_{0}$ and every $j \geq k$,
there is a family of \textit{positive constants} $\{c_{j,k}(\widetilde{Q}_{j,\beta}): \beta \in \mathcal{A}_{j}(\operatorname{S},\epsilon)\}$
such that
\begin{equation}
\label{eqq.positive_constants}
\mathfrak{m}^{j,k}(\widetilde{Q}_{j,\beta} \cap S)=c_{j,k}(\widetilde{Q}_{j,\beta})\mathfrak{m}^{j,j}(\widetilde{Q}_{j,\beta} \cap S) \quad \text{for all} \quad \beta \in \mathcal{A}_{j}(\operatorname{S},\epsilon).
\end{equation}
Furthermore, by \eqref{eqq.important_lower_bound} for each $k \in \mathbb{N}_{0}$, $i \in \{0,...,k\}$ and every $j \geq k$
\begin{equation}
\label{eq.twosided}
\epsilon^{i}c_{j,k}(\widetilde{Q}_{j,\beta}) \le c_{j,k-i}(\widetilde{Q}_{j,\beta}) \le c_{j,k}(\widetilde{Q}_{j,\beta}) \quad \text{for all}  \quad \beta \in \mathcal{A}_{j}(\operatorname{S},\epsilon).
\end{equation}

\textit{Step 3.}
Using the right-hand side of \eqref{eq.uppergrowth} and \eqref{eqq.intersection'}, one can get
$\sup_{j \geq k}\mathfrak{m}^{j,k}(B) < \infty$ for every closed ball $B \subset \operatorname{X}$.
Hence, by Proposition \ref{Prop.loccomp} and Lemma \ref{Lm.weakconv}, given $k \in \mathbb{N}_{0}$, there is a subsequence $\{\mathfrak{m}^{j_{s},k}\}$ and a (Borel regular) measure $\mathfrak{m}_{k}$ on $\operatorname{X}$
such that $\mathfrak{m}^{j_{s},k} \rightharpoonup \mathfrak{m}_{k}$, $s \to \infty$.
In fact, by the standard diagonal arguments we conclude that there is a strictly increasing sequence $\{j_{l}\}_{l=1}^{\infty} \subset \mathbb{N}$
such that $\mathfrak{m}^{j_{l},k} \rightharpoonup \mathfrak{m}_{k}$, $l \to \infty$ for every $k \in \mathbb{N}_{0}$
(in the case $j_{l} < k$ we put formally $\mathfrak{m}^{j_{l},k}:=\mathfrak{m}^{k,k}$).

In the next steps we show that the sequence $\{\mathfrak{m}_{k}\}:=\{\mathfrak{m}_{k}\}_{k=0}^{\infty}$ satisfies conditions $(\textbf{M}1) - (\textbf{M}5)$.

\textit{Step 4.}
Given $k \in \mathbb{N}_{0}$, by \eqref{weight2} it follows that $\operatorname{supp}\mathfrak{m}^{j_{l},k}$ is an $\epsilon^{j_{l}}$-separated subset of $S$
for all large enough $l \in \mathbb{N}$. Hence, by the construction of $\mathfrak{m}_{k}$, $k \in \mathbb{N}_{0}$,
we conclude that $\operatorname{supp}\mathfrak{m}_{k}=S$ for every $k \in \mathbb{N}_{0}$.
This implies that condition $(\textbf{M}1)$ is satisfied.

\textit{Step 5.}
We fix arbitrary $k,i \in \mathbb{N}_{0}$. By \eqref{eqq.positive_constants}, \eqref{eq.twosided} we obtain for each $\varphi \in C_{c}(\operatorname{X})$ and all large enough $l \in \mathbb{N}$,
$\epsilon^{i \theta }\int_{\operatorname{X}}\varphi(x)\,d\mathfrak{m}^{j_{l},k+i}(x) \le \int_{\operatorname{X}}\varphi(x)\,d\mathfrak{m}^{j_{l},k}(x) \le \int_{\operatorname{X}}\varphi(x)\,d\mathfrak{m}^{j_{l},k+i}(x).$
Hence, passing to the limit as $l \to \infty$, we get
$$
\epsilon^{i \theta}\int_{\operatorname{X}}\varphi(x)\,d\mathfrak{m}_{k+i}(x) \le \int_{\operatorname{X}}\varphi(x)\,d\mathfrak{m}_{k}(x)
\le \int_{\operatorname{X}}\varphi(x)\,d\mathfrak{m}_{k+i}(x) \quad  \text{for all} \quad \varphi \in C_{c}(\operatorname{X}).
$$
As a result, using the Borel regularity of the measures $\mathfrak{m}_{k}$, $k \in \mathbb{N}_{0}$ and the Radon--Nikod\'ym theorem, we see
that $(\textbf{M}4)$ is satisfied with $C_{3}=1$.

\textit{Step 6.} Given $k \in \mathbb{N}_{0}$, we fix an arbitrary closed ball $B_{r}(x)$ with $x \in \operatorname{X}$ and $r \in (0,\epsilon^{k}]$.
If $\widetilde{Q}_{k(r),\alpha} \cap B_{2r}(x) \neq \emptyset$ for some $\alpha \in \mathcal{A}_{k(r)}(\operatorname{S},\epsilon)$, then by  \eqref{eqq.pseudo_dyadic_4}
we have $\widetilde{Q}_{k(r),\alpha}  \subset B_{cr}(x)$ with $c=4\epsilon^{k(r)}/r+2$.
Furthermore, since $k(r) \geq k$, by \eqref{eq.uppergrowth} we have $\mathfrak{m}^{j_{l},k}(\widetilde{Q}_{k(r),\alpha}) \le h_{k(r),\alpha}$ for each
$\alpha \in \mathcal{A}_{k(r)}(\operatorname{S},\epsilon)$ and all large enough $l \in \mathbb{N}$. Finally, by the construction we have $\mathfrak{m}^{j_{l},k}(\partial\widetilde{Q}_{k(r),\alpha})=0$
for each $\alpha \in \mathcal{A}_{k(r)}(\operatorname{S},\epsilon)$ and all large enough $l \in \mathbb{N}$.
As a result, we apply Proposition \ref{Prop.weakconv} with $G=\operatorname{int}B_{2r}(x)$, then take into account the above observations, and, finally, use
the uniformly locally doubling property of the measure $\mu$. This gives
\begin{equation}
\begin{split}
&\mathfrak{m}_{k}(B_{r}(x)) \le \mathfrak{m}_{k}(\operatorname{int}B_{2r}(x)) \le \varliminf\limits_{l \to \infty}\mathfrak{m}^{j_{l},k}(\operatorname{int}B_{2r}(x))\\
&\le \varliminf\limits_{l \to \infty}\sum\{\mathfrak{m}^{j_{l},k}(\widetilde{Q}_{k(r),\alpha}):\widetilde{Q}_{k(r),\alpha}  \cap B_{2r}(x) \neq \emptyset\} \le C \frac{\mu(B_{cr}(x))}{r^{\theta}} \le C \frac{\mu(B_{r}(x))}{r^{\theta}}.
\end{split}
\end{equation}
Hence, condition $(\textbf{M}2)$ is satisfied.

\textit{Step 7.}
To verify condition $(\textbf{M}3)$, it is sufficient to show that there is a constant $C > 0$ such that (we put $B_{k}(x):=B_{\epsilon^{k}}(x)$ for brevity)
\begin{equation}
\label{eqq.417''}
\mathfrak{m}_{k}(B_{k}(x)) \geq C \frac{\mu(B_{k}(x))}{\epsilon^{k\theta}} \quad \text{for all} \quad k \in \mathbb{N}_{0} \quad \text{and all} \quad x \in S.
\end{equation}
Indeed, assume that we have already proved \eqref{eqq.417''}. Then, given $k \in \mathbb{N}_{0}$ and $r \in [\epsilon^{k},1]$, we note that $k(r) \le k$ according to notation \eqref{eqq.important_index_notation}.
Hence, using  (\textbf{M.4}) and, taking into account
the uniformly locally doubling property of $\mu$, we get, for each $x \in S$, the required estimate
\begin{equation}
\notag
\mathfrak{m}_{k}(B_{r}(x)) \geq \mathfrak{m}_{k(r)}(B_{k(r)+1}(x)) \geq \frac{1}{\epsilon}\mathfrak{m}_{k(r)+1}(B_{k(r)+1}(x)) \geq C \frac{\mu(B_{k(r)+1}(x))}{\epsilon^{(k(r)+1)\theta}} \geq C \frac{\mu(B_{r}(x))}{r^{\theta}}.
\end{equation}

To prove \eqref{eqq.417''} we fix an arbitrary $k \in \mathbb{N}_{0}$ and $x \in S$. Using the subadditivity property of the set function $\mathcal{H}_{\theta,\epsilon^{k}}$ and \eqref{eqq.intersection'}
we find $\widetilde{Q}_{k,\alpha}$, $\alpha \in \mathcal{A}_{k}(\operatorname{S},\epsilon)$ such that
\begin{equation}
\label{eqq.544''}
\mathcal{H}_{\theta,\epsilon^{k}}(\widetilde{Q}_{k,\alpha} \cap S) \geq \frac{1}{C_{D}(1,0)}\mathcal{H}_{\theta,\epsilon^{k}}(B_{\epsilon^{k}}(x) \cap S).
\end{equation}
Note that, for each $j \geq k$ and any $z_{j,\beta} \preceq z_{k,\alpha}$ with $\beta \in \mathcal{A}_{j}(\operatorname{S},\epsilon)$,
there is a \textit{maximal} among all numbers $s \in \{k,...,j\}$ for which there exists $\gamma \in \mathcal{A}_{s}(\operatorname{S},\epsilon)$ such that
$z_{j,\beta} \preceq z_{s,\gamma} \preceq z_{k,\alpha}$ and $\mathfrak{m}^{j,k}(\widetilde{Q}_{s,\gamma} \cap S)=h_{s,\gamma}$.
Thus, there exists a disjoint finite family $\{\widetilde{Q}_{s_{i},\gamma_{i}}\}_{i=1}^{N}$
with $k \le s_{i} \le j$, $i \in \{1,...,N\}$ such that $\mathfrak{m}^{j,k}(\widetilde{Q}_{k,\alpha} \cap S)
\geq \sum_{i=1}^{N}\mathfrak{m}^{j,k}(\widetilde{Q}_{s_{i},\gamma_{i}} \cap S)=\sum_{i=1}^{N} h_{s_{i},\gamma_{i}}$.
At the same time, by \eqref{eqq.pseudo_dyadic_1}, \eqref{eqq.pseudo_dyadic_4} it is clear that $\widetilde{Q}_{k,\alpha} \cap S \subset \cup_{i=1}^{N} \operatorname{cl}\widetilde{Q}_{s_{i},\gamma_{i}} \subset \cup_{i=1}^{N} B_{2\epsilon^{s_{i}}}(z_{s_{i},\gamma_{i}})$
and $\frac{1}{8}B_{\epsilon^{s_{i}}}(z_{s_{i},\gamma_{i}}) \subset \widetilde{Q}_{s_{i},\gamma_{i}}$ for all $i \in \{1,...,N\}$. Consequently, using the uniformly locally doubling property of the measure $\mu$ in combination with Proposition \ref{Prop.finite_intersection}, we obtain
\begin{equation}
\label{eqq.545''}
\begin{split}
&\mathfrak{m}^{j,k}(\widetilde{Q}_{k,\alpha} \cap S) \geq \sum_{i=1}^{N}
\frac{\mu(\frac{1}{8}B_{\epsilon^{s_{i}}}(z_{s_{i},\gamma_{i}}))}{\epsilon^{s_{i}\theta}} \geq C\sum_{i=1}^{N} \frac{\mu(4B_{\epsilon^{s_{i}}}(z_{s_{i},\gamma_{i}}))}{\epsilon^{s_{i}\theta}}
\\
&\geq C\sum_{i=1}^{N}\epsilon^{-s_{i}\theta}\sum\limits_{B \in \mathcal{B}_{s_{i}}(\operatorname{X},\epsilon)}\{\mu(B): B \cap B_{2\epsilon^{s_{i}}}(z_{s_{i},\gamma_{i}}) \neq \emptyset\}
\geq C\mathcal{H}_{\theta,\epsilon^{k}}(\widetilde{Q}_{k,\alpha} \cap S).
\end{split}
\end{equation}
Since the closed ball $B_{5\epsilon^{k}}(x)$ is a compact set and $\widetilde{Q}_{k,\alpha} \subset B_{5\epsilon^{k}}(x)$, we apply
Proposition \ref{Prop.weakconv} with $F=B_{5\epsilon^{k}}(x)$, then combine \eqref{eqq.544''}, \eqref{eqq.545''}, and, finally, take into account Definition \ref{Def.content}. This gives us the crucial estimate
\begin{equation}
\label{eqq.432'''}
\mathfrak{m}_{k}(B_{5\epsilon^{k}}(x)) \geq \varlimsup\limits_{l \to \infty}\mathfrak{m}^{j_{l},k}(B_{5\epsilon^{k}}(x)) \geq \varlimsup\limits_{l \to \infty}\mathfrak{m}^{j_{l},k}(\widetilde{Q}_{k,\alpha}) \geq C\frac{\mu(B_{\epsilon^{k}}(x))}{\epsilon^{k\theta}}.
\end{equation}
As a result, using \eqref{eqq.432'''} and Theorem \ref{Th.doubling_type} we arrive at \eqref{eqq.417''} completing the proof of (\textbf{M}3).

\textit{Step 8.} By (\textbf{M}4) and (\textbf{M}3), which were verified at steps 5 and 7, respectively, we have $\mathfrak{m}_{0}(B_{k}(x)) \geq \epsilon^{k\theta}\mathfrak{m}_{k}(B_{k}(x)) \geq C\mu(B_{k}(x))$
for all $x \in S$.
By Definition \ref{Def.weaknoncol} this implies that the measure $\mathfrak{m}_{0}$ is weakly noncollapsed.
We fix an arbitrary Borel set $E \subset S$ and recall notation \eqref{eq.lower_and_upper_densities}. Throughout this step we set $c=4/\epsilon+2$ and apply Lemma \ref{prop.density}.
This gives us the existence of a set $E' \subset E$ with $\mathfrak{m}_{0}(E \setminus E')=0$ such that, for each point $x \in E'$, one can find a strictly decreasing to zero sequence $\{r_{l}(x)\}$ such that
(we recall \eqref{eq.doubling.sequence})
\begin{equation}
\notag
\overline{D}^{\mathfrak{m}_{0}}_{E}(x) \geq \overline{\operatorname{D}}(x):=\varliminf\limits_{l \to \infty}\frac{\mathfrak{m}_{0}(B_{r_{l}(x)}(x)\cap E)}{\mathfrak{m}_{0}(B_{r_{l}(x)}(x))} \geq \frac{1}{2\underline{\operatorname{C}}_{\mathfrak{m}_{0}}(5c)} \quad \text{and} \quad \varlimsup\limits_{l \to \infty}\frac{\mathfrak{m}_{0}(B_{cr_{l}(x)}(x))}{\mathfrak{m}_{0}(B_{r_{l}(x)}(x))} \le \underline{\operatorname{C}}_{\mathfrak{m}_{0}}(5c).
\end{equation}
Furthermore, we fix $\underline{x} \in E'$ and $\varepsilon \in (0,1)$. We recall notation \eqref{eqq.important_index_notation} and
put $r_{l}:=r_{l}(\underline{x})$, $k_{l}:=k_{\epsilon}(r_{l})$ for all $l \in \mathbb{N}_{0}$.
Clearly, there is $L=L(\underline{x},\varepsilon) \in \mathbb{N}$ such that, for all $l \geq L$,
\begin{equation}
\label{eqq.418'}
\frac{\mathfrak{m}_{0}(B_{r_{l}}(\underline{x})\cap E)}{\mathfrak{m}_{0}(B_{r_{l}}(\underline{x}))} > (1-\frac{\varepsilon}{8})\overline{\operatorname{D}}(\underline{x}) \quad \text{and} \quad \frac{\mathfrak{m}_{0}(B_{cr_{l}}(\underline{x}))}{\mathfrak{m}_{0}(B_{r_{l}}(\underline{x}))} \le 2\underline{\operatorname{C}}_{\mathfrak{m}_{0}}(5c).
\end{equation}

Using Borel regularity of the measure $\mathfrak{m}_{0}$, given $l \in \mathbb{N}$, we find
an open set $\Omega_{l} \subset B_{2r_{l}}(\underline{x})$ containing $B_{r_{l}}(\underline{x})\cap E$ and
a compact set $K_{l} \subset B_{r_{l}}(\underline{x})\cap E$ such that
\begin{equation}
\label{eq.opencompact}
|\mathfrak{m}_{0}(\Omega_{l})-\mathfrak{m}_{0}(K_{l})| \le \frac{\varepsilon}{8}\overline{\operatorname{D}}(\underline{x})\mathfrak{m}_{0}(B_{r_{l}}(\underline{x})).
\end{equation}
Since $\sigma_{l}:=\operatorname{dist}(K_{l},\operatorname{X} \setminus \Omega_{l}) > 0$ for the $\frac{\sigma_{l}}{2}$-neighborhood $U_{\frac{\sigma_{l}}{2}}(K_{l})$ of the set $K_{l}$, we get
\begin{equation}
\label{eq.neighborhincl}
K_{l} \subset U_{\frac{\sigma_{l}}{2}}(K_{l}) \subset \operatorname{cl}U_{\frac{\sigma_{l}}{2}}(K_{l}) \subset \Omega_{l} \subset B_{2r_{l}}(\underline{x}).
\end{equation}

To verify (\textbf{M}5) it is sufficient to establish existence of $C(\underline{x}) > 0$ independent on $\varepsilon$ and $l$ such that for each $l \geq L$ there is
$N=N(\underline{x},l,\varepsilon) \geq k_{l}$ such that  for any $j \geq N$
\begin{equation}
\label{eqq.the_crucial_observation}
\mathfrak{m}^{j,k_{l}}(\widetilde{Q}_{k_{l},\beta(j)} \cap U_{\frac{\sigma_{l}}{2}}(K_{l})) \geq C(\underline{x})\mathfrak{m}^{j,k_{l}}(\widetilde{Q}_{k_{l},\beta(j)}) \quad
\text{for some} \quad \beta(j) \in \mathcal{A}_{k_{l}}(\operatorname{S},\epsilon).
\end{equation}
Indeed, suppose that we have already proved \eqref{eqq.the_crucial_observation}. Then, given $l \geq L$, we use \eqref{eqq.pseudo_dyadic_1}, \eqref{eq.neighborhincl} and apply Proposition \ref{Prop.weakconv} with
$F = \operatorname{cl}(U_{\frac{\sigma}{2}}(K_{l}))$, $G=\operatorname{int}B_{\frac{r_{l}}{8}}(z_{k_{l},\beta})$. This gives
\begin{equation}
\begin{split}
&\mathfrak{m}_{k_{l}}(\Omega_{l}) \geq \mathfrak{m}_{k_{l}}(\operatorname{cl} U_{\frac{\sigma_{l}}{2}}(K_{l})) \geq \varlimsup\limits_{j \to \infty}\mathfrak{m}^{j,k_{l}}(\operatorname{cl} U_{\frac{\sigma_{l}}{2}}(K_{l})) \geq
\varliminf\limits_{j \to \infty}
\mathfrak{m}^{j,k_{l}}(U_{\frac{\sigma_{l}}{2}}(K_{l}))\\
&\geq \varliminf\limits_{j \to \infty}\mathfrak{m}^{j,k_{l}}(\widetilde{Q}_{k_{l},\beta(j)} \cap U_{\frac{\sigma_{l}}{2}}(K_{l})) \geq C(\underline{x})\varliminf\limits_{j \to \infty}\mathfrak{m}^{j,k_{l}}(\widetilde{Q}_{k_{l},\beta(j)})\\
&\geq C(\underline{x})\varliminf\limits_{j \to \infty}\mathfrak{m}^{j,k_{l}}(\operatorname{int}\frac{1}{8}B_{k_{l}}(z_{k_{l},\beta(j)})) \geq C(\underline{x})\mathfrak{m}_{k_{l}}(\operatorname{int}\frac{1}{8}B_{k_{l}}(z_{k_{l},\beta(j)})).
\end{split}
\end{equation}
Since $\widetilde{Q}_{k_{l},\beta(j)} \cap B_{2r_{l}}(\underline{x}) \neq \emptyset$, from \eqref{eqq.pseudo_dyadic_4} we have $B_{r_{l}}(\underline{x}) \subset 6B_{k_{l}}(z_{k_{l},\beta(j)})$. The crucial fact is that $z_{k_{l},\beta(j)} \in S$ and we can apply Theorem \ref{Th.doubling_type}
with $c=6$ and $y=z_{k_{l},\beta(j)}$. As a result, 
\begin{equation}
\label{eq.2.26}
\mathfrak{m}_{k_{l}}(\Omega_{l}) \geq C(\underline{x})\mathfrak{m}_{k_{l}}(\operatorname{int}B_{\frac{r_{l}}{8}}(z_{k_{l},\beta(j)})) \geq C
\mathfrak{m}_{k_{l}}(6B_{k_{l}}(z_{k_{l},\beta(j)})) \geq C\mathfrak{m}_{k_{l}}(B_{r_{l}}(\underline{x})),
\end{equation}
where the constant $C > 0$ in the right-hand side of \eqref{eq.2.26} does not depend on $l$ and $\varepsilon$.
Finally, since $\Omega_{l} \supset E \cap B_{r_{l}}(\underline{x})$ was chosen arbitrarily and since the measure $\mathfrak{m}_{k_{l}}$ is Borel regular, we obtain
the required estimate $\mathfrak{m}_{k_{l}}(E \cap B_{r_{l}}(\underline{x})) \geq C\mathfrak{m}_{k_{l}}(B_{r_{l}}(\underline{x}))$ with the
constant $C > 0$  independent on $l \in \mathbb{N}$.
Since $l \geq L$ was chosen arbitrarily, this verifies $(\textbf{M}5)$.

To prove \eqref{eqq.the_crucial_observation} we proceed as follows. We fix $l \geq L$ and apply Proposition \ref{Prop.weakconv} with $G=U_{\frac{\sigma}{2}}(K_{l})$, $F=B_{2r_{l}}(\underline{x})$ and use \eqref{eqq.418'}--\eqref{eq.neighborhincl}. This gives (recall  that $\varepsilon \in (0,1)$)
\begin{equation}
\notag
\begin{split}
&\varliminf\limits_{j \to \infty}\mathfrak{m}^{j,0}(U_{\frac{\sigma}{2}}(K_{l})) \geq \mathfrak{m}_{0}(U_{\frac{\sigma}{2}}(K_{l}))  \geq (1-\frac{\varepsilon}{2})\overline{\operatorname{D}}(\underline{x})\mathfrak{m}_{0}(B_{r_{l}}(\underline{x}))\\
& \geq C\mathfrak{m}_{0}(B_{cr_{l}}(\underline{x})) \geq C\varlimsup\limits_{j \to \infty}\mathfrak{m}^{j,0}(B_{cr_{l}}(\underline{x})),
\end{split}
\end{equation}
where $C=\overline{\operatorname{D}}(\underline{x})/4\underline{\operatorname{C}}_{\mathfrak{m}_{0}}(5c)$. Hence, there exists $N=N(\underline{x},l,\varepsilon) \geq k_{l}$ such that, for all $j \geq N$,
\begin{equation}
\label{eqq.552''}
\begin{split}
&\mathfrak{m}^{j,0}(U_{\frac{\sigma_{l}}{2}}(K_{l})) \geq \frac{\overline{\operatorname{D}}(\underline{x})}{5\underline{\operatorname{C}}_{\mathfrak{m}_{0}}(5c)}\mathfrak{m}^{j,0}(B_{cr_{l}}(\underline{x})).
\end{split}
\end{equation}
By \eqref{eqq.intersection'} and \eqref{eq.neighborhincl}
we see that there are at most $C_{D}(2,0)$ generalized dyadic cubes $\widetilde{Q}_{k_{l},\beta} \cap S$ in $\operatorname{S}$
whose closures have nonempty intersections with $U_{\frac{\sigma_{l}}{2}}(K_{l})$. Furthermore, any such a
cube is contained with its closure in $B_{cr_{l}}(\underline{x})$. Hence, using \eqref{eqq.552''} we conclude that for each $j \geq N$
there exists  $\beta(j) \in \mathcal{A}_{k_{l}}(\operatorname{S},\epsilon)$ such that the following inequality
\begin{equation}
\label{eq.2.23}
\begin{split}
&\mathfrak{m}^{j,0}(\widetilde{Q}_{k_{l},\beta(j)} \cap U_{\frac{\sigma_{l}}{2}}(K_{l})) \geq C(\underline{x})\mathfrak{m}^{j,0}(B_{cr_{l}}(\underline{x}))
\geq C(\underline{x})\mathfrak{m}^{j,0}(\widetilde{Q}_{k_{l},\beta(j)})
\end{split}
\end{equation}
holds with the constant $C(\underline{x}):=(5C_{D}(2,0)\underline{\operatorname{C}}_{\mathfrak{m}_{0}}(5c))^{-1}\overline{\operatorname{D}}(\underline{x}).$
As a result, taking into account \eqref{eqq.positive_constants} we deduce from \eqref{eq.2.23} the required estimate
\begin{equation}
\notag
\begin{split}
&\mathfrak{m}^{j,k_{l}}(\widetilde{Q}_{k_{l},\beta(j)} \cap U_{\frac{\sigma_{l}}{2}}(K_{l}))=\frac{\mathfrak{m}^{j,k_{l}}(\widetilde{Q}_{k_{l},\beta(j)})}{\mathfrak{m}^{j,0}(\widetilde{Q}_{k_{l},\beta(j)})} \mathfrak{m}^{j,0}(\widetilde{Q}_{k_{l},\beta(j)} \cap U_{\frac{\sigma_{l}}{2}}(K_{l}))\\ &\geq C(\underline{x})\frac{\mathfrak{m}^{j,k_{l}}(\widetilde{Q}_{k_{l},\beta(j)})}{\mathfrak{m}^{j,0}(\widetilde{Q}_{k_{l},\beta(j)})}\mathfrak{m}^{j,0}(\widetilde{Q}_{k_{l},\beta(j)})
=C(\underline{x})\mathfrak{m}^{j,k_{l}}(\widetilde{Q}_{k_{l},\beta(j)}).
\end{split}
\end{equation}
The proof is complete.
\end{proof}

\textit{Proof of Theorem \ref{Th.regularsequence}.} The required result is a combination of Theorems \ref{Th.necessary_regular_sequence} and \ref{Th.sufficiency_regular_sequence}.
\hfil$\Box$

\subsection{Some examples} In this subsection we show that, for some sets $S \in \mathcal{LCR}_{\theta}(\operatorname{X})$, $\theta > 0$, one can easily
built concrete examples of sequences $\{\mathfrak{m}_{k}\} \in \mathfrak{M}_{\theta}(S)$. For generic sets $S \in \mathcal{LCR}_{\theta}(\operatorname{X})$ with $\theta > 0$, finding
explicit examples of sequences $\{\mathfrak{m}_{k}\} \in \mathfrak{M}_{\theta}(S)$ is a quite sophisticated problem.
In \cite{T2} an explicit example of $\{\mathfrak{m}_{k}\} \in \mathfrak{M}_{1}(\Gamma)$ was constructed for the case when $\Gamma \subset \mathbb{R}^{2}$
is a simple planar rectifiable curve of positive length.  In \cite{TV1} an explicit example of $\{\mathfrak{m}_{k}\} \in \mathfrak{M}_{n-1}(K)$ was given
for the case of a single cusp $K$ in $\mathbb{R}^{n}$. Using Theorem \ref{Th.coincideness_of_classes_of_measures} we conclude that in fact the sequences of measures introduced in \cite{T2} and \cite{TV1} belong
to the more narrow classes $\mathfrak{M}^{str}_{1}(\Gamma)$ and $\mathfrak{M}^{str}_{n-1}(K)$, respectively.

\textbf{Example 5.1.} We recall Remark \ref{Rem.Hausdorff_is_a_Borel_measure}. Let $\underline{\theta} \in [0,\underline{Q}_{\mu}(R))$ for some $R > 0$
and let $S \in \mathcal{ADR}_{\underline{\theta}}(\operatorname{X})$. In this case, given $\theta \geq \underline{\theta}$, we put $\epsilon=\frac{1}{2}$
and define
\begin{equation}
\notag
\mathfrak{m}_{k}:=2^{k(\theta-\underline{\theta})}\mathcal{H}_{\underline{\theta}}\lfloor_{S}, \quad k \in \mathbb{N}_{0}.
\end{equation}
It is easy to verify that the sequence of measures $\{\mathfrak{m}_{k}\}:=\{\mathfrak{m}_{k}\}_{k=0}^{\infty}$ lies in $\mathfrak{M}_{\theta}^{str}(S)$. Indeed, conditions (\textbf{M}1)--(\textbf{M}4)
follows immediately by the construction. To verify (\textbf{M}5) one should repeat with minor technical modifications (keeping in mind \eqref{eqq.ADR}) the corresponding arguments from the proof of Theorem 6.2 in \cite{Mat}.

\hfill$\Box$

\textbf{Example 5.2.}
Let $N \in \mathbb{N}$ and $\{\underline{\theta}_{1},...,\underline{\theta}_{N}\} \subset [0,\underline{Q}_{\mu}(R))$. Given $i \in \{1,...,N\}$, let $S_{i} \in  \mathcal{ADR}_{\underline{\theta}_{i}}(\operatorname{X})$. We put $\overline{\theta}:=\max\{\underline{\theta}_{1},...,\underline{\theta}_{N}\}$.
Now we put $\epsilon=\frac{1}{2}$ and, given $\theta \geq \overline{\theta}$, we define
\begin{equation}
\notag
\mathfrak{m}_{k}:=\sum\limits_{i=1}^{N}2^{k(\theta-\underline{\theta}_{i})}\mathcal{H}_{\underline{\theta}_{i}}\lfloor_{S_{i}},  \quad k \in \mathbb{N}_{0}.
\end{equation}
Based on the previous example, we immediately have $\{\mathfrak{m}_{k}\}:=\{\mathfrak{m}_{k}\}_{k=0}^{\infty} \in \mathfrak{M}_{\theta}^{str}(S)$.

\hfill$\Box$

\section{Lebesgue points of functions}

\textit{Throughout this section we fix the following data}:

\begin{itemize}
\item[\(\rm (\textbf{D.6.1})\)] a parameter $p \in (1,\infty)$ and a metric measure space $\operatorname{X}=(\operatorname{X},\operatorname{d},\mu) \in \mathfrak{A}_{p}$;

\item[\(\rm (\textbf{D.6.2})\)] a parameter $\theta \in [0,p)$ and a closed set $S \in \mathcal{LCR}_{\theta}(\operatorname{X})$;

\item[\(\rm (\textbf{D.6.3})\)] a sequence of measures $\{\mathfrak{m}_{k}\}:=\{\mathfrak{m}_{k}\}_{k=0}^{\infty} \in \mathfrak{M}^{str}_{\theta}(S)$ with $\epsilon=\epsilon(\{\mathfrak{m}_{k}\}) \in (0,\frac{1}{10}]$.
\end{itemize}

In this section, for any $x \in \operatorname{X}$ and $k \in \mathbb{Z}$, we will use \textit{notation} $B_{k}(x):=B_{\epsilon^{k}}(x)$.

\begin{Def}
\label{Def.sc_good}
Given $c > 0$ and $\delta > 0$, we say that a family of closed balls $\mathcal{B}:=\{B_{r_{i}}(x_{i})\}_{i=1}^{N}$ with $N \in \mathbb{N}$ is $(S,c,\delta)$-nice if the following conditions hold:

\begin{itemize}
\item[\((\textbf{B1})\)] $B_{r_{i}}(x_{i}) \cap B_{r_{j}}(x_{j}) = \emptyset$ if $i,j \in \{1,...,N\}$ and $i \neq j$;

\item[\((\textbf{B2})\)] $\sup\{r_{i}:i=1,...,N\} \le \delta$;

\item[\((\textbf{B3})\)] $B_{cr_{i}}(x_{i}) \cap S \neq \emptyset$ for all $i \in \{0,...,N\}$.
\end{itemize}

Furthermore, we say that $\mathcal{B}$ is an $(S,c,\delta)$-Whitney family if it satisfies (\textbf{B1})-(\textbf{B3}) and

\begin{itemize}
\item[\((\textbf{B4})\)] $B \subset \operatorname{X} \setminus S$ for all $B \in \mathcal{B}$.
\end{itemize}

We will call $(S,c,1)$-nice families and $(S,c,1)$-Whitney families just $(S,c)$-nice families and $(S,c)$-Whitney families, respectively.
\end{Def}

\begin{Remark}
\label{Rem.31'}
Given $\delta \in (0,1]$, $c \geq 1$ every $(S,c,\delta)$-Whitney family is an $(S,c',\delta')$-Whitney family and every $(S,c,\delta)$-nice family is an $(S,c',\delta')$-nice family
for any $\delta' \in [\delta,1]$ and $c' \geq c$.
\end{Remark}
We recall notation \eqref{eqq.important_index_notation} and, given a ball $B=B_{r}(x)$, we put $k(B):=k(r(B))$. Furthermore, we recall the notation given
at the beginning of Section 5.

\begin{Prop}
\label{Prop.niceball_estimate}
Let $c \geq 1$ and $c' \geq c+1$. If a closed ball $B=B_{r}(x)$ in $\operatorname{X}$ is such that $r \in (0,1]$ and $cB \cap S \neq \emptyset$, then
\begin{equation}
\frac{\mu(B)}{\mathfrak{m}_{k(B)}(c'B)} \le \frac{(C_{\mu}(c'))^{\log_{2}2c'+1}}{\epsilon^{\theta}}\frac{C_{\{\mathfrak{m}_{k}\},3}}{C_{\{\mathfrak{m}_{k}\},2}} (r(B))^{\theta}.
\end{equation}
\end{Prop}

\begin{proof}
Note that there exists a ball $\underline{B} \subset c'B$
such that $r(\underline{B})=r(B)$ and the center $\underline{x}$ of $\underline{B}$ belongs to $S$.
In this case, we have $B \subset 2c'\underline{B}$. Furthermore,
according to our notation $\epsilon^{k(B)+1} < r(B) \le \epsilon^{k(B)}$.
Hence, by \eqref{M.3}, \eqref{M.4} and the uniformly locally doubling property of $\mu$ (we let $[c]:=\max\{k\in\mathbb{Z}:k \le c\}$),
\begin{equation}
\begin{split}
&\frac{\mu(B)}{\mathfrak{m}_{k(B)}(2cB)} \le \frac{\mu(2c'\underline{B})}{\mathfrak{m}_{k(B)}(\underline{B})}
\le \frac{(C_{\mu}(c'))^{[\log_{2}2c']+1}}{\epsilon^{\theta}}\frac{C_{\{\mathfrak{m}_{k}\},3}}{C_{\{\mathfrak{m}_{k}\},2}}(r(B))^{\theta}.
\end{split}
\end{equation}
This completes the proof.
\end{proof}

In this and subsequent sections we will need a Brudnyi--Shvartsman functional on ``small scales''. More precisely we recall \eqref{eqq.tricky_average} and formulate the following concept.
\begin{Def}
\label{Def.Br_Shv}
Given $\delta \in (0,1]$ and $c > 1$, we introduce the $\delta$-scale Brudnyi--Shvartsman functional on $L_{1}^{loc}(\{\mathfrak{m}_{k}\})$ (it takes values in $[0,+\infty]$) by letting
\begin{equation}
\begin{split}
\label{eq.main2''}
\operatorname{BSN}^{\delta}_{p,\{\mathfrak{m}_{k}\},c}(f):=\|f|L_{p}(\mathfrak{m}_{0})\|+
\sup&\Bigl(\sum\limits_{B \in \mathcal{B}^{\delta}} \frac{\mu(B)}{(r(B))^{p}}\Bigl(\widetilde{\mathcal{E}}_{\mathfrak{m}_{k(B)}}(f,cB)\Bigr)^{p}\Bigr)^{\frac{1}{p}},
\end{split}
\end{equation}
where the supremum is taken over all families $(S,c,\delta)$-nice families $\mathcal{B}^{\delta}$.
\end{Def}

\begin{Remark}
Keeping in mind \eqref{eq.main2} and notation used in Theorems \ref{Th.SecondMain}, \ref{Th.FourthMain},
in the case $\delta = 1$ we will write
$\operatorname{BSN}_{p,\{\mathfrak{m}_{k}\},c}(f)$ instead of $\operatorname{BSN}^{1}_{p,\{\mathfrak{m}_{k}\},c}(f)$.
\end{Remark}

\begin{Lm}
\label{Lm.largescale_estimate}
Let $\delta \in (0,1]$ and $c \geq 1$. Then there is a constant $C > 0$ depending only on $\delta$, $\mathcal{C}_{\{\mathfrak{m}_{k}\}}$, $c$, $\epsilon$, $\theta$ and $C_{\mu}(2c)$ such that if $\mathcal{B}_{\delta}$ is an arbitrary $(S,c)$-nice family of balls such that $r(B) \geq \delta$
for all $B \in \mathcal{B}_{\delta}$, then, for each $f \in L_{p}(\mathfrak{m}_{0})$,
\begin{equation}
\sum\limits_{B \in \mathcal{B}_{\delta}} \frac{\mu(B)}{(r(B))^{p}}\Bigl(\mathcal{E}_{\mathfrak{m}_{k(B)}}(f,2cB)\Bigr)^{p}
\le C \int\limits_{S}|f(x)|^{p}\,d\mathfrak{m}_{0}(x).
\end{equation}
\end{Lm}

\begin{proof}
We fix $f \in L_{p}(\mathfrak{m}_{0})$ and an $(S,c)$-nice family $\mathcal{B}_{\delta}$ such that $r(B) \geq \delta$
for all $B \in \mathcal{B}_{\delta}$. Let $\underline{k}$ be the largest integer $k$ satisfying the inequality $\epsilon^{k} \geq \delta$.  
Below we will write explicitly all intermediate constants to indicate their dependence on $\underline{k}$ (and hence on $\delta$).
By Proposition \ref{Prop.finite_intersection} we obtain (we take into account that $N_{\mu}(\epsilon^{k},C) \le N_{\mu}(1,C)$ for each $k \in \{0,...,\underline{k}\}$ and $C > 0$)
\begin{equation}
\label{eq.3.covering_multiplicity}
\mathcal{M}(\{2cB:B \in \mathcal{B}_{\delta}\}) \le \sum\limits_{k=0}^{\underline{k}}\mathcal{M}(\{2cB:B \in \mathcal{B}_{\delta}(k,\epsilon)\}) \le \underline{k}N_{\mu}(1,\frac{4c}{\epsilon}).
\end{equation}
Given $k \in \{0,...,\underline{k}\}$, an application of Proposition \ref{Prop.rough_estimate} with $\mathfrak{m}=\mathfrak{m}_{k}$
gives
\begin{equation}
\notag
\begin{split}
&\Bigl(\mathcal{E}_{\mathfrak{m}_{k}}(f,2cB)\Bigr)^{p} \le
2^{p} \fint\limits_{2cB}|f(z)|^{p}\,d\mathfrak{m}_{k}(z) \quad \text{for any} \quad B \in \mathcal{B}_{\delta}(k,\epsilon).
\end{split}
\end{equation}
Hence, by Proposition \ref{Prop.niceball_estimate} with $c' =2c$ and \eqref{M.4} we get, for any $k \in \{0,...,\underline{k}\}$ and $B \in \mathcal{B}_{\delta}(k,\epsilon)$,
\begin{equation}
\label{eq.3.rough_estimate}
\begin{split}
&\mu(B)\Bigl(\mathcal{E}_{\mathfrak{m}_{k}}(f,2cB)\Bigr)^{p} \le 2^{p}\frac{C_{\{\mathfrak{m}_{k}\},3}}{\epsilon^{\underline{k}\theta}}\frac{\mu(B)}{\mathfrak{m}_{k}(2cB)}\int\limits_{2cB}|f(z)|^{p}\,d\mathfrak{m}_{0}(z)\\
&\le
2^{p}\frac{(C_{\mu}(2c))^{\log_{2}4c+1}}{\epsilon^{(\underline{k}+1)\theta}}\frac{(C_{\{\mathfrak{m}_{k}\},3})^{2}}{C_{\{\mathfrak{m}_{k}\},2}}(r(B))^{\theta}\int\limits_{2cB}|f(z)|^{p}\,d\mathfrak{m}_{0}(z).
\end{split}
\end{equation}
Note that the right-hand side of \eqref{eq.3.rough_estimate} depends on $\underline{k}$ but does not depend on $k \in \{0,...,\underline{k}\}$.
Consequently, using Proposition \ref{Prop.covering_multiplicity} and  \eqref{eq.3.covering_multiplicity}, \eqref{eq.3.rough_estimate} we obtain (recall that $\theta \in [0,p)$)
\begin{equation}
\notag
\begin{split}
&\sum\limits_{B \in \mathcal{B}_{\delta}} \frac{\mu(B)}{(r(B))^{p}}\Bigl(\mathcal{E}_{\mathfrak{m}_{k(B)}}(f,2cB)\Bigr)^{p}
\le \sum\limits_{B \in \mathcal{B}_{\delta}} \frac{2^{p}}{\delta^{p-\theta}}\frac{(C_{\mu}(2c))^{\log_{2}4c+1}}{\epsilon^{(\underline{k}+1)\theta}}\frac{(C_{\{\mathfrak{m}_{k}\},3})^{2}}{C_{\{\mathfrak{m}_{k}\},2}}
\int\limits_{2cB}|f(z)|^{p}\,d\mathfrak{m}_{0}(z)\\
&\le \underline{k}N_{\mu}(1,\frac{4c}{\epsilon})\frac{2^{p}}{\delta^{p}}\frac{(C_{\mu}(2c))^{\log_{2}4c+1}}{\epsilon^{\theta}}
\frac{(C_{\{\mathfrak{m}_{k}\},3})^{2}}{C_{\{\mathfrak{m}_{k}\},2}}\int\limits_{S}|f(x)|^{p}\,d\mathfrak{m}_{0}(x).
\end{split}
\end{equation}

The proof is complete.

\end{proof}

It is natural to ask whether the finiteness of $\operatorname{BSN}^{\delta}_{p,\{\mathfrak{m}_{k}\},c}(f)$ for small $\delta$ implies that of $\operatorname{BSN}_{p,\{\mathfrak{m}_{k}\},c}(f)$. Fortunately, we have the affirmative answer.

\begin{Lm}
\label{Lm.different_scales}
$\operatorname{BSN}_{p,\{\mathfrak{m}_{k}\},c}(f) < +\infty$ if and only $\operatorname{BSN}^{\delta}_{p,\{\mathfrak{m}_{k}\},c}(f) < +\infty$ for some $\delta \in (0,1]$.
\end{Lm}

\begin{proof}
The necessity follows from Remark \ref{Rem.31'}.
To prove the sufficiency, given an $(S,c)$-nice family $\mathcal{B}$, we split it into two subfamilies. More precisely, we put
$\mathcal{B}^{\delta}:=\{B \in \mathcal{B}:r(B) \le \delta\}$ and $\mathcal{B}_{\delta}:=\mathcal{B} \setminus \mathcal{B}^{\delta}$.
Now the claim follows from Lemma \ref{Lm.largescale_estimate}.
\end{proof}

We start with the following lemma (we will use notation $B_{k}(x):=B_{\epsilon^{k}}(x)$, $k \in \mathbb{N}_{0}$, $x \in \operatorname{X}$).

\begin{Lm}
\label{Lm.quasiLebesgue}
Let $f \in L_{1}^{loc}(\{\mathfrak{m}_{k}\})$ be such that $\operatorname{BSN}^{\delta}_{p,\{\mathfrak{m}_{k}\},c}(f) < +\infty$ for some $c > 1$ and $\delta \in (0,1]$. Then there exists a Borel function $\overline{f}: S \to \mathbb{R}$ and a Borel set $\underline{S} \subset S$ with $\mathcal{H}_{\theta}(S \setminus \underline{S})=0$ such that
\begin{equation}
\label{eq.quasiLebesgue}
\varlimsup\limits_{k \to \infty}\fint\limits_{B_{k}(x)}|\overline{f}(x)-f(y)|\,d\mathfrak{m}_{k}(y) = 0 \quad \text{for all}  \quad x \in S.
\end{equation}
\end{Lm}
\begin{proof}
We fix $\varepsilon \in (0,\frac{p-\theta}{2p})$ and split the proof into several steps.

\textit{Step 1.} Consider the function
\begin{equation}
\label{eq.r.function}
R(x):=\varlimsup\limits_{r \to 0}\sum\limits_{\epsilon^{k}< r}\Bigl(\mathcal{E}_{\mathfrak{m}_{k}}(f,B_{k}(x))\Bigr)^{p}, \quad x \in \operatorname{X}.
\end{equation}
It is clear that, given $\delta' \in (0,\delta]$,
\begin{equation}
R(x) \le \varlimsup\limits_{r \to 0}\sum\limits_{\epsilon^{k}< r}\frac{\epsilon^{k\varepsilon p}}{\epsilon^{k\varepsilon p}}\Bigl(\mathcal{E}_{\mathfrak{m}_{k}}(f,B_{k}(x))\Bigr)^{p} \le C(\delta')^{\varepsilon p}
\sup\limits_{\epsilon^{k}<\delta'}\frac{1}{\epsilon^{k\varepsilon p}}\Bigl(\mathcal{E}_{\mathfrak{m}_{k}}(f,B_{k}(x))\Bigr)^{p}.
\end{equation}

Given $t > 0$, we introduce the $t$-superlevel set of $R$ by letting $E_{t}:=\{x \in S:R(x) \geq t\}$.
Our aim is to show that
\begin{equation}
\label{eq.levels}
\mathcal{H}_{\theta}(E_{t})=0 \quad \text{for all} \quad t > 0.
\end{equation}

Now we fix arbitrary $t > 0$ and $\delta' \in (0,\delta]$. For each $x \in E_{t}$ we find $k_{x} \in \mathbb{N}_{0}$ such that
\begin{equation}
\notag
\sup\limits_{\epsilon^{k} < \delta'}\epsilon^{-k\varepsilon}\mathcal{E}_{\mathfrak{m}_{k}}(f,B_{k}(x)) \le 2\epsilon^{-k_{x}\varepsilon}\mathcal{E}_{\mathfrak{m}_{k_{x}}}(f,B_{k_{x}}(x)).
\end{equation}
Clearly, the family of balls $\mathcal{B}:\{B_{k_{x}}(x):x \in E_{t}\}$ is a covering of $E_{t}$. Using the $5B$-covering lemma, we find a
disjoint subfamily $\widetilde{\mathcal{B}} \subset \mathcal{B}$ such that $E_{t} \subset \bigcup \{5B: B \in \widetilde{\mathcal{B}}\}$.
Hence, we have
\begin{equation}
\label{eq.hausdroffrombelow}
\sum\{\frac{\mu(B)}{(r(B))^{\theta}}:B \in \widetilde{\mathcal{B}}\} \geq C\sum\{\frac{\mu(5B)}{(r(5B))^{\theta}}:B \in \widetilde{\mathcal{B}}\} \geq C\mathcal{H}_{\theta,5\delta'}(E_{t}).
\end{equation}
Note that any $(S,c,\delta')$-nice family is also an $(S,c,\delta)$-nice family. Furthermore, by Theorem \ref{Th.doubling_type}
and Remark \ref{Rem.best_approx} it is easy to see that $\mathcal{E}_{\mathfrak{m}_{k(B)}}(f,B) \le C \mathcal{E}_{\mathfrak{m}_{k(B)}}(f,cB)$ for all
$B \in \widetilde{\mathcal{B}}$. As a result, we obtain
\begin{equation}
\begin{split}
&t \mathcal{H}_{\theta,5\delta'}(E_{t}) \le C \sum\limits_{B \in \widetilde{\mathcal{B}}}\frac{\mu(B)}{(r(B))^{\theta}}\frac{(\delta')^{\varepsilon p}}{(r(B))^{\varepsilon p}}\Bigl(\mathcal{E}_{\mathfrak{m}_{k(B)}}(f,B)\Bigr)^{p}\\
&\le C (\delta')^{\varepsilon p}\sum\limits_{B \in \widetilde{\mathcal{B}}}\frac{\mu(B)}{(r(B))^{p}}\Bigl(\mathcal{E}_{\mathfrak{m}_{k(B)}}(f,cB)\Bigr)^{p} \le
C (\delta')^{\varepsilon p} \Bigl(\operatorname{BSN}^{\delta}_{p,\{\mathfrak{m}_{k}\},c}(f)\Bigr)^{p}.
\end{split}
\end{equation}
Passing to the limit as $\delta' \to 0$ and taking into account the arbitrariness of $t > 0$ we get \eqref{eq.levels}.

\textit{Step 2.} Fix $\delta' \in (0,\delta]$. If $l,k \in \mathbb{N}_{0}$
are such that $l > k$ and $2^{-k} < \delta$, then by Remark \ref{Rem.best_approx} and Theorem \ref{Th.doubling_type}, it is easy to see that
\begin{equation}
\label{eq.3.8}
\begin{split}
&\fint\limits_{B_{l}(x)}\fint\limits_{B_{k}(x)}|f(y)-f(z)|\,d\mathfrak{m}_{l}(y)\,d\mathfrak{m}_{k}(z)\\
&\le \sum\limits_{i=k}^{l-1}\fint\limits_{B_{i}(x)}\fint\limits_{B_{i+1}(x)}|f(y)-f(z)|\,d\mathfrak{m}_{i}(z)\,d\mathfrak{m}_{i+1}(y) \le C\sum\limits_{i=k}^{l}\mathcal{E}_{\mathfrak{m}_{i}}(f,B_{i}(x)).
\end{split}
\end{equation}
Consider the set $\underline{S}:=S \setminus \cup_{t > 0}E_{t}$. Since $R(x)=0$ for all $x \in \underline{S}$,
it follows from \eqref{eq.r.function}, \eqref{eq.3.8} that if $x \in \underline{S}$ then
$\{\fint_{B_{k}(x)}f(y)\,d\mathfrak{m}_{k}(y)\}_{k=1}^{\infty}$ is a Cauchy sequence. Hence, for every $x \in \underline{S}$, there exists a finite limit
$\overline{f}(x):=\lim_{l \to \infty}\fint_{B_{l}(x)}f(z)\,d\mathfrak{m}_{l}(z).$
An application of Fatou's lemma leads to the required estimate
\begin{equation}
\label{eq.almostLebesgue}
\begin{split}
&\varlimsup\limits_{k \to \infty}\fint\limits_{B_{k}(x)}|\overline{f}(x)-f(y)|\,d\mathfrak{m}_{k}(y) \le \varlimsup\limits_{k \to \infty}\varliminf\limits_{l \to \infty}\fint\limits_{B_{k}(x)}\Bigl|\fint\limits_{B_{l}(x)}f(z)\,d\mathfrak{m}_{l}(z)-f(y)\Bigr|\,d\mathfrak{m}_{k}(y)\\
&\le \varlimsup\limits_{k \to \infty}\varliminf\limits_{l \to \infty}\fint\limits_{B_{l}(x)}\fint\limits_{B_{k}(x)}|f(y)-f(z)|\,d\mathfrak{m}_{k}(y)\,d\mathfrak{m}_{k}(z)\\
&\le \varlimsup\limits_{k \to \infty}\sum\limits_{i=k}^{\infty}\mathcal{E}_{\mathfrak{m}_{i}}(f,B_{i}(x)) \le R(x)=0 \quad \text{for all} \quad x \in \underline{S}.
\end{split}
\end{equation}
From \eqref{eq.levels} we obviously have $\mathcal{H}_{\theta}(S \setminus \underline{S}) = 0$ completing the proof.
\end{proof}
We recall the following classical result (see Corollary 3.3.51 in \cite{HKST}).
\begin{Prop}
\label{Prop.Lebesguepoints}
Let $\mathfrak{m}$ be a measure on $\operatorname{X}$. Given $p \in [1,\infty)$, for each $f \in L_{p}(\mathfrak{m})$, for every $\varepsilon > 0$, there exists
an open set $O \subset \operatorname{X}$ such that $\mathfrak{m}(O) < \varepsilon$ and a continuous function $f_{\varepsilon}$ such that $f_{\varepsilon}(x)=f(x)$
for $\mathfrak{m}$-a.e. point $x \in \operatorname{X} \setminus O$.
\end{Prop}
Now we are ready to state the main result of this section. We recall Definition \ref{Def.multiweight_Lebesgue}.

\begin{Th}
\label{Th.trueLebesgue}
Let $f \in L_{1}^{loc}(\{\mathfrak{m}_{k}\})$ be such that $\operatorname{BSN}^{\delta}_{p,\{\mathfrak{m}_{k}\},c}(f) < +\infty$ for some $c > 1$ and $\delta \in (0,1]$. Then
\begin{equation}
\label{eq.leb.points}
\mathfrak{m}_{0}(S \setminus (\mathfrak{R}_{\{\mathfrak{m}_{k}\},\epsilon}(f))=0.
\end{equation}
\end{Th}

\begin{proof}
Let $\overline{f}$ and $\underline{S}$ be the same as in Lemma \ref{Lm.quasiLebesgue}. By Proposition \ref{Prop.absolute_continuity_measure} we have $\mathfrak{m}_{0}(S \setminus \underline{S})=\mathcal{H}_{\theta}(S \setminus \underline{S})=0$. Hence, in order to establish \eqref{eq.leb.points}
it is sufficient to show that
\begin{equation}
\label{eq.trueLebesgue}
f(x)=\overline{f}(x) \quad \text{for $\mathfrak{m}_{0}$--a.e.} \quad x \in \underline{S}.
\end{equation}

We apply Proposition \ref{Prop.Lebesguepoints} with $\mathfrak{m}=\mathfrak{m}_{0}$.
This gives, for each $\varepsilon >0$, the existence of a function $f_{\varepsilon} \in C(\operatorname{X})$ and an open set $O_{\varepsilon} \subset \operatorname{X}$ such that
$\mathfrak{m}_{0}(O_{\varepsilon}) < \varepsilon$ and $f_{\varepsilon}(x)=f(x)$ for $\mathfrak{m}_{0}$--a.e. point $x \in S \setminus O_{\varepsilon}$.
We recall \eqref{eq.lower_and_upper_densities} and put
$S_{\varepsilon}:=\{x \in S \setminus O_{\varepsilon}: f(x)=f_{\varepsilon}(x) \text{ and } \overline{D}^{\{\mathfrak{m}_{k}\}}_{S \setminus O_{\varepsilon}}(x,\epsilon) > 0\}.$
Taking into account
(\textbf{D.6.3}) and \eqref{M5}, we obtain, for each sequence $\varepsilon_{n} \downarrow 0$, $n \to \infty$,
\begin{equation}
\label{eq.varepsiloncover}
\mathfrak{m}_{0}\Bigl(S \setminus \cup_{n \in \mathbb{N}} S_{\varepsilon_{n}}\Bigr)=0.
\end{equation}

Fix a small enough $\varepsilon > 0$ and a point $\underline{x} \in S_{\varepsilon} \cap \underline{S}$.
By the Chebyshev inequality, for any fixed $\sigma > 0$, we have
\begin{equation}
\begin{split}
&(\mathfrak{m}_{k}(B_{k}(\underline{x})))^{-1}\mathfrak{m}_{k}(\{y \in B_{k}(x):|f(y)-\overline{f}(\underline{x})| \geq \sigma\})\\
&\le \frac{1}{\sigma}\fint\limits_{B_{k}(x)}|f(y)-\overline{f}(\underline{x})|\,d\mathfrak{m}_{k}(y) \to 0, \quad k \to \infty.
\end{split}
\end{equation}
We set $c(\underline{x}):=\overline{D}^{\{\mathfrak{m}_{k}\}}_{S \setminus O_{\varepsilon}}(x,\epsilon)$ for brevity. Hence, given $\sigma > 0$, there exists a sufficiently large
number $k=k(\underline{x},\sigma) \in \mathbb{N}$ such that
\begin{equation}
\notag
\mathfrak{m}_{k}(\{y \in B_{k}(x):|f(y)-\overline{f}(\underline{x})| < \frac{\sigma}{2} \text{ and } |f(y)-f(\underline{x})| < \frac{\sigma}{2}\}) \geq \frac{c(\underline{x})}{2}\mathfrak{m}_{k}(B_{k}(\underline{x})).
\end{equation}
As a result, by the triangle inequality, $|\overline{f}(\underline{x})-f(\underline{x})| < \sigma$.
Since, given $\underline{x} \in S_{\varepsilon}$, one can chose $\sigma > 0$ arbitrarily small, we obtain
\begin{equation}
\label{eq.epsilonLebesgue}
\overline{f}(\underline{x})=f(\underline{x}) \quad \text{for all} \quad \underline{x} \in S_{\varepsilon}.
\end{equation}
Finally, taking into account that $\varepsilon > 0$ can be chosen arbitrary small and combining \eqref{eq.varepsiloncover}, \eqref{eq.epsilonLebesgue}
we deduce \eqref{eq.trueLebesgue} and complete the proof.
\end{proof}

\section{Extension operator}

\textit{Throughout this section we fix the following data}:

\begin{itemize}
\item[\(\rm (\textbf{D.7.1})\)] a parameter $p \in (1,\infty)$ and a metric measure space $\operatorname{X}=(\operatorname{X},\operatorname{d},\mu) \in \mathfrak{A}_{p}$;

\item[\(\rm (\textbf{D.7.2})\)] a parameter $\theta \in [0,p)$ and a closed set $S \in \mathcal{LCR}_{\theta}(\operatorname{X})$;

\item[\(\rm (\textbf{D.7.3})\)] a sequence of measures $\{\mathfrak{m}_{k}\} \in \mathfrak{M}^{str}_{\theta}(S)$ with parameter $\epsilon=\epsilon(\{\mathfrak{m}_{k}\}) \in (0,\frac{1}{10}]$.
\end{itemize}
In this section, we put $B_{k}(x):=B_{\epsilon^{k}}(x)$ for each $k \in \mathbb{Z}$ and all $x \in \operatorname{X}$.
Furthermore, we recall notation \eqref{eqq.centers_of_dyadic_cubes} and fix a sequence $\{Z_{k}(\operatorname{X},\epsilon)\}:=\{Z_{k}(\operatorname{X},\epsilon)\}_{k \in \mathbb{Z}}$. We recall \eqref{eq.lattice} and put
\begin{equation}
\label{eqq.important_balls}
\widetilde{B}_{k,\alpha}:=B_{2\epsilon^{k}}(z_{k,\alpha}) \quad \text{for each} \quad k \in \mathbb{Z}, \quad \text{for every} \quad \alpha \in \mathcal{A}_{k}(\operatorname{X},\epsilon).
\end{equation}
Given $k \in \mathbb{Z}$, the \textit{$k$th neighborhood of $S$} and \textit{the $k$th layer of $S$}, respectively, are defined by
\begin{equation}
\label{eq.5.1'}
U_{k}(S):=\{x \in \operatorname{X}: \operatorname{dist}(x,S) < 5 \epsilon^{k}\}, \quad V_{k}(S):=U_{k-1}(S) \setminus U_{k}(S).
\end{equation}
The advantages of such layers are clear from the following elementary proposition.

\begin{Prop}
\label{Prop.5_useful}
Let $k,k' \in \mathbb{Z}$ be such that $|k-k'| \geq 2$. Then, for any ball $B=\widetilde{B}_{k,\alpha}$ with $B \cap V_{k}(S) \neq \emptyset$
and any ball $B'=\widetilde{B}_{k',\alpha'}$ with $B' \cap V_{k'}(S) \neq \emptyset$, we have $B \cap B' = \emptyset$.
\end{Prop}

\begin{proof} We fix arbitrary balls $B$, $B'$ satisfying the assumptions of the lemma. If $B \cap B' \neq \emptyset$, then by the triangle inequality we
get $\operatorname{dist}(V_{k}(S),V_{k'}(S)) \le 4(\epsilon^{k}+\epsilon^{k'})$.
On the other hand, since $\epsilon \le \frac{1}{10}$, we get $\operatorname{dist}(V_{k}(S),V_{k'}(S)) \geq 4\epsilon^{\min\{k,k'\}}+5\epsilon^{\max\{k,k'\}}$.
This contradiction gives the required claim.

\end{proof}

A useful property of our space $\operatorname{X}$ is the following simple and known
result about partitions of unity (see Lemma 2.4 in \cite{GT} for the details):

\begin{Lm}
\label{Lm.partitionunity}
There exists a constant $C> 0$ depending only on $C_{\mu}(10)$ such that
\begin{equation}
\label{eq.partunity1}
0 \le \varphi_{k,\alpha} \le \chi_{\widetilde{B}_{k,\alpha}} \quad \text{and} \quad \operatorname{lip}\varphi_{k,\alpha} \le \frac{C}{\epsilon^{k}}\chi_{\widetilde{B}_{k,\alpha}}
\quad \text{for each} \quad k \in \mathbb{N}_{0}, \quad \text{for all} \quad \alpha \in \mathcal{A}_{k}(\operatorname{X},\epsilon),
\end{equation}
and, furthermore,
\begin{equation}
\label{eq.partunity2}
\sum_{\alpha \in \mathcal{A}_{k}(\operatorname{X},\epsilon)}\varphi_{k,\alpha}=1 \quad \text{for each} \quad k \in \mathbb{N}_{0}.
\end{equation}
\end{Lm}

Now we establish a simple combinatorial result which is a folklore. Nevertheless, we present the details
for completeness.

\begin{Lm}
\label{Lm.combinatorial_result}
There exists a constant $\underline{N} \in \mathbb{N}$ depending only on $C_{\mu}(10)$ such that for each $k \in \mathbb{N}_{0}$ the family $\widetilde{\mathcal{B}}_{k}:=\{\widetilde{B}_{k,\alpha}:\alpha \in \mathcal{A}_{k}(\operatorname{X},\epsilon)\}$
can be decomposed into $N \le \underline{N}$ disjoint subfamilies $\{\widetilde{\mathcal{B}}^{i}_{k}\}_{i=1}^{N}$.
\end{Lm}

\begin{proof}
We fix $k \in \mathbb{N}_{0}$ and put $E(0)=Z_{k}(\operatorname{X},\epsilon)$. By $Z(1)$ we denote a maximal $5\epsilon^{k}$-separated subset of $Z_{k}(\operatorname{X},\epsilon)$ and put $E(1):=Z_{k}(\operatorname{X},\epsilon) \setminus Z(1)$.
Arguing by induction, suppose that for some $i \in \mathbb{N}$ we have already built sets $Z(1),...,Z(i)$ and $E(1),...,E(i)$ in such a way that $E(i')=Z_{k}(\operatorname{X},\epsilon) \setminus \cup_{l=1}^{i'}Z(l)$, $i' \in \{1,...,i\}$.
Let $Z(i+1)$ be a maximal $5\epsilon^{k}$-separated subset of $E(i)$ and $E(i+1):=E(i) \setminus Z(i+1)$. We put $\underline{N}:=\lceil N_{\mu}(1,24)\rceil$ (where $N_{\mu}(R,c)$ is the same as in Proposition \ref{Prop.metric_doubling}).
We show that $E(i)=\emptyset$ for each $i > \underline{N}$.
Indeed, assume that there is a number $i > \underline{N}$ such that $E(i) \neq \emptyset$ and fix $\underline{x}(i) \in E(i)$. Given $i' \in \{1,...,i\}$, by the maximality of $Z(i')$,
and the obvious inclusion $E_{i} \subset E_{i'-1}$ it follows that there is a point $\underline{x}(i') \in Z(i')$ such that $\operatorname{d}(\underline{x}(i'),\underline{x}(i)) < 5\epsilon^{k}$.
Hence, we have $B_{\epsilon^{k}/4}(\underline{x}(i')) \subset B_{6\epsilon^{k}}(\underline{x}(i))$. As a result, since $i'$ can be chosen arbitrarily, there exists
a family $\mathcal{F}:=\{B_{\epsilon^{k}/4}(\underline{x}(i')):i' \in \{1,...,i\}\}$ of pairwise disjoint balls contained in the ball $B_{6\epsilon^{k}}(\underline{x}(i))$ such that $\#\mathcal{F}=i$.
Combining this observation with Proposition \ref{Prop.metric_doubling}, we get a contradiction. Hence, $\#\{i \in \mathbb{N}:E(i) \neq \emptyset\} \le \underline{N}$. It remains to note
that, for each $i \in \{1,...,\underline{N}\}$, for every $z,z' \in Z(i)$, we have $B_{2\epsilon^{k}}(z) \cap B_{2\epsilon^{k}}(z') = \emptyset$.
\end{proof}

Given $k \in \mathbb{Z}$, we set
\begin{equation}
\label{eq.5.2'}
\mathcal{A}_{k}(S):=\{\alpha \in \mathcal{A}_{k}(\operatorname{X},\epsilon):\widetilde{B}_{k,\alpha} \cap U_{k-1}(S) \neq \emptyset\}.
\end{equation}

\begin{Remark}
\label{Rem.important_inclusion}
Since $\epsilon \in (0,\frac{1}{10}]$, by \eqref{eq.5.1'} and \eqref{eq.5.2'} it is easy to see that, for each $k \in \mathbb{Z}$ and any $\alpha \in \mathcal{A}_{k}(S)$, there exists a point $\underline{x} \in S$ such that $B_{\epsilon^{k}}(\underline{x}) \subset \frac{3}{\epsilon}\widetilde{B}_{k,\alpha}=B_{\frac{6}{\epsilon}}(z_{k,\alpha})$.
\end{Remark}

The following result is an immediate consequence of \eqref{eq.5.1'}--\eqref{eq.5.2'}.

\begin{Prop}
\label{Prop.supportsinclusions}
For each $k \in \mathbb{N}_{0}$,  $\chi_{U_{k-1}(S)}(x) \le \sum_{\alpha \in \mathcal{A}_{k}(S)} \varphi_{k,\alpha}(x) \le \chi_{U_{k-2}(S)}(x)$, $x \in \operatorname{X}$.
\end{Prop}

Now, keeping in mind Remark \ref{Rem.important_inclusion}, given an element $f \in L^{loc}_{1}(\{\mathfrak{m}_{k}\})$, we define for each $k \in \mathbb{N}_{0}$ the special family of numbers.
More precisely,  we put
\begin{equation}
\label{eq.5.1}
f_{k,\alpha}:=
\begin{cases}
\fint\limits_{\frac{3}{\epsilon}\widetilde{B}_{k,\alpha}}f(x)\,d\mathfrak{m}_{k}(x), \quad \text{if} \quad \alpha \in \mathcal{A}_{k}(S);\\
0, \quad \text{if} \quad \alpha \in \mathcal{A}_{k}(\operatorname{X},\epsilon) \setminus \mathcal{A}_{k}(S).
\end{cases}
\end{equation}

The following simple proposition will be useful in what follows. We recall notation \eqref{eqq.tricky_average}.

\begin{Prop}
\label{Prop.5.2''}
There exists a constant $C > 0$ such that, for each $f \in L_{1}^{loc}(\{\mathfrak{m}_{k}\})$, for each $k \in \mathbb{N}_{0}$, the following inequality
\begin{equation}
\label{eqq.5.2}
|f_{k,\alpha}-f_{k',\beta}| \le C \widetilde{\mathcal{E}}_{\mathfrak{m}_{k}}(f,\frac{3}{\epsilon}\widetilde{B}_{k,\alpha})
\end{equation}
holds for each $k' \in \{k,k+1\}$ and any $\alpha \in \mathcal{A}_{k}(S)$, $\beta \in \mathcal{A}_{k'}(S)$
for which $\widetilde{B}_{k,\alpha} \cap \widetilde{B}_{k',\beta} \neq \emptyset$.
\end{Prop}

\begin{proof}
We fix $f \in L_{1}^{loc}(\{\mathfrak{m}_{k}\})$, numbers $k \in \mathbb{N}_{0}$, $k' \in \{k,k+1\}$, and indexes $\alpha \in \mathcal{A}_{k}(S)$, $\beta \in \mathcal{A}_{k'}(S)$
such that $\widetilde{B}_{k,\alpha} \cap \widetilde{B}_{k',\beta} \neq \emptyset$.
By \eqref{eqq.important_balls} we have $\frac{3}{\epsilon}\widetilde{B}_{k',\beta} \subset (\frac{3}{\epsilon}+4)\widetilde{B}_{k,\alpha} \subset \frac{6}{\epsilon}\widetilde{B}_{k,\alpha}.$
Hence, using \eqref{eq.5.1}, \eqref{M.4}, Theorem \ref{Th.doubling_type} and Remark \ref{Rem.best_approx} we get the required estimate
\begin{equation}
\label{eqq.5.3}
\begin{split}
&|f_{k,\alpha}-f_{k',\beta}| \le \fint\limits_{\frac{3}{\epsilon}\widetilde{B}_{k,\alpha}}\fint\limits_{\frac{3}{\epsilon}\widetilde{B}_{k',\beta}}|f(y)-f(y')|\,d\mathfrak{m}_{k}(y)d\mathfrak{m}_{k'}(y') \\
&\le C \fint\limits_{\frac{6}{\epsilon}\widetilde{B}_{k,\alpha}}\fint\limits_{\frac{6}{\epsilon}\widetilde{B}_{k,\alpha}}|f(y)-f(y')|\,d\mathfrak{m}_{k}(y)d\mathfrak{m}_{k}(y')
\le C \widetilde{\mathcal{E}}_{\mathfrak{m}_{k}}(f,\frac{3}{\epsilon}\widetilde{B}_{k,\alpha}).
\end{split}
\end{equation}
The proof is complete.
\end{proof}

Given an element $f \in L^{loc}_{1}(\{\mathfrak{m}_{k}\})$, for each $k \in \mathbb{N}_{0}$, we put
\begin{equation}
\label{eq.aproxsequence}
f_{k}(x):=\sum\limits_{\alpha \in \mathcal{A}_{k}(\operatorname{X},\epsilon)} \varphi_{k,\alpha}(x)f_{k,\alpha}=\sum\limits_{\alpha \in \mathcal{A}_{k}(S)} \varphi_{k,\alpha}(x)f_{k,\alpha}, \quad x \in \operatorname{X}.
\end{equation}

Having at our disposal Propositions \ref{Prop.supportsinclusions}, \ref{Prop.5.2''}, we get nice pointwise estimates of the local Lipschitz constants of the functions $f_{k}$, $k \in \mathbb{N}_{0}$.
We recall \eqref{eqq.local_Lipschitz_constant}.

\begin{Prop}
\label{Prop.lipestimate1}
There exists a constant $C > 0$ such that, for each $f \in L_{1}^{loc}(\{\mathfrak{m}_{k}\})$, for every $k \in \mathbb{N}_{0}$ and every $\underline{x} \in U_{k-1}(S)$, the following inequality
\begin{equation}
\operatorname{lip}f_{k}(\underline{x}) \le \frac{C}{\epsilon^{k}}
\widetilde{\mathcal{E}}_{\mathfrak{m}_{k}}(f,\frac{3}{\epsilon}\widetilde{B}_{k,\underline{\alpha}})
\end{equation}
holds for any index $\underline{\alpha} \in \mathcal{A}_{k}(S)$ satisfying the condition $\widetilde{B}_{k,\underline{\alpha}} \ni \underline{x}$.
\end{Prop}

\begin{proof} We fix $k \in \mathbb{N}_{0}$ and $\underline{x} \in U_{k-1}(S)$. We also fix an arbitrary index $\underline{\alpha} \in \mathcal{A}_{k}(S)$ such that
$\underline{x} \in \widetilde{B}_{k,\underline{\alpha}}$.
By (3) in Proposition \ref{Prop.lipprop} and \eqref{eq.partunity2} we get $\operatorname{lip}f_{k}(\underline{x})=\operatorname{lip}\Bigl(f_{k}-\sum_{\alpha \in \mathcal{A}_{k}(\operatorname{X},\epsilon)}\varphi_{k,\alpha}f_{k,\underline{\alpha}}\Bigr)(\underline{x}).$
By \eqref{eq.5.2'} it follows that, given $\alpha \in \mathcal{A}_{k}(\operatorname{X},\epsilon)$, $\varphi_{k,\alpha}(\underline{x}) \neq 0$ implies $\alpha \in \mathcal{A}_{k}(S)$.
Hence, we use \eqref{eq.aproxsequence} together with \eqref{eq.partunity1} and, finally, take into account Propositions \ref{Prop.finite_intersection}, \ref{Prop.5.2''}. This gives
\begin{equation}
\label{eq.5.4}
\begin{split}
&\operatorname{lip}f_{k}(\underline{x}) \le \sum\limits_{\alpha \in \mathcal{A}_{k}(\operatorname{X},\epsilon)}\operatorname{lip}\varphi_{k,\alpha}(\underline{x})|f_{k,\underline{\alpha}}-f_{k,\alpha}| = \sum\limits_{\alpha \in \mathcal{A}_{k}(S)}\operatorname{lip}\varphi_{k,\alpha}(\underline{x})|f_{k,\underline{\alpha}}-f_{k,\alpha}|\\
&\le C\sum\limits_{\alpha \in \mathcal{A}_{k}(S)}\chi_{\widetilde{B}_{k,\alpha}}(\underline{x})\frac{1}{\epsilon^{k}}|f_{k,\underline{\alpha}}-f_{k,\alpha}| \le \frac{C}{\epsilon^{k}}
\widetilde{\mathcal{E}}_{\mathfrak{m}_{k}}(f,\frac{3}{\epsilon}\widetilde{B}_{k,\underline{\alpha}}).
\end{split}
\end{equation}
The proof is complete.
\end{proof}

This pointwise estimate leads to the nice estimate in the $L_{p}$-norm. Recall that, given a Borel set $E \subset \operatorname{X}$ and an element $f \in L_{p}^{loc}(\operatorname{X})$,
we put $\|f|L_{p}(E)\|:=\|f|L_{p}(E,\mu)\|$.

\begin{Lm}
\label{Lm.lipestimate1}
There exists a constant $C > 0$ such that, for each $f \in L_{1}^{loc}(\{\mathfrak{m}_{k}\})$, for each $k \in \mathbb{N}_{0}$ and any Borel set $E \subset U_{k-1}(S)$,
\begin{equation}
\|\operatorname{lip}f_{k}|L_{p}(E)\|^{p} \le C
\sum\limits_{E \cap \widetilde{B}_{k,\alpha} \neq \emptyset}\frac{\mu(\widetilde{B}_{k,\alpha})}{\epsilon^{kp}}\Bigl(\widetilde{\mathcal{E}}_{\mathfrak{m}_{k}}(f,\frac{3}{\epsilon}\widetilde{B}_{k,\alpha})\Bigr)^{p}.
\end{equation}
\end{Lm}

\begin{proof}
By Proposition \ref{Prop.lipestimate1} for any ball $\widetilde{B}_{k,\alpha}$ with $E \cap \widetilde{B}_{k,\alpha} \neq \emptyset$ we clearly have
\begin{equation}
\notag
\int\limits_{E \cap \widetilde{B}_{k,\alpha}}\Bigl(\operatorname{lip}f_{k}(x)\Bigr)^{p}\,d\mu(x) \le C\frac{\mu(E \cap \widetilde{B}_{k,\alpha})}{\epsilon^{kp}}
\Bigl(\widetilde{\mathcal{E}}_{\mathfrak{m}_{k}}(f,\frac{3}{\epsilon}\widetilde{B}_{k,\alpha})\Bigr)^{p}.
\end{equation}
This observation in combination with \eqref{eq.5.2'} gives
\begin{equation}
\begin{split}
&\int\limits_{E}\Bigl(\operatorname{lip}f_{k}(x)\Bigr)^{p}\,d\mu(x) \le
\sum\limits_{\alpha \in \mathcal{A}_{k}(S)}\int\limits_{E \cap \widetilde{B}_{k,\alpha}}\Bigl(\operatorname{lip}f_{k}(x)\Bigr)^{p}\,d\mu(x)\\
&\le C
\sum\limits_{E \cap \widetilde{B}_{k,\alpha} \neq \emptyset}\frac{\mu(\widetilde{B}_{k,\alpha})}{\epsilon^{kp}}\Bigl(\widetilde{\mathcal{E}}_{\mathfrak{m}_{k}}(f,\frac{3}{\epsilon}\widetilde{B}_{k,\alpha})\Bigr)^{p}.
\end{split}
\end{equation}
This completes the proof.
 \end{proof}

To construct our extension operator, for any given $f \in L^{loc}_{1}(\{\mathfrak{m}_{k}\})$, we built some special sequence $\{f^{j}\}_{j \in \mathbb{N}}$. Informally speaking,
the graph of each $f^{j}$ looks like a stairway composed of elementary steps $\operatorname{St}_{i}[f]$, $i=1,...,j$. More precisely,
given $f \in L^{loc}_{1}(\{\mathfrak{m}_{k}\})$, we set $f_{0}:=0$. Furthermore, arguing inductively,
for each $i \in \mathbb{N}$, we define the \textit{elementary $i$th step of $f$} by
\begin{equation}
\label{eqq.415'}
\operatorname{St}_{i}[f](x):=\sum\limits_{\substack{\alpha \in \mathcal{A}_{i}(S)}}\varphi_{i,\alpha}(x)(f_{i,\alpha}-f_{i-1}(x)), \quad x \in \operatorname{X}.
\end{equation}

\begin{Remark}
\label{Rem.support_steps}
In view of Proposition \ref{Prop.supportsinclusions} it is clear that $\operatorname{supp}\operatorname{St}_{i}[f] \subset U_{i-2}(S)$ for all $i \in \mathbb{N}$.
\end{Remark}
Finally, we define the \textit{special approximating sequence} by letting, for each $j \in \mathbb{N}$,
\begin{equation}
\label{eq.specaproxsequence}
f^{j}(x):=\sum\limits_{i=1}^{j}\operatorname{St}_{i}[f](x), \qquad x \in \operatorname{X}.
\end{equation}

\begin{Prop}
\label{Prop.stac_sequence}
For each point $x \in \operatorname{X} \setminus S$, there exists $j(x) \in \mathbb{N}$ such that
$f^{j}(x)=f^{j(x)}(x)$ for all $j \geq j(x)$.
\end{Prop}

\begin{proof}
For each $x \in \operatorname{X} \setminus S$, from Remark \ref{Rem.support_steps} and \eqref{eq.specaproxsequence} we have $f^{j}(x)=f^{j+1}(x)$ provided that $x \in \operatorname{X} \setminus U_{j-1}(S)$.
Since $U_{j+1}(S) \subset U_{j}(S)$ for all $j \in \mathbb{N}$, the claim follows.
\end{proof}

Now we are ready to present our extension operator.

\begin{Def}
\label{Def.Extension_operator}
The operator $\operatorname{Ext}_{S,\{\mathfrak{m}_{k}\}}:L_{1}^{loc}(\{\mathfrak{m}_{k}\}) \to L_{0}(\mathfrak{m}_{0}) \cap \mathfrak{B}(\operatorname{X} \setminus S)$ is defined
by
\begin{equation}
\label{eq.Extension_operator}
\operatorname{Ext}_{S,\{\mathfrak{m}_{k}\}}(f):=\chi_{S}f+\chi_{\operatorname{X}\setminus S}\lim\limits_{j \to \infty}f^{j}, \quad f \in L_{1}^{loc}(\{\mathfrak{m}_{k}\}),
\end{equation}
where by $\chi_{S}f$ we mean an $\mathfrak{m}_{0}$-equivalence class and $\chi_{\operatorname{X}\setminus S}\lim_{j \to \infty}f^{j}$ denotes the pointwise limit of the sequence $\{f^{j}\}$
on the set $\operatorname{X}\setminus S$.
\end{Def}

\begin{Remark}
By Proposition \ref{Prop.stac_sequence}, the operator $\operatorname{Ext}_{S,\{\mathfrak{m}_{k}\}}$ is well defined and linear.
\end{Remark}


The main reason for introducing the sequence $\{f^{j}\}$ in that way is the presence of remarkable pointwise properties of
steps $\operatorname{St}_{i}[f]$, $i \in \mathbb{N}$. More precisely, the following result holds.

\begin{Prop}
\label{Prop.5.4''}
Let $f \in L_{1}^{loc}(\{\mathfrak{m}_{k}\})$. Let $i \in \mathbb{N}$, $\underline{x} \in U_{i-2}(S)$ and
$\underline{\alpha} \in \mathcal{A}_{i-1}(S)$ be such that $\underline{x} \in \widetilde{B}_{i-1,\underline{\alpha}}$. Then
\begin{equation}
\widetilde{\operatorname{St}}_{i}[f](\underline{x}):= \sum\limits_{\substack{\alpha \in \mathcal{A}_{i}(S)}}\chi_{\widetilde{B}_{i,\alpha}}(\underline{x})|f_{i,\alpha}-f_{i-1}(\underline{x})|
\le C \widetilde{\mathcal{E}}_{\mathfrak{m}_{i-1}}(f,\frac{3}{\epsilon}\widetilde{B}_{i-1,\underline{\alpha}}),
\end{equation}
where the constant $C > 0$ depends neither on $f$ nor on $i,\underline{x}$ and $\underline{\alpha}$.
\end{Prop}

\begin{proof}
By \eqref{eq.partunity2} and \eqref{eq.aproxsequence} we have
\begin{equation}
\notag
\widetilde{\operatorname{St}}_{i}[f](\underline{x}) \le \sum\limits_{\alpha \in \mathcal{A}_{i}(S)}\sum\limits_{\alpha' \in \mathcal{A}_{i-1}(\operatorname{X},\epsilon)}\chi_{\widetilde{B}_{i,\alpha}}(\underline{x})
\varphi_{i-1,\alpha'}(\underline{x})|f_{i,\alpha}-f_{i-1,\alpha'}|.
\end{equation}
Using the triangle inequality we have $|f_{i,\alpha}-f_{i-1,\alpha'}| \le |f_{i,\alpha}-f_{i-1,\underline{\alpha}}|+|f_{i-1,\underline{\alpha}}-f_{i-1,\alpha'}|$. Hence, using \eqref{eq.partunity1}, Proposition \ref{Prop.5.2''}, and taking into account Proposition \ref{Prop.finite_intersection}, we get
\begin{equation}
\notag
\widetilde{\operatorname{St}}_{i}[f](\underline{x}) \le
\sum\limits_{i'=i-1}^{i}\sum\limits_{\alpha' \in \mathcal{A}_{i'}(\operatorname{X},\epsilon)}\chi_{\widetilde{B}_{i',\alpha'}}(\underline{x})|f_{i',\alpha'}-f_{i-1,\underline{\alpha}}|
\le C \widetilde{\mathcal{E}}_{\mathfrak{m}_{i-1}}(f,\frac{3}{\epsilon}\widetilde{B}_{i-1,\underline{\alpha}}).
\end{equation}
The proof is complete.
\end{proof}
This proposition leads to useful estimates of $L_{p}$-norms of steps and their local Lipschitz constants.

\begin{Lm}
\label{Lm.steps_estimate}
There exists $C > 0$ such that, for each $i \in \mathbb{N}$,
the following properties hold:

\begin{itemize}
\item[\(\rm 1)\)] for each Borel set $E \subset U_{i-2}(S)$ and any measure $\nu$ on $\operatorname{X}$,
\begin{equation}
\notag
\|\operatorname{St}_{i}[f]|L_{p}(E,\nu)\|^{p} \le C \sum\limits_{\substack{\widetilde{B}_{i-1,\alpha} \cap E \neq \emptyset}}\nu(\widetilde{B}_{i-1,\alpha})\Bigl(\widetilde{\mathcal{E}}_{\mathfrak{m}_{i-1}}(f,\frac{3}{\epsilon}\widetilde{B}_{i-1,\alpha})\Bigr)^{p};
\end{equation}

\item[\(\rm 2)\)] for each Borel set $E \subset U_{i-2}(S)$,
\begin{equation}
\notag
\|\operatorname{lip}(\operatorname{St}_{i}[f])|L_{p}(E)\|^{p} \le C \sum\limits_{\substack{\widetilde{B}_{i-1,\alpha} \cap E \neq \emptyset}}\frac{\mu(\widetilde{B}_{i-1,\alpha})}{\epsilon^{(i-1)p}}\Bigl(\widetilde{\mathcal{E}}_{\mathfrak{m}_{i-1}}(f,\frac{3}{\epsilon}\widetilde{B}_{i-1,\alpha})\Bigr)^{p}.
\end{equation}

\end{itemize}
\end{Lm}

\begin{proof}
We fix $i \in \mathbb{N}$ and a Borel set $E \subset U_{i-2}(S)$.

To prove the first claim, we note that by \eqref{eq.partunity1} we have $|\operatorname{St}_{i}[f](x)| \le \widetilde{\operatorname{St}}_{i}[f](x)$ for all $x \in \operatorname{X}$.
Hence, an application of Proposition \ref{Prop.5.4''} gives
\begin{equation}
\label{eq.417''}
\begin{split}
&\|\operatorname{St}_{i}[f]|L_{p}(E,\nu)\|^{p} \le \|\widetilde{\operatorname{St}}_{i}[f]|L_{p}(E,\nu)\|^{p} \le
\sum\limits_{\alpha \in \mathcal{A}_{i-1}(S)}\int\limits_{\widetilde{B}_{i-1,\alpha} \cap E}
\Bigl(\widetilde{\operatorname{St}}_{i}[f](x)\Bigr)^{p}\,d\nu(x)\\
&\le C \sum\limits_{\substack{\widetilde{B}_{i-1,\alpha} \cap E \neq \emptyset}}\nu(\widetilde{B}_{i-1,\alpha})\Bigl(\widetilde{\mathcal{E}}_{\mathfrak{m}_{i-1}}(f,\frac{3}{\epsilon}\widetilde{B}_{i-1,\alpha})\Bigr)^{p}.
\end{split}
\end{equation}

To prove the second claim, we note that by \eqref{eqq.415'}, Proposition \ref{Prop.lipprop}, and \eqref{eq.partunity1}, \eqref{eq.partunity2} we have
\begin{equation}
\begin{split}
\notag
&\operatorname{lip}(\operatorname{St}_{i}[f])(x) \le \sum\limits_{\substack{\alpha \in \mathcal{A}_{i}(S)}}\operatorname{lip}\varphi_{i,\alpha}(x)|f_{i,\alpha}-f_{i-1}(x)|+
\sum\limits_{\substack{\alpha \in \mathcal{A}_{i}(S)}}\varphi_{i,\alpha}(x)\operatorname{lip}f_{i-1}(x)\\
&\le \frac{C}{\epsilon^{i}}\widetilde{\operatorname{St}}_{i}[f](x)+\operatorname{lip}f_{i-1}(x) \quad \text{for all} \quad x \in \operatorname{X}.
\end{split}
\end{equation}
Hence, using Lemma \ref{Lm.lipestimate1} with $k=i-1$ and \eqref{eq.417''} with $\nu=\mu$ we obtain the required estimate
\begin{equation}
\begin{split}
&\|\operatorname{lip}(\operatorname{St}_{i}[f])|L_{p}(E)\|^{p} \le
C \sum\limits_{\substack{\widetilde{B}_{i-1,\alpha} \cap E \neq \emptyset}}\frac{\mu(\widetilde{B}_{i-1,\alpha})}{\epsilon^{(i-1)p}}\Bigl(\widetilde{\mathcal{E}}_{\mathfrak{m}_{i-1}}(f,\frac{3}{\epsilon}\widetilde{B}_{i-1,\alpha})\Bigr)^{p}.
\end{split}
\end{equation}

The proof is complete.

\end{proof}

The crucial observation is given in the following lemma.

\begin{Lm}
\label{Lm.derivappseq}
There is a constant $C > 0$ such that, for each $f \in L^{loc}_{1}(\{\mathfrak{m}_{k}\})$, for every $j \in \mathbb{N}$,
\begin{equation}\
\label{eqq.5.9}
\|\operatorname{lip}f^{j}|L_{p}(U_{0}(S))\|^{p} \le  C \sum\limits_{i=1}^{j}
\sum\limits_{\substack{\widetilde{B}_{i,\alpha} \cap \widehat{V}_{i}(S) \neq \emptyset}}\frac{\mu(\widetilde{B}_{i,\alpha})}{\epsilon^{ip}}\Bigl(\widetilde{\mathcal{E}}_{\mathfrak{m}_{i}}(f,\frac{3}{\epsilon}\widetilde{B}_{i,\alpha})\Bigr)^{p},
\end{equation}
where $\widehat{V}_{i}(S):=V_{i}(S)$ if $j \geq 2$, $i\in \{1,...,j-1\}$ and $\widehat{V}_{j}(S):=U_{j-1}(S)$.
\end{Lm}

\begin{proof}
Fix for a moment $j \in \mathbb{N}$, $j \geq 2$ and
$i \in \{1,...,j-1\}$. Given $\underline{x} \in \widehat{V}_{i}(S)$, by \eqref{eq.5.1'} and Remark \ref{Rem.support_steps} it is clear that
$\operatorname{St}_{i'}[f](\underline{x})=0$ for all $i' \geq i+2$. Hence,
using \eqref{eq.specaproxsequence} and Proposition \ref{Prop.lipprop} we obtain, $\operatorname{lip}f^{j}(\underline{x}) \le \operatorname{lip}f_{i}(\underline{x})+\operatorname{lip}\operatorname{St}_{i+1}[f](\underline{x})$.
Thus, applying Lemma \ref{Lm.lipestimate1} and Lemma \ref{Lm.steps_estimate} with $E=\widehat{V}_{i}(S)$, we deduce
\begin{equation}
\begin{split}
\label{eqq.4.20}
&\|\operatorname{lip}f^{j}|L_{p}(\widehat{V}_{i}(S))\|^{p}  \le \|\operatorname{lip}f_{i}|L_{p}(\widehat{V}_{i}(S))\|^{p}
 + \|\operatorname{lip}\operatorname{St}_{i+1}[f]|L_{p}(\widehat{V}_{i}(S))\|^{p}  \\
&\le C \sum\limits_{\substack{\widetilde{B}_{i,\alpha} \cap \widehat{V}_{i}(S) \neq \emptyset}}\frac{\mu(\widetilde{B}_{i,\alpha})}{\epsilon^{ip}}\Bigl(\widetilde{\mathcal{E}}_{\mathfrak{m}_{i}}(f,\frac{3}{\epsilon}\widetilde{B}_{i,\alpha})\Bigr)^{p}.
\end{split}
\end{equation}

On the other hand, given $j \in \mathbb{N}$, an application of Lemma  \ref{Lm.lipestimate1} with $E=U_{j-1}(S)$ gives
\begin{equation}
\label{eqq.4.21}
\begin{split}
&\|\operatorname{lip}f^{j}|L_{p}(U_{j-1}(S))\|^{p}=\|\operatorname{lip}f_{j}|L_{p}(U_{j-1}(S))\|^{p} \\
&\le C
\sum\limits_{\substack{\widetilde{B}_{j,\alpha} \cap U_{j-1}(S) \neq \emptyset}}
\frac{\mu(\widetilde{B}_{j,\alpha})}{\epsilon^{jp}}\Bigl(\widetilde{\mathcal{E}}_{\mathfrak{m}_{j}}(f,\frac{3}{\epsilon}\widetilde{B}_{j,\alpha})\Bigr)^{p}.
\end{split}
\end{equation}

Summing \eqref{eqq.4.20} over all $i \in \{1,...,j-1\}$ and then taking into account \eqref{eqq.4.21}, we arrive at \eqref{eqq.5.9} completing the proof.

\end{proof}

\begin{Lm}
\label{Lm.derivappseq'}
There exists a constant $C > 0$ such that, for each $f \in L^{loc}_{1}(\{\mathfrak{m}_{k}\})$,
\begin{equation}\
\label{eq.5.9''}
\|f_{1}|L_{p}(\operatorname{X})\|^{p}+\|\operatorname{lip}f_{1}|L_{p}(\operatorname{X})\|^{p} \le  C\|f|L_{p}(\mathfrak{m}_{0})\|^{p}.
\end{equation}
\end{Lm}

\begin{proof}
Combining the first inequality in \eqref{eq.partunity1} with \eqref{eq.5.1}, \eqref{eq.aproxsequence} and using H\"older's inequality, we get
\begin{equation}
\label{eq.423''}
|f_{1}(x)|^{p} \le C \sum\limits_{\alpha \in \mathcal{A}_{1}(S)}\chi_{\widetilde{B}_{1,\alpha}}(x)
\Bigl(\fint\limits_{\frac{3}{\epsilon}\widetilde{B}_{1,\alpha}}|f(y)|\,d\mathfrak{m}_{1}(y)\Bigr)^{p}, \quad x \in \operatorname{X}.
\end{equation}
Similarly, taking into account the second inequality in \eqref{eq.partunity1}, we have
\begin{equation}
\label{eq.424''}
(\operatorname{lip}f_{1}(x))^{p} \le C \sum\limits_{\alpha \in \mathcal{A}_{1}(S)}\chi_{\widetilde{B}_{1,\alpha}}(x)
\Bigl(\fint\limits_{\frac{3}{\epsilon}\widetilde{B}_{1,\alpha}}|f(y)|\,d\mathfrak{m}_{1}(y)\Bigr)^{p}, \quad x \in \operatorname{X}.
\end{equation}
Combining \eqref{eq.423''} with \eqref{eq.424''} and using H\"older's inequality, we get
\begin{equation}
\label{eqq.427}
\begin{split}
&\|f_{1}|L_{p}(\mu)\|^{p}+\|\operatorname{lip}f_{1}|L_{p}(\mu)\|^{p} \le C \sum\limits_{\alpha \in \mathcal{A}_{1}(S)}\mu(\widetilde{B}_{1,\alpha})\fint\limits_{\frac{3}{\epsilon}\widetilde{B}_{1,\alpha}}|f(y)|^{p}\,d\mathfrak{m}_{1}(y).
\end{split}
\end{equation}
By \eqref{eq.5.2'} we have $(\frac{6}{\epsilon}-1) B_{1,\alpha} \cap S \neq \emptyset$ for all $\alpha \in \mathcal{A}_{1}(S)$.
Hence, applying Proposition \ref{Prop.niceball_estimate} with $c=\frac{6}{\epsilon}-1$, $c' = \frac{6}{\epsilon}$ and taking into account the uniformly locally
doubling property of $\mu$, we get $\mu(\widetilde{B}_{1,\alpha})/\mathfrak{m}_{1}(\frac{3}{\epsilon}\widetilde{B}_{1,\alpha}) \le C$
for all $\alpha \in \mathcal{A}_{1}(S)$.
As a result, using this observation, Propositions \ref{Prop.covering_multiplicity}, \ref{Prop.finite_intersection} and taking into account \eqref{M.4}, we obtain
\begin{equation}
\label{eqq.427'''''}
\begin{split}
&\sum\limits_{\alpha \in \mathcal{A}_{1}(S)}\mu(\widetilde{B}_{1,\alpha})\fint\limits_{\frac{3}{\epsilon}\widetilde{B}_{1,\alpha}}|f(y)|^{p}\,d\mathfrak{m}_{1}(y) \le C \sum\limits_{\alpha \in \mathcal{A}_{1}(S)}\int\limits_{\frac{3}{\epsilon}\widetilde{B}_{1,\alpha}}|f(y)|^{p}\,d\mathfrak{m}_{1}(y)\\
& \le C \|f|L_{p}(\mathfrak{m}_{1})\|^{p} \le C \|f|L_{p}(\mathfrak{m}_{0})\|^{p}.
\end{split}
\end{equation}
Combining \eqref{eqq.427} and \eqref{eqq.427'''''} we get the required estimate and complete the proof.
\end{proof}

We recall Definitions \ref{Def.sc_good} and \ref{Def.Br_Shv}. We also recall \eqref{eqq.tricky_average} and write, for brevity, $k(B):=k(r(B))$. Now we introduce a new useful functional.
\begin{Def}
\label{Def.kostyl}
Given $f \in L_{1}^{loc}(\{\mathfrak{m}_{k}\})$, we put
\begin{equation}
\begin{split}
\label{eq.kostyl}
&\operatorname{N}_{p,\{\mathfrak{m}_{k}\},c}(f):=\varliminf\limits_{\delta \to 0}\operatorname{BSN}^{\delta}_{p,\{\mathfrak{m}_{k}\},c}(f)
+\sup \Bigl(\sum\limits_{B \in \mathcal{B}} \frac{\mu(B)}{(r(B))^{p}}\Bigl(\widetilde{\mathcal{E}}_{\mathfrak{m}_{k(B)}}(f,cB)\Bigr)^{p}\Bigr)^{\frac{1}{p}},
\end{split}
\end{equation}
where the supremum in the second term is taken over all $(S,c)$-Whitney families $\mathcal{B}$.
\end{Def}

Now we present a keystone estimate for the local Lipschitz constants of functions $f^{j}$, $j \in \mathbb{N}$.

\begin{Th}
\label{Th.derivappseq}
Given $c \geq \frac{3}{\epsilon}$, there exists a constant $C > 0$ such that
\begin{equation}
\label{eq.5.9}
\varliminf\limits_{j \to \infty}\|\operatorname{lip}f^{j}|L_{p}(\operatorname{X})\|^{p}  \le  C \operatorname{N}_{p,\{\mathfrak{m}_{k}\},c}(f) \quad \text{for all}
\quad f \in L^{loc}_{1}(\{\mathfrak{m}_{k}\}).
\end{equation}
\end{Th}

\begin{proof}

Without loss of generality wee may assume that $\operatorname{N}_{p,\{\mathfrak{m}_{k}\},c}(f) < +\infty$, since otherwise the inequality is trivial. We split the proof into several steps.

\textit{Step 1.}
We claim that for each $j \in \mathbb{N}$, $j \geq 2$ there is an $(S,c)$-Whitney family of balls $\mathcal{B}^{j}_{1}(S)$ such that (we put $k(B):=k(r(B))$, as usual, and recall \eqref{eqq.tricky_average})
\begin{equation}
\label{eq.sum_splitting}
\sum\limits_{B \in \mathcal{B}^{j}_{1}(S)}\frac{\mu(B)}{(r(B))^{p}}\Bigl(\widetilde{\mathcal{E}}_{\mathfrak{m}_{k(B)}}(f,cB)\Bigr)^{p} \geq \frac{1}{2\underline{N}}\sum\limits_{i=1}^{j-1}
\sum\limits_{\substack{\widetilde{B}_{i,\alpha} \cap V_{i}(S) \neq \emptyset}}\frac{\mu(\widetilde{B}_{i,\alpha})}{\epsilon^{ip}}\Bigl(\widetilde{\mathcal{E}}_{\mathfrak{m}_{i}}(f,c\widetilde{B}_{i,\alpha})\Bigr)^{p},
\end{equation}
where the constant $\underline{N}$ is the same as in Lemma \ref{Lm.combinatorial_result}.
Indeed, we split the sum in the right-hand side of \eqref{eq.sum_splitting} into sums
over odd and even $i \in \{1,....,j-1\}$, respectively. Without loss of generality, we may assume that the sum over odd indices is not smaller than that over even indices.
Next, for each odd $i \in \{1,....,j-1\}$ we apply Lemma \ref{Lm.combinatorial_result} and split the family $\{\widetilde{B}_{i,\alpha}:\widetilde{B}_{i,\alpha} \cap V_{i}\}$ into
at most $\underline{N}$ disjoint subfamilies. For each odd $i \in \{1,....,j-1\}$ we choose the subfamily which maximizes the corresponding sum and denote it by $\mathcal{G}_{i}$.
By Proposition \ref{Prop.5_useful} we have $\mathcal{G}_{i} \cap \mathcal{G}_{i'} = \emptyset$ if $i \neq i'$.
Finally, we set $\mathcal{B}^{j}_{1}(S):=\cup \{\mathcal{G}_{i}\}$, where the union is taken over all odd $i \in \{1,...,j-1\}$. This clearly
gives \eqref{eq.sum_splitting}.
On the other hand, it is clear that
\begin{equation}
\label{eq.sum_splitting'}
\sum\limits_{B \in \mathcal{B}^{j}_{1}(S)}\frac{\mu(B)}{(r(B))^{p}}\Bigl(\widetilde{\mathcal{E}}_{\mathfrak{m}_{k(B)}}(f,cB)\Bigr)^{p} \le
\sup \sum\limits_{B \in \mathcal{B}} \frac{\mu(B)}{(r(B))^{p}}\Bigl(\widetilde{\mathcal{E}}_{\mathfrak{m}_{k(B)}}(f,cB)\Bigr)^{p},
\end{equation}
where the supremum in \eqref{eq.sum_splitting'} is taken over all $(S,c)$-Whitney families $\mathcal{B}$.

\textit{Step 2.} Given $j \in \mathbb{N}$, by Lemma \ref{Lm.combinatorial_result} there is a disjoint $(S,c)$-nice family $\mathcal{B}^{j}_{2}(S)$ such that
\begin{equation}
\label{eq.sum_splitting''}
\sum\limits_{B \in \mathcal{B}^{j}_{2}(S)}\frac{\mu(B)}{(r(B))^{p}}\Bigl(\widetilde{\mathcal{E}}_{\mathfrak{m}_{k(B)}}(f,cB)\Bigr)^{p} \geq \frac{1}{\underline{N}}
\sum\limits_{\substack{\alpha \in \mathcal{A}_{j}(S)}}\frac{\mu(B_{j,\alpha})}{\epsilon^{jp}}\Bigl(\widetilde{\mathcal{E}}_{\mathfrak{m}_{j}}(f,c\widetilde{B}_{j,\alpha})\Bigr)^{p}.
\end{equation}
By Definitions \ref{Def.sc_good} and \ref{Def.Br_Shv} we have
\begin{equation}
\label{eq.sum_splitting'''}
\sum\limits_{B \in \mathcal{B}^{j}_{2}(S)}\frac{\mu(B)}{(r(B))^{p}}\Bigl(\widetilde{\mathcal{E}}_{\mathfrak{m}_{k(B)}}(f,cB)\Bigr)^{p}
\le \operatorname{BSN}^{2\epsilon^{j}}_{p,\{\mathfrak{m}_{k}\},c}(f).
\end{equation}

\textit{Step 3.}
Using Lemma \ref{Lm.derivappseq} and \eqref{eq.sum_splitting}, \eqref{eq.sum_splitting''} we obtain, for any large enough $j \in \mathbb{N}$,
\begin{equation}
\begin{split}
\notag
&\int\limits_{U_{0}(S)}\Bigl(\operatorname{lip}f^{j}(x)\bigr)^{p}\,d\mu(x) \le C \sum\limits_{B \in \mathcal{B}^{j}_{1}(S) \cup \mathcal{B}^{j}_{2}(S)}\frac{\mu(B)}{(r(B))^{p}}\Bigl(\widetilde{\mathcal{E}}_{\mathfrak{m}_{k(B)}}(f,cB)\Bigr)^{p}.\\
\end{split}
\end{equation}
Combining this inequality with \eqref{eq.sum_splitting'}, \eqref{eq.sum_splitting'''} and \eqref{eq.kostyl} we deduce
\begin{equation}
\begin{split}
\label{eqq.435}
\varliminf\limits_{j \to \infty}\|\operatorname{lip}f^{j}|L_{p}(U_{0}(S))\|^{p}  \le C\operatorname{N}_{p,\{\mathfrak{m}_{k}\},c}(f).
\end{split}
\end{equation}

\textit{Step 4.} By Remark \ref{Rem.support_steps} and \eqref{eq.specaproxsequence} it follows that $f_{1}(x)=f^{j}(x)$ for each $j \in \mathbb{N}$ and all $x \in \operatorname{X} \setminus U_{0}(S)$. Hence, using
Lemma \ref{Lm.derivappseq'} we get
\begin{equation}
\label{eqq.436}
\varliminf\limits_{j \to \infty}\|\operatorname{lip}f^{j}|L_{p}(\operatorname{X} \setminus U_{0}(S))\|^{p}=
\|\operatorname{lip}f_{1}|L_{p}(\operatorname{X} \setminus U_{0}(S))\|^{p}
\le C \|f|L_{p}(\mathfrak{m}_{0})\|^{p}.
\end{equation}

\textit{Step 5.} Combining \eqref{eqq.435}, \eqref{eqq.436} and taking into account \eqref{eq.main2''}, \eqref{eq.kostyl} we obtain the required inequality
\eqref{eq.5.9} and complete the proof.

\end{proof}

The finiteness of $\operatorname{N}_{p,\{\mathfrak{m}_{k}\},c}(f)$ allows one to establish interesting convergence properties of the sequence $\{f^{j}\}$.
More precisely, the following assertion holds.

\begin{Th}
\label{Th.converge1}
If $\operatorname{BSN}^{\delta}_{p,\{\mathfrak{m}_{k}\},c}(f) < +\infty$ for some $c \geq \frac{3}{\epsilon}$
and $\delta \in (0,1]$, then:

\begin{itemize}
\item[\(\rm (i)\)] $\{f^{j}\}$ converges $\mathfrak{m}_{0}$-a.e. to $f$ and converges everywhere to $\operatorname{Ext}_{S,\{\mathfrak{m}_{k}\}}(f)$ on $\operatorname{X} \setminus S$;

\item[\(\rm (ii)\)] $\|f^{j}-\operatorname{Ext}_{S,\{\mathfrak{m}_{k}\}}(f)|L_{p}(\operatorname{X})\| \to 0$, $j \to \infty$;

\item[\(\rm (iii)\)] for each $k \in \mathbb{N}$, $\|f^{j}-f|L_{p}(\mathfrak{m}_{k})\| \to 0$, $j \to \infty$.
\end{itemize}
\end{Th}

\begin{proof}
By Proposition \ref{Prop.stac_sequence} and Definition \ref{Def.Extension_operator}, the sequence $\{f^{j}\}$ converges everywhere on $\operatorname{X} \setminus S$ to the function
$\operatorname{Ext}_{S,\{\mathfrak{m}_{k}\}}(f)$.
Now we recall Definition \ref{Def.multiweight_Lebesgue} and fix an arbitrary point $\underline{x} \in \mathfrak{R}_{\{\mathfrak{m}_{k}\},\epsilon}(f)$.
Since $\mathfrak{R}_{\{\mathfrak{m}_{k}\},\epsilon}(f) \subset S$, by \eqref{eq.aproxsequence}, \eqref{eq.specaproxsequence} it follows that $f^{j}(\underline{x})=f_{j}(\underline{x})$ for all $j \in \mathbb{N}$.
Hence, by \eqref{eq.partunity1}, \eqref{eq.5.1} and Theorem \ref{Th.doubling_type}, for every $j \geq 2$,
\begin{equation}
\begin{split}
&|f(\underline{x})-f^{j}(\underline{x})| \le \sum\limits_{\alpha \in \mathcal{A}_{j}(S)}\varphi_{j,\alpha}(\underline{x})|f(\underline{x})-f_{j,\alpha}|\\
&\le
\sum\limits_{\alpha \in \mathcal{A}_{j}(S)}\chi_{\widetilde{B}_{j,\alpha}}(\underline{x})\fint\limits_{\frac{3}{\epsilon}\widetilde{B}_{j,\alpha}}|f(\underline{x})-f(y)|\,d\mathfrak{m}_{j}(y) \le C \fint\limits_{B_{j-2}(\underline{x})}|f(\underline{x})-f(y)|\,d\mathfrak{m}_{j}(y).
\end{split}
\end{equation}
Hence, using \eqref{M.4} we get
\begin{equation}
\notag
|f(\underline{x})-f^{j}(\underline{x})| \le C \fint\limits_{B_{j-2}(\underline{x})}|f(\underline{x})-f(y)|\,d\mathfrak{m}_{j-2}(y) \to 0, \quad j \to \infty.
\end{equation}
Combining this observation with Theorem \ref{Th.trueLebesgue} we arrive at assertion (i).

Now let us verify (ii). Given $i \geq 2$, by Lemma
\ref{Lm.combinatorial_result} there exists an $(S,\frac{3}{\epsilon},2\epsilon^{i-1}$)-nice family $\mathcal{B}$ such that
\begin{equation}
\label{eqq.737}
\sum\limits_{\substack{\widetilde{B}_{i-1,\alpha} \cap U_{i-2}(S) \neq \emptyset}}\mu(\widetilde{B}_{i-1,\alpha})\Bigl(\widetilde{\mathcal{E}}_{\mathfrak{m}_{i-1}}(f,\frac{3}{\epsilon}\widetilde{B}_{i-1,\alpha})\Bigr)^{p} \le C
\sum\limits_{\substack{B \in \mathcal{B}}}\mu(B)\Bigl(\widetilde{\mathcal{E}}_{\mathfrak{m}_{i-1}}(f,\frac{3}{\epsilon}B)\Bigr)^{p}.
\end{equation}
Given $i \geq 2$, by Remark \ref{Rem.support_steps} we have $\|\operatorname{St}_{i}[f]|L_{p}(\operatorname{X})\| = \|\operatorname{St}_{i}[f]|L_{p}(U_{i-2}(S))\|$.
Hence, we apply Lemma \ref{Lm.steps_estimate} with $\nu=\mu$, $E=U_{i-2}(S)$, use \eqref{eqq.737} and take into account Definition \ref{Def.Br_Shv} and Remark \ref{Rem.31'}. As a result, given $i \geq 2$, we get
\begin{equation}
\notag
\|\operatorname{St}_{i}[f]|L_{p}(\operatorname{X})\|  \le C \epsilon^{i} \operatorname{BSN}^{2\epsilon^{i-1}}_{p,\{\mathfrak{m}_{k}\},c}(f).
\end{equation}
Thus, by the triangle inequality, for each $l \in \mathbb{N}$ with $\epsilon^{l-1} \le \delta$ and any $m > l$,
\begin{equation}
\notag
\|f^{l}-f^{m}|L_{p}(\operatorname{X})\| \le \sum\limits_{i=l+1}^{m}\|\operatorname{St}_{i}[f]|L_{p}(\operatorname{X})\| \le C \epsilon^{l} \operatorname{BSN}^{\delta}_{p,\{\mathfrak{m}_{k}\},c}(f).
\end{equation}
Consequently, $\|f^{l}-f^{m}|L_{p}(\operatorname{X})\| \to 0$, $l,m \to \infty$.  Since $L_{p}(\operatorname{X})$ is complete, there exists $F \in L_{p}(\operatorname{X})$
such that $\|F-f^{j}|L_{p}(\operatorname{X})\| \to 0$, $j \to \infty$. The classical arguments give the existence of a subsequence $\{f^{j_{s}}\}$ converging
$\mu$-a.e. to $F$. On the other hand, the measure $\mu$ is absolutely continuous with respect to the measure $\mathfrak{m}_{0}$, and, consequently,
by the already proved assertion (i) we get $F=\operatorname{Ext}_{S,\{\mathfrak{m}_{k}\}}(f)$ in $\mu$-a.e.\ sense. This proves the claim.

To establish (iii) we fix $k \in \mathbb{N}_{0}$. Given $i \geq 2$, by Lemma \ref{Lm.combinatorial_result} there exists
an $(S,\frac{3}{\epsilon},2\epsilon^{i-1}$)-nice family $\mathcal{B}$ such that
\begin{equation}
\label{eqq.738}
\sum\limits_{\substack{\widetilde{B}_{i-1,\alpha} \cap U_{i-2}(S) \neq \emptyset}}
\mathfrak{m}_{k}(\widetilde{B}_{i-1,\alpha})\Bigl(\widetilde{\mathcal{E}}_{\mathfrak{m}_{i-1}}(f,\frac{3}{\epsilon}\widetilde{B}_{i-1,\alpha})\Bigr)^{p} \le C \sum\limits_{\substack{B \in \mathcal{B}}}\mathfrak{m}_{k}(B)\Bigl(\widetilde{\mathcal{E}}_{\mathfrak{m}_{i-1}}(f,\frac{3}{\epsilon}B)\Bigr)^{p}.
\end{equation}
For each $i \in \mathbb{N}$, $i \geq 2$ we apply Lemma \ref{Lm.steps_estimate} with $\nu=\mathfrak{m}_{k}$, $E=U_{i-2}(S)$, use \eqref{eqq.738} and take into account
Remark \ref{Rem.support_steps}. This gives
\begin{equation}
\notag
\|\operatorname{St}_{i}[f]|L_{p}(\mathfrak{m}_{k})\|^{p} \le C \sum\limits_{\substack{B \in \mathcal{B}}}\mathfrak{m}_{k}(B)\Bigl(\widetilde{\mathcal{E}}_{\mathfrak{m}_{i-1}}(f,\frac{3}{\epsilon}B)\Bigr)^{p}.
\end{equation}
Thus, by \eqref{M.2}, Definition \ref{Def.Br_Shv} and Remark \ref{Rem.31'}, for each $i \geq 2$,
\begin{equation}
\begin{split}
\notag
&\|\operatorname{St}_{i}[f]|L_{p}(\mathfrak{m}_{k})\|^{p} \le C \epsilon^{(p-\theta) i} \Bigl(\operatorname{BSN}^{2\epsilon^{i-1}}_{p,\{\mathfrak{m}_{k}\},c}(f)\Bigr)^{p}.
\end{split}
\end{equation}
Since $\theta \in [0,p)$, given $\delta \in (0,1]$,
an application of the triangle inequality gives, for all large enough $l \in \mathbb{N}$ and all $m > l$,
\begin{equation}
\notag
\|f^{l}-f^{m}|L_{p}(\mathfrak{m}_{k})\| \le \sum\limits_{i=l+1}^{m}\|\operatorname{St}_{i}(f)|L_{p}(\mathfrak{m}_{k})\| \le C \epsilon^{\frac{(p-\theta)l}{p}} \operatorname{BSN}^{\delta}_{p,\{\mathfrak{m}_{k}\},c}(f).
\end{equation}
Since $L_{p}(\mathfrak{m}_{k})$ is complete, there is $h \in L_{p}(\mathfrak{m}_{k})$
such that $\|h-f^{j}|L_{p}(\mathfrak{m}_{k})\| \to 0$, $j \to \infty$. Thus, there is a subsequence $\{f^{j_{s}}\}$ converging
$\mathfrak{m}_{k}$-a.e.\ to $h$.
Combining this fact with the above assertion (i) and \eqref{M.4} we get $h=f$ in $\mathfrak{m}_{k}$-a.e.\ sense and conclude the proof of  (iii).

The proof is complete.
\end{proof}

While the finiteness of $\operatorname{BSN}^{\delta}_{p,\{\mathfrak{m}_{k}\},c}(f)$ for small $\delta > 0$ is sufficient to control
convergence of the sequence $\{f^{j}\}$ in $L_{p}(\operatorname{X})$-sense, it is still not sufficiently powerful to obtain an appropriate estimate of $L_{p}(\operatorname{X})$-norm
of the limit.

\begin{Th}
\label{Th.lpnorm_estimate}
Given $c \geq \frac{3}{\epsilon}$, there exists a constant $C > 0$ such that
\begin{equation}
\|\operatorname{Ext}_{S,\{\mathfrak{m}_{k}\}}(f)|L_{p}(\operatorname{X})\| \le C \operatorname{BSN}_{p,\{\mathfrak{m}_{k}\},c}(f) \quad \text{for all} \quad f \in L^{loc}_{1}(\{\mathfrak{m}_{k}\}).
\end{equation}
\end{Th}

\begin{proof}

The arguments used in the proof of assertion (ii) of Theorem \ref{Th.converge1} give, for any $m \geq 2$,
\begin{equation}
\notag
\|f^{m}-f^{1}|L_{p}(\operatorname{X})\| \le C \operatorname{BSN}_{p,\{\mathfrak{m}_{k}\},c}(f).
\end{equation}
On the other hand, by Lemma \ref{Lm.derivappseq'} we have
\begin{equation}
\notag
\|f^{1}|L_{p}(\operatorname{X})\| \le C \|f|L_{p}(\mathfrak{m}_{0})\|.
\end{equation}
As a result, by the triangle inequality and assertion $(ii)$ of Theorem \ref{Th.converge1} we obtain
\begin{equation}
\notag
\|\operatorname{Ext}_{S,\{\mathfrak{m}_{k}\}}(f)|L_{p}(\operatorname{X})\|=\lim\limits_{l \to \infty}\|f^{l}|L_{p}(\operatorname{X})\| \le C \operatorname{BSN}_{p,\{\mathfrak{m}_{k}\},c}(f).
\end{equation}
The proof is complete.
\end{proof}
Unfortunately, it is difficult to estimate $\|\operatorname{Ext}_{S,\{\mathfrak{m}_{k}\}}(f)|L_{p}(\operatorname{X})\|$ from above in terms of $\operatorname{N}_{p,\{\mathfrak{m}_{k}\},c}(f)$
with some constant $C > 0$ independent on $f$. On the other hand, we have a weaker result which, however, will be sufficient for our purposes.

\begin{Ca}
\label{Ca.weaklp_estimate}
Given $c \geq \frac{3}{\epsilon}$, for each $f \in L_{p}(\mathfrak{m}_{0})$, there exists a constant $C_{f} > 0$ such that
\begin{equation}
\|\operatorname{Ext}_{S,\{\mathfrak{m}_{k}\}}(f)|L_{p}(\operatorname{X})\| \le C_{f} \operatorname{N}_{p,\{\mathfrak{m}_{k}\},c}(f).
\end{equation}
\end{Ca}

\begin{proof}
If $\operatorname{N}_{p,\{\mathfrak{m}_{k}\},c}(f)=+\infty$ one can put $C_{f}=1$.
If $f \in L_{p}(\mathfrak{m}_{0})$ and $\operatorname{N}_{p,\{\mathfrak{m}_{k}\},c}(f) < +\infty$, then by \eqref{eq.kostyl} there is $\delta=\delta(f) \in (0,1)$
such that $\operatorname{BSN}^{\delta}_{p,\{\mathfrak{m}_{k}\},c}(f) < +\infty$. Hence, by Lemma \ref{Lm.different_scales} we have $\operatorname{BSN}_{p,\{\mathfrak{m}_{k}\},c}(f) < +\infty$.
This fact in combination with Theorem \ref{Th.lpnorm_estimate} proves the claim.
\end{proof}

Now we are ready to prove the key result of this section. We recall Definition \ref{Def.Cheeger_Sobolev}.

\begin{Th}
\label{Th.firstext}
If $\operatorname{N}_{p,\{\mathfrak{m}_{k}\},c}(f) < +\infty$ for some $c \geq \frac{3}{\epsilon}$, then $\operatorname{Ext}_{S,\{\mathfrak{m}_{k}\}}(f) \in W_{p}^{1}(\operatorname{X})$.
Furthermore, there exists a constant $C > 0$ such that
\begin{equation}
\label{eq.reversetrace}
\operatorname{Ch}_{p}(\operatorname{Ext}_{S,\{\mathfrak{m}_{k}\}}(f)) \le C \operatorname{N}_{p,\{\mathfrak{m}_{k}\},c}(f) \quad \text{for all} \quad f \in L^{loc}_{1}(\{\mathfrak{m}_{k}\}).
\end{equation}
\end{Th}

\begin{proof}

If $\operatorname{N}_{p,\{\mathfrak{m}_{k},c\}}(f)=+\infty$, then inequality \eqref{eq.reversetrace} is obvious.
If $\operatorname{N}_{p,\{\mathfrak{m}_{k}\},c}(f)<+\infty$, then
by Theorem \ref{Th.converge1} and Corollary \ref{Ca.weaklp_estimate} we have $\operatorname{Ext}_{S,\{\mathfrak{m}_{k}\}}(f) \in L_{p}(\operatorname{X})$
and $f^{j} \to \operatorname{Ext}_{S,\{\mathfrak{m}_{k}\}}(f)$, $j \to \infty$ in $L_{p}(\operatorname{X})$-sense. Furthermore, by Theorem \ref{Th.derivappseq}
\begin{equation}
\notag
\operatorname{Ch}_{p}(\operatorname{Ext}_{S,\{\mathfrak{m}_{k}\}}(f)) \le \varliminf\limits_{j \to \infty}\|\operatorname{lip}f^{j}|L_{p}(\operatorname{X})\| \le C \operatorname{N}_{p,\{\mathfrak{m}_{k}\},c}(f).
\end{equation}
By Definition \ref{Def.Cheeger_Sobolev} this implies that $F \in W_{p}^{1}(\operatorname{X})$ and \eqref{eq.reversetrace} holds.

The proof is complete.

\end{proof}

\section{Comparison of different trace functionals}

The aim of this section is to compare the functionals $\operatorname{CN}_{p,\{\mathfrak{m}_{k}\}}$, $\operatorname{BN}_{p,\{\mathfrak{m}_{k}\},\sigma}$, $\operatorname{BSN}_{p,\{\mathfrak{m}_{k}\},c}$ and
$\operatorname{N}_{p,\{\mathfrak{m}_{k}\},c}$, respectively.
Recall that initially these functionals are defined on the space $L_{1}^{loc}(\{\mathfrak{m}_{k}\})$ and take
their values in $[0,+\infty]$.
The following \textit{data are assumed to be fixed}:

\begin{itemize}
\item[\((\rm \textbf{D.8.1})\)] a parameter $p \in (1,\infty)$ and a metric measure space $\operatorname{X}=(\operatorname{X},\operatorname{d},\mu) \in \mathfrak{A}_{p}$;

\item[\((\rm \textbf{D.8.2})\)] a parameter $\theta \in [0,p)$ and a closed set $S \in \mathcal{LCR}_{\theta}(\operatorname{X})$;

\item[\((\rm \textbf{D.8.3})\)] a sequence of measures $\{\mathfrak{m}_{k}\} \in \mathfrak{M}_{\theta}(S)$ with parameter $\epsilon=\epsilon(\{\mathfrak{m}_{k}\}) \in (0,\frac{1}{10}]$.
\end{itemize}

\textit{The first keystone result} of this section is straightforward. We recall Definitions \ref{Def.Br_Shv}, \ref{Def.kostyl}.

\begin{Th}
\label{Th.comparison1}
For each $c \geq 1$, $\operatorname{N}_{p,\{\mathfrak{m}_{k}\},c}(f) \le 2 \operatorname{BSN}_{p,\{\mathfrak{m}_{k}\},c}(f)$ for all $f \in L_{1}^{loc}(\{\mathfrak{m}_{k}\}).$
\end{Th}

\begin{proof}
In the case $\operatorname{BSN}_{p,\{\mathfrak{m}_{k}\},c}(f)=+\infty$ the inequality is trivial. We fix $f \in L_{1}^{loc}(\{\mathfrak{m}_{k}\})$ such that $\operatorname{BSN}_{p,\{\mathfrak{m}_{k}\},c}(f) < +\infty$.
Since each $(S,c)$-Whitney family is an $(S,c)$-nice family, we have
\begin{equation}
\notag
\sup\Bigl(\sum\limits_{B \in \mathcal{B}} \frac{\mu(B)}{(r(B))^{p}}\Bigl(\widetilde{\mathcal{E}}_{\mathfrak{m}_{k(B)}}(f,cB)\Bigr)^{p}\Bigr)^{\frac{1}{p}} \le
\operatorname{BSN}_{p,\{\mathfrak{m}_{k}\},c}(f),
\end{equation}
where the supremum in the left-hand side is taken over all $(S,c)$-Whitney families $\mathcal{B}$.
On the other hand, by Remark \ref{Rem.31'} we have $\varliminf_{\delta \to 0}\operatorname{BSN}^{\delta}_{p,\{\mathfrak{m}_{k}\},c}(f) \le \operatorname{BSN}_{p,\{\mathfrak{m}_{k}\},c}(f)$.
Combining these observations, we get the required estimate and complete the proof.
\end{proof}

To go further we recall notation adopted at the beginning of Section 5. Furthermore, we set $B_{k}(x):=B_{\epsilon^{k}}(x)$ for all $k \in \mathbb{Z}$ and $x \in \operatorname{X}$.

\begin{Prop}
\label{Prop.41}
Let $c \geq 1$ and let $\underline{k} \in \mathbb{N}_{0}$ be the smallest $k' \in \mathbb{N}_{0}$ satisfying $\epsilon^{k'} \le \frac{1}{2c}$.
Then there exists a constant $C > 0$ depending on $p$, $\theta$, $\underline{k}$, $c$ and $\mathcal{C}_{\{\mathfrak{m}_{k}\}}$ such that
if $k \geq \underline{k}$, $r \in (\epsilon^{k+1},\epsilon^{k}]$ and $B_{r}(y)$ is a closed ball such that $S \cap B_{cr}(y) \neq \emptyset$ and $B_{r}(x) \subset B_{cr}(y)$ for some $x \in S \cap B_{cr}(y)$, then
\begin{equation}
\label{eq.42}
\mathcal{E}_{\mathfrak{m}_{k}}(f,B_{cr}(y)) \le C \inf\limits_{z \in  B_{cr}(y)} \mathcal{E}_{\mathfrak{m}_{k-\underline{k}}}(f,B_{k-\underline{k}}(z)) \quad \text{for all} \quad f \in L_{1}^{loc}(\{\mathfrak{m}_{k}\}).
\end{equation}
\end{Prop}

\begin{proof}
We fix a ball $B_{r}(x) \subset B_{cr}(y)$ with $x \in S \cap cB_{r}(y)$. Since $\frac{r}{\epsilon} \geq \epsilon^{k}$, we have
\begin{equation}
\label{eq.53''}
B_{k-\underline{k}}(z) \subset (2c+\frac{1}{\epsilon^{\underline{k}+1}})B_{r}(x) \quad \text{and} \quad cB_{r}(y) \subset  B_{k-\underline{k}}(z) \quad \text{for all $z \in cB_{r}(y)$}.
\end{equation}
Since $x \in S$, by \eqref{eq.53''} and Theorem \ref{Th.doubling_type} we obtain
\begin{equation}
\notag
\sup\limits_{z \in  B_{cr}(y)}\mathfrak{m}_{k}(B_{k-\underline{k}}(z)) \le \mathfrak{m}_{k}((2c+\frac{1}{\epsilon^{\underline{k}+1}})B_{r}(x)) \le C \mathfrak{m}_{k}(B_{r}(x)) \le C \mathfrak{m}_{k}(cB_{r}(y)).
\end{equation}
Hence, by the second inclusion in \eqref{eq.53''}, Remark \ref{Rem.best_approx}  and \eqref{M.4}, for each $z \in B_{cr}(y)$, we have
\begin{equation}
\label{eq.43}
\begin{split}
&\mathcal{E}_{\mathfrak{m}_{k}}(f,cB_{r}(y)) \le \Bigl(\frac{1}{\mathfrak{m}_{k}(cB_{r}(y))}\Bigr)^{2}\int\limits_{B_{k-\underline{k}}(z)}\int\limits_{B_{k-\underline{k}}(z)}|f(v)-f(w)|\,d\mathfrak{m}_{k}(v)\,d\mathfrak{m}_{k}(w)\\
&\le C \fint\limits_{B_{k-\underline{k}}(z)}\fint\limits_{B_{k-\underline{k}}(z)}|f(v)-f(w)|\,d\mathfrak{m}_{k-\underline{k}}(v)\,d\mathfrak{m}_{k-\underline{k}}(w)
\le C \mathcal{E}_{\mathfrak{m}_{k-\underline{k}}}(f,B_{k-\underline{k}}(z)).
\end{split}
\end{equation}
Since $z \in B_{cr}(y)$ was chosen arbitrarily, the claim follows.
\end{proof}

We recall that in Theorem \ref{Th.SecondMain} we defined $\operatorname{CN}_{p,\{\mathfrak{m}_{k}\}}(f):=\mathcal{CN}_{p,\{\mathfrak{m}_{k}\}}(f)+\|f|L_{p}(\mathfrak{m}_{0})\|$
for $f \in L_{1}^{loc}(\{\mathfrak{m}_{k}\})$.
The following assertion is \textit{the second keystone result} of this section.

\begin{Th}
\label{Theore.53}
Given $c \geq 1$, there exists a constant $C > 0$ such that
\begin{equation}
\label{eq.52'}
\operatorname{BSN}_{p,\{\mathfrak{m}_{k}\},c}(f) \le C \operatorname{CN}_{p,\{\mathfrak{m}_{k}\}}(f) \quad \text{for all} \quad f \in L_{1}^{loc}(\{\mathfrak{m}_{k}\}).
\end{equation}
Furthermore, for each $\delta \in (0,\frac{1}{4c}]$, there is a constant $C > 0$ (depending on $\delta$) such that
\begin{equation}
\label{eq.53'}
\Bigl|\operatorname{BSN}^{\delta}_{p,\{\mathfrak{m}_{k}\},c}(f)-\|f|L_{p}(\mathfrak{m}_{0})\|\Bigr| \le C \|f^{\sharp}_{\{\mathfrak{m}_{k}\}}|L_{p}(U_{(c+1)\delta}(S))\|.
\end{equation}
\end{Th}

\begin{proof}
We put $\overline{\delta}:=\frac{1}{4c}$. Let $\underline{k}$ be the smallest $k \in \mathbb{N}_{0}$ sutisfying $\epsilon^{k} \le \overline{\delta}$.

We start with the second claim. We fix $\delta \in (0,\overline{\delta}]$. Without loss of generality we may assume that $f^{\sharp}_{\{\mathfrak{m}_{k}\}} \in L_{p}(U_{(c+1)\delta}(S))$ because otherwise inequality \eqref{eq.53'} is trivial.
Let $\mathcal{B}^{\delta}$ be an arbitrary $(S,c,\delta)$-nice family of closed balls. Since $cB \cap S \neq \emptyset$ for all $B \in \mathcal{B}^{\delta}$ we get
\begin{equation}
\label{eq.54'}
B \subset U_{(c+1)\delta}(S) \quad \text{for all} \quad B \in \mathcal{B}^{\delta}.
\end{equation}
We recall notation \eqref{eq.strata}. Given $k \geq \underline{k}$ and $B \in \mathcal{B}^{\delta}(k,\epsilon)$, we apply Proposition \ref{Prop.41} with $c$ replaced by $2c$.
This gives $\mathcal{E}_{\mathfrak{m}_{k}}(f,2cB) \le C \inf\limits_{z \in B}\mathcal{E}_{\mathfrak{m}_{k-\underline{k}}}(f,B_{k-\underline{k}}(z)).$
Consequently,  we have
\begin{equation}
\label{eq.55'}
\frac{\mu(B)}{(r(B))^{p}}\Bigl(\mathcal{E}_{\mathfrak{m}_{k(B)}}(f,2cB)\Bigr)^{p} \le
C \int\limits_{B}(f^{\sharp}_{\{\mathfrak{m}_{k}\}}(y))^{p}\,d\mu(y).
\end{equation}
Combining  \eqref{eq.54'}, \eqref{eq.55'} and taking into account that $\mathcal{B}^{\delta}$ is a disjoint family of balls we obtain
\begin{equation}
\notag
\begin{split}
&\sum\limits_{B \in \mathcal{B}^{\delta}}\frac{\mu(B)}{(r(B))^{p}}\Bigl(\mathcal{E}_{\mathfrak{m}_{k(B)}}(f,2cB)\Bigr)^{p} \le C
\sum\limits_{k \geq \underline{k}}\sum\limits_{B \in \mathcal{B}^{\delta}(k,\epsilon)} \int\limits_{B}(f^{\sharp}_{\{\mathfrak{m}_{k}\}}(y))^{p}\,d\mu(y)\\
&\le C \int\limits_{U_{(c+1)\delta}(S)}(f^{\sharp}_{\{\mathfrak{m}_{k}\}}(x))^{p}\,d\mu(x).
\end{split}
\end{equation}
Since $\mathcal{B}^{\delta}$ was chosen arbitrarily, the claim follows from \eqref{eq.main2''}.

To prove \eqref{eq.52'}, given an $(S,c)$-nice family of closed balls $\mathcal{B}$, we split it into two subfamilies. The first one $\mathcal{B}^{1}$ consists of the balls of radii
greater or equal than $\overline{\delta}$ and the second one is $\mathcal{B}^{2}:=\mathcal{B} \setminus \mathcal{B}^{1}$. Applying Lemma \ref{Lm.largescale_estimate}
we obtain
\begin{equation}
\label{eq.56'}
\begin{split}
&\sum\limits_{B \in \mathcal{B}^{1}}\frac{\mu(B)}{(r(B))^{p}}\Bigl(\mathcal{E}_{\mathfrak{m}_{k(B)}}(f,2cB)\Bigr)^{p} \le C \int\limits_{S}|f(y)|^{p}\,d\mathfrak{m}_{0}(y).
\end{split}
\end{equation}
On the other hand, using \eqref{eq.53'} just proved and taking into account \eqref{eq.main1} we get
\begin{equation}
\label{eq.57'}
\sum\limits_{B \in \mathcal{B}^{2}}\frac{\mu(B)}{(r(B))^{p}}\Bigl(\mathcal{E}_{\mathfrak{m}_{k(B)}}(f,2cB)\Bigr)^{p} \le C \Bigl(\mathcal{CN}_{p,\{\mathfrak{m}_{k}\}}(f)\Bigr)^{p}.
\end{equation}
It remains to combine \eqref{eq.56'}, \eqref{eq.57'} and take into account that $\mathcal{B}$ was chosen arbitrarily.
This proves the first claim.

The proof is complete.
\end{proof}


The following lemma is an important ingredient for comparing the functionals $\operatorname{N}_{p,\{\mathfrak{m}_{k}\},c}$ and $\operatorname{BN}_{p,\{\mathfrak{m}_{k}\},\sigma}$.
Recall that in Theorem \ref{Th.SecondMain} we set $\operatorname{BN}_{p,\{\mathfrak{m}_{k}\},\sigma}(f):=\|f|L_{p}(\mathfrak{m}_{0})\|+\mathcal{BN}_{p,\{\mathfrak{m}_{k}\},\sigma}(f)$ for $f \in  L_{1}^{loc}(\{\mathfrak{m}_{k}\})$.

\begin{Lm}
\label{Lm.admis_estimate}
For each $c \geq 1$, for every $\sigma \in (0,\frac{\epsilon^{2}}{4c})$, there exists $C > 0$
(depending on $\sigma$) such that, for each $(S,c)$-Whitney family of closed balls $\mathcal{B}$,
\begin{equation}
\label{eq.admis_estimate}
\begin{split}
&\Sigma:=\sum\limits_{B \in \mathcal{B}}\frac{\mu(B)}{(r(B))^{p}}\Bigl(\mathcal{E}_{\mathfrak{m}_{k(B)}}(f,2cB)\Bigr)^{p}
\le C \Bigl(\operatorname{BN}_{p,\{\mathfrak{m}_{k}\},\sigma}(f)\Bigr)^{p} \quad \text{for all} \quad  f \in  L_{1}^{loc}(\{\mathfrak{m}_{k}\}).
\end{split}
\end{equation}
\end{Lm}

\begin{proof} We put $\overline{\sigma}:=\frac{\epsilon^{2}}{4c}$ and let $\underline{k}$ be the smallest $k \in \mathbb{N}_{0}$
satisfying $\epsilon^{k} \le \frac{1}{4c}$.
Since $\mathcal{B}$ is an $(S,c)$-Whitney family we have $cB \cap S \neq \emptyset$ for all $B \in \mathcal{B}$.
We recall notation \eqref{eq.strata}. Given $k \geq \underline{k}$ and $B \in \mathcal{B}(k,\epsilon)$,
it is easy to see that $B \subset B_{(2c+1)r(B)}(x) \subset B_{\epsilon^{k-\underline{k}}}(x)$ for all $x \in 2cB \cap S$.
If $B \in \mathcal{B}(k,\epsilon)$ for some $k \geq \underline{k}$, then $r(B) \geq \frac{\epsilon^{2}}{4c}\epsilon^{k-\underline{k}}$.
Hence, for each $k \geq \underline{k}$, and every $\sigma \in (0,\overline{\sigma}]$, we have
\begin{equation}
\label{eq.513'}
2cB \cap S \subset S_{k-\underline{k}}(\sigma) \quad \text{for all} \quad B \in \mathcal{B}(k,\epsilon).
\end{equation}
Given $k \geq \underline{k}$ and $B \in \mathcal{B}(k,\epsilon)$, by Proposition \ref{Prop.41} (applied with $2c$ instead of $c$) we have
\begin{equation}
\label{eqq.812}
\mathcal{E}_{\mathfrak{m}_{k}}(f,2cB) \le C \inf\limits_{y \in 2cB}\mathcal{E}_{\mathfrak{m}_{k-\underline{k}}}(f,B_{k-\underline{k}}(y)).
\end{equation}
By \eqref{M.3} and \eqref{M.4} we have $\mu(B)/(r(B))^{p} \le C \mathfrak{m}_{k-\underline{k}}(B) \le C \mathfrak{m}_{k-\underline{k}}(2cB \cap S)$.
This estimate in combination with \eqref{eqq.812} leads to
\begin{equation}
\notag
\frac{\mu(B)}{(r(B))^{p}}\Bigl(\mathcal{E}_{\mathfrak{m}_{k}}(f,2cB)\Bigr)^{p} \le C \int\limits_{2cB \cap S} \Bigl(\mathcal{E}_{\mathfrak{m}_{k-\underline{k}}}(f,B_{k-\underline{k}}(y))\Bigr)^{p}\,d\mathfrak{m}_{k-\underline{k}}(y).
\end{equation}
As a result, we use \eqref{eq.513'} and take into account Propositions \ref{Prop.covering_multiplicity}, \ref{Prop.finite_intersection}. This gives
\begin{equation}
\label{eq.514'}
\begin{split}
&\sum\limits_{k=\underline{k}}^{\infty}\sum\limits_{B \in \mathcal{B}(k,\epsilon)}\frac{\mu(B)}{(r(B))^{p}}\Bigl(\mathcal{E}_{\mathfrak{m}_{k}}(f,2cB)\Bigr)^{p}\\
&\le C \sum\limits_{k=\underline{k}}^{\infty}\int\limits_{S_{k-\underline{k}}(\sigma)} \Bigl(\mathcal{E}_{\mathfrak{m}_{k-\underline{k}}}(f,B_{k-\underline{k}}(y))\Bigr)^{p}\,d\mathfrak{m}_{k-\underline{k}}(y)
\le C \Bigl(\mathcal{BN}_{p,\{\mathfrak{m}_{k}\},\sigma}(f)\Bigr)^{p} .
\end{split}
\end{equation}
On the other hand, by Lemma \ref{Lm.largescale_estimate} we clearly have
\begin{equation}
\label{eq.515'}
\sum\limits_{k=0}^{\underline{k}}\sum\limits_{B \in \mathcal{B}(k,\epsilon)}\frac{\mu(B)}{(r(B))^{p}}\Bigl(\mathcal{E}_{\mathfrak{m}_{k}}(f,2cB))\Bigr)^{p}
 \le C \int\limits_{S}|f(z)|^{p}\,d\mathfrak{m}_{0}(z).
\end{equation}
Combining \eqref{eq.514'} and \eqref{eq.515'}, we obtain \eqref{eq.admis_estimate} and complete the proof.
\end{proof}

Now we are ready to formulate and prove the \textit{third keystone result of this section.}

\begin{Th}
\label{Th.Functional1}
Given $c \geq 1$, there exists a constant $C > 0$ such that, for each $\sigma \in (0,\frac{\epsilon^{2}}{4c})$, the following inequality
\begin{equation}
\label{eq.Functional1}
\operatorname{N}_{p,\{\mathfrak{m}_{k}\},c}(f) \le C \operatorname{BN}_{p,\{\mathfrak{m}_{k}\},\sigma}(f)
\end{equation}
holds for any function $f \in L_{1}^{loc}(\{\mathfrak{m}_{k}\})$ with $\operatorname{CN}_{p,\{\mathfrak{m}_{k}\}}(f) < +\infty$.
\end{Th}

\begin{proof}
We fix an arbitrary $\sigma \in (0,\frac{\epsilon^{2}}{4c})$.
If $\operatorname{CN}_{p,\{\mathfrak{m}_{k}\}}(f) < +\infty$, then by Theorem \ref{Theore.53} we have
\begin{equation}
\begin{split}
\label{eq.517'}
&\varliminf\limits_{\delta \to 0}\operatorname{BSN}^{\delta}_{p,\{\mathfrak{m}_{k}\},c}(f) \le C \Bigl(\varliminf\limits_{\delta \to 0}\|f^{\sharp}_{\{\mathfrak{m}_{k}\}}|L_{p}(U_{(c+1)\delta}(S))\|+\|f|L_{p}(\mathfrak{m}_{0})\|\Bigr)\\
&=C \Bigl(\|f^{\sharp}_{\{\mathfrak{m}_{k}\}}|L_{p}(S)\|+\|f|L_{p}(\mathfrak{m}_{0})\|\Bigr).
\end{split}
\end{equation}
On the other hand, by Lemma \ref{Lm.admis_estimate}, we get
\begin{equation}
\label{eq.518'}
\sup \sum\limits_{B \in \mathcal{B}}\frac{\mu(B)}{(r(B))^{p}}\Bigl(\mathcal{E}_{\mathfrak{m}_{k(B)}}(f,2cB)\Bigr)^{p}
\le C \Bigl(\operatorname{BN}_{p,\{\mathfrak{m}_{k}\},\sigma}(f)\Bigr)^{p},
\end{equation}
where the supremum is taken over all $(S,c)$-Whitney families $\mathcal{B}$.
Collecting estimates \eqref{eq.517'} and \eqref{eq.518'} we obtain \eqref{eq.Functional1} and complete the proof.
\end{proof}

Finally, \textit{the fourth keystone result} of this section reads as follows. We recall \eqref{eq.main1}, \eqref{eq.main3}.
\begin{Th}
\label{Th.Functional2}
For each $\sigma \in (0,1)$ there exists a constant $C > 0$ such that
\begin{equation}
\label{eq.Functional2}
\mathcal{BN}_{p,\{\mathfrak{m}_{k}\},\sigma}(f) \le C  \mathcal{CN}_{p,\{\mathfrak{m}_{k}\}}(f) \quad \text{for all} \quad f \in  L_{1}^{loc}(\{\mathfrak{m}_{k}\}).
\end{equation}
\end{Th}

\begin{proof}
We recall notation \eqref{eqq.notation_porous} and \eqref{eq.lattice}. We put $S_{k}(\sigma):=S_{\epsilon^{k}}(\sigma)$ and $\mathcal{B}_{k}:=\mathcal{B}_{k}(\operatorname{X},\epsilon)$, $k \in \mathbb{N}_{0}$ for brevity.
We set $B_{k}(x):=B_{\epsilon^{k}}(x)$ for $k \in \mathbb{N}_{0}$, $x \in \operatorname{X}$, as usual. For each $k \in \mathbb{N}_{0}$, for any ball $B \in \mathcal{B}_{k}$ with nonempty intersection with $S_{k}(\sigma)$, we fix a point $x_{B} \in B \cap S_{k}(\sigma)$ and a ball $B'=B'(B) \subset B_{k}(x_{B}) \setminus S$ with $r(B') \geq \sigma \epsilon^{k}$.

Since $\epsilon^{-1} \geq 10$, given $k \in \mathbb{N}$ and $B \in \mathcal{B}_{k}$, for any point $z \in B'(B)$ we have $B_{k-1}(z) \supset 2B_{k}(y)$
and $2B_{k}(y) \supset B$ for all $y \in B \cap S_{k}(\sigma)$.
Hence, by Proposition \ref{Prop.41} (applied with $x=y$, $c=2$ and $\underline{k}=1$)
\begin{equation}
\notag
\mathcal{E}_{\mathfrak{m}_{k}}(f,2B_{k}(y)) \le C \mathcal{E}_{\mathfrak{m}_{k-1}}(f,B_{k-1}(z)) \quad \text{for all} \quad y \in B \cap S_{k}(\sigma) \quad \text{and all} \quad z \in B'.
\end{equation}
On the other hand, by Remark \ref{Rem.best_approx} and Theorem \ref{Th.doubling_type} we have $\mathcal{E}_{\mathfrak{m}_{k}}(f,B_{k}(y)) \le C \mathcal{E}_{\mathfrak{m}_{k}}(f,2B_{k}(y))$
with $C > 0$ independent on $f$, $k$ and $y$.
These observations together with the uniformly locally doubling property of $\mu$  and \eqref{M.2} give
\begin{equation}
\label{eq.521'}
\begin{split}
&\epsilon^{k(\theta-p)}\int\limits_{B \cap S_{k}(\sigma)}\Bigl(\mathcal{E}_{\mathfrak{m}_{k}}(f,B_{k}(y))\Bigr)^{p}\,d\mathfrak{m}_{k}(y)
\le \epsilon^{k(\theta-p)}\mathfrak{m}_{k}(B)
\inf\limits_{z \in \frac{1}{2}B'}\Bigl(\mathcal{E}_{\mathfrak{m}_{k-1}}(f,B_{k-1}(z))\Bigr)^{p}\\
&\le C \mu(B')\inf\limits_{z \in \frac{1}{2}B'}\Bigl(f^{\sharp}_{\{\mathfrak{m}_{k}\}}(z)\Bigr)^{p}  \le C \int\limits_{\frac{1}{2}B'}\Bigl(f^{\sharp}_{\{\mathfrak{m}_{k}\}}(z)\Bigr)^{p}\,d\mu(z).
\end{split}
\end{equation}

Given $k \in \mathbb{N}$, the key property of balls $B'(B)$, $B \in \mathcal{B}_{k}$ is that (because $\epsilon \in (0,\frac{1}{10}]$)
\begin{equation}
\label{eq.522'}
\frac{1}{2}B'(B) \subset U_{2\epsilon^{k}}(S) \setminus  U_{\epsilon^{k+1}}(S) \quad \text{for any} \quad B \in \mathcal{B}_{k} \quad \text{with} \quad B \cap S_{k}(\sigma) \neq \emptyset.
\end{equation}
Furthermore, since $\frac{1}{2}B'(B) \subset 3B$ for all $B \in \mathcal{B}_{k}$ and the family $\{\frac{1}{2}B: B \in\mathcal{B}_{k}\}$ is disjoint
by Proposition \ref{Prop.finite_intersection}, we have
\begin{equation}
\label{eq.523'}
\mathcal{M}(\{\frac{1}{2}B'(B):B \in \mathcal{B}_{k} \text{ and } B \cap S_{k}(\sigma) \neq \emptyset\}) \le C.
\end{equation}
Hence, using \eqref{eq.521'}--\eqref{eq.523'} and Proposition \ref{Prop.covering_multiplicity} we obtain, for each $k \in \mathbb{N}$,
\begin{equation}
\notag
\begin{split}
&\epsilon^{k(\theta-p)}\int\limits_{S_{k}(\sigma)}\Bigl(\mathcal{E}_{\mathfrak{m}_{k}}(f,B_{k}(x))\Bigr)^{p}\,d\mathfrak{m}_{k}(x)
\le \epsilon^{k(\theta-p)}\sum\limits_{B \in \mathcal{B}_{k}}\int\limits_{B \cap S_{k}(\sigma)}\Bigl(\mathcal{E}_{\mathfrak{m}_{k}}(f,B_{k}(x))\Bigr)^{p}\,d\mathfrak{m}_{k}(x)\\
&\le
C \sum\limits_{\substack{B \in \mathcal{B}_{k}\\ B \cap S_{k}(\sigma) \neq \emptyset}} \int\limits_{\frac{1}{2}B'(B)}\Bigl(f^{\sharp}_{\{\mathfrak{m}_{k}\}}(y)\Bigr)^{p}\,d\mu(y) \le C \int\limits_{U_{k-1}(S)\setminus U_{k+1}(S)}\Bigl(f^{\sharp}_{\{\mathfrak{m}_{k}\}}(y)\Bigr)^{p}\,d\mu(y).
\end{split}
\end{equation}
As a result, we have
\begin{equation}
\label{eq.527'}
\sum\limits_{k=1}^{\infty}\epsilon^{k(\theta-p)}\int\limits_{S_{k}(\sigma)}\Bigl(\mathcal{E}_{\mathfrak{m}_{k}}(f,B_{k}(x))\Bigr)^{p}\,d\mathfrak{m}_{k}(x)
\le C \int\limits_{U_{0}(S)}\Bigl(f^{\sharp}_{\{\mathfrak{m}_{k}\}}(y)\Bigr)^{p}\,d\mu(y).
\end{equation}
Combining \eqref{eq.527'} with \eqref{eq.main1} and \eqref{eq.main3} we get \eqref{eq.Functional2}, which completes the proof.

\end{proof}

\section{Trace inequalities for the Riesz potentials}

The aim of this section is to establish the Hedberg--Wolff-type inequality. Of course, the results of this section are not
surprising for  experts. Nevertheless, the author has not succeeded in finding a precise reference in the literature.
We present the details for completeness.

Throughout the whole section \textit{we fix the following data}:

\begin{itemize}
\item[\(\rm (\textbf{D.9.1})\)] a parameter $q \in [1,\infty)$ and a metric measure space $\operatorname{X}=(\operatorname{X},\operatorname{d},\mu) \in \mathfrak{A}_{q}$;

\item[\(\rm (\textbf{D.9.2})\)] a measure $\mathfrak{m}$ on $\operatorname{X}$, a parameter $\epsilon \in (0,\frac{1}{10}]$, and a number $\underline{k} \in \mathbb{Z}$.

\end{itemize}

Having at our disposal Proposition \ref{Prop.mms_dyadic} we fix a family $\{Q_{k,\alpha}\}$ of generalized dyadic cubes in $\operatorname{X}$ and introduce the \textit{essential part} of $\operatorname{X}$
by letting $\underline{\operatorname{X}}:= \cap_{k=\underline{k}}^{\infty}\cup_{\alpha \in \mathcal{A}_{k}(\operatorname{X},\epsilon)}Q_{k,\alpha}$.
Furthermore, given a generalized dyadic cube $Q_{k,\alpha}$ in $\operatorname{X}$, we put
\begin{equation}
\label{eq.modif_dyadic}
\widehat{Q}_{k,\alpha}:=\bigcup\{\operatorname{cl}Q_{k,\alpha'}:\operatorname{cl}Q_{k,\alpha'} \cap 5B_{\epsilon^{k}}(z_{k,\alpha}) \neq \emptyset\}.
\end{equation}
From Propositions \ref{Prop.finite_intersection}, \ref{Prop.mms_dyadic} we easily obtain the following assertion.
\begin{Prop}
\label{Prop.5.1}
There exists a constant $C > 0$ such that
\begin{equation}
\notag
\mathcal{M}(\{\widehat{Q}_{k,\alpha}:\alpha \in \mathcal{A}_{k}(\operatorname{X},\epsilon)\}) \le C \quad \text{for all} \quad k \geq \underline{k}.
\end{equation}
\end{Prop}
By \eqref{eq.modif_dyadic} and (DQ2), (DQ3) of Proposition \ref{Prop.mms_dyadic} we immediately get the following result.
\begin{Prop}
\label{Prop.5.2}
For each $k \geq \underline{k}$ and $\alpha \in \mathcal{A}_{k}(\operatorname{X},\epsilon)$
\begin{equation}
\notag
\bigcup \{\widehat{Q}_{j,\beta}:Q_{j,\beta} \subset Q_{k,\alpha}\} \subset \widehat{Q}_{k,\alpha} \quad \text{for all} \quad j \geq k.
\end{equation}
\end{Prop}

Given a point $x \in \underline{\operatorname{X}}$ and a number $k \geq \underline{k}$,
there is a unique $\alpha(x) \in \mathcal{A}_{k}(\operatorname{X},\epsilon)$ for which $x \in Q_{k,\alpha(x)}$. In what follows, we set
\begin{equation}
\label{eq.5.4''}
Q_{k,\alpha}(x):=Q_{k,\alpha(x)}, \qquad \widehat{Q}_{k,\alpha}(x):=\widehat{Q}_{k,\alpha(x)}.
\end{equation}

Given a Borel set $E \subset \operatorname{X}$ with $\mu(E) > 0$, we put
\begin{equation}
\label{eq.specialcoefficients}
a_{\mathfrak{m}}(E):=\frac{\mathfrak{m}(E)}{\mu(E)}\operatorname{diam}E.
\end{equation}
Given $R \in (0,\epsilon^{\underline{k}}]$, the \textit{restricted Riesz potential} of $\mathfrak{m}$ is the mapping $I^{R}[\mathfrak{m}]:\operatorname{X} \to [0,+\infty]$ defined as
\begin{equation}
\label{eq.Riesz}
I^{R}[\mathfrak{m}](x):=\sum_{\epsilon^{k} \le R}a_{\mathfrak{m}}(B_{\epsilon^{k}}(x)), \quad x \in \operatorname{X}.
\end{equation}
We also introduce the \textit{restricted dyadic Riesz potential} of $\mathfrak{m}$ by the formula
\begin{equation}
\label{eq.Riesz_d}
\widehat{I}^{R}[\mathfrak{m}](x):=
\begin{cases}
\sum_{\epsilon^{k} \le R}a_{\mathfrak{m}}(\widehat{Q}_{k,\alpha}(x)), \quad x \in \underline{\operatorname{X}};\\
0, \quad x \in \operatorname{X} \setminus \underline{\operatorname{X}}.
\end{cases}
\end{equation}
Given $p \in (1,\infty)$, we set $p':=\frac{p}{p-1}$. Given a Borel set $E\subset \operatorname{X}$ and a parameter $R \in (0,\epsilon^{\underline{k}}]$, the \textit{restricted energy} and the \textit{restricted dyadic energy} of the measure $\mathfrak{m}$ are defined by
\begin{equation}
\label{eq.energy}
\mathfrak{E}^{R}_{p}[\mathfrak{m}](E):= \int\limits_{E} \Bigl(I^{R}[\mathfrak{m}](x)\Bigr)^{p'}\,d\mu(x), \quad
\widehat{\mathfrak{E}}^{R}_{p}[\mathfrak{m}](E):= \int\limits_{E} \Bigl(\widehat{I}^{R}[\mathfrak{m}](x)\Bigr)^{p'}\,d\mu(x).
\end{equation}
By (DQ4) of Proposition \ref{Prop.mms_dyadic} it is clear that
\begin{equation}
\label{eqq.9.incl}
\widehat{Q}_{k,\alpha} \subset 9 B_{\epsilon^{k}}(z_{k,\alpha}).
\end{equation}
Hence, by the uniformly locally doubling property of $\mu$, there is a constant $C > 0$ such that, for each $R \in (0,\epsilon^{\underline{k}}]$,
\begin{equation}
\label{eqq.88}
I^{R}[\mathfrak{m}](x) \le C\widehat{I}^{R}[\mathfrak{m}](x) \quad \text{for all} \quad x \in \underline{\operatorname{X}}.
\end{equation}

The following elementary observation will be important in the proof of Theorem \ref{Th.Wolf2} below.

\begin{Prop}
\label{Prop.5.4}
Let $j \geq \underline{k}$ and $\beta \in \mathcal{A}_{j}(\operatorname{X},\epsilon)$. Then
\begin{equation}
\notag
\sum\limits_{k'=k}^{j}\sum_{\substack{Q_{k',\alpha} \supset Q_{j,\beta}}}a_{\mathfrak{m}}(\widehat{Q}_{k',\alpha})
\le \inf\limits_{x \in Q_{j,\beta}} \widehat{I}^{\epsilon^{k}}[\mathfrak{m}](x), \quad k \in \{\underline{k},...,j\}.
\end{equation}
\end{Prop}

\begin{proof}
Indeed, by \eqref{eq.5.4''}
\begin{equation}
\sum\limits_{k'=k}^{j}\sum_{\substack{Q_{k',\alpha} \supset Q_{j,\beta}}}a_{\mathfrak{m}}(\widehat{Q}_{k',\alpha}) = \sum\limits_{k'=k}^{j}a_{\mathfrak{m}}(\widehat{Q}_{k',\alpha}(x)) \quad \text{for all} \quad x \in Q_{j,\beta}.
\end{equation}
Combining this observation with \eqref{eq.Riesz_d} we obtain the required estimate.
\end{proof}

The following three assertions will be crucial for us.
In fact, the corresponding proofs are based on the ideas similar to those used in the proof of Proposition 2.2 from \cite{CasOrtVer}.
However, we should be careful because of the lack of the Euclidean structure. We present the details.

\begin{Lm}
\label{Lm.5.main}
Let $\mathfrak{m}$ be a Borel measure and $p \in (1,\infty)$.
Then there exists a constant $C > 0$ such that, for each $R \in (0,\epsilon^{\underline{k}}]$,
\begin{equation}
\widehat{\mathfrak{E}}^{R}_{p}[\mathfrak{m}](E) \le
p'\sum_{\epsilon^{k} \le R}\sum_{\alpha \in \mathcal{A}_{k}(\operatorname{X},\epsilon)}a_{\mathfrak{m}}(\widehat{Q}_{k,\alpha})\int\limits_{Q_{k,\alpha} \cap E}\Bigl(\widehat{I}^{\epsilon^{k}}[\mathfrak{m}](x)\Bigr)^{p'-1} \,d\mu(x)
\end{equation}
for any Borel set $E \subset \operatorname{X}$.
\end{Lm}

\begin{proof} First of all, we claim that
\begin{equation}
\begin{split}
\label{eq.2.6}
&\Bigl(\widehat{I}^{R}[\mathfrak{m}](x)\Bigr)^{p'} \le p'\sum_{\epsilon^{k} \le R}a_{\mathfrak{m}}(\widehat{Q}_{k,\alpha}(x))
\Bigl(\widehat{I}^{\epsilon^{k}}[\mathfrak{m}](x)\Bigr)^{p'-1}, \quad x \in \underline{\operatorname{X}}.
\end{split}
\end{equation}
Indeed, in the case $\widehat{I}^{R}[\mathfrak{m}](x)=+\infty$, the inequality is obvious.
Assume that $\widehat{I}^{R}[\mathfrak{m}](x) < +\infty$.
Recall that, for any $s \geq 1$, the elementary inequality $\beta^{s}-\alpha^{s} \le s(\beta-\alpha)\beta^{s-1}$
holds for all real numbers $0 \le \alpha \le \beta$.
Hence, given $x \in \underline{\operatorname{X}}$, if $k \in \mathbb{Z}$ is such that $\epsilon^{k} \le R$, then
\begin{align*}
&\Bigl(\sum\limits_{j \geq k}a_{\mathfrak{m}}(\widehat{Q}_{j,\alpha}(x))\Bigr)^{p'}-
\Bigl(\sum\limits_{j \geq k+1}a_{\mathfrak{m}}(\widehat{Q}_{j,\alpha}(x))\Bigr)^{p'} \le p'(a_{\mathfrak{m}}(\widehat{Q}_{k,\alpha}(x)))
\Bigl(\sum\limits_{j \geq k}a_{\mathfrak{m}}(\widehat{Q}_{j,\alpha}(x))\Bigr)^{p'-1}.
\end{align*}
Clearly, if $\widehat{I}^{R}[\mathfrak{m}](x) < +\infty$, then we have
$\widehat{I}^{\epsilon^{k}}[\mathfrak{m}](x) \to 0$, $k \to \infty$. Thus, the standard telescopic-type arguments together with the above inequality lead
to \eqref{eq.2.6}.
As a result, combining  \eqref{eq.Riesz_d}, \eqref{eq.energy}, \eqref{eq.2.6} and taking into account (DQ5) of Proposition \ref{Prop.mms_dyadic} we arrive at the required estimate and complete the proof.
\end{proof}

Now we can estimate the restricted dyadic energy from above in the case $p \geq 2$.

\begin{Lm}
\label{Th.Wolf1}
Let $p \in [2,\infty)$. Then there exists a constant $C > 0$ such that, for each $R \in (0,\epsilon^{\underline{k}}]$ and any Borel set $E \subset \operatorname{X}$,
\begin{equation}
\label{eq.est}
\widehat{\mathfrak{E}}^{R}_{p}[\mathfrak{m}](E) \le
C \sum\limits_{\epsilon^{k} \le R} \sum\limits_{\substack{Q_{k,\alpha} \cap E \neq \emptyset}}\epsilon^{k} \mathfrak{m}(\widehat{Q}_{k,\alpha})
\Bigl(a_{\mathfrak{m}}(\widehat{Q}_{k,\alpha})\Bigr)^{p'-1}.
\end{equation}
\end{Lm}

\begin{proof}
We fix $R \in (0,\epsilon^{\underline{k}}]$ and a Borel set $E \subset \operatorname{X}$.
We have $p' \in (1,2]$ and hence $p'-1 \le 1$.
We change the sum and the integral, use (DQ1), (DQ2), (DQ5) of Proposition \ref{Prop.mms_dyadic}, and take into account \eqref{eq.specialcoefficients}, \eqref{eq.Riesz_d}, \eqref{eqq.9.incl}. Thus, given $k \geq \underline{k}$ and $\alpha \in \mathcal{A}_{k}(\operatorname{X},\epsilon)$, we obtain
\begin{equation}
\notag
\begin{split}
&\frac{1}{\mu(Q_{k,\alpha})}\int\limits_{Q_{k,\alpha} \cap E}\Bigl(\widehat{I}^{\epsilon^{k}}[\mathfrak{m}](x)\Bigr)^{p'-1}\,d\mu(x)
\le \Bigl(\fint\limits_{Q_{k,\alpha}}\sum\limits_{j \geq k} a_{\mathfrak{m}}(\widehat{Q}_{j,\beta}(x))\,d\mu(x)\Bigr)^{p'-1}\\
&\le  9^{p'-1} \Bigl(\frac{1}{\mu(Q_{k,\alpha})}\sum\limits_{j \geq k}\epsilon^{j}\sum_{Q_{j,\beta} \subset Q_{k,\alpha}}\mathfrak{m}(\widehat{Q}_{j,\beta})\Bigr)^{p'-1}.
\end{split}
\end{equation}
Given $j \geq k$, using Propositions \ref{Prop.5.1} and \ref{Prop.5.2}, we have $\sum_{Q_{j,\beta} \subset Q_{k,\alpha}}\mathfrak{m}(\widehat{Q}_{j,\beta}) \le C\mathfrak{m}(\widehat{Q}_{k,\alpha})$. Hence, by the above
inequality, by (DQ4) of Proposition \ref{Prop.mms_dyadic}, and the uniformly locally doubling property of the measure $\mu$, we derive
\begin{equation}
\begin{split}
\label{eq.2.9}
&\frac{1}{\mu(Q_{k,\alpha})}\int\limits_{Q_{k,\alpha} \cap E}\Bigl(\widehat{I}^{\epsilon^{k}}[\mathfrak{m}](x)\Bigr)^{p'-1}\,d\mu(x)
\le C \Bigl(\epsilon^{k}\frac{\mathfrak{m}(\widehat{Q}_{k,\alpha})}{\mu(Q_{k,\alpha})}\Bigr)^{p'-1} \le C \Bigl(a_{\mathfrak{m}}(\widehat{Q}_{k,\alpha})\Bigr)^{p'-1}.
\end{split}
\end{equation}
Finally, using Lemma \ref{Lm.5.main} and \eqref{eq.2.9}, we obtain  \eqref{eq.est} and complete the proof.
\end{proof}

Now we are ready to establish the keystone estimate.

\begin{Th}
\label{Th.Wolf2}
Let $p \in (1,\infty)$. Then there exists a constant $C > 0$ such that, for each $k \geq \underline{k}$,
\begin{equation}
\label{eq.est}
\widehat{\mathfrak{E}}^{\epsilon^{k}}_{p}[\mathfrak{m}](Q_{k,\alpha}) \le
C \sum\limits_{j = k}^{\infty} \sum\limits_{Q_{j,\beta} \subset Q_{k,\alpha}}\epsilon^{j}
\mathfrak{m}(\widehat{Q}_{j,\beta})\Bigl(a_{\mathfrak{m}}(\widehat{Q}_{j,\beta})\Bigr)^{p'-1} \quad \text{for all} \quad \alpha \in \mathcal{A}_{k}(\operatorname{X},\epsilon).
\end{equation}
\end{Th}
\begin{proof}
In the case $p \in [2,\infty)$ the assertion of the theorem follows from Lemma \ref{Th.Wolf1}.

Consider the case $p \in (1,2)$. We fix $k \geq \underline{k}$, $\alpha \in \mathcal{A}_{k}(\operatorname{X},\epsilon)$ and argue by induction.
More precisely, the \textit{base of induction} is that \eqref{eq.est} holds for $p' \in (1,l]$ with $l=2$.
We are going to show that \eqref{eq.est} holds for any $p' > 1$.
Assume that \eqref{eq.est} is proved for $p' \in (1,l]$ for some $l \in \mathbb{N} \cap [2,\infty)$ and show that
\eqref{eq.est} holds for all $p' \in (1,l+1]$.
We use Lemma \ref{Lm.5.main}, and then take into account that $p'-1 \in (1,l]$.
As a result, we obtain
\begin{equation}
\notag
\begin{split}
&\widehat{\mathfrak{E}}^{\epsilon^{k}}_{p}[\mathfrak{m}](Q_{k,\alpha})  \le
p'\sum_{k'=k}^{\infty}\sum_{\substack{Q_{k',\alpha'} \subset Q_{k,\alpha}}}a_{\mathfrak{m}}(\widehat{Q}_{k',\alpha'})\int\limits_{Q_{k',\alpha'}}\Bigl(\widehat{I}^{\epsilon^{k'}}[\mathfrak{m}](x)\Bigr)^{p'-1} \,d\mu(x)\\
&\le C
\sum_{k'=k}^{\infty}\sum_{\substack{Q_{k',\alpha'} \subset Q_{k,\alpha}}}a_{\mathfrak{m}}(\widehat{Q}_{k',\alpha'})
\sum\limits_{j=k'}^{\infty} \sum\limits_{Q_{j,\beta} \subset Q_{k',\alpha'}}\epsilon^{j}\mathfrak{m}(\widehat{Q}_{j,\beta})\Bigl(a_{\mathfrak{m}}(\widehat{Q}_{j,\beta})\Bigr)^{p'-2}.
\end{split}
\end{equation}
Changing the order of summation, we get
\begin{equation}
\label{eq.2.10'}
\begin{split}
&\widehat{\mathfrak{E}}^{\epsilon^{k}}_{p}[\mathfrak{m}](Q_{k,\alpha}) \le C
\sum\limits_{j=k}^{\infty} \sum\limits_{Q_{j,\beta} \subset Q_{k,\alpha}}\epsilon^{j}\mathfrak{m}(\widehat{Q}_{j,\beta})\Bigl(a_{\mathfrak{m}}(\widehat{Q}_{j,\beta})\Bigr)^{p'-2}
\sum\limits_{k'=k}^{j}\sum_{\substack{Q_{k',\alpha'} \supset Q_{j,\beta}}}a_{\mathfrak{m}}(\widehat{Q}_{k',\alpha'}).
\end{split}
\end{equation}
Hence, by Proposition \ref{Prop.5.4}, \eqref{eq.specialcoefficients}, \eqref{eqq.9.incl}, and the uniformly locally doubling property of $\mu$,
\begin{equation}
\label{eq.2.13}
\begin{split}
&\widehat{\mathfrak{E}}^{\epsilon^{k}}_{p}[\mathfrak{m}](Q_{k,\alpha}) \le C\sum_{j=k}^{\infty}\sum\limits_{Q_{j,\beta} \subset Q_{k,\alpha}}\int\limits_{Q_{j,\beta}}
\frac{\epsilon^{j}\mathfrak{m}(\widehat{Q}_{j,\beta})}{\mu(Q_{j,\beta})}(a_{\mathfrak{m}}(\widehat{Q}_{j,\beta}))^{p'-2}\widehat{I}^{\epsilon^{k}}[\mathfrak{m}](x)\,d\mu(x)\\
&\le C \int\limits_{Q_{k,\alpha}}\sum_{j=k}^{\infty}
\Bigl(a_{\mathfrak{m}}(\widehat{Q}_{j,\beta}(x))\Bigr)^{p'-1}\Bigl(\widehat{I}^{\epsilon^{k}}[\mathfrak{m}](x)\Bigr)\,d\mu(x).
\end{split}
\end{equation}
An application of H\"older's inequality for sums with exponents $q=\frac{p'-1}{p'-2}$ and $q'=p'-1$ gives
\begin{equation}
\label{eq.2.14}
\begin{split}
&\sum_{j=k}^{\infty}\Bigl(a_{\mathfrak{m}}(\widehat{Q}_{j,\beta}(x))\Bigr)^{p'-1}=\sum_{j=k}^{\infty}\Bigl(a_{\mathfrak{m}}(\widehat{Q}_{j,\beta}(x))\Bigr)^{\frac{1}{p'-1}}
\Bigl(a_{\mathfrak{m}}(\widehat{Q}_{j,\beta}(x))\Bigr)^{p'-2+\frac{p'-2}{p'-1}}\\
&\le C\Bigl(\widehat{I}^{\epsilon^{k}}[\mathfrak{m}](x)\Bigr)^{\frac{1}{p'-1}}
\Bigl(\sum_{j=k}^{\infty}\Bigl(a_{\mathfrak{m}}(\widehat{Q}_{j,\beta}(x))\Bigr)^{p'}\Bigr)^{\frac{p'-2}{p'-1}}.
\end{split}
\end{equation}
Now we plug \eqref{eq.2.14} into \eqref{eq.2.13} and apply H\"older's inequality for integrals with exponents $p'-1$ and $\frac{p'-1}{p'-2}$.
This gives
\begin{equation}
\notag
\widehat{\mathfrak{E}}^{\epsilon^{k}}_{p}[\mathfrak{m}](Q_{k,\alpha}) \le C \Bigl(\widehat{\mathfrak{E}}^{\epsilon^{k}}_{p}[\mathfrak{m}](Q_{k,\alpha})\Bigr)^{\frac{1}{p'-1}}
\Bigl(\int\limits_{Q_{k,\alpha}}\sum_{j=k}^{\infty}\Bigl(a_{\mathfrak{m}}(\widehat{Q}_{j,\beta}(x))\Bigr)^{p'}\,d\mu(x)\Bigr)^{\frac{p'-2}{p'-1}}.
\end{equation}
As a result, if $\widehat{\mathfrak{E}}^{\epsilon^{k}}_{p}[\mathfrak{m}](Q_{k,\alpha}) < +\infty$, then we have the required inequality
\begin{equation}
\begin{split}
\label{e}
&\widehat{\mathfrak{E}}^{\epsilon^{k}}_{p}[\mathfrak{m}](Q_{k,\alpha}) \le C
\int\limits_{Q_{k,\alpha}}\sum_{j=k}^{\infty}\Bigl(a_{\mathfrak{m}}(\widehat{Q}_{j,\beta}(x))\Bigr)^{p'}\,d\mu(x)\\
&\le C\sum_{j=k}^{\infty}\sum\limits_{Q_{j,\beta} \subset Q_{k,\alpha}}\epsilon^{j}\mathfrak{m}(\widehat{Q}_{j,\beta})(a_{\mathfrak{m}}(\widehat{Q}_{j,\beta}))^{p'-1}.
\end{split}
\end{equation}

To remove the assumption  $\widehat{\mathfrak{E}}^{\epsilon^{k}}_{p}[\mathfrak{m}](Q_{k,\alpha}) < +\infty$ we proceed as follows. Given $l \in \mathbb{N}$, we consider the $l$th truncation
of the restricted Riesz potential $\widehat{I}^{\epsilon^{l},\epsilon^{k}}[\mathfrak{m}]$ obtained by summing in \eqref{eq.Riesz_d} over all $l \le k' \le k$ only (note that
in this case we should use the index $k'$ instead of $k$ in \eqref{eq.Riesz_d}).
Clearly, the corresponding truncations of the restricted dyadic energies $\widehat{\mathfrak{E}}^{\epsilon^{l},\epsilon^{k}}_{p}[\mathfrak{m}](Q_{k,\alpha})$ are finite for all $l \in \mathbb{N}$.
Repeating with minor changes the above arguments we obtain, for any fixed $l \in \mathbb{N}$, the analog of \eqref{e} with $\widehat{\mathfrak{E}}^{\epsilon^{k}}_{p}[\mathfrak{m}](Q_{k,\alpha})$
replaced by $\widehat{\mathfrak{E}}^{\epsilon^{l},\epsilon^{k}}_{p}[\mathfrak{m}](Q_{k,\alpha})$. Passing to the limit as $l$ goes to infinity we complete the proof.
\end{proof}

Given $R > 0$, the \textit{restricted Wolff potential} is a mapping $\mathcal{W}^{R}_{p}[\mathfrak{m}]:\operatorname{X} \to [0,+\infty]$ defined by 
\begin{equation}
\begin{split}
\label{eq.Wolf}
&\mathcal{W}^{R}_{p}[\mathfrak{m}](x):= \sum_{\epsilon^{k} \le R}\Bigl(\epsilon^{kp}\frac{\mathfrak{m}(B_{\epsilon^{k}}(x))}{\mu(B_{\epsilon^{k}}(x))}\Bigr)^{p'-1}, \quad x \in \operatorname{X}.
\end{split}
\end{equation}

Now we are ready to establish an analog of the Hedberg--Wolff-type inequality.

\begin{Ca}
\label{Ca.Wolff}
There are constants $c_{1},c_{2} > 1$ depending on $\epsilon$ only such that the following holds.
Given $p \in (1,\infty)$, there is a constant $C > 0$ such that, for each $R \in (0,\epsilon^{\underline{k}}]$, for every Borel set $E \subset \operatorname{X}$,
\begin{equation}
\label{eq.5.23}
\mathfrak{E}^{R}_{p}[\mathfrak{m}](E) \le C \int\limits_{U_{c_{2}R}(E)}\mathcal{W}^{c_{1}R}_{p}[\mathfrak{m}](y)\,d\mathfrak{m}(y),
\end{equation}
where $U_{c_{2}R}(E):=\{y \in \operatorname{X}: \inf_{x \in E}\operatorname{d}(y,x) < c_{2}R\}$.
\end{Ca}

\begin{proof}
We fix $R \in (0,\epsilon^{\underline{k}}]$, recall notation \eqref{eqq.important_index_notation}, and put $k:=k_{\epsilon}(R)$. By \eqref{eq.energy}, \eqref{eqq.88}
and (DQ5) of Proposition \ref{Prop.mms_dyadic},
\begin{equation}
\label{eq.5.24}
\mathfrak{E}^{R}_{p}[\mathfrak{m}](E) \le C\sum\limits_{Q_{k,\alpha} \cap E \neq \emptyset}\widehat{\mathfrak{E}}^{\epsilon^{k}}_{p}[\mathfrak{m}](Q_{k,\alpha}).
\end{equation}
From (DQ4) of Proposition \ref{Prop.mms_dyadic} and \eqref{eq.modif_dyadic}, given $j \geq k$, it is clear that, for any generalized dyadic cube $Q_{j,\beta}$ in $\operatorname{X}$, we have $B_{\frac{\epsilon^{j}}{2}}(z_{j,\beta}) \subset \widehat{Q}_{j,\beta} \subset B_{18\epsilon^{j}}(y) \subset B_{27\epsilon^{j}}(z_{j,\beta})$ for all $y \in \widehat{Q}_{j,\beta}$. Hence,
using the uniformly locally doubling property of the measure $\mu$, we get
\begin{equation}
\notag
\epsilon^{j}\mathfrak{m}(\widehat{Q}_{j,\beta})\Bigl(a_{\mathfrak{m}}(\widehat{Q}_{j,\beta})\Bigr)^{p'-1} \le
C\int\limits_{\widehat{Q}_{j,\beta}}\Bigl(\epsilon^{jp}\frac{\mathfrak{m}(B_{18\epsilon^{j}}(y))}{\mu(B_{18\epsilon^{j}}(y))}\Bigr)^{p'-1}\,d\mathfrak{m}(y).
\end{equation}
Combining this estimate with Theorem \ref{Th.Wolf2}, we deduce
\begin{equation}
\label{eq.5.25''}
\widehat{\mathfrak{E}}^{\epsilon^{k}}_{p}[\mathfrak{m}](Q_{k,\alpha}) \le
C \int\limits_{\operatorname{X}} \sum\limits_{j=k}^{\infty}\sum\limits_{Q_{j,\beta} \subset Q_{k,\alpha}}\chi_{\widehat{Q}_{j,\beta}}(y)
\Bigl(\epsilon^{jp}\frac{\mathfrak{m}(B_{18\epsilon^{j}}(y))}{\mu(B_{18\epsilon^{j}}(y))}\Bigr)^{p'-1}\,d\mathfrak{m}(y).
\end{equation}
For each $j \geq k$, we put $k(j):=k_{\epsilon}(18\epsilon^{j})$. By the uniformly locally doubling property of $\mu$ it is easy to see that
$\epsilon^{j}\frac{\mathfrak{m}(B_{18\epsilon^{j}}(y))}{\mu(B_{18\epsilon^{j}}(y))} \le C \epsilon^{k(j)}\frac{\mathfrak{m}(B_{\epsilon^{k(j)}}(y))}{\mu(B_{\epsilon^{k(j)}}(y))}$.
Hence, letting $c_{1}=\frac{18\epsilon^{k}}{R}$, and taking into account \eqref{eq.Wolf} and Proposition \ref{Prop.5.1}, we continue \eqref{eq.5.25''}. This gives
\begin{equation}
\label{eq.5.25}
\begin{split}
\widehat{\mathfrak{E}}^{\epsilon^{k}}_{p}[\mathfrak{m}](Q_{k,\alpha})\le C\int\limits_{\widehat{Q}_{k,\alpha}}\mathcal{W}^{c_{1}R}_{p}[\mathfrak{m}](y)\,d\mathfrak{m}(y).
\end{split}
\end{equation}

We put $c_{2}=\frac{11\epsilon^{k}}{R}$. By \eqref{eqq.9.incl} we have $\widehat{Q}_{k,\alpha} \subset U_{c_{2}R}(E)$ provided that $Q_{k,\alpha} \cap E \neq \emptyset$.  Hence, combining \eqref{eq.5.24}, \eqref{eq.5.25}, (DQ5) of Proposition \ref{Prop.mms_dyadic}, and taking into account Propositions \ref{Prop.5.1}, \ref{Prop.covering_multiplicity}, we obtain \eqref{eq.5.23} and complete the proof.
\end{proof}

\section{Proofs of the main results}

In this section we prove Theorems \ref{Th.SecondMain}--\ref{Th.FifthMain}.
We \textit{fix the following data}:

\begin{itemize}
\item[\(\rm (\textbf{D.10.1})\)] a parameter $p \in (1,\infty)$, an m.m.s.\ $\operatorname{X}=(\operatorname{X},\operatorname{d},\mu) \in \mathfrak{A}_{p}$,
and a parameter $q \in (1,p)$ such that $\operatorname{X} \in \mathfrak{A}_{q}$ (we recall Proposition \ref{KeithZhong});

\item[\(\rm (\textbf{D.10.2})\)] a parameter $\theta \in [0,q)$  and a set $S \in \mathcal{LCR}_{\theta}(\operatorname{X})$;

\item[\(\rm (\textbf{D.10.3})\)] a sequence of measure $\{\mathfrak{m}_{k}\} \in \mathfrak{M}_{\theta}(S)$ with parameter $\epsilon = \epsilon(\{\mathfrak{m}_{k}\}) \in (0,\frac{1}{10}]$.

\end{itemize}

We recall again notation \eqref{eqq.important_index_notation} and \eqref{eq.Maximal_function}.
\begin{Lm}
\label{Lm.9.1}
Let $\alpha > \frac{\theta}{p}$. Then, for each $R > 0$, there is a constant $C > 0$ such that
\begin{equation}
\label{eqq.9.1'''}
\|M^{R}_{q,\alpha}(g)|L_{p}(\mathfrak{m}_{0})\| \le C \|g|L_{p}(\operatorname{X})\| \quad \text{for all} \quad g \in L_{p}(\operatorname{X}).
\end{equation}
\end{Lm}

\begin{proof}
Given $g \in L_{p}(\operatorname{X})$, let $\operatorname{X}_{0}(g):=\{x \in \operatorname{X}:\lim_{\widetilde{R} \to 0}M^{\widetilde{R}}_{p,\alpha}(g)=0\}$.
By Proposition \ref{Prop.Evans} and H\"older's inequality,
$\mathcal{H}_{\alpha p}(\operatorname{X} \setminus \operatorname{X}_{0}(g))=0$.
Since $\alpha p > \theta$, by  \eqref{eqq.definition_Hausdorff_measure} we get $\mathcal{H}_{\theta}(\operatorname{X} \setminus \operatorname{X}_{0}(g))=0$.
By Proposition \ref{Prop.absolute_continuity_measure}, we have $\mathfrak{m}_{0}(\operatorname{X} \setminus \operatorname{X}_{0}(g))=0$. Using this observation and H\"older's inequality,
we see that $M^{R}_{q,\alpha}(g)(x) < +\infty$ for all $x \in \operatorname{X}_{0}(g)$.
Given $x \in \operatorname{X}_{0}(g)$, we fix $r_{x} \in (0,R]$ such that
$$
(r_{x})^{\alpha}\Bigl(\fint_{B_{r_{x}}(x)}|g(y)|^{q}\,d\mu(y)\Bigr)^{\frac{1}{q}} > \frac{1}{2}M^{R}_{q,\alpha}(g)(x).
$$
We put $\mathcal{G}:=\{B_{r_{x}}(x):x \in \operatorname{X}_{0}(g)\}$.
Given $k \in \mathbb{N}_{0}$, we put $E_{k}:=\{x \in \operatorname{X}_{0}(g): r_{x} \in (\frac{R}{2^{k+1}},\frac{R}{2^{k}}]\}$ and $\mathcal{G}_{k}:=\{B_{r_{x}}(x):x \in E_{k}\}$.
Obviously, $E_{k} \cap E_{j} = \emptyset$ for $k \neq j$ and $\cup_{k \in \mathbb{N}_{0}}E_{k} = \operatorname{X}_{0}(g)$. Hence,
\begin{equation}
\label{eqq.9.2'''}
\int\limits_{\operatorname{X}}(M^{R}_{q,\alpha}(g)(x))^{p}\,d\mathfrak{m}_{0}(x)=\sum\limits_{k=0}^{\infty}\int\limits_{E_{k}}(M^{R}_{q,\alpha}(g)(x))^{p}\,d\mathfrak{m}_{0}(x).
\end{equation}
Given $k \in \mathbb{N}_{0}$, using the $5B$-covering lemma we find a disjoint family of balls $\widetilde{\mathcal{G}}_{k} \subset \mathcal{G}_{k}$ such that $E_{k} \subset \cup\{5B: B \in \widetilde{\mathcal{G}}_{k}\}$.
Using the uniformly locally doubling property of $\mu$, we have
\begin{equation}
\notag
\begin{split}
&\int\limits_{E_{k}}(M^{R}_{q,\alpha}(g)(x))^{p}\,d\mathfrak{m}_{0}(x) \le 2^{p}\sum\limits_{B \in \widetilde{\mathcal{G}}_{k}}\int\limits_{5B \cap E_{k}}(r_{x})^{\alpha p}\Bigl(\fint\limits_{B_{r_{x}}(x)}|g(y)|^{q}\,d\mu(y)\Bigr)^{\frac{p}{q}}\,d\mathfrak{m}_{0}(x)\\
&\le C\sum\limits_{B \in \widetilde{\mathcal{G}}_{k}}\mathfrak{m}_{0}(5B)\Bigl(\frac{R}{2^{k}}\Bigr)^{\alpha p}\Bigl(\fint\limits_{6B}|g(y)|^{q}\,d\mu(y)\Bigr)^{\frac{p}{q}}.
\end{split}
\end{equation}
By Theorem \ref{Th.doubling_type} and \eqref{M.2} we have $\mathfrak{m}_{0}(5B) \le C2^{k\theta}\mu(B)$ for all $B \in \widetilde{\mathcal{G}}_{k}$. Hence, using
H\"older's inequality and Propositions \ref{Prop.covering_multiplicity}, \ref{Prop.finite_intersection},
we deduce
\begin{equation}
\label{eqq.9.3'''}
\begin{split}
&\int\limits_{E_{k}}(M^{R}_{q,\alpha}(g)(x))^{p}\,d\mathfrak{m}_{0}(x)
\le C 2^{k(\theta-\alpha p)}\sum\limits_{B \in \widetilde{\mathcal{G}}_{k}}\mu(5B)\fint\limits_{6B}|g(y)|^{p}\,d\mu(y)\\
&\le C 2^{k(\theta-\alpha p)}\sum\limits_{B \in \widetilde{\mathcal{G}}_{k}}\int\limits_{6B}|g(y)|^{p}\,d\mu(y) \le C 2^{k(\theta-\alpha p)}\int\limits_{\operatorname{X}}|g(y)|^{p}\,d\mu(y).
\end{split}
\end{equation}
Since $\alpha p > \theta$, a combination of \eqref{eqq.9.2'''} and \eqref{eqq.9.3'''}  gives \eqref{eqq.9.1'''}.

\end{proof}

The following result is crucial for our analysis. We recall \eqref{eqq.notation_average}, and
\textit{throughout this section}, we put $F_{G}:=F_{G,\mu}$. Furthermore, we set $B_{k}(x):=B_{\epsilon^{k}}(x)$, as usual.

\begin{Th}
\label{Th.main_est}
For each $c \geq 1$ and $\widetilde{p} \in (q,p]$, there exist $C > 0$, $\widetilde{c} \geq c$ such that, if $B_{k}(x') \subset cB_{k}(x)$  for some $k \in \mathbb{N}_{0}$, $x' \in S$ and $x \in \operatorname{X}$, then
\begin{equation}
\label{eq.7.1}
\fint\limits_{cB_{k}(x)}|F|_{S}^{\mathfrak{m}_{0}}(y)-F_{cB_{k}(x)}|\,d\mathfrak{m}_{k}(y)
\le C \epsilon^{k}\Bigl(\fint\limits_{\widetilde{c}B_{k}(x)}(|DF|_{p}(y))^{\widetilde{p}}\,d\mu(y)\Bigr)^{\frac{1}{\widetilde{p}}} \quad \text{for all} \quad F \in W_{p}^{1}(\operatorname{X}).
\end{equation}
\end{Th}

\begin{proof} We fix $c \geq 1$. We also fix $k \in \mathbb{N}_{0}$, $x' \in S$, $x \in \operatorname{X}$ such that $B_{k}(x') \subset cB_{k}(x)$. To the end of the proof we put $B:=B_{k}(x)$.

\textit{Step 1.} By Definition \ref{Def.m_trace_space}, for $\mathfrak{m}_{0}$-a.e.\ $y \in cB \cap S$, we have
\begin{equation}
\label{eqq.92}
\begin{split}
&|F|_{S}^{\mathfrak{m}_{0}}(y)-F_{cB}|=\lim\limits_{i \to \infty}|F_{B_{i}(y)}-F_{cB}| \le |F_{B_{k}(y)}-F_{cB}|+\sum\limits_{i=k}^{\infty}|F_{B_{i}(y)}-F_{B_{i+1}(y)}|.\\
\end{split}
\end{equation}
Combing Remark \ref{Rem.best_approx} with Proposition \ref{Prop.Sobolev_Poincare} we obtain, for $\mathfrak{m}_{0}$-a.e. $y \in cB \cap S$,
\begin{equation}
\label{eqq.93'}
\sum\limits_{i=k}^{\infty}|F_{B_{i}(y)}-F_{B_{i+1}(y)}| \le C \sum\limits_{i=k}^{\infty}\epsilon^{i}\Bigl(\fint\limits_{\lambda B_{i}(y)}(|DF|_{p}(v))^{q}\,d\mu(v)\Bigr)^{\frac{1}{q}}.
\end{equation}
We recall \eqref{eq.Maximal_function}. Using the uniformly locally doubling property of $\mu$ it is easy to see that, for each $i \geq k$,
\begin{equation}
\label{eq.7.2}
\begin{split}
&\Bigl(\fint\limits_{\lambda B_{i}(y)}(|DF|_{p}(v))^{q}\,d\mu(v)\Bigr)^{\frac{1}{q}} \le C \inf\limits_{z \in B_{i}(y)}\Bigl(\fint\limits_{(\lambda+1)B_{i}(z)}(|DF|_{p}(v))^{q}\,d\mu(v)\Bigr)^{\frac{1}{q}}\\
&\le C \fint\limits_{B_{i}(y)}M^{(\lambda+1)\epsilon^{i}}_{q,0}(|DF|_{p})(z)\,d\mu(z) \quad \text{for all} \quad y \in cB.
\end{split}
\end{equation}
As a result, using estimates \eqref{eqq.93'}, \eqref{eq.7.2} and taking into account \eqref{eqq.weighted_measure}, \eqref{eq.Riesz}, we obtain
\begin{equation}
\label{eqq.94}
\fint\limits_{cB}\sum\limits_{i=k}^{\infty}|F_{B_{i}(y)}-F_{B_{i+1}(y)}|\,d\mathfrak{m}_{k}(y) \le C \Bigl(\fint\limits_{cB}I^{\epsilon^{k}}[M^{(\lambda+1)\epsilon^{k}}_{q,0}(|DF|_{p})\mu](y)\,d\mathfrak{m}_{k}(y)\Bigr).
\end{equation}
On the other hand, $B_{k}(y) \subset (c+1)B$ for all $y \in cB$. Hence, using Remark \ref{Rem.best_approx},
Proposition \ref{Prop.Sobolev_Poincare} and H\"older's inequality, we obtain
\begin{equation}
\label{eqq.95'}
|F_{B_{k}(y)}-F_{cB}|  \le C \epsilon^{k}\Bigl(\fint\limits_{\lambda(c+1)B}(|DF|_{p}(z))^{\widetilde{p}}\,d\mu(z)\Bigr)^{\frac{1}{\widetilde{p}}} \quad \text{for all} \quad y \in cB.
\end{equation}

\textit{Step 2.} By the standard duality arguments it is clear that, given a constant $\underline{C} > 0$, a parameter $R > 0$, and measures $\nu$, $\sigma$ on $\operatorname{X}$,
\begin{equation}
\notag
\|I^{R}[g\nu]|L_{1}(\sigma)\|  \le \underline{C} \|g|L_{\widetilde{p}}(\nu)\| \quad \text{for all} \quad g \in L_{\widetilde{p}}(\nu)
\end{equation}
if and only if (we set $\widetilde{p}':=\frac{\widetilde{p}}{\widetilde{p}-1}$, as usual)
\begin{equation}
\notag
\|I^{R}[h\sigma]|L_{\widetilde{p}'}(\nu)\|  \le \underline{C} \|h|L_{\infty}(\sigma)\| \quad \text{for all} \quad h \in L_{\infty}(\sigma).
\end{equation}
Note that in \eqref{eqq.94} we work only with the part of the measure $\mu$ concentrated in $(c+2+\lambda)B$.
Now we set $\widetilde{c}:=\max\{c+2+\lambda,\lambda(c+1)\}$ and apply the duality arguments given above with the measures $\sigma=\mathfrak{m}_{k}\lfloor_{cB}$, $\nu=\mu\lfloor_{\widetilde{c}B}$ and with $g=M^{2\epsilon^{k}}_{q}(|DF|_{p}\chi_{\widetilde{c}B})$. This gives
\begin{equation}
\label{eq.7.4}
\int\limits_{cB}I^{\epsilon^{k}}[M^{(\lambda+1)\epsilon^{k}}_{q,0}(|DF|_{p}\chi_{\widetilde{c}B})\mu](y)\,d\mathfrak{m}_{k}(y) \le \underline{C} \Bigl(\int\limits_{\widetilde{c}B}(M^{(\lambda+1)\epsilon^{k}}_{q,0}(|DF|_{p}\chi_{\widetilde{c}B}))^{\widetilde{p}}\,d\mu(y)\Bigr)^{\frac{1}{\widetilde{p}}}
\end{equation}
with the constant $\underline{C}:=\Bigl(\mathfrak{E}^{\epsilon^{k}}_{\widetilde{p}}[\mathfrak{m}_{k}](\widetilde{c}B)\Bigr)^{\frac{1}{\widetilde{p}'}}.$

\textit{Step 3.}  By Corollary \ref{Ca.Wolff} we have $\Bigl(\underline{C}\Bigr)^{\widetilde{p}'} \le C \int_{(c_{2}+\widetilde{c})B}\mathcal{W}^{c_{1}\epsilon^{k}}_{\widetilde{p}}[\mathfrak{m}_{k}](y)d\mathfrak{m}_{k}(y)$.
Hence, using \eqref{eq.Wolf}, \eqref{M.2} and Theorem \ref{Th.doubling_type} we obtain
\begin{equation}
\label{eq.7.7}
\begin{split}
&\Bigl(\underline{C}\Bigr)^{\widetilde{p}'} \le C \int\limits_{(c_{2}+\widetilde{c})B}\sum\limits_{\epsilon^{i} \le c_{1}\epsilon^{k}}(\epsilon^{i(\widetilde{p}-\theta)})^{\widetilde{p}'-1}\,d\mathfrak{m}_{k}(y) \le C \mu(B)\epsilon^{k(\widetilde{p}-\theta)\frac{\widetilde{p}'}{\widetilde{p}}-\theta}.
\end{split}
\end{equation}

\textit{Step 4.} Since $B_{k}(x') \subset cB$, we have $B \subset (c+1)B_{k}(x')$. Hence, using
\eqref{M.3} and the uniformly locally doubling property of $\mu$, we get
\begin{equation}
\label{eq.7.8}
\frac{1}{\mathfrak{m}_{k}(cB)} \le \frac{1}{\mathfrak{m}_{k}(B_{k}(x'))}
\le C \frac{\epsilon^{k\theta}}{\mu(B_{k}(x'))} \le C \frac{\epsilon^{k\theta}}{\mu((c+1)B_{k}(x'))} \le C \frac{\epsilon^{k\theta}}{\mu(B)}.
\end{equation}

\textit{Step 5.} Combining \eqref{eqq.92} with \eqref{eqq.94}--\eqref{eq.7.8} and taking into account Proposition \ref{Prop.maximal_function} we get \eqref{eq.7.1}, completing the proof.
\end{proof}

Now we recall \eqref{eq.main1} and establish the following powerful estimate.

\begin{Ca}
\label{Ca.7.2}
There exists a constant $C > 0$ such that
\begin{equation}
\label{eq.79}
\mathcal{CN}_{p,\{\mathfrak{m}_{k}\}}(F|_{S}^{\mathfrak{m}_{0}}) \le C \||DF|_{p}|L_{p}(\operatorname{X})\| \quad \text{for all} \quad F \in W_{p}^{1}(\operatorname{X})
\end{equation}
and, furthermore,
\begin{equation}
\label{eq.79'}
\|F|_{S}^{\mathfrak{m}_{0}}|L_{p}(\mathfrak{m}_{0})\| \le C \|F|W^{1}_{p}(\operatorname{X})\| \quad \text{for all} \quad F \in W_{p}^{1}(\operatorname{X}).
\end{equation}
\end{Ca}

\begin{proof}
We fix $\widetilde{p} \in (q,p)$ and put $f:=F|_{S}^{\mathfrak{m}_{0}}$ for brevity. Using \eqref{eqq.best_approximation_constant}
and applying Theorem \ref{Th.main_est} with $c=2$ we get,
for each ball $B_{k}(x)$ with $k \in \mathbb{N}_{0}$ and $B_{k}(x) \cap S \neq \emptyset$,
\begin{equation}
\label{eqq.1015}
\begin{split}
&\mathcal{E}_{\mathfrak{m}_{k}}(f,2B_{k}(x))\le \fint\limits_{2B_{k}(x)}|f(y)-F_{2B_{k}(x)}|\,d\mathfrak{m}_{k}(y) \le C \epsilon^{k}\Bigl(\fint\limits_{\widetilde{c}B_{k}(x)}(|DF|_{p})^{\widetilde{p}}\,d\mu(y)\Bigr)^{\frac{1}{\widetilde{p}}}.
\end{split}
\end{equation}
By Theorem \ref{Th.doubling_type} it is easy to see that
$f^{\sharp}_{\{\mathfrak{m}_{k}\}}(x) \le C\sup_{k \in \mathbb{N}_{0}}\epsilon^{-k}\mathcal{E}_{\mathfrak{m}_{k}}(f,2B_{k}(x))$
for all $x \in \operatorname{X}$. Hence, using \eqref{eqq.1015} we obtain
\begin{equation}
\notag
f^{\sharp}_{\{\mathfrak{m}_{k}\}}(x) \le C M^{\widetilde{c}}_{\widetilde{p},0}(|DF|_{p})(x) \quad \text{for all} \quad x \in \operatorname{X}.
\end{equation}
As a result, an application of Proposition \ref{Prop.maximal_function} gives \eqref{eq.79}.

We recall notation \eqref{eq.lattice} and fix a family $\mathcal{B}:=\mathcal{B}_{0}(\operatorname{X},\epsilon)$. By H\"older's inequality, 
\begin{equation}
\label{eqq.915'''}
\begin{split}
&\int\limits_{S}|f(x)|^{p}\,d\mathfrak{m}_{0}(x) \le C \sum\limits_{\substack{B \in \mathcal{B}\\ B \cap S \neq \emptyset}}\int\limits_{B}|f(x)|^{p}\,d\mathfrak{m}_{0}(x)\\
&\le C \sum\limits_{\substack{B \in \mathcal{B}\\B \cap S \neq \emptyset}}\Bigl(\int\limits_{B}|f(x)-F_{2B}|^{p}\,d\mathfrak{m}_{0}(x)+\frac{\mathfrak{m}_{0}(B)}{\mu(2B)}\int\limits_{2B}|F(x)|^{p}\,d\mu(x)\Bigr).
\end{split}
\end{equation}
Given $B \in \mathcal{B}$ with $B \cap S \neq \emptyset$, by the triangle inequality
\begin{equation}
\label{eqq.916'''}
\begin{split}
&\int\limits_{B}|f(x)-F_{2B}|^{p}\,d\mathfrak{m}_{0}(x) \le \int\limits_{B}|F_{B_{0}(x)}-F_{2B}|^{p}\,d\mathfrak{m}_{0}(x)+\int\limits_{B}|f(x)-F_{B_{0}(x)}|^{p}\,d\mathfrak{m}_{0}(x).\\
\end{split}
\end{equation}
By Remark \ref{Rem.best_approx}, $|F_{B_{0}(x)}-F_{2B}| \le C\mathcal{E}_{\mu}(F,2B)$ for all $x \in B \cap S$.
Hence, by Proposition \ref{Prop.Sobolev_Poincare},
\begin{equation}
\label{eqq.916'''}
\begin{split}
\int\limits_{B}|F_{B_{0}(x)}-F_{2B}|^{p}\,d\mathfrak{m}_{0}(x) \le C\mathfrak{m}_{0}(B) \fint\limits_{2\lambda B}\Bigl((|DF|_{p})(y)\Bigr)^{p}\,d\mu(y).
\end{split}
\end{equation}
By Definition \ref{Def.m_trace_space}, given $\delta \in (0,1/2)$, we have $|f(x)-F_{B_{0}(x)}| \le \sum_{i=0}^{\infty}\frac{\epsilon^{i\delta}}{\epsilon^{i\delta}}|F_{B_{i}(x)}-F_{B_{i+1}(x)}|$
for $\mathfrak{m}_{0}$-a.e. $x \in S$. Using H\"older's inequality for sums, Remark \ref{Rem.best_approx}, Proposition \ref{Prop.Sobolev_Poincare} and taking into account \eqref{eq.Maximal_function}, we get, for $\mathfrak{m}_{0}$-a.e. $x \in S$,
\begin{equation}
\label{eqq.917'''}
\begin{split}
&|f(x)-F_{B_{0}(x)}|^{p} \le C\sum_{i=0}^{\infty}\epsilon^{-ip\delta}|F_{B_{i}(x)}-F_{B_{i+1}(x)}|^{p}\\
&\le
C \sum_{i=0}^{\infty}\epsilon^{i(p-p\delta)}\Bigl(\fint\limits_{\lambda B_{i}(x)}(|DF|_{p}(y))^{q}\,d\mu(y)\Bigr)^{\frac{p}{q}} \le C \Bigl(M^{\lambda}_{q,1-2\delta}(|DF|_{p})(x)\Bigr)^{p}.
\end{split}
\end{equation}
By \eqref{M.2}, $\mathfrak{m}_{0}(B) \le C\mu(B)$ for all $B \in \mathcal{B}$.
We combine estimates \eqref{eqq.915'''}--\eqref{eqq.917'''}, and take into account
Propositions \ref{Prop.covering_multiplicity}, \ref{Prop.finite_intersection}. Finally, choosing $\delta > 0$ so small that $p-2p\delta > \theta$ we apply Lemma \ref{Lm.9.1} with $\alpha=1-2\delta$. This gives \eqref{eq.79'}
and completes the proof.
\end{proof}

We should warn the reader that the following result is not a consequence of Theorem \ref{Th.main_est}. Indeed,
at this moment it has not yet been proved that $\operatorname{Ext}_{S,\{\mathfrak{m}_{k}\}}$ is the right inverse of $\operatorname{Tr}|^{\mathfrak{m}_{0}}_{S}$.

\begin{Th}
\label{Th.trace_of_extension}
Assume that $\{\mathfrak{m}_{k}\} \in \mathfrak{M}^{str}_{\theta}(S)$.
Then, for each $c \geq \frac{3}{\epsilon}$, there exists a constant $C > 0$ such that, for every
$f \in L_{1}^{loc}(\{\mathfrak{m}_{k}\})$ satisfying $\operatorname{N}_{p,\{\mathfrak{m}_{k}\},c}(f) < +\infty$, for each $k \in \mathbb{N}_{0}$, and any
ball $B=B_{\epsilon^{k}}(x)$ with $x \in S$,
\begin{equation}
\label{eq.7.9}
\fint\limits_{B}|f(y)-F_{B}|\,d\mathfrak{m}_{k}(y)
\le C \epsilon^{k}\Bigl(\fint\limits_{B}(|DF|_{p}(y))^{p}\,d\mu(y)\Bigr)^{\frac{1}{p}}
\end{equation}
where $F:=\operatorname{Ext}_{S,\{\mathfrak{m}_{k}\}}(f)$.
\end{Th}

\begin{proof}
We fix $f \in L_{1}^{loc}(\{\mathfrak{m}_{k}\})$ with $\operatorname{N}_{p,\{\mathfrak{m}_{k}\},c}(f) < +\infty$.
By Theorem \ref{Th.firstext} we have $F=\operatorname{Ext}_{S,\{\mathfrak{m}_{k}\}}(f) \in W_{p}^{1}(\operatorname{X})$. Hence, the right-hand side of inequality \eqref{eq.7.9} makes sense.
We recall the concept of a special approximating sequence introduced in \eqref{eq.specaproxsequence}.
By Theorem \ref{Th.derivappseq} there is a constant $C > 0$ independent on $f$ and a subsequence $\{f^{j_{s}}\}$ of the sequence $\{f^{j}\}$ such that for all
large enough $s \in \mathbb{N}$ we have $\|\operatorname{lip}f^{j_{s}}|L_{p}(\operatorname{X})\| \le C \operatorname{N}_{p,\{\mathfrak{m}_{k}\},c}(f)$.
By Definition \ref{Def.Cheeger_Sobolev} this implies that $\operatorname{Ch}_{p}(f^{j_{s}})<+\infty$ and, moreover, by Remark \ref{Rem.minimal_upper_gradient}
$\||Df^{j_{s}}|_{p}|L_{p}(\operatorname{X})\| \le C \operatorname{N}_{p,\{\mathfrak{m}_{k}\},c}(f)$ for all large enough $s \in \mathbb{N}$
with constant $C > 0$ independent on $f$ and $s$. At the same time, by Theorem \ref{Th.converge1} and Corollary \ref{Ca.weaklp_estimate}
we have $\sup_{j \in \mathbb{N}}\|f^{j}|L_{p}(\operatorname{X})\| \le C_{f}\operatorname{N}_{p,\{\mathfrak{m}_{k}\},c}(f) < +\infty$.
Hence, the sequence $\{f^{j_{s}}\}$ is bounded in $W_{p}^{1}(\operatorname{X})$. In view of Proposition \ref{Prop.reflexivity}, there exists
a weakly convergent subsequence of the sequence $\{f^{j_{s}}\}$. By Mazur's lemma there is an increasing sequence $\{N_{l}\}_{l \in \mathbb{N}}$ and a sequence of convex combinations
$\widetilde{f}^{N_{l}}:=\sum_{i=0}^{M_{l}}\lambda_{N_{l}}^{i}f^{N_{l}+i}$ with $\lambda^{i}_{N_{l}} \geq 0$ and $\sum_{i=0}^{M_{l}}\lambda^{i}_{N_{l}}=1$
such that $\|F-\widetilde{f}^{N_{l}}|W_{p}^{1}(\operatorname{X})\| \to 0$ as $l \to \infty.$
Since, given $l \in \mathbb{N}$, the function $\widetilde{f}^{N_{l}}$ is Lipschitz, the $\mathfrak{m}_{0}$-trace of $\widetilde{f}^{N_{l}}$ to $S$ is an $\mathfrak{m}_{0}$-equivalence class of the pointwise
restriction of $\widetilde{f}^{N_{l}}$ to $S$. Hence, by Theorem \ref{Th.main_est} we get
\begin{equation}
\label{eq.7.11}
\fint\limits_{B}|\widetilde{f}^{N_{l}}|_{S}^{\mathfrak{m}_{0}}(y)-\widetilde{f}^{N_{l}}_{B}|\,d\mathfrak{m}_{k}(y) \le \epsilon^{k}\Bigl(\fint\limits_{B}(|D\widetilde{f}^{N_{l}}|_{p}(y))^{p}\,d\mu(y)\Bigr)^{\frac{1}{p}}.
\end{equation}
By Theorem \ref{Th.converge1} the sequence $\{f^{j}\}$ converges both in $L_{p}(\operatorname{X})$ and $L_{p}(\mathfrak{m}_{k})$ to $F$ and $f$, respectively. Clearly, the same holds true for
the sequence $\{\widetilde{f}^{N_{l}}\}$. As a result, passing to the limit in \eqref{eq.7.11} we get the required estimate.

\end{proof}

Now we are ready to show that our extension operator $\operatorname{Ext}_{S,\{\mathfrak{m}_{k}\}}$ is the right inverse of the $\mathfrak{m}_{0}$-trace
operator $\operatorname{Tr}|_{S}^{\mathfrak{m}_{0}}$. We recall Definition \ref{Def.multiweight_Lebesgue}.
\begin{Ca}
\label{Ca.ext_inverse_trace}
Assume that
$\{\mathfrak{m}_{k}\} \in \mathfrak{M}^{str}_{\theta}(S)$.
If $\operatorname{N}_{p,\{\mathfrak{m}_{k}\},c}(f) < +\infty$ for some $c \geq \frac{3}{\epsilon}$, then
\begin{equation}
\label{eq.1022}
\lim\limits_{k \to \infty}\fint\limits_{B_{k}(x)}|f(x)-\operatorname{Ext}_{S,\{\mathfrak{m}_{k}\}}(f)|\,d\mu(y) = 0 \quad \text{for} \quad \mathcal{H}_{p}-a.e. \quad x \in \mathcal{R}_{\{\mathfrak{m}_{k}\},\epsilon}(f).
\end{equation}
In particular, $\operatorname{Tr}|_{S}^{\mathfrak{m}_{0}} \circ \operatorname{Ext}_{S,\{\mathfrak{m}_{k}\}}(f)=f$.
\end{Ca}

\begin{proof}
We put $F:=\operatorname{Ext}_{S,\{\mathfrak{m}_{k}\}}(f)$. It is clear that, given $k \in \mathbb{N}_{0}$ and $x \in S$, we have
\begin{equation}
\begin{split}
\label{eq.714}
&\fint\limits_{B_{k}(x)}|f(x)-F(y)|\,d\mu(y) \le \Bigl|f(x)-\fint\limits_{B_{k}(x)}f(y)\,d\mathfrak{m}_{k}(y)\Bigr|\\
&+\fint\limits_{B_{k}(x)}\fint\limits_{B_{k}(x)}|f(y)-F(y)|\,d\mathfrak{m}_{k}(y)\,d\mu(y)=:\sum\limits_{i=1}^{2}R_{k}^{i}(x).
\end{split}
\end{equation}
By Definition \ref{Def.multiweight_Lebesgue} we have $\lim_{k \to \infty}R_{k}^{1}(x)=0$ for all $x \in \mathcal{R}_{\{\mathfrak{m}_{k}\},\epsilon}(f)$.
On the other hand, combining Remark \ref{Rem.best_approx}, Propositions \ref{Prop.Evans}, and Theorem \ref{Th.trace_of_extension}, we obtain
\begin{equation}
\label{eq.716}
\Bigl(R_{k}^{2}(x)\Bigr)^{p} \le C \epsilon^{kp}\fint\limits_{B_{k}(x)}(|DF|_{p}(y))^{p}\,d\mu(y) \to 0, \quad k \to \infty \quad \text{for} \quad \mathcal{H}_{p}-\text{a.e.} \quad x \in S.
\end{equation}
Combining \eqref{eq.714} and \eqref{eq.716} we get \eqref{eq.1022}. Finally, to prove the second claim it is sufficient to use Theorem \ref{Th.trueLebesgue}, Proposition \ref{Prop.absolute_continuity_measure}
and take into account that the inequality $p > \theta$ implies that $\mathcal{H}_{\theta}\lfloor_{S}$ is absolutely continuous with respect $\mathcal{H}_{p}\lfloor_{S}$.
\end{proof}

We recall the concept of a $p$-sharp representative introduced in Section 2.3.

\begin{Prop}
\label{Prop.from_my_previous_paper}
Given $F \in W_{p}^{1}(\operatorname{X})$, let $\overline{F}$ be an arbitrary $p$-sharp representative of $F$. Then $C_{p}(S\setminus \mathcal{R}_{\{\mathfrak{m}_{k}\},\epsilon}(\overline{F}))=0$.
\end{Prop}

\begin{proof}
One should repeat almost verbatim the proof of Lemma 4.3 in \cite{TV1} using (in appropriate places) Propositions \ref{Prop.Capacity}, \ref{Prop.Evans} and Theorem \ref{Th.main_est} of this paper instead of Propositions 2.4, 3.1 and Theorem 3.1 from \cite{TV1}.
\end{proof}


\textit{Proofs of Theorems \ref{Th.SecondMain}--\ref{Th.FifthMain}.}
We fix $c \geq \frac{3}{\epsilon}$, $\sigma \in (0,\frac{\epsilon^{2}}{4c})$ and split the proof into several steps.

\textit{Step 1.} We recall Definition \ref{Def.kostyl} and fix $f \in L_{1}^{loc}(\{\mathfrak{m}_{k}\})$ such that $\operatorname{N}_{p,\{\mathfrak{m}_{k}\},c}(f) < +\infty$.
By Theorem \ref{Th.firstext} we have $F:=\operatorname{Ext}_{S,\{\mathfrak{m}_{k}\}}(f) \in W^{1}_{p}(\operatorname{X})$.
By Corollary \ref{Ca.ext_inverse_trace} we conclude that $f \in W_{p}^{1}(\operatorname{X})|_{S}^{\mathfrak{m}_{0}}$ and $f=\operatorname{Tr}|_{S}^{\mathfrak{m}_{0}}(F)$.
Hence, by Corollary \ref{Ca.7.2}, we get $\operatorname{CN}_{p,\{\mathfrak{m}_{k}\}}(f) < +\infty$. By Theorem \ref{Theore.53}, this implies that $\operatorname{BSN}_{p,\{\mathfrak{m}_{k}\},c}(f) < +\infty$
for any $c \geq 1$. Combining these observations with Theorem \ref{Th.comparison1} we prove equivalence of $(i)$, $(ii)$ and $(iii)$ in Theorem \ref{Th.SecondMain}.
An application of Theorems \ref{Th.Functional1}, \ref{Th.Functional2} verifies the equivalence of $(i)-(iv)$ in Theorem \ref{Th.SecondMain}.

\textit{Step 2.} By Theorems \ref{Th.lpnorm_estimate}, \ref{Th.firstext}, \ref{Th.comparison1} and \ref{Theore.53}, we have (we put $F=\operatorname{Ext}_{S,\{\mathfrak{m}_{k}\}}(f)$)
\begin{equation}
\label{eq.718}
\begin{split}
&\|f|W_{p}^{1}(\operatorname{X})|_{S}^{\mathfrak{m}_{0}}\| \le \|F|W_{p}^{1}(\operatorname{X})\| \le C \operatorname{BSN}_{p,\{\mathfrak{m}_{k}\},c}(f)
\le C \operatorname{CN}_{p,\{\mathfrak{m}_{k}\}}(f).
\end{split}
\end{equation}
From Remark \ref{Rem.minimal_upper_gradient}, Theorem \ref{Th.firstext}, and Corollary \ref{Ca.7.2}, we have
\begin{equation}
\label{eqq.1026}
C^{-1}\mathcal{CN}_{p,\{\mathfrak{m}_{k}\}}(f) \le \||DF|_{p}|L_{p}(\operatorname{X})\| \le C \operatorname{N}_{p,\{\mathfrak{m}_{k}\},c}(f).
\end{equation}
By \eqref{eq.718}, \eqref{eqq.1026} and Theorem \ref{Th.comparison1} we get $\operatorname{N}_{p,\{\mathfrak{m}_{k}\},c}(f) \approx \operatorname{BSN}_{p,\{\mathfrak{m}_{k}\},c}(f) \approx \operatorname{CN}_{p,\{\mathfrak{m}_{k}\}}(f)$.
Finally, combining this fact with \eqref{eq.718} and Theorems  \ref{Th.Functional1}, \ref{Th.Functional2} we get \eqref{eq.717}.
As a result, we complete the proof of Theorem \ref{Th.SecondMain}, and, furthermore, we prove assertion (1) in Theorem \ref{Th.FifthMain}.

\textit{Step 3.} Now we prove the first claim in Theorem \ref{Th.FourthMain}. If $f \in W_{p}^{1}(\operatorname{X})|_{S}$, then by Definitions \ref{Def.sharptrace} and \ref{Def.m_trace_space} it is clear that
$\operatorname{I}_{\mathfrak{m}_{0}}(f) \in W_{p}^{1}(\operatorname{X})|^{\mathfrak{m}_{0}}_{S}$. Hence, condition $(A)$ in Theorem \ref{Th.FourthMain} holds. Furthermore, by Proposition \ref{Prop.from_my_previous_paper}
condition $(B)$ in Theorem \ref{Th.FourthMain} holds. Conversely, assume that a function $f \in \mathfrak{B}(S)$ is such that
conditions $(A)$, $(B)$ in Theorem \ref{Th.FourthMain} hold true. By Theorem \ref{Th.SecondMain}, Definition \ref{Def.multiweight_Lebesgue} and Corollary \ref{Ca.ext_inverse_trace}
we get $f \in W_{p}^{1}(\operatorname{X})|_{S}$.

\textit{Step 4.}
By Definitions \ref{Def.sharptrace} and \ref{Def.m_trace_space} it is clear that $\operatorname{I}_{\mathfrak{m}_{0}}: W_{p}^{1}(\operatorname{X})|_{S} \to W_{p}^{1}(\operatorname{X})|^{\mathfrak{m}_{0}}_{S}$ is a continuous embedding. By
Theorem \ref{Th.SecondMain}, Definition \ref{Def.multiweight_Lebesgue}, Corollary \ref{Ca.ext_inverse_trace} and Proposition \ref{Prop.from_my_previous_paper} it follows that for each $[f] \in W_{p}^{1}(\operatorname{X})|^{\mathfrak{m}_{0}}_{S}$ there
is a representative $f \in W_{p}^{1}(\operatorname{X})|_{S}$. Hence, $\operatorname{I}_{\mathfrak{m}_{0}}$ is a surjection.
Now we show that the mapping $\operatorname{I}_{\mathfrak{m}_{0}}$ is injective on $W_{p}^{1}(\operatorname{X})|_{S}$. Assume that, given $f \in W_{p}^{1}(\operatorname{X})|_{S}$, we have
$\operatorname{I}_{\mathfrak{m}_{0}}(f)(x)=0$ for $\mathfrak{m}_{0}$-a.e. $x \in \operatorname{X}$. Hence, by $(B)$ in Theorem \ref{Th.FourthMain}
this implies that $f(x)=0$ everywhere on $S$ except a set of $p$-capacity zero, i.e., $f=0$ in the sense of $W_{p}^{1}(\operatorname{X})|_{S}$.
As a result, from Definitions \ref{Def.sharptrace} and \ref{Def.m_trace_space}
it follows easily that $\operatorname{I}_{\mathfrak{m}_{0}}$ is an isometric isomorphism. This verifies (3) in Theorem \ref{Th.FifthMain}.

\textit{Step 5.} We put $\overline{\operatorname{Ext}}_{S,\{\mathfrak{m}_{k}\},p}:=
\operatorname{Ext}_{S,\{\mathfrak{m}_{k}\}}\circ\operatorname{I}_{\mathfrak{m}_{0}}$. This gives a well-defined linear operator $\overline{\operatorname{Ext}}_{S,\{\mathfrak{m}_{k}\},p}:W_{p}^{1}(\operatorname{X})|_{S} \to W_{p}^{1}(\operatorname{X})$. Furthermore, by Corollary \ref{Ca.ext_inverse_trace}, Proposition \ref{Prop.from_my_previous_paper} we see that $\operatorname{Tr}|_{S}=(\operatorname{I}_{\mathfrak{m}_{0}})^{-1} \circ \operatorname{Tr}|^{\mathfrak{m}_{0}}_{S}$ and get the commutativity of the diagram in Theorem \ref{Th.FifthMain}. Since $\operatorname{I}_{\mathfrak{m}_{0}}:W_{p}^{1}(\operatorname{X})|_{S} \to W_{p}^{1}(\operatorname{X})|^{\mathfrak{m}_{0}}_{S}$ is an isometric isomorphism, we get
$\|\overline{\operatorname{Ext}}_{S,\{\mathfrak{m}_{k}\},p}\| = \|\operatorname{Ext}_{S,\{\mathfrak{m}_{k}\}}\|$.
Combining this fact with \eqref{eq.717} we get \eqref{eq.main2'} and complete the proof of Theorems \ref{Th.FourthMain}, \ref{Th.FifthMain}.

\section{Examples}

In this section we show that many results related to Problems \ref{TraceProblem} and \ref{MeasureTraceProblem}
available in the literature are particular cases of Theorems \ref{Th.SecondMain}, \ref{Th.FourthMain} and \ref{Th.FifthMain} of the present paper.
In addition, we present a model Example 11.4 which does not fall into the scope of the previously known investigations.
In fact, in some examples we present only sketches of the corresponding proofs leaving the routine verifications to the reader.

\textit{Example 11.1.} First of all, we note that in the particular case $\operatorname{X}=(\mathbb{R}^{n},\|\cdot\|_{2},\mathcal{L}^{n})$ by Theorem \ref{Th.FourthMain} and assertion (2) of Theorem \ref{Th.FifthMain}
we get a clarification of the results obtained in \cite{TV1}. Indeed, in contrast to Theorem \ref{Th.FourthMain}, the criterion presented
in Theorem 2.1 of the paper \cite{TV1} was based on a more subtle Besov-type norm. Furthermore, characterizations via the Brudnyi--Shvartsman-type functionals
were not considered in \cite{TV1}.

For the next examples we recall notation \eqref{eqq.notation_porous} and Definition \ref{Def.Ahlfors_regular_space}.

\textit{Example 11.2.} Let $p \in (1,\infty)$ and $\operatorname{X} \in \mathfrak{A}_{p}$. Assume, in addition, that $\operatorname{X}$ is an Ahlfors $Q$-regular space for some $Q > 0$.
Let $\theta \in (0,\min\{p,Q\})$ and $S \in \mathcal{ADR}_{\theta}(\operatorname{X})$. Since $\theta > 0$, we have $\mu(S) = 0$.
By Theorem 5.3 from \cite{JJKR} there exists $\sigma > 0$ such that for each $\epsilon \in (0,1]$ we have $S_{\epsilon^{k}}(\sigma)=S$ for all $k \in \mathbb{N}_{0}$, i.e.,
the set $S$ is porous. We put $\mathfrak{m}_{k}=\mathcal{H}_{\theta}\lfloor_{S}$ for all $k \in \mathbb{N}_{0}$. In accordance with Example 5.1 we have
$\{\mathfrak{m}_{k}\}:=\{\mathfrak{m}_{k}\}_{k=0}^{\infty} \in \mathfrak{M}^{str}_{\theta}(S)$. As a result,
using the equivalence between $(i)$ and $(iv)$ in Theorem \ref{Th.SecondMain} and assertion (1) of Theorem \ref{Th.FifthMain} we get the following criterion by E.~Saksman and T.~Soto \cite{SakSot} (see Theorems 1.5 and 1.7 therein):

\textit{
A function $f \in L_{p}(\mathcal{H}_{\theta}\lfloor_{S})$ belongs to the space $W^{1}_{p}(\operatorname{X})|^{\mathcal{H}_{\theta}}_{S}$ if and only if the Besov seminorm of $f$ is finite, i.e.,
$\|f|B^{1-\frac{\theta}{p}}_{p,p}(S)\|^{p}:=\sum_{k=0}^{\infty}2^{k(p-\theta)}\int_{S}(\mathcal{E}_{\mathcal{H}_{\theta}}(f,B_{2^{-k}}(x)))^{p}\,d\mathcal{H}_{\theta}(x) < + \infty$.
Furthermore, $\|F|W_{p}^{1}(\operatorname{X})|_{S}^{\mathcal{H}_{\theta}}\| \approx \|f|L_{p}(\mathcal{H}_{\theta}\lfloor_{S})\|+\|f|B^{1-\frac{\theta}{p}}_{p,p}(S)\|$ with the equivalence constants
independent on $f$. Moreover, there exists
an $\mathcal{H}_{\theta}$-extension operator $\operatorname{Ext}_{S} \in \mathcal{L}( W_{p}^{1}(\operatorname{X})|_{S}^{\mathcal{H}_{\theta}}, W_{p}^{1}(\operatorname{X}))$.
}

\textit{Example 11.3.} Let $p \in (1,\infty)$, $\operatorname{X}=(\operatorname{X},\operatorname{d},\mu) \in \mathfrak{A}_{p}$ and $S \in \mathcal{ADR}_{0}(\operatorname{X})$.
We recall Example 5.1 and note that by letting $\mathfrak{m}_{k}=\mu\lfloor_{S}$, $k \in \mathbb{N}_{0}$ we obtain $\{\mathfrak{m}_{k}\} \in \mathfrak{M}^{str}_{0}(S)$.
We recall a combinatorial result, which is a slight modification of Theorem 2.6 in \cite{Shv1}.

\begin{Prop}
\label{Prop.10.1}
Let $\mathcal{B}$ be an $(S,c)$-Whitney family of balls.
Then there exist constants $c_{1},c_{2} > 0$, $\tau \in (0,1)$, and a family $\mathcal{U}:=\{U(B):B \in \mathcal{B}\}$ of Borel subsets of $S$ such that
$U(B) \subset c_{1}B$ and $\mu(U(B)) \geq \tau \mu(B)$ for all $B \in \mathcal{B}$, $\mathcal{M}(\{U(B):B \in \mathcal{B}\}) \le c_{2}$.
\end{Prop}
Based on Proposition \ref{Prop.10.1} one can repeat, with minor modifications, the arguments used in Example 6.1 of \cite{TV1} and deduce, for each $\sigma \in (0,1)$, the existence of a constant $C > 0$ such that,
for every $f \in L^{loc}_{1}(\mu\lfloor_{S})$,
\begin{equation}
\label{eqq.10.1}
\sum\limits_{k=1}^{\infty}
2^{k(p-\theta)}\int\limits_{S_{k}(\sigma)}\Bigl(\mathcal{E}_{\mu\lfloor_{S}}(f,B_{2^{-k}}(x))\Bigr)^{p}\,d\mu\lfloor_{S}(x) \le C (\|f|L_{p}(\mu\lfloor_{S})\|+\|f^{\sharp}_{\mu\lfloor_{S}}|L_{p}(\mu\lfloor_{S})\|).
\end{equation}
Hence, using the equivalence between $(i)$ and $(iv)$ in Theorem \ref{Th.SecondMain} and assertion (1) of Theorem \ref{Th.FifthMain}
we arrive at Shvartsman's criterion \cite{Shv1}:

\textit{
A function $f \in L_{p}(\mu\lfloor_{S})$ belongs to the space $W_{p}^{1}(\operatorname{X})|^{\mu}_{S}$ if and only if
$f^{\sharp}_{S,\mu} \in L_{p}(\mu\lfloor_{S})$. Furthermore, $\|f|W_{p}^{1}(\operatorname{X})|^{\mu}_{S}\| \approx \|f|L_{p}(\mu\lfloor_{S})\|+\|f^{\sharp}_{\mu\lfloor_{S}}|L_{p}(\mu\lfloor_{S})\|$
with the equivalence constants independent of $f$. Moreover, there exists a $\mu$-extension operator $\operatorname{Ext}_{S} \in \mathcal{L}(W_{p}^{1}(\operatorname{X})|^{\mu}_{S}, W_{p}^{1}(\operatorname{X}))$.}

Now we present a model example, which exhibits an interesting effect arising in the case when we describe the trace spaces to closed sets composed of pieces of different dimensions.

\textit{Example 11.4.} Let $p \in (1,\infty)$ and $\operatorname{X}=(\operatorname{X},\operatorname{d},\mu) \in \mathfrak{A}_{p}$. Assume that $\operatorname{X}$ is Ahlfors $Q$-regular
for some $Q > 1$. Let $\underline{B}=B_{R}(x)$ be a closed ball of radius $R > 0$ centered at $x \in \operatorname{X}$ and $\gamma: [0,1] \to \operatorname{X}$ be a rectifiable curve such that
$\Gamma:=\gamma([0,1]) \in \mathcal{ADR}_{Q-1}(\operatorname{X})$, $\Gamma \cap \underline{B} = \{\underline{x}\}$ for some point $\underline{x} \in \partial \underline{B}$, and, furthermore,
for some $\kappa > 0$, $\operatorname{dist}(\gamma(t),\underline{B}) \geq \kappa l(\gamma([0,t]))$ for all $t \in [0,1]$. We put $S:=\underline{B} \cup \Gamma$ and recall Example 5.2.
For each $k \in \mathbb{N}_{0}$ we set $\mathfrak{m}_{k}:=2^{k(Q-1)}\mu\lfloor_{\underline{B}}+\mathcal{H}_{Q-1}\lfloor_{\Gamma}$. Hence, we get $\{\mathfrak{m}_{k}\}:=\{\mathfrak{m}_{k}\}_{k=0}^{\infty} \in \mathfrak{M}^{str}_{Q-1}(S)$.

Given $k \in \mathbb{N}_{0}$, we introduce the $k$th gluing functional by letting for each $f \in L_{1}^{loc}(\{\mathfrak{m}_{k}\})$,
\begin{equation}
\label{eqq.k_gluing}
\operatorname{gl}_{k}(f,\underline{x}):=\fint\limits_{\underline{B} \cap B_{k}(\underline{x})}\fint\limits_{\Gamma \cap B_{k}(\underline{x})}|f(y)-f(z)|\,d\mu(y)\,d\mathcal{H}_{Q-1}(z).
\end{equation}

We put $\underline{k}:=\min\{k \in \mathbb{N}:2^{k} > \frac{1}{\kappa}\}$ and $\alpha:=1-\frac{Q-1}{p}$. Using Remark \ref{Rem.best_approx} and letting $B_{k}(x):=B_{2^{-k}}(x)$, $S_{k}(\sigma):=S_{2^{-k}}(\sigma)$, it is easy to see that, given $\sigma \in (0,1)$,
\begin{equation}
\begin{split}
\label{eqq.10.3}
&\Bigl(\sum\limits_{k=1}^{\underline{k}+1}
2^{k\alpha p}\int\limits_{S_{k}(\sigma)}\Bigl(\mathcal{E}_{\mathfrak{m}_{k}}(f,B_{k}(x))\Bigr)^{p}\,d\mathfrak{m}_{k}(x)\Bigr)^{\frac{1}{p}} \le C \Bigl(\|f|L_{p}(\mu\lfloor_{\underline{B}})\|+
\|f|L_{p}(\mathcal{H}_{Q-1}\lfloor_{\Gamma})\|\Bigr).
\end{split}
\end{equation}
Using arguments given in Examples 11.2, 11.3 and taking into account that for each $k > \underline{k}$ we have $S \cap B_{k}(x) = \underline{B} \cap B_{k}(x)$ for all $x \in \underline{B} \setminus B_{k-\underline{k}}(\underline{x})$
and $S \cap B_{k}(x)  = \Gamma \cap B_{k}(x)$ for all $x \in \Gamma \setminus B_{k-\underline{k}}(\underline{x})$, one can show that, given a small enough $\sigma \in (0,1)$ (depending on $\Gamma$),
\begin{equation}
\label{eqq.10.4}
\begin{split}
&\sum\limits_{k=\underline{k}+2}^{\infty}
2^{k\alpha p}\int\limits_{S_{k}(\sigma) \setminus B_{k-\underline{k}}(\underline{x})}\Bigl(\mathcal{E}_{\mathfrak{m}_{k}}(f,B_{k}(x))\Bigr)^{p}\,d\mathfrak{m}_{k}(x) \le C \Bigl(\|f^{\sharp}_{\mu\lfloor_{\underline{B}}}|L_{p}(\mu\lfloor_{\underline{B}})\|^{p}+
\|f|B^{\alpha}_{p,p}(\Gamma)\|^{p}\Bigr).
\end{split}
\end{equation}

Since $1/\epsilon \geq 10$, given $x \in S_{k}(\sigma) \cap  B_{k-\underline{k}}(\underline{x})$ and $k > \underline{k}$, we have 
$B_{k}(x) \subset B_{k-\underline{k}-1}(\underline{x})$ and 
$B_{k-\underline{k}-1}(\underline{x}) \subset B_{k-\underline{k}-2}(x)$.
By Remark \ref{Rem.best_approx} and Theorem \ref{Th.doubling_type}, $\mathcal{E}_{\mathfrak{m}_{k}}(f,B_{k}(x)) \le C\operatorname{gl}_{k-\underline{k}-1}(f,\underline{x})$. Furthermore, it is easy to see that $\mathfrak{m}_{k}(B_{k-\underline{k}}(\underline{x})) \le C2^{-k}$.
As a result,
\begin{equation}
\label{eqq.10.5}
\begin{split}
&\int\limits_{S_{k}(\sigma) \cap  B_{k-\underline{k}}(\underline{x})}\Bigl(\mathcal{E}_{\mathfrak{m}_{k}}(f,B_{k}(x))\Bigr)^{p}\,d\mathfrak{m}_{k}(x) \le C 
2^{-k}\Bigl(\operatorname{gl}_{k-\underline{k}-1}(f,\underline{x})\Bigr)^{p}.
\end{split}
\end{equation}
On the other hand, given $x \in S_{k}(\sigma) \cap  B_{k-\underline{k}}(\underline{x})$ and $k > \underline{k}+1$, we have $B_{k-\underline{k}}(\underline{x}) \subset B_{k-\underline{k}-1}(x)$ and 
$B_{k-\underline{k}-1}(x) \subset B_{k-\underline{k}-2}(\underline{x})$. Hence, by Theorem \ref{Th.doubling_type} and Remark \ref{Rem.best_approx} it is easy to get 
$\operatorname{gl}_{k-\underline{k}-1}(f,\underline{x}) \le C\mathcal{E}_{\mathfrak{m}_{k}}(f,B_{k-\underline{k}-2}(x))$.
If $\sigma \in (0,1)$ is small enough, then one can show that $\mathfrak{m}_{k}(S_{k}(\sigma) \cap B_{k}(\underline{x})) \approx 2^{-k}$. As a result,
\begin{equation}
\label{eqq.10.6}
\begin{split}
&\Bigl(\operatorname{gl}_{k-\underline{k}}(f,\underline{x})\Bigr)^{p} \le C 2^{k}\int\limits_{S_{k}(\sigma) \cap  B_{k}(\underline{x})}\Bigl(\mathcal{E}_{\mathfrak{m}_{k}}(f,B_{k-\underline{k}-1}(x))\Bigr)^{p}\,d\mathfrak{m}_{k}(x).
\end{split}
\end{equation}

At the same, time it follows from Examples 11.2, 11.3 that there is $C > 0$ such that for all $F \in W_{p}^{1}(\mathbb{R}^{n})$ with $f=F|_{S}^{\mathfrak{m}_{0}}$
\begin{equation}
\label{eqq.11.7}
\|f|L_{p}(\mu\lfloor_{\underline{B}})\|+
\|f|L_{p}(\mathcal{H}_{Q-1}\lfloor_{\Gamma})\|+\|f^{\sharp}_{\mu\lfloor_{\underline{B}}}|L_{p}(\mu\lfloor_{\underline{B}})\|+
\|f|B^{\alpha}_{p,p}(\Gamma)\| \le C \|F|W_{p}^{1}(\mathbb{R}^{n})\|.
\end{equation}

Finally, combining \eqref{eqq.10.3}--\eqref{eqq.11.7}, it is easy to deduce from Theorem \ref{Th.SecondMain} and assertion (1) of Theorem \ref{Th.FifthMain} the following result.

\textit{A function $f \in L_{p}(\mu\lfloor_{\underline{B}}) \cap L_{p}(\mathcal{H}_{\theta}\lfloor_{\Gamma})$
belongs to $W_{p}^{1}(\operatorname{X})|^{\mathfrak{m}_{0}}_{S}$ if and only if
$f^{\sharp}_{\mu\lfloor_{\underline{B}}} \in L_{p}(\mu\lfloor_{\underline{B}})$, $f \in B^{1-\frac{Q-1}{p}}_{p,p}(\Gamma)$ and $\Bigl(\mathcal{GL}(f,\underline{x})\Bigr)^{p}:=\sum_{k=1}^{\infty}
2^{k(p-Q)}(\operatorname{gl}_{k}(f,\underline{x}))^{p} < +\infty$.
Furthermore,
\begin{equation}
\begin{split}
&\|f|W_{p}^{1}(\operatorname{X})|^{\mathfrak{m}_{0}}_{S}\| \approx \|f|L_{p}(\mu\lfloor_{\underline{B}})\|+\|f|L_{p}(\mathcal{H}_{\theta}\lfloor_{\Gamma})\|\\
&+(\|f^{\sharp}_{\mu\lfloor_{\underline{B}}}|L_{p}(\mu\lfloor_{\underline{B}})\|+
\|f|B^{1-\frac{Q-1}{p}}_{p,p}(\Gamma)\|+\mathcal{GL}(f,\underline{x}).
\end{split}
\end{equation}
with the equivalence constants independent of $f$. Moreover, there exists an $\mathfrak{m}_{0}$-extension operator $\operatorname{Ext}_{S,\{\mathfrak{m}_{k}\}} \in \mathcal{L}(W_{p}^{1}(\operatorname{X})|^{\mathfrak{m}_{0}}_{S}, W_{p}^{1}(\operatorname{X}))$.}

We conclude by presenting a natural generalization of Theorem 1.2 from \cite{Shv2}.

\textit{Example 11.5.} Let $\operatorname{X}=(\operatorname{X},\operatorname{d},\mu)$ be a geodesic metric measure space. Assume that $\operatorname{X}$ is Ahlfors $Q$-regular for some $Q \geq 1$.
Furthermore, assume that $p \in (Q,\infty)$ and $\operatorname{X} \in \mathfrak{A}_{p}$. We also fix a parameter $\theta \in [Q,p)$ and a nonempty closed set $S \subset \operatorname{X}$. For simplicity we assume that $S \subset B_{1}(\underline{x})$ for some $\underline{x} \in \operatorname{X}$.
According to Example 4.9, we have $S \in \mathcal{LCR}_{\theta}(\operatorname{X})$.

Combining results from \cite{AGS2} with Theorem 9.1.15 from \cite{HKST} we obtain that, for each
$F \in W_{p}^{1}(\operatorname{X})$,  there exists a \textit{continuous representative} $\overline{F}$ such that, for each closed ball $B$,
\begin{equation}
\label{eqq.11.8}
\sup_{x \in B}|\overline{F}(x)-F_{B}| \le C r(B)\Bigl(\fint\limits_{B}(|DF|_{p}(y))^{p}\,d\mu(y)\Bigr)^{\frac{1}{p}}.
\end{equation}
Furthermore, it is clear that, given $F \in W_{p}^{1}(\operatorname{X})$, the $p$-sharp trace of $F|_{S}$ is well defined
and coincides with the pointwise restriction of $\overline{F}$ to the set $S$.

By Remark \ref{Rem.best_approx} it is easy to see that if a ball $B_{r}(x)$ with $r \in (0,1]$, $x \in \operatorname{X}$ is such that $B_{cr}(x) \cap S \neq \emptyset$ for some $c \geq 1$, then,
for each sequence $\{\mathfrak{m}_{k}\} \in \mathfrak{M}^{str}_{\theta}(S)$,
\begin{equation}
\label{eqq.11.9}
\mathcal{E}_{\mathfrak{m}_{k(r)}}(f,B_{2cr}(x)) \le \sup_{y,z \in B_{2cr}(x) \cap S}|f(y)-f(z)| \quad \text{for all} \quad f \in C(S).
\end{equation}

We define the modified \textit{Brudnyi--Shvartsman-type} functional on $C(S)$ (with values in $[0,+\infty]$) as follows. Given a function $f \in C(\operatorname{X})$, we put
\begin{equation}
\notag
\widetilde{\mathcal{BSN}}_{p}(f):=\sup \Bigl(\sum\limits_{i=1}^{N}\frac{\mu(B_{r_{i}})(x_{i})}{r_{i}^{p}}\sup\limits_{y,z \in B_{60r_{i}}(x_{i})\cap S}|f(y)-f(z)|^{p}\Bigr)^{\frac{1}{p}},
\end{equation}
where the supremum is taken over all finite disjoint families of closed balls $\{B_{i}\}_{i=1}^{N}$ in $\operatorname{X}$ with radii $\le 1$. In the particular
case $\operatorname{X}=(\mathbb{R}^{n},\|\cdot\|_{\infty},\mathcal{L}^{n})$, our functional is very similar to that used in \cite{Shv2}.
The only difference is that in \cite{Shv2} the corresponding coefficient of dilatation of balls was $11$ rather than $60$.
Using \eqref{eqq.11.8} and Proposition \ref{Prop.finite_intersection} it is easy to see that
\begin{equation}
\label{eqq.11.10}
\widetilde{\mathcal{BSN}}_{p}(F|_{S})+\inf\limits_{z \in S}|F|_{S}(z)| \le C \|F|W_{p}^{1}(\operatorname{X})\| \quad \text{for all} \quad F \in W_{p}^{1}(\operatorname{X}).
\end{equation}
On the other hand,
if $x_{m} \in S$ is a minimum point and $x_{M} \in S$ is a maximum point,
then it is easy to see that $|f(x_{m})-f(x_{M})| \le (\mu(B_{1}(\underline{x})))^{-\frac{1}{p}} \widetilde{\mathcal{BSN}}_{p}(f)$. As a result, by \eqref{eqq.11.9}
\begin{equation}
\label{eqq.11.11}
\operatorname{BSN}_{p,\{\mathfrak{m}_{k}\},30}(f) \le \widetilde{\mathcal{BSN}}_{p}(f)+\sup_{x \in S}|f(x)| \le C(\widetilde{\mathcal{BSN}}_{p}(f)+\inf_{z \in S}|f(z)|), \quad  f \in C(S).
\end{equation}
Finally, we apply Theorem \ref{Th.FourthMain} with $\epsilon=\frac{1}{10}$, $c=30$ and use \eqref{eqq.11.10}, \eqref{eqq.11.11}. This leads to the following criterion.

\textit{A function $f \in C(S)$ belongs to $W_{p}^{1}(\operatorname{X})|_{S}$ if and only if $\widetilde{\mathcal{BSN}}_{p}(f)<+\infty$. Furthermore,
$\|f|W_{p}^{1}(\operatorname{X})|_{S}\| \approx \widetilde{\mathcal{BSN}}_{p}(f)+\inf_{z \in S}|f(z)|$, with the equivalence constants independent on $f$. Moreover, there
exists a $p$-sharp extension operator $\operatorname{Ext}_{S,p} \in \mathcal{L}(W_{p}^{1}(\operatorname{X})|_{S}, W_{p}^{1}(\operatorname{X}))$}.

\end{document}